\documentclass[11pt]{article}
\usepackage{pifont}
\usepackage{amsfonts,amsmath,amssymb,amsthm}
\usepackage{mathrsfs, amscd}
\usepackage{amsbsy}
\usepackage{xy}
\xyoption{all}
\usepackage{color}

\newcommand{\mput}{\multiput}
\newcommand{\bcen}{\begin{center}}
\newcommand{\ecen}{\end{center}}

\def\dz{\delta}
\def\lz{\lambda}

\def\A{{\cal A}}

\def\C{{\cal C}}

\def\G{{\cal G}}

\def\K{{\cal K}}
\def\U{{\cal U}}
\def\V{{\cal V}}
\def\M{{\cal M}}

\def\N{{\cal N}}
\def\T{{{\cal T}}}

\def\f{\textsf{f}}
\def\tl{\mbox l}
\def\tr{\mbox r}
\def\tm{\mbox m}
\def\td{\mbox d}
\def\te{\mbox e}

\def\l{\underline{l}}
\def\n{{\underline{n}}}
\def\m{{\underline{m}}}

\def\tt{{t\times t}}

\def\rad{\mbox{rad}}

\def\wt{\widetilde}

\def\Hom{\mbox{Hom}}

\def\Hom{\mbox{Hom}}

\def\dim{\mbox{dim}}
\def\mod{\mbox{mod}}

\def\Im{\mbox{Im}}

\def\rank{\mbox{rank}}
\def\mf{\mathfrak}

\def\mf{\mathfrak}

\def\IM{\hbox{\rm I \hskip -4.9pt M}}
\def\IE{\hbox{\rm I\hskip -2pt E}}
\def\IF{\hbox{\rm I\hskip -2pt F}}
\def\id{\hbox{\rm 1\hskip -2.2pt l}}

\begin{document}

\begin{center}{\Large\bf Homogeneity implies Tameness}

\end{center}

\begin{center} Zhang Yingbo, Beijing Normal University\end{center}
\begin{center} Xu Yunge, Hubei University\end{center}

\bigskip
\centerline{\bf Introduction}
\bigskip

Throughout the paper we always assume that $k$ is an algebraically
closed field, all the rings or algebras contain identities. We write
our maps on either the left or the right,
but always compose them as if they were written on the right.

In 1977 Drozd showed in \cite{D1} that  a finite-dimensional algebra
over an algebraically closed field  is either of tame representation
type or of wild representation type, which has been one of
foundations of representation theory of finite dimensional algebras.

\medskip

{\bf Definition 1} \cite{D1, CB1, DS}\, A finite-dimensional
$k$-algebra $\Lambda$ is of tame representation type, if for any
positive integer $d$, there are a finite number of localizations
$R_i=k[x,\phi_i(x)^{-1}]$ of $k[x]$ and $\Lambda$-$R_i$-bimodules
$T_i$ which are free  as right $R_i$-modules, such that almost all
(except finitely many) indecomposable $\Lambda$-modules of dimension
at most $d$ are isomorphic to
$$T_i\otimes_{R_i}R_i/(x-\lz)^{m},$$ for some  $\lz\in k,$
$\phi_i(\lz)\ne 0$, and some positive integer $m$.

\medskip

{\bf Definition 2} \cite{D1, CB1}\, A finite-dimensional
$k$-algebra $\Lambda$ is  of wild representation type if there is a
finitely generated $\Lambda$-$k\langle x,y\rangle$-bi-module $T$,
which is free as a right $k\langle x,y\rangle$-module, such that the
functor
$$T\otimes_{k\langle x,y\rangle}-: k\langle
x,y\rangle\mbox{-mod}\rightarrow \Lambda\mbox{-mod}$$ preserves
indecomposability and isomorphism classes.

\medskip

The proof of the Drozd's Tame-Wild Theorem is highly indirect. The
argument relies on the notion of a bocs (the brief for ``bi-module
of co-algebra structure"), introduced first by Rojter in \cite{Ro}.
In 1988, Crawley-Boevey formalized the theory of bocses and showed
that for a tame algebra $\Lambda$, and for each dimension $d$, all
but finitely many isomorphism classes of indecomposable
$\Lambda$-modules of dimension $d$ are isomorphic to their
Auslander-Reiten translations and hence belong to homogeneous tubes.

Since then, many authors have tried to prove the
converse of Crawley-Boevey's Theorem. They expected to find
infinitely many non-isomorphic indecomposable representations
$\{M_i\mid i\in I\}$ of the same dimension in the representation
category of a layered bocs such that $M_i\ncong DTr(M_i)$. However,
four authors constructed a  strongly homogeneous wild layered locs
$\mf{B}$ in \cite{BCLZ} in 2000 such that each representation of
$\mf{B}$ is homogeneous (i.e., $DTr(M)\cong M$). It shows that the
converse of the Crawley-Boevey Theorem does not hold true for
layered bocses in general. But it is still open for finite dimensional
$k$-algebras.

Let $\Lambda$ be  a finite-dimensional algebra over an algebraically
closed field. We say that an indecomposable $\Lambda$-module $M$ is
{\it homogeneous} if $DTr(M)\cong M$. The category $\mod\Lambda$ is
said to be {\it homogeneous} provided that for each dimension $d$
all but finitely many isomorphism classes of indecomposable
$\Lambda$-modules of dimension $d$ are homogeneous. Our purpose in
this paper is to prove the following theorem

\medskip

{\bf Main Theorem 3} Let $\Lambda$ be  a finite-dimensional
algebra over an algebraically closed field. Then $\Lambda$ is of
tame representation type if and only if $\mod\Lambda$ is
homogeneous.

\medskip

The necessity of the theorem is proved by Crawley-Boevey. Our proof
for the sufficiency relies on the notions of matrix bi-module
problem and its associated bocs, as well as their reduction
techniques. The key of the argument is to find a full subcategory of
representation category of a bipartite matrix bi-module problem
which admits infinitely many non-homogeneous representations of dimension $d$, and the
fact that matrix bi-module problems associated to finite-dimensional algebras
are bipartite.

\bigskip
\bigskip
\bigskip

\centerline{\bf 1 Matrix Bi-module Problems}

\bigskip

Matrix problems, which include special cases the representation
theory of finite-dimensional algebras, subspace problems and
projective modules, get their importance in studying questions about
representation type. In this section we will introduce the notions
of matrix bi-module problems over some minimal algebras and their associated bi-co-module problems,
which seems to be more convenient for calculational purposes. We will unify some
established notions for matrix problems, such as linear matrix
problem [S], bocs \cite{Ro,D1,CB1}, differential graded category \cite{Ro}
and so on, to formalize the matrix methods and approach to
representation theory.


\bigskip
\bigskip

\noindent{\bf \S 1.1 Matrix bi-module problems}

\bigskip

In the present subsection, we will construct a $k$-algebra $\Delta$ of
a minimal algebra $R$, then define matrix bi-module problems over
$\Delta$. The concepts and the results are proposed by S. Liu.

\medskip

Let ${\cal T}$ be a vertex set whose elements are divided into two
disjoint families: the subset $\T_0$ of trivial vertices and the
subset $\T_1$ of non-trivial vertices.  To each $X \in \T_1$, we
associate an indeterminate $x$ and a fixed non-zero polynomial
$\phi_{_X}(x)$ in $k[x]$. Given any $X\in {\cal T}$, we define a
$k$-algebra $R_X$ with identity $1_X$ by $R_X=k$ if $X$ is trivial;
and otherwise, $R_X=k[x, \phi_{_X}(x)^{-1}]$ the localization of
$k[x]$ at $\phi_{_X}(x)$, and $x$ is said to be the {\it
parameter associated to $X\in {\cal T}_1$}. Now we call the
$k$-algebra $R=\Pi_{X\in {\cal T}}R_X$ a {\it minimal algebra} over
$\cal T$ with a set of orthogonal primitive idempotents $\{1_{_X}\mid X\in\T\}$.

\medskip

Define a tensor product of $p\geqslant 1$ copies of $R$ over $k$:
$$\begin{array}{c}
R^{\otimes p}=R\otimes_k\cdots\otimes_k R=\sum_{(X_1,\cdots,X_p)
\in\Pi^p\T}R_{X_1}
\otimes_k\cdots\otimes_kR_{X_{p}}.\end{array}\eqno{(1.1\mbox{-}1)}$$
There exists a natural left and right $R$-module structure on $R^{\otimes p}$:
for any $\alpha=r_1\otimes_k r_2\otimes_k\cdots\otimes_k
r_{p-1}\otimes_kr_p\in R^{\otimes p}$, $s\in R$, the left module action is given by:
$s\otimes_{_R}\alpha=(s\otimes_{_R}r_1)\otimes_kr_2\otimes_k\cdots\otimes_kr_p$;
and the right one by:
$\alpha\otimes_{_R}s=r_1\otimes_k\cdots\otimes_kr_{p-1}\otimes_k(r_p\otimes_{_R}s)$.
If $r_i\in R_{X_i}, s\in R_Y$, then $s\otimes_{_{R}}\alpha=0$ for $Y\ne X_1$ and
$\alpha\otimes_{_R}s=0$ for $Y\ne X_p$. Thus $R^{\otimes p}$ can be viewed as an
{\it $R$-$R$-bi-module}, or simply an $R^{\otimes 2}$-module, with the module action:
$$(r\otimes_ks)\otimes_{_{R^{\otimes 2}}}\alpha=r\otimes_{_R}\alpha\otimes_{_R}s=
\alpha\otimes_{_{R^{\otimes 2}}}(r\otimes_ks),\quad \forall\, r,s\in R.\eqno{(1.1\mbox{-}2)}$$

Consider the direct sum of $R^{\otimes p}$, which is still an $R^{\otimes 2}$-module:
$$\begin{array}{c}\Delta=\oplus_{p=1}^\infty R^{\otimes p},\quad \mbox{let}\,\,
\bar\Delta=\oplus_{p=2}^\infty R^{\otimes p},\, \Delta=R\oplus\bar\Delta.\end{array}\eqno{(1.1\mbox{-}3)}$$
We define a multiplication on $R^{\otimes 2}$-module $\Delta$, given by $\Delta\times\Delta
\rightarrow\Delta\otimes_{_R}\Delta\subseteq\Delta$:
$$\begin{array}{c}\Delta^{\otimes p}\otimes_{_R}\Delta^{\otimes q}\subseteq\Delta^{\otimes (p+q-1)},\quad\alpha\otimes_{_R}\beta=r_1\otimes_k
\cdots\otimes_k(r_ps_1)\otimes_ks_2\cdots\otimes_ks_q,
\end{array}\eqno{(1.1\mbox{-}4)}$$
where $\beta=s_1\otimes_k\cdots\otimes_ks_q$, and if $r_i\in R_{X_i},s_j\in R_{Y_j}$,
$\alpha\otimes_{_R}\beta=0$ for $X_p\ne Y_1$. Thus we obtain {\it an associative non-commutative
$k$-algebra $(\Delta,\otimes_{_R},1_R)$}
with the set of orthogonal primitive idempotents $\{1_X\mid X\in\T\}$.

Moreover $\Delta\otimes_{_R}\Delta$ can be viewed as an $R^{\otimes 3}$-module:
for any $\alpha,\beta\in\Delta, r,s,w\in R$,
$$\begin{array}{c}(r\otimes_ks\otimes_kw)
\otimes_{_{R^{\otimes 3}}}(\alpha\otimes_{_R}\beta)=(\alpha\otimes_k\beta)
\otimes_{_{R^{\otimes 3}}}(r\otimes_ks\otimes_kw)\\
=r\otimes_{_R}\alpha\otimes_{_R}s\otimes_{_R}\beta\otimes_{_R}w.
\end{array}\eqno{(1.1\mbox{-}5)}$$

Denote by $\IM_{m\times n}(\Delta)$ the set of
{\it matrices over $\Delta$} of size $m\times n$; and
by $\mathbb T_n(\Delta), \mathbb N_n(\Delta)$, $\mathbb D(\Delta)$ the set of
$n\times n$-upper triangular, strictly upper triangular, and diagonal $\Delta$-matrices respectively.
The product of two $\Delta$-matrices is the {\it usual matrix product}.
If $H=(h_{ij})\in$IM$_{m\times n}(R)$, $U=(u_{ij})\in$IM$_{m\times n}(R\otimes_kR)$, $\alpha\in\Delta$, define
$$\begin{array}{c} H\otimes_{_R}\alpha=(h_{ij}\otimes_{_R}\alpha)\in\IM_{m\times n}(\Delta),\\
\alpha\otimes_{_R}H=(\alpha\otimes_{_R} h_{ij})\in\IM_{m\times n}(\Delta);\\
U\otimes_{_{R^{\otimes 2}}}\alpha=\big(\alpha\otimes_{_{R^{\otimes 2}}}u_{ij}
\big)=\alpha\otimes_{_{R^{\otimes 2}}}U\in\IM_{m\times n}(\Delta).
\end{array}\eqno{(1.1\mbox{-}6)}$$
Since for any $u,v\in R^{\otimes 2},\delta\in R^{\otimes 3}$, $\big((\alpha\otimes_{R^{\otimes 2}}u)
\otimes_R(\beta\otimes_{R^{\otimes 2}}v)\big)
\otimes_{R^{\otimes 3}}\delta=(\alpha\otimes_R\beta)
\otimes_{R^{\otimes 3}}\big((u\otimes_R v)\otimes_{R^{\otimes 3}}\delta\big)$,
if $V=(v_{jl})\in$IM$_{n\times r}(R\otimes_kR)$, $\beta\in\Delta$, we have
$$\begin{array}{c}(\alpha\otimes_{_{R^{\otimes 2}}}U)(\beta\otimes_{_{R^{\otimes 2}}}V)
=(\alpha\otimes_{_R}\beta)\otimes_{_{R^{\otimes 3}}}(UV).
\end{array}\eqno{(1.1\mbox{-}7)}$$

An $R$-$R$-bi-module
${\mathcal S}_1$ is said to be a {\it quasi-free bi-module} finitely generated by
$U_1,\cdots,U_m$, if the morphism
$$(R_{X_1}\otimes_kR_{Y_1})\oplus \cdots \oplus (R_{X_m}\otimes_kR_{Y_m})\rightarrow
\mathcal S_1,\quad 1_{X_i}\otimes_k 1_{Y_i}\mapsto U_i$$
is an isomorphism. In this case, $\{U_1,\ldots,U_m\}$ is called an
{\it $R$-$R$-quasi-free basis} of $\mathcal S_1$.

Let $R$-$R$-bi-module
$\mathcal S_{p}=R^{\otimes (p+1)}\otimes_{R^{\otimes 2}}\mathcal S_1$ be given by
$(r\otimes_ks)\otimes_{R^{\otimes 2}}(\alpha\otimes_{R^{\otimes 2}}U)=(r\alpha s)\otimes_{R^{\otimes 2}}U
=\alpha\otimes_{R^{\otimes 2}}((r\otimes_ks)\otimes_{R^{\otimes 2}}U)$,
for $r,s\in R,\alpha\in R^{p+1},U\in\mathcal S$, since it is valid on the basis.
Moreover, $\mathcal S=\sum_{p=1}^\infty\mathcal S_p=\bar\Delta\otimes_{R^{\otimes 2}}\mathcal S_1$
is an $R$-$R$-bi-module, and $\mathcal S_p$ is said to be {\it index $p$}.
\medskip

{\bf Definition 1.1.1}\, Let $T=\{1,2,\cdots,t\}$ be a set of integers,
and let $\sim$ be an equivalent relation on
$T$, such that there is a one-to-one correspondence between the set
$T/\sim$ of equivalence classes and the set $\T$ of vertices of a
minimal algebra $R$. We may write $\T=T/\sim$.

\medskip

{\bf Definition 1.1.2}\, (i)\, Define an $R$-$R$-bi-module
$\K_0=\{$diag$(s_{11},\cdots,s_{tt})\mid s_{ii}\in R_X,\forall\, i\in X;$ and $s_{ii}=s_{jj},
\forall\, i\sim j\}$, which is isomorphic to $R$ as algebras.
Set $E_X\in\mathbb D_{t}(R_X)$ with the element $s_{ii}=1_X$ if $i\in X$ and $s_{ii}=0$ if $i\notin X$, then
$\{E_X\mid X\in\T\}$ is called a quasi-free $R$-basis of $\K_0$.

(ii)\, Define a quasi-free $R$-$R$-bi-module
$\K_1\subseteq\mathbb N(R\otimes_k R)$ with an $R$-$R$-quasi-basis:
$$\begin{array}{c}\V=\cup_{(X,Y)\in\T\times\T}\V_{XY}=\{V_1,V_2,\cdots,V_m\},
\quad \V_{XY}\subset\mathbb N_t(R_{X}\otimes_k R_{Y}).\end{array}$$
Write $1_XV1_Y=V$ for $V\in\V_{XY}$.

(iii)\, Suppose $\K=\K_0
\oplus(\bar\Delta\otimes_{_{R^{\otimes 2}}}\K_1)$ possesses an algebra structure, where
multiplication $\tm: \K\times\K\rightarrow \K$ is the usual matrix product over $\Delta$; unit
$\te: R\cong \K_0\hookrightarrow\K$.

\medskip

The $k$ algebra $\K$ is said to be {\it finitely generated in index
$(0,1)$ over $\Delta$}, because $V_i V_j\in R^{\otimes 3}\otimes_{R^{\otimes 2}}\K_1$.
$E_XV_j=V_j$ if $ 1_XV_j=V_j$, or $0$ otherwise, and similarly for
$V_jE_X$, thus $\{E_X\mid X\in\T\}$ is a set of orthogonal
primitive idempotents of $\K$, $E=\te (1_R)=\sum_{X\in\T}E_X$ is the {\it identity matrix}.
$\tm_{pq}:\K_p\times\K_q\rightarrow \K_{p+q}$, since $R^{\otimes (p+1)}\otimes_R
R^{\otimes (q+1)}\simeq R^{\otimes(p+q+1)}$.

\medskip

Let $ T=\{1,2,\ldots,t\}$ and $ T'=\{1,2,\ldots,t'\}$ be two sets
of integers. An {\it order} on $ T\times T'$ is defined as follows:
$(i,j)\preccurlyeq (i',j')$ provided that $i>i'$, or $i=i'$ but
$j\leqslant j'$. Thus we obtain an order on the index set of entries of a matrix in
$\IM_{t\times t'}(\Delta)$.
Let $M=(\lambda_{ij})\in \IM_{t\times t'}(\Delta)$, $\lambda_{pq}$ is
said to be the {\it leading entry} of
$M$ if $\lambda_{pq}\ne 0$, and any $\lambda_{ij}\ne 0$ implies that
$(p,q)\preccurlyeq (i,j)$. Let $\bar{M}=(C_{ij})$ be a partitioned
matrix over $\Delta$, one defines similarly the {\it leading
block} of $\bar{M}$.  In both cases, the pair $(p,q)$ is
called the {\it leading position} of $M$ resp. $\bar{M}$.

Let ${\cal S}$ be a subspace of
$\IM_t(k)$. An ordered basis $\U=\{U_1,U_2,\ldots,U_r\}$ with the
leading positions $(p_1,q_1), \ldots,(p_r,q_r)$ respectively is
called a {\it normalized basis} of ${\cal S}$ provided that
$$\begin{array}{l}\mbox{(i)\, the leading entry of}\,\, U_i\,\, {\mbox{is}}\, 1;\\
\mbox{(ii)\, the}\,\, (p_i,q_i)\mbox{-entry of}\,\, U_j\,\, \mbox{is}\,\, 0\,\, \mbox{for}\,\,j\ne i;\\
\mbox{(iii)\,}\, U_i\preccurlyeq U_j\,\, \mbox{if and only if}\,\, (p_i,q_i)\preccurlyeq
(p_j,q_j).\end{array}\eqno{(1.1\mbox{-}8)}$$
The basis  $\U$ is a linear ordered set.
It is easy to see that ${\cal S}$ has a normalized basis by Linear algebra.
In fact, taken $t^2$ variables $x_{ij}$ under the order of matrix indices defined as above. Then ${\cal
S}$ will be the solution space of some  system of linear equations
$$\begin{array}{c}
\sum_{(i,j)\in\T\times\T} a_{ij}^lx_{ij}=0, \quad a^l_{ij}\in k,
\quad 1\leqslant l\leqslant s\end{array}
$$
for some positive integer $s$. Reducing the coefficient matrix to the simplest echelon
form, we assume that $x_{p_1q_1}, x_{p_2q_2},$ $\ldots,x_{p_rq_r}$
are all the free variables. Evaluated $x_{p_iq_i}$ the column vector
whose $(p_i,q_i)$-entry is $1$ and $(p,q)\prec (p_i,q_i)$-entry is $0$ for $i=1,2,\ldots,r$,
we obtain a normalized basis of ${\cal S}$.

\medskip

{\bf Definition 1.1.3}\, (i)\, Define a quasi-free $R$-$R$-bi-module
$\M_1\subseteq\IM_t(R\otimes_k R)$, such that $E_X\M_1E_Y$ has a
normalized basis $\A_{XY}\subseteq\IM_t(k1_X\otimes_k1_Yk)\simeq\IM_t(k)$, write $1_XA1_Y=A$
for $A\in\A_{XY}$. Thus there is a {\it normalized basis}:
$$\begin{array}{c}\A=\cup_{(X,Y)\in\T\times\T}\A_{XY}=\{A_1,A_2,\cdots,A_n\}.\end{array}$$

(ii)\, Let $\M=\bar\Delta\otimes_{_{R^{\otimes 2}}}\M_1$, and the
algebra $\K$ be given by definition 1.1.2.
Define an $\K$-$\K$-bi-module structure on $\M$, with a left module action
$\tl: \K\times\M\rightarrow \M$ and a right one
$\tr: \M\times\K\rightarrow \M$ given by the usual
matrix product.

\medskip

$\tl_{pq}:\K_p\times\M_q\rightarrow\M_{p+q}$
and $\tr_{pq}: \M_p\times\K_q\rightarrow\M_{p+q}$.
The $\K$-$\K$-bi-module $\M$ is said to be {\it finitely generated in index
$(0,1)$} with $\M_0=\{0\}$.

\medskip

{\bf Definition 1.1.4}\, Let $H=\sum_{X\in\T}H_X\in\IM_t(R)$ be a matrix,
$H_X=(h_{ij})_{t\times t}\in E_X\IM(R)E_X$ with $h_{ij}\in
R_X$ for $i,j\in X$, and $h_{ij}=0$ otherwise. Suppose $H$ yields a
derivation $\td: \K\rightarrow \M, U\mapsto UH-HU$, such that
$\td(\K_0)=\{0\};\td(\K_1)\subseteq \M_1$. Clearly $\td_p:\K_p\rightarrow \M_p$.

\medskip

{\bf Definition 1.1.5}\, A quadruple $\mf{A}=(R,\K,\M,H)$ is called a
{\it matrix bi-module problem} provided

(i)\, $R$ is a minimal algebra with a vertex set $\T$;

(ii)\, $\K$ is an algebra given by Definition 1.1.2;

(iii)\, $\M$ is a $\K$-$\K$-bi-module given by Definition 1.1.3;

(iv)\, There is a derivation $\td: \K\rightarrow \M$
given by Definition 1.1.4.

\smallskip

In particular, if $\M=0$, $\mf{A}$ is said to be a {\it minimal matrix
bi-module problem}.

\bigskip
\bigskip

\noindent{\bf \S 1.2 Bi-co-module problems and Bocses}

\bigskip

We will define a notion of bi-co-module problems associated to matrix bi-module problems,
which is the transition into bocses. The concepts and the proofs are proposed by Y. Han.

\medskip

Since $\K_1$ and $\M_1$ are both quasi-free $R$-$R$-bi-modules, we have their
$R$-$R$-dual structures $\C_1$ and $\N_1$ with quasi-free $R$-$R$-quasi-basis $\V^*$ and
$\A^*$ respectively:
$$\begin{array}{c}\C_1=\Hom_{_{R^{\otimes 2}}}(\K_1,R^{\otimes 2}),\quad \V^*=\{v_1,v_2,\cdots,v_m\};\\
\N_1=\Hom_{_{R^{\otimes 2}}}(\M_1,R^{\otimes 2}),\quad \A^*=\{a_1,a_2,\cdots,a_n\}.\end{array}\eqno{(1.2\mbox{-}1)}$$
Write $v:X\mapsto Y$ (resp. $a:X\mapsto Y$) provided $1_XV1_Y=V$ (resp. $1_XA1_Y=A$).

\medskip

{\bf Definition 1.2.1}\, Let $\K$ be a $k$-algebra as in Definition 1.1.2.
We define a quasi-free $R$-module $\C_0=\sum_{X\in\T} R_Xe_{_X}\simeq R$ with an $R$-quasi-basis
$\{e_{_X}\}_{X\in \T}$; and a quasi-free $R$-$R$-bi-module $\C_1$ with an $R$-$R$-quasi-basis
$\V^\ast$ defined by the first formula of (1.2-1). Write $\C=\C_0\oplus\C_1$,
define a co-algebra structure $\varepsilon: \C\rightarrow R,\, e_{_X}\mapsto 1_{_X},v_j\mapsto 0$
and $\mu:\C\mapsto\C\otimes_{R}\C$ dual to $(\tm_{00},\tm_{01},\tm_{10},\tm_{11})$:
$$\begin{array}{c}\mu=\left(\begin{array}{c}\mu_{00}\\
\mu_{10}+\mu_{01}+\mu_{11}\end{array}\right):\left(\begin{array}{c}\C_0\\ C_1\end{array}\right)
\rightarrow \left(\begin{array}{c}\C_0\otimes_{_R}\C_0,\\
\C_1\otimes_{_R}\C_0\oplus\C_0\otimes_{_R}\C_1\oplus\C_1\otimes_{_R}\C_1\end{array}\right)\\
\mu_{00}(e_{_X})=e_{_X}\otimes_{_R}e_{_X};\,\,\mu_{10}(v_l)=v_l\otimes_{_R}e_{t(v_l)};
\,\,\mu_{01}(v_l)=e_{s(v_l)}\otimes_{_R}v_l;\\
\mu_{11}(v_l)=\sum_{i,j}\gamma_{ijl}\otimes_{_{R^{\otimes 3}}}(v_i\otimes_{_R} v_j),\,\,\mbox{if}\,\,
V_iV_j=\sum_l\gamma_{ijl}\otimes_{_{R^{\otimes 2}}}V_l.
\quad\end{array}$$

{\bf Definition 1.2.2}\, Let $\M$ be a $\K$-$\K$-bi-module of Definition 1.1.2, we define
a quasi free $R$-$R$-bi-module $\N_1$ with a quasi-basis $\A^\ast$ given
by the second formula of (1.2-1). Write $\N=\N_1$, suppose  $V_iA_j=\sum_l\eta_{ijl}\otimes_{_{R^{\otimes 2}}}V_l$;
$A_iV_j=\sum_l\sigma_{ijl}\otimes_{_{R^{\otimes 2}}}V_l$,
then $\N$ has a $\C$-$\C$-bi-co-module structure with the left and right co-module actions
dual to $(\tl_{01},\tl_{11})$ and $(\tr_{10},\tr_{01})$ respectively:
$$\begin{array}{c}\iota=(\iota_0+\iota_1): \N\rightarrow \C\otimes_{_R}\N=\C_0\otimes_{_R}\N\oplus\C_1\otimes_{_R}\N,\\
\iota_0(a_l)=e_{s(a_l)}\otimes_{_R}a_l,\,\iota_1(a_l)=\sum_{i,j}\eta_{ijl}
\otimes_{_{R^{\otimes 3}}}(v_i\otimes_{_R} a_j);\\
\tau=(\tau_0+\tau_1):\N\rightarrow\N\otimes_{_R}\C=\N\otimes_{_R}\C_0\oplus\N\otimes_{_R}\C_1,\\
\tau_0(a_l)=a_l\otimes_{_R}e_{t(a_l)},\,
\tau_1(a_l)=\sum_{i,j}\sigma_{ijl}\otimes_{_{R^{\otimes 3}}}(a_i\otimes_{_R} v_j).\end{array}$$

{\bf Definition 1.2.3}\, Let $H$ be a matrix over $R$ in Definition 1.1.4,
and $\C,\N$ be given in 1.2.1-1.2.2. The map
$\partial=(\partial_0,\partial_1): \N\rightarrow\C=\C_0\oplus\C_1$ with
$\partial_0=0$, and $\partial_1(a_l)=
\sum_i\zeta_{il}\otimes_{_{R^{\otimes 2}}}v_i$ if $\td(v_i)=\sum_l\zeta_{il}\otimes_{_{R^{\otimes 2}}}A_l$
is a co-derivation dual to $(\td_0,\td_1)$, such that
$\mu\partial=(\id\otimes\partial)\iota+(\partial\otimes\id)\tau$.

\medskip

{\bf Definition 1.2.4}\, Let $\mf A=(R,\K,\M,H)$ be a matrix bi-module problem.
A quadruple $\mf{C}=(R,\C,\N,\partial)$ is
said to be a {\it bi-co-module problem associated to $\mf A$} provided

(i) $R$ is a minimal algebra with a vertex set $\T$;

(ii)\, $\C$ is a co-algebra given in Definition 1.2.1;

(iii)\, $\N$ is a $\C$-$\C$-bi-co-module
given by Definition 1.2.2;

(iv)\, $\partial: \N\rightarrow\C$ is a co-derivation given by Definition 1.2.3.

\medskip

Now we construct the {\it bocs associated to a matrix bi-module problem $\mf A$}
based on the bi-co-module problem $\mf C$ associated to $\mf A$ presented by Roiter in \cite{Ro}.

Write $\N^{\otimes p}=\N\otimes_{_R} \cdots \otimes_{_R}\N$ with $p$ copies
of $\N$ and $\N^{\otimes 0}=R$. Define a tensor algebra $\Gamma$ of $\N$ over $R$,
whose multiplication is given by the natural isomorphisms:
$$\Gamma=\oplus_{p=0}^{\infty}
\N^{\otimes p}; \quad \N^{\otimes p}\otimes_{_R}\N^{\otimes
q}\simeq\N^{\otimes (p+q)}
$$
Let $\Xi= \Gamma \otimes_{R} \C
\otimes_{R} \Gamma$ be a $\Gamma$-$\Gamma$-bi-module of co-algebra
structure induced by $R\hookrightarrow\Gamma$, and denoted by $(\Xi,
{\mu}_{\Xi}, {\varepsilon}_{\Xi})$.
Define three $R$-$R$-bi-module maps:
$$\begin{array}{c}\kappa_1 : \N \stackrel{{\iota}}{\rightarrow}
\C \otimes_{_R} \N\stackrel{\cong}{\rightarrow} R
\otimes_{_R} \C\otimes_{_R} \N \hookrightarrow \Gamma
\otimes_{_R} \C \otimes_{_R} \Gamma,\\
\kappa_2 : \N \stackrel{{\tau}}{\rightarrow}
\N \otimes_{_R} \C \stackrel{\cong}{\rightarrow} \N
\otimes_{_R} \C \otimes_{_R} R \hookrightarrow \Gamma
\otimes_{_R} \C \otimes_{_R} \Gamma,\\
\kappa_3 : \N
\stackrel{{\partial}}{\rightarrow} \C
\stackrel{\cong}{\rightarrow} R \otimes_{_R} \C \otimes_{_R}
R \hookrightarrow \Gamma \otimes_{_R} \C \otimes_R
\Gamma.\end{array}$$

\medskip

{\bf Lemma 1.2.5} [CB1]\,  $\Im (\kappa_1 - \kappa_2 +
\kappa_3)$ is a $\Gamma$-co-ideal in $\Xi$. Thus $\Omega := \Xi / \Im
(\kappa_1 - \kappa_2 + \kappa_3)$ is a $\Gamma$-$\Gamma$-bi-module
of co-algebra structure.

\smallskip

{\bf Proof}\, Recall the law of bi-co-module: $({\mu} \otimes\id)
{\iota} = (\id\otimes {\iota}) {\iota}, (\id\otimes
{\mu}) {\tau} = ({\tau} \otimes\id){\tau} , (\id
\otimes {\tau}){\iota} = ({\iota} \otimes\id){\tau}
, (\id\otimes {\partial}){\iota} -{\mu}{\partial} +
({\partial} \otimes\id){\tau} = 0$. Thus, for any $b \in
\N$, we have
$$\begin{array}{l} \mu_{\Xi}(\kappa_1 - \kappa_2 + \kappa_3)(b)\end{array}$$
$$\begin{array}{ll}  = & \mu_{\Xi}(1_\Gamma\otimes {\iota}(b)) - \mu_{\Xi}({\tau}(b)
\otimes 1_\Gamma) + \mu_{\Xi}(\id\otimes {\partial}(b) \otimes1_\Gamma)
\\ = & ({\mu} \otimes\id){\iota}(b) - (\id\otimes
{\mu}){\tau}(b) + {\mu}({\partial}(b)) \\
= & (\id\otimes {\iota}){\iota}(b) - ({\tau} \otimes\id
){\tau}(b) + (\id\otimes{\partial}){\iota}(b)+
({\partial} \otimes\id)
{\tau}(b) \\
= & (\id\otimes {\iota} - \id\otimes{\tau} + \id\otimes
{\partial}){\iota}(b) + ({\iota} \otimes \id -{\tau}
\otimes\id + {\partial} \otimes\id) {\tau}(b)\\
= & u_{(1)} \otimes ({\iota} - {\tau} +
{\partial})(b_{(1)}) + ({\iota}- {\tau} +
{\partial})(b_{(2)}) \otimes u_{(2)} \\
= & u_{(1)} \otimes (\kappa_1 - \kappa_2 + \kappa_3)(b_{(1)}) +
(\kappa_1 - \kappa_2 + \kappa_3)(b_{(2)}) \otimes u_{(2)}
\\ \in & \Xi \otimes \Im(\kappa_1 - \kappa_2 + \kappa_3) + \Im(\kappa_1 -
\kappa_2 + \kappa_3) \otimes \Xi \end{array}$$
where ${\iota}(b) := u_{(1)} \otimes b_{(1)} , {\tau}(b) :=
b_{(2)} \otimes u_{(2)} $, and each term in each step is viewed as
an element in $\Xi \otimes_{\Gamma} \Xi$ naturally. The proof is
completed.

\medskip

Recall from \cite[3.4 Definition]{CB1}, that $\mf{B}=(\Gamma,\Omega)$ defined as above is a
bocs with a layer
$$
L=(R; \omega; a_1,a_2,\ldots,a_n; v_1,v_2,\ldots,v_m).
$$
Denote by $\varepsilon_{_\Omega}$ and $\mu_{_\Omega}$ the induced co-unit and co-multiplication,
then $\bar\Omega=$ker$\varepsilon_{_\Omega}$ is a $\Gamma$-$\Gamma$-bi-module
freely generated by $v_1, v_2,\ldots,v_m$, and $\Omega=\Gamma\oplus\bar{\Omega}$ as bi-modules.

From this, we use the imbedding: $\C_0\oplus \C_1\oplus\N\otimes\C_1\oplus\C_1\otimes\N\hookrightarrow
\Gamma\otimes_R(\C_0\oplus \C_1\oplus\N\otimes_R\C_1\oplus\C_1\otimes_R\N)\otimes_R\Gamma\subset\Omega$;
and the isomorphism: $\bar\Omega\otimes_R\bar\Omega\simeq\bar\Omega\otimes_\Gamma\bar\Omega$.
The group-like $\omega: R\rightarrow\Omega, 1_X\mapsto e_{_X}$ is an $R$-$R$-bi-module map.
Since $\dz_1(a_i)=\tau_0(a_i)-\iota_0(a_i)$, and
$(\iota_0(a_i)+{\iota_1}(a_i))-(\tau_0(a_i)+{\tau_1}(a_i))+
{\partial_1}(a_i)= ({\iota}-{\tau}+{\partial})(a_i)=
(\kappa_1-\kappa_2+\kappa_3)(a_i)=0$ in $\Omega$, the pair of the differentials determined by $\omega$ is given by
$$\begin{array}{c}\dz_1: \Gamma\rightarrow \bar{\Omega}, 1_X\mapsto 0,
a_i\mapsto{\iota_1}(a_i)-{\tau_1}(a_i)+
{\partial_1}(a_i),\,\,\,\forall\,\,a_i\in \A^*(X,Y);\\
\dz_2:\bar\Omega\mapsto\bar\Omega\otimes_\Gamma\bar\Omega, v_j\mapsto\mu_{11}(v_j),
\,\,\,\forall\,\,v_j\in \V^*(X,Y).\end{array} \eqno{(1.2\mbox{-}2)}$$

Recall from \cite{CB1} that a representation of a layered bocs $\mf{B}$ is a
left $\Gamma$-module $P$ of finite dimension $\underline d$, which means that
$P$ consists of a set of vector spaces and a set of linear maps:
$$\begin{array}{c}\{P_X=k^{d_X},P(x):P_X\mapsto P_X\mid X\in\T\};\\
\{P(a_i):k^{m_{_{X_i}}}\rightarrow k^{m_{_{Y_i}}}\mid\, a_i:X_i\mapsto Y_i,\,\,
1\leqslant i\leqslant n\}.\end{array}\eqno{(1.2\mbox{-}3)} $$
A morphism from $P$ to $Q$ is given by a
$\Gamma$-map $f:\Omega\otimes_{\Gamma}P\mapsto Q$, which is equivalent to write
$f=(f_X;f(v_j))_{X\in\T;j=1,\cdots,m}$, such that
$$\begin{array}{c}P(a_l)f_{_{Y_l}}-f_{_{X_l}}Q(a_l)=
\sum_{i,j}\eta_{ijl}
\otimes_{_{R^{\otimes 3}}}(f(v_i)\otimes_{_R}Q(a_j))\\
-\sum_{i,j}\sigma_{ijl}\otimes_{_{R^{\otimes 3}}}(P(a_i)\otimes_{_R} f(v_j))
+\sum_i\zeta_{il}\otimes_{_{R^{\otimes 2}}}f(v_i)\end{array}\eqno{(1.2\mbox{-}4)}$$
for all $a_l\in{\cal A}^\ast, 1\leqslant l\leqslant n$, by substituting $P(x)$ for $x$, see \cite{BK}. It is clear that
$$\begin{array}{c}
\Hom_{\Gamma}(\bar{\Omega}\otimes_{\Gamma}P, Q)\simeq
\oplus_{j=1}^m
\Hom_{\Gamma}(\Gamma1_{s(v_j)}\otimes_k1_{t(v_j)}P, Q)\simeq
\oplus_{j=1}^m \Hom_k(1_{t(v_j)}P, 1_{s(v_j)}Q).
\end{array}\eqno{(1.2\mbox{-}5)}$$

Denote by {\it $R(\mf{B})$ the representation category of $\mf{B}$}.

\medskip
\bigskip

\noindent{\bf \S 1.3  The representation category of a matrix bi-module problem}

\bigskip

This sub-section is devoted to defining the representation category
of a matrix bi-module problem. Which is relatively complicated, but
extremely useful for the proof of the main theorem.

\medskip

{\bf  Definition 1.3.1} Let $ J(\lambda)=J_d(\lambda)^{e_d}\oplus
J_{d-1}(\lambda)^{e_{d-1}}\oplus \cdots\oplus J_1(\lambda)^{e_1}$,
with $e_i$ non-negative integers, be a Jordan matrix. Set
$$m_j=e_d+e_{d-1}+\cdots+e_j.$$ The
following partitioned matrix $W(\lambda)$ similar to $J(\lambda)$ is called
a {\it Weyr matrix of eigenvalue $\lambda$}:
$$
W(\lz)=\left(
\begin{array}{cccccc}
\lz I_{m_1} & W_{12} &0&  \cdots & 0&0 \\
  & \lz I_{m_2} & W_{23}&\cdots &0&0 \\
  && \lz I_{m_3} &\cdots & 0 &0\\
&&&\ddots&\vdots&\vdots\\
  & & & & \lz I_{m_{d-1}}& W_{d-1,d} \\
 & & & && \lz I_{m_d}
\end{array}
\right)_{d\times d}, \ \ W_{j,j+1}=\left(
\begin{array}{c}
I_{m_{j+1}}\\ 0
\end{array}
\right)_{m_j\times m_{j+1}}.
$$ A direct sum $W=W(\lambda_1)\oplus W(\lambda_2)\oplus
\cdots \oplus W(\lambda_s)$ with
distinct eigenvalues $\lambda_i$ is said to be a {\it Weyr matrix}.
We may define {\it an  order $<$ on the base field} $k$, so that
each Weyr matrix has a unique form.
Similarly, let $Z_{ij}$ be a set of vertices, $S=k1_{Z_{ij}}\oplus k[z,\phi(z)^{-1}]$ be a
minimal algebra,
$\bar W\simeq \oplus J_{ij}(\lambda_{i})^{e_{ij}}1_{Z_{ij}}\oplus (z)^{\delta}$
with $\delta=0$ or $1$ is said to be a{\it Weyr matrix over $S$}.
If $S=\oplus_{j} k1_{Z_j}$ trivial, then $\bar W=\oplus_j I_{n_j}1_{Z_j}$ is still called a Weyr matrix over $S$.

Form now on, we assume that $\mf A=(R,\K,\M,H)$ is a matrix bi-module problem
with $T=\{1,2,\cdots,t\}$ and its partition $\T$. Let $X\in \T_1$. A Weyr matrix $W$ over $k$ is called {\it
$R_X$-regular} if its eigenvalues are $R_X$-regular, i.e.
$\phi_X(\lambda)\ne 0,$ for all the eigenvalues $\lambda$.
If $X\in\T_0$, an identity matrix $I$ is called $R_X$-regular.

A vector of non-negative integers is called a {\it size vector} over $T$,
if $\m=(m_1,m_2,\ldots,m_t)$ with $m_i=m_j,\forall\,i\sim j.$
$m=\sum_{i=1}^tm_i$ is said to be the {\it size}
of $\m$. Given a size vector $\m$, the vector
$\underline{d}=(m_{_X})_{_{X\in \T}}$ with $m_{_X}=m_i,
\forall\,i\in X$ is called a dimension vector {\it over $\T$} determined by
$\m$.

\medskip

{\bf Definition 1.3.2}\, Let $\mf A=(R,\K,\M,H)$ be a matrix bi-module problem,
let $S$ be a minimal algebra, $\Sigma=\oplus_{p=1}^\infty S^{\otimes p}$, see Formula (1.1-3).

(i)\, Write $H_X=(h_{ij}(x)1_X)_\tt$ with $h_{ij}(x)\in k[x]$ for $X\in\T_1$, and $x=1$ for $X\in\T_0$.
Let $\bar W_X$ be a Weyr
matrix of size $m_{_X}$ over $S$. We define an
$\m\times \m$-partitioned matrix:
$$H_X(\bar W_X)=(B_{ij})_{\tt},\quad B_{ij}=\left\{\begin{array}{cc}
h_{ij}(\bar W_X)_{m_i\times m_j},&i,j\in X,
\\(0)_{m_i\times m_j},&i\notin X\,\, \mbox{or}\,\, j\notin X,\end{array}\right.$$

(ii)\, Let $\m=(m_1,\cdots,m_t)$ and
$\n=(n_1,\cdots, n_t)$ be two size vectors over $\T$, and let
$F\in\IM_{m_{_X}\times n_{_X}}(S^{\otimes p}),p=1,2$,
with an $R_X$-$R_X$-module structure. The star product
$\ast$ of $F_X$ and $E_X$ is defined to be a diagonal $\m\times \n$-partitioned matrix:
$$\begin{array}{c}F_{_X}\ast E_X=(B_{ij})_\tt,\quad B_{ii}=
\left\{\begin{array}{cc}F_{_X},& i\in X,\\
0& i\notin X,\end{array}\right.\end{array}$$

(iii)\, Let $U=V_j$ or $A_i$ for some $j=1,\cdots,m, i=1,\cdots,n$, and $U=(\alpha_{ij})\in E_X$IM$_{t}(R\otimes_k R)E_Y$;
let $C\in\IM_{m_{_X}\times n_{_Y}}(S^{\otimes p})$ with $p=2$, and possibly $p=1$ for $U=A_1$;
$\{\bar W_X\in\IM_{m_{_X}}(S)\mid X\in\T\}$,
$\{W'_Y\in\IM_{n_{_Y}}(S)\mid Y\in\T\}$ be two sets of regular Weyr matrices. Suppose
there is an $R_X$-$R_Y$-bi-module structure on $\IM_{m_{_X}\times n_{_Y}}(S^{\otimes p})$:
$$\begin{array}{c}
(x\otimes_ky)\otimes_{_{R^{\otimes 2}}}C=W_XC W'_Y.\end{array}\eqno(1.3\mbox{-}1)$$
The star product
$\ast$ of $C$ and $U$ is defined to be a $(\m\times \n)$-partitioned matrix:
$$\begin{array}{c}C\ast U=(B_{ij})_\tt,\quad\quad B_{ij}=\left\{\begin{array}{cl}
C\otimes_{_{R^{\otimes 2}}}\alpha_{ij}, &i\in X,j\in Y;\\ 0_{m_{i}\times n_j},
&i\notin X,\,\mbox{or}\, j\notin Y.\end{array}\right.\end{array}$$

\smallskip

{\bf Lemma 1.3.3}\, Let $\mf A=(R,\K,\M,H)$ be a matrix bi-module problem.

(i)\, If $C\in\IM_{m_{_X}\times n_{_Y}}(S^{\otimes 2})$,
$\zeta_{il}\in R_X\otimes_k R_Y$ are given in Definition 1.2.3, then by the
usual product of $\Sigma$-matrices:
$$\begin{array}{cc}(C\ast V_i)H_Y(\bar W_Y)-H_X(\bar W_X)(C\ast V_i)=
\sum_{l}(C\otimes_{R^{\otimes 2}}\zeta_{il})\ast A_l.\end{array}$$

(ii)\, If $F_X\in\IM_{m_{_X}\times n_{_X}}(S^{\otimes p}),p=1,2$;
$C\in\IM_{m_{_X}\times n_{_Y}}(S^{\otimes q})$, where $q=2$, and possibly $q=1$
for $U=A_1$, then by the usual multiplication of matrices over $\Sigma$:
$$\begin{array}{c}(F_{_X}\ast E_X)(C\ast U)=\left\{\begin{array}{cc}(F_{_X}C)\ast U,&1_XU=U;\\
0,&1_XU=0.\end{array}\right.\end{array}$$
Similarly, $(C\ast U)(F_{_X}\ast E_X)=(CF_{_X})\ast U$ if $U1_X=U$, or $0$ if $U1_X=0$.

\smallskip

(iii)\, Let $U\in E_X\IM_t(R\otimes_kR)E_Y, V\in
E_Y\IM_t(R\otimes_kR)E_Z$ and $UV=\sum_{l=1}^n\epsilon_l\otimes_{_{R^{\otimes 2}}}G_l$ with
$\epsilon_l\in R_X\otimes_kR_Y\otimes_kR_Z$ and $G_l\in E_X\IM_t(R\otimes_kR)E_Z$.
Let $\m, \n, \l$ be size vectors over $\T$, the $R\otimes_kR$-module
structures given by Formula (1.3-1) yield an $R^{\otimes 3}$-module structure on
$\oplus_{(X,Y,Z)\, \in \T\times \T}\,\IM_{m_{_X}n_{_Y}}(S^{\otimes p})
\otimes_R\IM_{n_{_Y}\times l_{_Z}}(S^{\otimes q})$ for $p,q=2$, and possibly $p=1$ for
$U=A_1$, or $q=1$ for $V=A_1$. Then by the usual $\Sigma$-matrix product:
$$(C\ast U)(D\ast V)={\sum}_{l=1}^n
\big((C\otimes_RD)\otimes_{_{R^{\otimes 3}}}\epsilon_l\big)\ast G_l.$$

\smallskip

{\bf Proof}\, (i)\, Write $H_X=(\gamma_{pq}),\gamma\in R_X; H_Y=(\delta_{pq}),
\delta\in R_Y; V_i=(\alpha_{pq}),\alpha\in R_X\otimes_k R_Y$,
The left side$=(\sum_l(C\otimes_{R^\otimes 2}\alpha_{pl})\otimes_{R}\delta_{lq})
-(\sum_l\gamma_{pl}\otimes_R(C\otimes_{R^\otimes 2}\alpha_{lq}))=
(C\otimes_{R^{\otimes 2}}\sum_l(\alpha_{pl}\delta_{lq}-\gamma_{pl}\alpha_{lq}))=
C\ast(V_iH_Y-H_XV_i)=C\ast\td(V_i)=C\ast(\sum_{l}\zeta_{il}A_l)=$the right side.

(ii)\, Write $U=(\alpha_{pq}),\alpha\in R\otimes_kR$, the left
side$=(F_X1_X(C\otimes_{R^{\otimes 2}}\alpha_{pq}))
=((F_XC)\otimes_{R^{\otimes 2}}\alpha_{pq})=$the right side.

(iii)\, Write $U=(\alpha_{pq}),\alpha\in R_X\otimes_kR_Y,V=(\beta_{pq}),
\beta\in R_Y\otimes_kR_Z$. The left
side$=(\sum_l(C\otimes_{R^{\otimes 2}}\alpha_{pl})(D\otimes_{R^{\otimes 2}}\beta_{lq}))
=((CD)\otimes_{R^{\otimes 3}}(\sum_l\alpha_{pl}\otimes_R\beta_{lq}))
=$the right side. The proof is finished.

\medskip

{\bf Definition 1.3.4}\, Let $\mf{A}=(R,\K,\M,H)$ be a matrix
bi-module problem, and $\m$ a size
vector over $\T$. Then a representation $\bar P$ of $\mf{A}$ can be written as an
$\m\times\m$-partitioned matrix over $k$:
$$\begin{array}{c}\bar P=\sum_{X\in\T} H_X(W_X)+\sum_{i=1}^n
\bar P(a_i)\ast A_i,\end{array}$$
where $W_X\in\IM_{m_{_X}}(k)$ is regular for any $X\in\T$,
$\bar P(a_i)\in\IM_{m_{_{X_i}}\times m_{_{Y_i}}}(k)$ see Formula (1.2-3).
Taken $S=k=\Sigma$, the first summand is defined in 1.3.2 (i), and the second one in (iii) .

\medskip

{\bf Definition 1.3.5}\, Let $\m,\bar P$ be given above, $\n$
be a size vector over $\T$, and $\bar Q$ a representation over $\mf A$. A morphism $\bar f:\bar P\rightarrow \bar Q$
can be written as an $\m\times \n$-partitioned matrix by Definition 1.3.2 (ii) and (iii) for $S=k=\Sigma$:
$$\begin{array}{c}\bar f=\sum_{X\in\T}\bar f_{_X}\ast E_X+\sum_{j=1}^m\bar f(v_j)\ast V_j,\quad \end{array}$$
where $\bar f_{_X}\in\IM_{m_{_X}\times n_{_X}}(k),\bar f(v_j)\in\IM_{m_{s(v_j)}\times n_{t(v_j)}}(k)$,
such that $\bar P\bar f=\bar f\bar Q$, where the multiplication is given according to Lemma 1.3.3 (i)-(iii).

\medskip

If $\bar f':\bar Q\rightarrow\bar U$ is also a morphisms over $\mf A$. Then
$\bar f\bar f':\bar P\rightarrow\bar U$ calculated according to Lemma 1.3.3 (ii)-(iii) is a morphism.
In fact, $(\bar f\bar f')\bar P=\bar f(\bar Q\bar f')=(\bar U\bar f)\bar f'=\bar U(\bar f\bar f')$.
We denote by {\it $R(\mf{A})$ the category of representations of the
matrix bi-module problem $\mf{A}$}.

\medskip

{\bf Theorem 1.3.6}\, Let $\mf A$ be a matrix bi-module problem,
and $\mf B$ the associated bocs.
Then the categories $R(\mf A)$ and $R({\mf{B}})$ are equivalent.

\smallskip

{\bf Proof}\, Without loss of generality, we may assume that $\{P(x)=W_X\mid X\in\T\}$
is a set of regular Weyr matrices. Then $\bar P$ in Definition 1.3.4 and $P$
in Formula (1.2-3) are one-to-on correspondent; $\bar f$ of Definition 1.3.5 and
$f$ in Formula (1.2-5) are one-to-on correspondent.
Moreover, $\bar P\bar f=\bar f\bar Q$ if and only if $f$ satisfying Formula (1.2-4), for the proof
of this assertion, we refer to
Formula (1.4-2) and Theorem 1.4.2 below.

\medskip

Thanks to Theorem 1.3.6, we will denote by $P,f$ in both $R(\mf A)$ and $R(\mf B)$ in a unified manner.

\bigskip
\bigskip

\noindent{\bf 1.4. Formal Products and Formal Equations}

\bigskip

\medskip

Now we introduce a notion of ``formal equation", which will build a
nice connection between matrix bi-module problems and associated bocses.

\medskip

Let $\mf{A}=(R,\K,\M,H)$ be a matrix bi-module
problem, with an associated bi-co-module problem $\mf{C}=(R,\C,\N,\partial)$.
Recall that $\{E_{_X}\}$ and $\{e_{_X}\}$
are dual bases of $(\K_0,\C_0)$; $\{V_1,\cdots,V_m\}$ and $\{v_1,\cdots,v_m\}$
are those of $(\K_1,\C_1)$; and $\{A_1,\dots,A_n\}$ and $\{a_1,\cdots,a_n\}$
of $(\M_1,\N_1)$. Then, taken $S=R,\Sigma=\Delta,\m=(1,\cdots,1)=\n$, in Definition 1.3.2 (ii)-(iii):
$$\begin{array}{lll} \Upsilon &=&\sum_{X\in \T}e_{_X}\ast E_{_X}\\
\Pi &=&\sum_{j=1}^mv_j\ast V_j\\
\Theta &=&\sum_{i=1}^na_i\ast A_i
\end{array}\eqno{(1.4\mbox{-}1)}$$
are called the {\it formal products} of $(\K_0,\C_0), (\K_1,\C_1)$
and $(\M_1,\N_1)$ respectively.

\medskip

{\bf Lemma 1.4.1}\, With the notations above, and
$\delta$ the differential in the associated bocs $\mf B$. Then the matrix multiplication yields:
$$\begin{array}{c}\big(\sum_{i=1}^mv_i\ast V_i\big)
\big(\sum_{j=1}^mv_j\ast V_j\big)=\sum_{l=1}^m
\mu_{11}(v_l)\ast V_l;\\
\big(\sum_{i=1}^na_i\ast A_i\big) \big(\sum_{j=1}^mv_j\ast
V_j\big)=\sum_{l=1}^n
\tau_1(a_l)\ast A_l;\\
\big(\sum_{j=1}^mv_j\ast V_j\big) \big(\sum_{i=1}^na_i\ast
A_i\big)=\sum_{l=1}^n
\iota_1(a_l) \ast A_l;\\
\big(\sum_{j=1}^mv_j\ast V_j\big) H-H\big(\sum_{j=1}^mv_j\ast
V_j\big)=\sum_{l=1}^n
\partial_1(a_l)\ast A_l;\\
\big(\sum_{l=1}^na_l\ast A_l\big) \big(\sum_{_{X\in \T}} e_{_X}\ast
E_{_X}\big)- \big(\sum_{_{X\in  \T}} e_{_X}\ast E_{_X}\big)
\big(\sum_{l=1}^na_i\ast A_l\big)=\sum_{l=1}^n\delta(a_l)\ast A_l.
\end{array}$$

\smallskip

{\bf Proof}\, \ding{172} We first prove the second equality, the proofs of the first and the third
are similar. By Lemma 1.3.3 (iii) for $S=R,p=q=2$, the left side
$=\sum_{l=1}^n \big(\sum_{i,j}\sigma_{ijl}\otimes_{_{R^{\otimes 3}}}
(a_i\otimes_Rv_i)\big)\ast A_l =$ the right side.
\ding{173} For the fourth equality, by Lemma 1.3.3 (i) the left side
$=\sum_{l=1}^n\big(\sum_{j=1}^m\zeta_{lj}\otimes_{_{R^{\otimes 2}}}v_j\big)\ast A_l=$
the right side. \ding{174} For the last one, by Lemma 1.3.3 (ii), $p=1,q=2$, the left side
$=\sum_{l=1}^n (a_l\otimes_{_R}e_{Y_l}-e_{X_l}\otimes_{_R} a_l)\ast A_l=$ the right side.
The proof of the lemma is completed.

\medskip

Denote by {\it $(\mf{A},\mf{B})$ the pair of a matrix bi-module problem and its
associated bocs}. Then the matrix equation
$(\Theta+H)(\Upsilon+\Pi)=(\Upsilon+\Pi)(\Theta+H)$, more precisely,
$$\begin{array}{cc}
&\big(\sum_{i=1}^n a_i\ast A_i+H\big) \big(\sum_{X\in \T}
e_X\ast E_X+\sum_{j=1}^mv_j\ast V_j\big)\\
=&\big(\sum_{X\in \T} e_X\ast E_X+\sum_{j=1}^mv_j\ast V_j\big)
\big(\sum_{i=1}^n a_i\ast
A_i+H\big)\end{array}\eqno{(1.4\mbox{-}2)}$$ is called the {\it
formal equation} of $(\mf A, \mf B)$ due to the following theorem.

\medskip

{\bf Theorem 1.4.2}\, The entry at the
leading position of $A_l$ in the formal equation is
$$\delta(a_l)=\iota_1(a_l)-\tau_1(a_l)+\partial_1(a_l).$$

\smallskip

{\bf Proof.} According to Formula (1.4-2) and Lemma 1.4.1:
$$\begin{array}{ll}&\sum_{l=1}^n\delta(a_l)\ast A_l=
\sum_{l=1}^n\big(a_le_{t(A_l)}-e_{s(A_l)}a_l\big)\ast A_l\\
=&\sum_{j,i}(v_j\ast V_j)(a_i\ast A_i)- \sum_{i,j}(a_i\ast
A_i)(v_j\ast V_j)+ \sum_{j}\big((v_j\ast V_j) H-
H(v_j\ast V_j)\big)\\
=&\sum_{l=1}^n\iota_1(a_l)\ast A_l -\sum_{l=1}^n\tau_1(a_l)
\ast A_l + \sum_{l=1}^n\partial_1(a_l)\ast A_l\\
=&\sum_{l=1}^n\big(\iota_1(a_l)-\tau_1(a_l)
+\partial_1(a_l)\big)\ast A_l.\end{array}$$
We obtain the expression at the leading position of
$A_l$ for $1\leqslant l\leqslant n$. The proof is finished.

\medskip

Moreover, the first formula of Lemma 1.4.1 gives:
$$\begin{array}{c}\big(\sum_{_{X\in \T}} e_{_X}\ast E_{_X}
+\sum_{i=1}^mv_i\ast V_i\big)\big(\sum_{_{X\in \T}}
e_{_X}\ast E_{_X}+\sum_{j=1}^mv_j\ast V_j\big)\\
=\sum_{_{X\in\T}} (e_{_X}\otimes_{_R}e_{_X})\ast E_{_X}+\sum_{l=1}^m
\mu(v_l)\ast V_l.\end{array}\eqno{(1.4\mbox{-}3)}$$

\smallskip

Now we define a special class of matrix
bi-module problems to end the sub-section.

\medskip

{\bf Definition 1.4.3}\, Let $\mf{A}=(R,\K,\M,H=0)$ be a matrix
bi-module problem with $R$ trivial. $\mf{A}$ is said to be {\it
bipartite} if $\T=\T'\dot{\cup}\T''$, $R=R'\times R''$, $\K=\K'\times \K''$ as direct
products of algebras, and $\M$ is a $\K'$-$\K''$-bi-module.

\medskip

Let $\Lambda$ be a finite-dimensional basic $k$-algebra,
$J=\rad(\Lambda)$ be the Jacobson radical of $\Lambda$ with the
nilpotent index $m$, and $S=\Lambda/J$. Suppose $\{e_1,\cdots,e_h\}$ is a complete set
of orthogonal primitive idempotents of $\Lambda$. Taken the pre-images of
$k$-bases of $e_i(J^i/J^{i+1})e_j$ under the canonical
projections $J^i\rightarrow J^i/J^{i+1}$ in turn for $i=m,
\cdots,1$, we obtain an ordered basis of $J$ under the length order,
see \cite[6.1]{CB1} for details. Then we construct the left regular
representation $\bar\Lambda$ of $\Lambda$ under the $k$-basis
$(a_{n},\cdots,a_2,a_1,e_1,\cdots,e_h)$ of $\Lambda$,
which yields a bipartite matrix bi-module problem $\mf A=(R,\K,\M,H=0)$ with
$$R=S\times S;\quad\K_0\oplus\K_1=\bar\Lambda\times\bar\Lambda;
\quad\M_1=\mbox{rad}(\bar\Lambda);\quad H=0.$$

{\bf Remark 1.4.4}\, A simple calculation shows that the {\it row indices} of the leading positions of
the base matrices in $\A$ are {\it pairwise different}, and the {\it column index} of the
leading position of $A\in\A_{XZ}$ equals $j_{_Z}=$max$\{j\in Z\}$ for any $X\in\T$, they are {\it concentrated} ,
the $j_{_Z}$-th column is said to be a {\it main column} over $Z$. Such a fact
is denoted by RDCC for short, which is not essential in the proof of the main theorem,
but makes it easier and more intuitive.

\medskip

{\bf Example 1.4.5}\,\cite{D1,R1} Let $Q= $ {\unitlength=1mm
\begin{picture}(20,4) \put(10,1){\circle*{0.6}}
\put(7,1){\circle{4}} \put(13,1){\circle{4}}
\put(9,2){\vector(0,1){0}} \put(11,2){\vector(0,1){0}}
\put(2,0){$a$} \put(16,0){$b$}
\end{picture}}
be a quiver, $I=\langle a^2,ba- ab, ab^2 , b^3\rangle$ be an ideal
of $k Q$, and $\Lambda=kQ/I$. Denote  the residue
classes of $e, a, b$ in $\Lambda$ still by $e, a,b $ respectively.
Moreover set $ c=b^2, d=ab$. Then an ordered $k$-basis $\{d, c, b,
a, e\}$ of $\Lambda$ yields a regular representation $\bar \Lambda$.
A matrix bi-module problem $\mf A$ follows by Theorem 1.4.4, with
its associated bocs $\mf B$.
The formal equation of the pair $(\mf A,\mf B)$ can be written as:
$$\small{
\left(
\begin{array}{ccccc}
e& 0 &  u_1    & u_2 & u_4 \\
  & e &  u_2    & 0 & u_3 \\
  &   &   e   & 0 & u_2 \\
  &   &       & e & u_1 \\
  &   &       &   & e
\end{array}
\right) \left(
\begin{array}{ccccc}
0 & 0 &  a    &  b & d \\
  & 0 &  b    & 0 & c \\
  &   &   0   & 0 & b \\
  &   &       & 0 & a \\
  &   &       &   & 0
\end{array}
\right)= \left(
\begin{array}{ccccc}
0 & 0 &  a    &  b & d \\
  & 0 &  b    & 0 & c \\
  &   &   0   & 0 & b \\
  &   &       & 0 & a \\
  &   &       &   & 0
\end{array}
\right) \left(
\begin{array}{ccccc}
f & 0 &  v_1    & v_2 & v_4 \\
  & f &  v_2    & 0 & v_3 \\
  &   &   f   & 0 & v_2 \\
  &   &       & f & v_1 \\
  &   &       &   & f
\end{array}
\right)}
$$ where $e=e_{_{X}}, f=e_{_{Y}}$ for simplicity. Denote by $A,B,C,D$
the $R$-$R$-quasi-basis of $\M_1$, and by $a,b,c,d$ the $R$-$R$-dual basis of $\N_1$.
From this we can obtain the associated bocs $\mf{B}$ with the
layer $L=(R;\omega; a,b,c,d; u_1,u_2,u_3,u_4,$
$v_1,v_2,v_3,v_4)$.

$$
\unitlength=1mm
\begin{picture}(80, 26)
\put(10,0){\circle*{1.00}} \put(10,20){\circle*{1.00}}
\qbezier(7,19)(5,11)(7,3) \put(7,3){\vector(1,-1){1}}
\qbezier(13,19)(15,11)(13,3) \put(13,3){\vector(-1,-1){1}}

\qbezier[10](9,20)(4.5,19)(5,21) \qbezier[10](9,20)(5.5,22)(5,21)
\qbezier[10](9,21)(6,23)(7,24) \qbezier[10](9,21)(9,25)(7,24)
\qbezier[10](11,20)(15.5,19)(15,21)
\qbezier[10](11,20)(14.5,22)(15,21)
\qbezier[10](11,21)(14,23)(13,24) \qbezier[10](11,21)(11,25)(13,24)

\qbezier[10](9,0)(4,1)(5,-1) \qbezier[10](9,0)(5,-2)(5,-1)
\qbezier[10](9,-1)(6,-3)(7,-4) \qbezier[10](9,-1)(9,-5)(7,-4)
\qbezier[10](11,0)(16,1)(15,-1) \qbezier[10](11,0)(15,-2)(15,-1)
\qbezier[10](11,-1)(14,-3)(13,-4) \qbezier[10](11,-1)(11,-5)(13,-4)

\put(17,20){$X$} \put(17,0){$Y$}

\put(9.00,19.00){\vector(0,-1){18}}
\put(11.00,19.00){\vector(0,-1){18}}

\put(4,10){\mbox{$a$}} \put(8,10){\mbox{$b$}}
\put(10.50,10){\mbox{$c$}} \put(14,10){\mbox{$d$}}

\put(35,8){\mbox{$ \left\{\begin{array}{l} \dz(a)=0,\\ \dz(b)=0, \\
\dz(c)=u_2b-bv_2, \\ \dz(d)=u_1b+u_2a-bv_1-av_2.
\end{array}\right.$}}
\end{picture}
$$

\bigskip
\bigskip
\bigskip
\centerline{\bf  2 Reductions for Matrix Bi-module Problems}
\bigskip

In the present section, we will define six reductions in terms of matrix bi-module problems
associated to those of bocses, and will give two additional ones.
Then we discuss defining systems of pairs, in order to construct the induced pairs in
series of reductions.

\bigskip
\bigskip

\noindent{\bf \S 2.1  Triangular properties and admissible bi-modules}

\bigskip

Let $\mf{A}=(R,\K,\M,\td)$ be a matrix bi-module problem. Then
for any $A_i\in \M_1\subseteq \IM_t(R\otimes_kR)$, the leading position of $V_j A_i$
(resp. $A_i V_j$) is strictly larger than that of $A_i$, since $\K_1\subseteq\mathbb N(R\otimes_kR)$. But
$\A=\{A_1,\cdots,A_n\}$ is an ordered set of normalized basis, the left and right module action
satisfies Formula (2.1-1) below, so called {\it triangular} property:
$$\tl(\K_1\times A_i),\, \tr(A_i\times\K_1),\subseteq
\oplus_{l=i+1}^n R^{\otimes 3}\otimes_{_{R^{\otimes 2}}}A_l.\eqno{(2.1\mbox{-}1)}$$
Since $(\K_1,\C_1)$ and $(\M_1,\N_1)$ are dual $R$-$R$-bi-modules in the associated bi-co-module problem
$\mf{C}=(R,\C,\N,\partial)$ of $\mf A$, the left and right co-module action also possesses the
{\it triangular} property:
$$\begin{array}{c}\iota_1(a_l)\in\, \oplus_{i=1}^{l-1}
\C_1\otimes_{_R}a_i,\quad \tau_1(a_l)\in\, \oplus_{i=1}^{l-1}a_i
\otimes_{_R} \C_1.\end{array}\eqno{(2.1\mbox{-}2)}$$
Define a $\K$-$\K$ sub-bi-module
and the corresponding $\K$-$\K$-quotient-bi-module:
$$\begin{array}{c}\M^{(h)}=\oplus_{i=h+1}^n\bar\Delta\otimes_{_{R^{\otimes 2}}} A_i\subseteq\M,
\quad \M^{[h]}=\M/\M^{(h)}.\end{array}$$
$\mf{A}^{[h]}=(R,\K,\M^{[h]},\bar\td)$, with $\bar\td$ induced from $\td$,
is said to be a {\it quotient} of $\mf A$, which itself might be no longer a
matrix bi-module problem. If $\Gamma^{(h)}$ is freely generated by $a_1,\cdots,a_{h}$, the associated bocs
$\mf{B}=(\Gamma,\Omega)$ has a sub-bocs:
$$\begin{array}{c}\mf B^{(h)}=(\Gamma^{(h)},\Gamma^{(h)}\otimes_R\Omega\otimes_R
\Gamma^{(h)}).\end{array}
$$

Note a simple fact: let $(\mf{A},\mf{C},\mf B)$ be a triple defined above, then
$$\begin{array}{ll}&\tl(\K_1\times\M_1),\tr(\M_1\times\K_1),\td(\K_1)\subseteq \M^{(h)}_1\,\, \mbox{in}\,\, \mf A\\
\Longleftrightarrow & \C_1\otimes_{_R}\N^{(h)}_1=0,\N^{(h)}_1\otimes_{_R}\C_1=0,\partial(\N^{(h)}_1)=0\,\, \mbox{in}\,\, \mf C\\
\Longleftrightarrow & \delta(\Gamma^{(h)})=0\,\,\mbox{in}\,\,\mf B.\end{array}$$
In fact, the condition in $\mf A$ is equivalent to $\eta_{jil}=0,\sigma_{ijl}=0, \zeta_{jl}=0$ for
$l=1,\cdots,h$ and any $i,j$, which is equivalent to the condition on $\mf C$ and $\mf B$.

\medskip

Let $R_X=k[x,\phi(x)^{-1}]$, $r$ a fixed positive integer, and
$\lambda_1, \ldots, \lambda_s\in k$ with
$\phi(\lambda_i)\ne 0$, write
$g(x)=(x-\lambda_1)\cdots (x-\lambda_s)$. Define a minimal algebra $S$
and a $R_X$-module $K$ over $S$:
$$\begin{array}{c}
S=\big(\prod_{i=1}^s\prod_{j=1}^r
k1_{Z_{ij}}\big)\times k[z,\phi(z)^{-1}g(z)^{-1}],\\
K=\big(\oplus_{i=1}^s\oplus_{j=1}^r\oplus_{q=1}^j
k1_{Z_{ijq}}\big)\oplus k[z,\phi(z)^{-1}g(z)^{-1}],\\
K(x)=\bar W: K\rightarrow K,\,\, \bar W\simeq \oplus_{i=1}^s
\oplus_{j=1}^r J_j(\lambda_j)1_{Z_{ij}}\oplus (z),\end{array}\eqno{(2.1{\mbox{-}}3)}$$
a Weyr matrix over $S$. Set $n_{_X}=\frac{1}{2}sr(r+1)+1$,
denote by $\{(i,j,l)\mid 1\leqslant l\leqslant j, 1\leqslant j\leqslant r,1\leqslant i\leqslant s\}\cup\{n_{_X}\}$,
the index set of the direct summands of $K$. The order on the set is defined by
$$\begin{array}{c}(i,j, l)\prec (i',j', l')\Longleftrightarrow i<i'; \,\mbox
{or}\,\, i=i', l<l';\,\mbox {or}\,\, i=i',l=l',
j>j'.\end{array}$$
there is a partition on the index set,
such that $Z_{ij}=\{(i,j,l)\mid l=1,\cdots,j\}$,
$Z=\{n_{_X}\}$. Suppose ${\textsf e}_{_{(i,j,l)}}$ and $\textsf f_{_{(i,j,l)}}$ are
$1\times n_{_X}$ and $n_{_X}\times 1$ matrices respectively,
with $1_{_{Z_{ij}}}$ at the $(i,j,l)$-th component and $0$ at others.
Let $m_{_X}$
square matrix $\f\otimes_{_S}\textsf e$ (resp. $\f\otimes_{k}\textsf e)$ stand for the
usual matrix tensor product over $S$ (resp. $k$). Let End$_{R_X}(K)$ be the
endomorphism ring of $K$ over $\Sigma=\sum_{p=1}^\infty S^{\otimes p}$.
Then the $S$-quasi basis of index $0$,
and $S$-$S$-quasi-basis of index $1$ are given respectively by
$$\begin{array}{c}\{F_{ij}=\sum_{l=1}^{j}\f_{_{(i,j,l)}}\otimes_{_S}\textsf e_{_{(i,j,l)}};
F_{m_{_X}m_{_X}}=\f_{_{m_X}}\otimes_{_S}\textsf e_{_{m_X}}\mid
1\leqslant j\leqslant r, 1\leqslant i\leqslant s\}.\end{array}$$
$$\begin{array}{c}
F_{ijj'l}=\sum_{h=1}^{j'-l+1}\f_{_{(i,j,h)}}\otimes_k\textsf e_{_{(i,j',l+h-1)}},
\left\{\begin{array}{c}l=1,\cdots,j',\,\,\mbox{if}\,\,j>j';\\
l=2,\cdots,j',\,\,\mbox{if}\,\,j=j';\\
l=j,\cdots,j',\,\,\mbox{if}\,\,j<j'.\end{array}\right.\end{array}
$$

Define a path algebra $R_{XY}$, a minimal algebra $S$ and a $R_{XY}$-module $K$ over $S$:
$$\begin{array}{c}R_{XY}:X\stackrel{a_1}\rightarrow Y,\quad S=\prod_{i=1}^3
S_{Z_i},\quad S_{Z_i}=k1_{Z_i}, i=1,2,3;\\
K_{X}=k1_{Z_2}\oplus k1_{Z_1},\quad K_{Y}=k1_{Z_3}\oplus k1_{Z_2},
\quad K(a_1)=\left(\begin{array}{cc} 0& 1_{Z_2}\\ 0&
0\end{array}\right): K_{X}\rightarrow K_{Y}.\end{array}\eqno{(2.1\mbox{-}4)}$$
Let End$_{R_{XY}}(K)$ be over $\Sigma$, the $S$-quasi basis
and $S$-$S$-quasi-basis are given respectively by
$$\begin{array}{c}F_{Z_1}=\f_{_{Z_{1(X,2)}}}\otimes_{_S}\textsf e_{_{Z_{1(X,2)}}},\quad
F_{Z_3}=\f_{_{Z_{3(Y,1)}}}\otimes_{_S}\textsf e_{_{Z_{3(Y,1)}}},\\
F_{Z_2}=(\f_{_{Z_{2(X,1)}}}\otimes_{_S}\textsf e_{_{Z_{2(X,1)}}},\f_{_{Z_{2(Y,2)}}}\otimes_{_S}\textsf e_{_{Z_{2(Y,2)}}});\\
F_{Z_2Z_1}=\f_{_{Z_{2(X,1)}}}\otimes_{k}\textsf e_{_{Z_{1(X,2)}}},\quad F_{Z_3Z_2}=\f_{_{Z_{3(Y,1)}}}\otimes_{k}\textsf e_{_{Z_{2(Y,2)}}}.
\end{array}$$

Define $R_{XY},S,K$ below, then
End$_{R_{XY}}(K)$ over $\Sigma$ has $S$-quasi basis $\{F_Z\}$ of index $0$, and the part of index $1$ in
End$_{R_{XY}}(K)$ is $0$. Let $R_X=k[x,\phi_{_X}(x)^{-1}]$, $R_Y=k1_Y$ or $k[y,\phi_{_Y}(y)]$:
$$\begin{array}{c}R_{XY}:\xymatrix{X\ar@(ul,dl)[]_x\ar[r]^{a_1}&Y},\,\,\mbox{or}\,\,\xymatrix{X\ar@(ul,dl)[]_x\ar[r]^{a_1}&Y\ar@(ur,dr)[]^y};
\quad S=\xymatrix{X\ar@(ul,dl)[]_z};\\
K_{X}=S,\, K_{Y}=S, \,\, K(a_1)=(1_Z),\, K(x)=(z),\,\,\mbox{or add}\,K(y)=(z)\,\mbox{if}\,Y\in\T_1;\\
F_{Z}=(\f_{_{Z_{X}}}\otimes_{_S}\textsf e_{_{Z_{X}}}, \f_{_{Z_{Y}}}\otimes_{_S}\textsf e_{_{Z_{Y}}}).\end{array}\eqno{(2.1\mbox{-}5)}$$

{\bf  Definition 2.1.1}\, Let $(\mf A,\mf B)$ be a pair, let
$R'$ be a minimal algebra with algebra $\Delta'=\sum_{p=1}^\infty {R'}^{\otimes p}$
in Formula (1.1-3), and $\underline d=(n_{_X}\mid X\in\T)$ a
dimension vector over $\T$. An $R'$-$\bar R$-bi-module $L$ (or an $\bar R$-module over $R'$)
is said to be admissible, if $L$ satisfies (a1)-(13) below.

{\bf (a1)}\, There are three cases:

\ding{172} $d_X=1$, or $0$ for any $X\in\T$. $\bar R=R$ or $R[a_1]$ with
$\dz(a_1)=0$;

\ding{173} $R_X=k[x,\phi(x)^{-1}],\bar R=R$, and
$R'=S\times\prod_{Y\ne X} R_Y$, $L=K\oplus(\oplus_{Y\in\T,Y\ne X}R_Y)$ with $S,K$ defined in Formula
(2.1-3);

\ding{174} $\bar R=R_{XY}\times\prod_{U\in\T,U\ne X,Y} R_{U}$
with $X,Y\in\T_0$, $\dz(a_1)=0$, and $R'=S\times\prod_{U\ne X,Y} R_Y$,
$L=K\oplus(\oplus_{U\in\T,U\ne X,Y}R_U)$ with $S,K$ defined in Formula
(2.1-4).

{\bf (a2)}\, Denote unified by $L=\oplus_{X\in\T}L_X$, $L_X=\oplus_{p=1}^{n_{_X}}R'_{Z_{(X,p)}}, X\in\T$.
Let $\textsf e_{_{Z_{(X,p)}}}$ be a $(1\times n_{_X})$-matrix row with the $p$-th
entry $1_{{Z_{(X,p)}}}$ and others zero. Let $L^\ast=\Hom_{R'}(L,R')$ be an $\bar R$-$R'$-bi-module,
$\f_{_{Z_{(X,p)}}}=\textsf e_{_{Z_{(X,p)}}}^\ast=\textsf e_{_{Z_{(X,p)}}}^T$ be
an $(n_{_X}\times 1)$-matrix column. And $\textsf e(\f)=\textsf e\f$.

{\bf (a3)}\, $\bar E=$End$_R(L)\subseteq\Pi_{X\in\T}\mathbb T_{n_{_X}}(\Delta')$, where $\bar E_0\simeq R'$;
$\bar E_1$ is a quasi-free $R'$-$R'$-bi-module; and
$\bar E$ is finitely generated in index $(0,1)$.
Forgotten the $\bar R$-$\bar R$-structure on $L^\ast\otimes L$, we assume
$$\begin{array}{c}\bar E_0\subseteq\Pi_{X\in\T}\mathbb D_{n_{_X}}(R')
\subseteq\Pi_{X\in\T}(L^\ast\otimes_{R'}L);\\
\bar E_1\subseteq\Pi_{X\in\T}\mathbb N_{n_{_X}}(R'\otimes_k R')
\subseteq\Pi_{X\in\T}(L^\ast\otimes_{k}L).\end{array}$$

{\bf Lemma 2.1.2\, }Let $D$ be a commutative algebra, and $\Lambda, \Sigma$ be commutative
$D$-algebras. Let $_\Lambda \mathcal{G}$ and $,\mathcal{S}_\Sigma$ be finitely generated
projective left $\Lambda$-module and right $\Sigma$-module respectively,
then there exists a $\Lambda\otimes_D\Sigma$-module isomorphism
$${\rm Hom}_{\Lambda}(\mathcal{G}, \Lambda)\otimes_D {\rm
Hom}_{\Sigma}(\mathcal{S}, \Sigma)\cong {\rm
Hom}_{\Lambda\otimes_D\Sigma}(\mathcal{G}\otimes_D\mathcal{S},
\Lambda\otimes_D\Sigma).$$

{\bf Proof}\, We first claim, that $\mathcal G\otimes_D\mathcal S$ is a projective
$\Lambda\otimes_D\Sigma$-module. In fact, suppose $\mathcal G\oplus\mathcal G'=\mathcal F_1$,
and $\mathcal S\oplus\mathcal S'=\mathcal F_2$ are free modules over $\Lambda$, $\Sigma$
respectively. Then $\mathcal G\otimes_D\mathcal S$ is a direct summand of the free
$\Lambda\otimes_D\otimes\Sigma$-module $\mathcal F_1\otimes\mathcal F_2$. Consequently,
${\rm Hom}_{\Lambda\otimes_D\Sigma}(\mathcal{G}\otimes_D\mathcal{S},
\Lambda\otimes_D\Sigma)$ is also a projective $\Lambda\otimes_D\Sigma$-module.
Consider the following commutative diagram
$$\begin{array}{c}\xymatrix{
{\rm Hom}_{\Lambda}(\mathcal{G}, \Lambda)\times {\rm
Hom}_{\Sigma}(\mathcal{S}, \Sigma)\ar[r] \ar[d]_{\psi}& {\rm
Hom}_{\Lambda}(\mathcal{G}, \Lambda)\otimes_D {\rm
Hom}_{\Sigma}(\mathcal{S}, \Sigma)\ar[ld]^{\wt{\psi}}\\
{\rm Hom}_{\Lambda\otimes_D\Sigma}(\mathcal{G}\otimes_D\mathcal{S},
\Lambda\otimes_D\Sigma) &} \end{array}$$
Let $f\in {\rm Hom}_\Lambda(\mathcal{G}, \Lambda)$ and $g\in {\rm
Hom}_\Sigma(\mathcal{S}, \Sigma)$. Since $f$ and $g$ are $D$-linear,
there exists a $\Lambda\otimes_D\Sigma$-linear map
$\psi: \mathcal{G}\otimes_D\mathcal{S}\to \Lambda\otimes_D\Sigma$, such that
$(\psi(f, g))(x\otimes y)=f(x)\otimes g(y),$
for $(x, y)\in \mathcal{G}\times \mathcal{S}$. Now
$\psi(fr,g)=\psi(f, r g)$ for $r\in D$. Thus there
exists a unique $(\Lambda\otimes_D\Sigma)$-linear map $\tilde\psi$ given by $f\otimes g\mapsto
\psi(f, g)$, which is clearly natural in
both $\mathcal{G}$ and $\mathcal{S}$. $\tilde\psi$ is an isomorphism
if $_\Lambda\mathcal G,\mathcal S_\Sigma$ are free, consequently so is
for $_\Lambda\mathcal G,\mathcal S_\Sigma$ being projective by \cite{J} p134, Proposition 3.4-3.5. This
completes the proof.

\medskip

{\bf Proposition 2.1.3}\, If $\bar R,R'$ are viewed as categories $A',B'$ respectively, and $L$ as a functor $\theta'$,
then $\theta'$ is an admissible functor given by Defined 4.3 of \cite{CB1}.

(We stress in particular, that the opposite construction is usually impossible.
Throughout the paper, we use the right module structure and upper triangular matrix, which is
opposite to the left module and lower triangular matrix used in \cite{CB1}.)

\smallskip

{\bf Proof}\, (i)\, Let $\theta'(X)=\oplus_{p=1}^{n_{_X}}Z_{(X,p)}$, (a1) implies (A1)-(A2) of Definition 4.3 in \cite{CB1}.

(ii)\, Let $\bar E^\ast_0=$Hom$_{R'}(\bar E_0,R')$, then $\bar E^\ast_0\simeq$Hom$_{R'}(R',R')\simeq R'$.
And by Lemma 1.2.2,
$$\begin{array}{ll}&\Hom_{R'}(L^\ast\otimes_{R'}L,R')\simeq\Hom_{R'\otimes_{R'}R'}(L^\ast\otimes_{R'}L,R'\otimes_{R'}R')\\
\simeq&\Hom_{R'}(L^\ast,R')\otimes_{R'}\Hom_{R'}(L,R')\simeq L\otimes_{R'}L^\ast.\end{array}$$
We may establish an equivalent relation $\sim$ on the elements of $L\otimes_{R'}L^\ast$: two elements are equivalent,
if and only if both of them acting on every $R'$-base matrix of $\bar E_0$ have the same value in $R'$.
In the case of Definition 2.1.1 (a1) \ding{172}, if $\bar R=R[a_1]$ with $R_{XY},S,K$ given in Formula (2.1-5),
$\textsf e_{_{Z_X}}\otimes_{R'}\f_{_{Z_X}}$ and
$\textsf e_{_{Z_Y}}\otimes_{R'}\f_{_{Z_Y}}$ acting on $F_Z$ equal $1_Z$, and $0$ on others;
while $\textsf e_{_{Z_X}}\otimes_{\bar R}\f_{_{Z_X}}=\textsf e_{_{Z_Y}}\otimes_{\bar R}\f_{_{Z_Y}}$
by carrying $a_1$ across the tensor product.  In (a1) \ding{173}, $\textsf e_{_{(ijl)}}\otimes_{_{R'}}\f_{_{(ijl)}}$
acting on $F_{ij}$ equal $1_{Z_{ij}}$, and $0$ on others;
while $\textsf e_{_{(ij1)}}\otimes_{_{\bar R}}\f_{_{(ij1)}}=\cdots=\textsf e_{_{(ijj)}}\otimes_{_{\bar R}}\f_{_{(ijj)}}$.
In (a1) \ding{174}, $\textsf e_{_{{Z_2}_{(X,1)}}}\otimes_{_{R'}}\f_{_{{Z_2}_{(X,1)}}}$ and
$\textsf e_{_{{Z_2}_{(Y,2)}}}\otimes_{_{R'}}\f_{_{{Z_2}_{(Y,2)}}}$ acting on $(F_{X,Z_2},F_{Y,Z_2})$ equal $1_{Z_2}$,
and $0$ on others; while $\textsf e_{_{{Z_2}_{(X,1)}}}\otimes_{_{\bar R}}\f_{_{{Z_2}_{(X,1)}}}=
\textsf e_{_{{Z_2}_{(Y,2)}}}\otimes_{_{\bar R}}\f_{_{{Z_2}_{(Y,2)}}}$. And in all the cases, we have a unique
$\textsf e_{_{Z_{Z}}}\otimes_{_{R'}}\f_{_{Z_Z}}$ acting on $\textsf e_{_{Z_{Z}}}\otimes_{_{R'}}\f_{_{Z_Z}}$ equals
$1_Z$, and $0$ on others. Moreover any element besides action on $\bar E_0$ is $0$.
Therefore $(L\otimes_{R'}L^\ast/\sim)\simeq L\otimes_{\bar R}L^\ast$, and
we obtain the $R'$-quasi-basis of $\bar E^\ast_0$ formed as
$\{\textsf e_{_{Z_{X,p}}}\otimes_{_{\bar R}}\f_{_{Z_{X,p}}}\}$ dual to that of $\bar E_0$.

(iii)\, Let $\bar E^\ast_1=$Hom$_{R'^{\otimes 2}}(\bar E_1,R'^{\otimes 2})$, by Lemma 2.1.2:
$$\Hom_{R'\otimes_{k}R'}(L^\ast\otimes_{k}L,R'\otimes_{k}R')
\simeq\Hom_{R'}(L^\ast,R')\otimes_{k}\Hom_{R'}(L,R')\simeq L\otimes_{k}L^\ast.\eqno(2.1\mbox{-}6)$$
We establish an equivalent relation $\sim$ on the elements of $L\otimes_{k}L^\ast$: two elements are equivalent,
if and only if both of them acting on every base matrix of $\bar E_1$ have the same value in $R'\otimes_kR'$.
In the case of Definition 2.1.1 (a1) \ding{172}, $\bar E_1=0$, so that $\bar E_1^\ast=0$.
In (a1) \ding{173}, $\textsf e_{i,j,h}\otimes_k\f_{i,j',l+h-1}$ acting on
$F_{ijj'l}$ equal $1_{Z_{ij}}\otimes_k1_{Z_{ij'}}$, and $0$ on others; while
$\textsf e_{_{(ij1)}}\otimes_{_{\bar R}}\f_{_{(ij'l)}}=\cdots=
\textsf e_{_{(ij,j'-l+1)}}\otimes_{_{\bar R}}\f_{_{(ij'j')}}, \forall\, j\geqslant j'$ and
$\textsf e_{_{(ij1)}}\otimes_{_{\bar R}}\f_{_{(ij'l)}}=\cdots=\textsf e_{_{(ij,j-l+1)}}\otimes_{_{\bar R}}\f_{_{(ij'j)}}
\forall\, j<j'$ by carrying $x$ across the tensor product. In (a1) \ding{174}, we have the unique element
$\textsf e_{_{{Z_2}_{(X,1)}}}\otimes_{_{R'}}\f_{_{{Z_1}_{(X,2)}}}$ acting on $F_{Z_2Z_1}$ equals $1_{Z_2}\otimes_k1_{Z_1}$,
and $0$ on $F_{Z_3Z_2}$; similar for $\textsf e_{_{{Z_3}_{(Y,1)}}}\otimes_{_{R'}}\f_{_{{Z_3}_{(Y,2)}}}$.
Moreover any element besides action on $\bar E$ is $0$. Therefore $(L\otimes_kL^\ast/\sim)\simeq L\otimes_{\bar R}L^\ast$, and
we obtain the $R'$-$R'$-quasi-basis of $\bar E^\ast_1$ formed as
$\{\textsf e_{_{Z_{X,p}}}\otimes_{_{\bar R}}\f_{_{Z_{X,q}}}\}$ dual to that of $\bar E_1$.
The picture below shows $\textsf e_{_{(ij1)}}\otimes_{_{\bar R}}\f_{_{(ij'l)}}$ in $\bar E_1^\ast$ of (a2)
as dotted arrows in the case of $s=1, r=3$:

\vskip 3mm

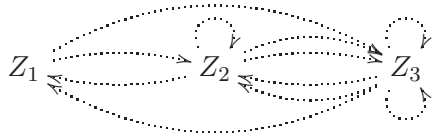
$$\xymatrix{
Z_1\ar@{.>}@/^2pc/[rrrr] \ar@{.>}@/^/[rr]   &&
Z_2\ar@{.>}@/^/[ll]\ar@{.>}@/^1pc/[rr]\ar@{.>}@/^/[rr]\ar@{.>}@(ul,ur) &&
Z_3\ar@{.>}@/^/[ll]\ar@{.>}@/^1pc/[ll] \ar@{.>}@/^2pc/[llll]\ar@{.>}@(ul,ur)\ar@{.>}@(dl,dr) }$$

\vskip 9mm

\noindent

Summary up (ii)-(iii), $B'\otimes_{A'}B'$ is determined by $\bar E^\ast_0\oplus\bar E^\ast_1$,
$J'$ by $\bar E^\ast_1$ is projective.
There is a natural map $(\bar E^\ast_0\oplus\bar E^\ast_1)\rightarrow R',
\textsf e_{_{Z_{(X,p)}}}\otimes_{_{\bar R}}\f_{_{Z_{Y,q}}}\mapsto\textsf e_{_{Z_{(X,p)}}}\f_{_{Z_{(Y,q)}}}$,
which equals $1_{Z_{(X,p)}}$ for $Z_{(X,p)}=Z_{(Y,q)}$, or $0$ otherwise.
Thus $J'$ is the kernel of the map $B'\otimes_{A'}B'\rightarrow B'$. (A3) follows.

(iv)\, (A4)-(A6) are easy. There is a natural
ordering on the basis of $\bar E^\ast_1$ transferred from that
of $\bar E_1$, (A7) follows. The lemma is proved.

\medskip

More generally, $\bar E_0$ has a $R'$-quasi-basis $\{F_Z=(F_{Z,X}\mid X\in\T)\mid Z\in\T'\}$, such that
$$\begin{array}{c}F_{Z,X}=\mbox{diag}(s_{_{Z_{(X,1)}}},\cdots,s_{_{Z_{(X,n_{_X})}}}),\quad
\left\{\begin{array}{ll}s_{_{Z_{(X,p)}}}=1_Z,&\mbox{for}\,\,Z_{(X,p)}=Z;\\
s_{_{Z_{(X,p)}}}=0,&\mbox{for}\,\,Z_{(X,p)}\ne Z.\end{array}\right.\end{array}$$
And $\bar E_1$ has a $R'$-$R'$-quasi basis $F_1,\cdots,F_l$, where for $i=1,\cdots,l$:
$$\begin{array}{c}F_{i}=(F_{i,X}\mid X\in\T),\quad F_{i,X}\in\mathbb N_{n_{_X}}
\big(R'\otimes_kR'\big).\end{array}$$
The multiplication of two base matrices is given by usual $\Delta'$-matrix product component wise.

IM$_{n_{_{X}}\times n_{_{Y}}}(R'\otimes_k R'),\forall\,X,Y\in\T$,
possesses an $\bar R$-$\bar R$-bi-module structure as follows: taken any $\f_{_{Z_{(X,p)}}}\otimes_k\textsf e_{_{Z_{(Y,q)}}}$,
$b, c\in \{x\mid X\in\T\}\cup\{a_1\}$ with $e(b)=X,s(c)=Y$,
$$\begin{array}{c}b\otimes_{\bar R}(\f_{_{Z_{(X,p)}}}\otimes_k\textsf e_{_{Z_{(Y,q)}}})\otimes_{\bar R}c
=L(b)(\f_{_{Z_{(X,p)}}}\otimes_k\textsf e_{_{Z_{(Y,q)}}})L(c).\end{array}$$

{\bf Construction 2.1.4}\, Let $\mf{A}=(R,\K,\M,\td)$ be a matrix
bi-module problem, and $\mf B$ the associated bocs. Suppose
$\bar R,R',L, \underline d$ are given in Definition 2.1.1.
Then there is an induced matrix bi-module problem $\mf A'=(R',\K',\M',H')$
in the following sense.

(i)\, The size vector of the matrices in $\K',\M'$ and $H'$ over $\T$ is
$\n$ determined by $\underline d$: $n_i=n_{_X}$ for $i\in X$. Then $t'=\sum_{i\in T}n_{i}$, the set of integers
$T'=\{1,\cdots,t'\}$.

(ii)\, $\K_0'\simeq R'$, an isomorphism $\bar E_0\stackrel{\nu_0}\rightarrow\K_0'$ gives the $R'$-quasi-basis
$\{E'_Z=\sum_{X\in\T}F_{Z,X}\ast E_X\in\mathbb D_{t'}(R')\mid Z'\in\T'\}$ of $\K_0'$,
where $\ast$ is given by Definition 1.3.2 (ii) for $S=R',p=1$.
$\K_1'=\K_{10}'\oplus\K_{11}'$. An isomorphism $\bar E\stackrel{\nu_1}\simeq\K_{10}$ gives the $R'$-$R'$-quasi basis
$\mathcal F'=\{\sum_{X\in\T}F_{i,X}\ast E_X \mid i=1,\cdots,l\}$ of $\K_{10}$ given by 1.3.2 (ii) for $S=R',p=2$;
$\K_{11}=(L^\ast\otimes_kL)\otimes_{R^{\otimes 2}}\K_1$ with a basis $\U'=\{(\f_{_{Z_{(X'_{j},p)}}}\otimes_k\textsf e_{_{Z_{(Y'_{j},q)}}})\ast
V_j\mid 1_{X'_{j}}V_j1_{Y'_{j}}=V_j,\,\forall\,p,q; j=1,\cdots,m\}\subseteq \mathbb N_{t'}(R'\otimes_k R')$ by 1.3.2 (iii) for $S=R',p=2$.
The basis of $\K_1'$ is $\V'=\mathcal F'\cup\U'$.

(iii)\, $\M'_1\simeq(L^\ast\otimes_kL)\otimes_{R^{\otimes 2}}\M_1$, with the normalized $R'$-$R'$-quasi-basis
$\A'=\{(\f_{_{Z_{(X_{i},p)}}}\otimes_k\textsf e_{_{Z_{(Y_{i},q)}}})\ast
A_i\mid 1_{X_{i}}A_i1_{Y_{i}}=A_i,\,\forall\,p,q\}$ given by
Definition 1.3.2 (iii) for $S=R',p=2$, where $1\leqslant i\leqslant n$ if
$\bar R=R$, and $1<i\leqslant n$ if $\bar R=R[a_1]$.

(iv)\, $H'=\sum_{X\in\T}H_X(L_X(x))+L(a_1)\ast A_1$, where
$L_X(x)=\bar W_X$, $H_X(\bar W_X)$ is defined in 1.3.2 (i);
and $\ast$ is given by 1.3.2 (iii) for $S=R',p=1$.

The multiplication is given by usual $\Delta'$-matrix product according to Lemma 1.3.3, as an example, we
calculate $\tr_{11}':\M_1'\otimes\K_1'\rightarrow \K_2'$ by substituting $\bar W_X$ for $x$:
$$\begin{array}{c}((\f_{X_l,p}\otimes_k\textsf e_{Z_{Y_l,q}})\ast A_l)(F_{h,Y_l}\ast E_{Y_l})
=((\f_{Z_{X_l,p}}\otimes_k\textsf e_{Z_{Y_l,q'}})F_{h,Y_l})\ast(A_lE_{Y_l});\\
((\f_{Z_{X_{i},p_1}}\otimes_k\textsf e_{Z_{Y_{i},q_1}})\ast A_i)((\f_{Z_{X'_{j},p_2}}
\otimes_k\textsf e_{Z_{Y_{j}',q_2}})\ast V_j)\\
=\sum_l((\f_{Z_{X_{i},p_1}}\otimes_k\textsf e_{Z_{Y_{i},q_1}})
(\f_{Z_{X'_j,p_2}}\otimes_k\textsf e_{Z_{Y'_j,q_2}})\otimes_{R^{\otimes 3}}\sigma_{ijl})\ast A_l,\end{array}$$
by 1.3.3 (ii) and (iii) for $p=2=q$. $\tm_{11}',\tl_{11}',\td_1'$ are similar. The proof is finished.

\medskip

{\bf Proposition 2.1.5}\, Let $(\mf A,\mf B)$ be a pair, and $\mf A'$ be given by
Construction 2.1.4. Then the associated bocs $\mf B'$ of $\mf A'$
is the induced bocs of $\mf B$ given by
Proposition 4.5 in \cite{CB1}. And $R(\mf A')\simeq R(\mf B')$.

\smallskip

{\bf Proof}\, Denote by $\mf C'=(R',\C',\N',\partial')$ the associated
bi-co-module problem of $\mf A'$.

(i) $\C'_0=\Hom_{R'}(\K_0',R')$. The isomorphism
$\bar E_0^\ast=\Hom_{R'}(\bar E_0,R')\stackrel{\nu_0^\ast}\rightarrow\Hom_{R'}(\K_0',R')=\C_0'$
gives the $R'$-quasi-basis $\{e'_{_Z}\mid Z\in\T'\}$ of $\C_0'$,
the image of the basis of $\bar E^\ast_0$ given in
the proof (ii) of Proposition 2.1.3 under $\nu_0^\ast$, which is $R'$-dual to that of $\K'_0$.

(ii)\, $\C_1'=\Hom_{R'^{\otimes 2}}(\K',R'^{\otimes 2})\simeq\C_{10}'\oplus\C_{11}'$.
$\bar E_1^\ast=\Hom_{R^{\otimes 2}}(\bar E_1,R^{\otimes 2})\stackrel{\nu_1^\ast}\rightarrow
\Hom_{R'^{\otimes 2}}(\K_{10}',R'^{\otimes 2})=\C_{10}'$ is an isomorphism,
the $R'$-$R'$-quasi-basis of $\C_{10}'$ is
the image of that in $\bar E^\ast_1$ given by the proof (iii) of 2.1.3 under $\nu_1^\ast$.
According to Lemma 2.1.2 and Formula (2.1-6):
$$\begin{array}{ll}&\C_{11}'=\Hom_{R'^{\otimes 2}}(\K_{11}',R'^{\otimes 2})
=\Hom_{R'^{\otimes 2}}((L^\ast\otimes_kL)\otimes_{R^{\otimes 2}}\K_1,R'^{\otimes 2})\\
\simeq&\Hom_{R'^{\otimes 2}\otimes_{R^{\otimes 2}}R^{\otimes 2}}((L^\ast\otimes_kL)\otimes_{R^{\otimes 2}}\K_1,
R'^{\otimes 2}\otimes_{R^{\otimes 2}}R^{\otimes 2})\\
\simeq&\Hom_{R'^{\otimes 2}}(L^\ast\otimes_kL,R'^{\otimes 2})\otimes_{R^{\otimes 2}}
(\K_1,R^{\otimes 2})\simeq (L\otimes_kL^\ast)\otimes_{R^{\otimes 2}}\C_1.\end{array}$$
Write ${\U'}^\ast=\{(\textsf e_{_{Z_{(X'_{j},p)}}}\otimes_{k}\f_{_{Z_{(Y'_{j},q)}}})\otimes_{_{R^{\otimes 2}}}v_j
\mid \forall\,p,q; 1\leqslant j\leqslant m\}$, which is $R'$-$R'$-dual to $\U'$ given in Construction 2.1.4 (ii).
The $R'$-$R'$-quasi basis ${\V'}^\ast$ of $\C_1'$ is given respectively by: ${\V'}^\ast={\U'}^\ast$
in the case of Definition 2.1.1 (a1) \ding{172};
${\V'}^\ast={\U'}^\ast\cup\{\textsf e_{(ij1)}\otimes_{_{\bar R}}\f_{(ij'l)}\mid i,j,l\}$ in (a1) \ding{173};
${\V'}^\ast={\U'}^\ast\cup\{\textsf e_{_{{Z_1}_{(X,2)}}}\otimes_{_{\bar R}}\f_{_{{Z_1}_{(X,2)}}},
\textsf e_{_{{Z_3}_{(Y,1)}}}\otimes_{_{\bar R}}\f_{_{{Z_3}_{(Y,1)}}}\}$ in (a1) \ding{174},
which is dual to the basis $\V'$ of $\K_1'$.

(iii)\, $\N_1'=\Hom_{R'^{\otimes 2}}(\M_1',R'^{\otimes 2})\simeq (L\otimes_kL^\ast)\otimes_{R^{\otimes 2}}\M_1$,
with $R'$-$R'$-quasi basis $\mathcal A'^\ast=\{(\textsf e_{_{Z_{(X_{i},p)}}}\otimes_k\f_{_{Z_{(Y_{i},q)}}})
\otimes_{_{R^{\otimes 2}}}a_i\mid \forall\,p,q\}$ dual to $\A'$, where
$1\leqslant i\leqslant n$ for $\bar R=R$, and $1<i\leqslant n$ for $\bar R=R[a_1]$.

(iv) The co-multiplication $\mu'$, the left, (resp.right) co-module action $\iota'$,
(resp.$\tau'$), and the co-derivation $\partial$
are dual to $\tm',\tl',\tr',\td'$ respectively. For example,
$\tau_{11}':\N_2'\rightarrow \N_1'\otimes_{R'}\C_1'$,
$$\begin{array}{c}\tau'_{11}((\textsf e_{_{Z_{X_l,p_1}}}\otimes_{k}\f_{_{Z_{Y_l,q_2}}})\otimes_{_{R^{\otimes 2}}}a_l)=
\sum_{q>q_2}\big((\textsf e_{_{Z_{X_l,p_1}}}\otimes_{_{R}}a_l\otimes_{_{R}}\f_{_{Z_{Y_l,q}}})
\otimes_{_{R'}}(\textsf e_{_{Z_{Y_l,q}}}\otimes_{_{\bar R}}\f_{_{Z_{Y_l,q_2}}})\big) \\
+\sum_{i,j}\big(\textsf e_{_{Z_{X_l,p_1}}}\otimes_{_{R}}a_i\otimes_{_{R}}(\sum_{X=e(a_i)=s(v_j);p,q}\f_{_{Z_{X,p}}}\otimes_{_{R}}
\textsf e_{_{Z_{X,q}}})\otimes_{_{R}}v_j\otimes_{_{R}}
\f_{_{Z_{Y_l,q_2}}}\big)\otimes_{R^{\otimes 3}}\sigma_{ijl}.\end{array}$$
Finally, we obtain the bocs $\mf B'$ of $\mf C'$, with a layered $L'$ and the differential $\dz'$ given in [CB1] 4.5,
the proof is finished.

\bigskip
\bigskip
\noindent{\bf 2.2  Reductions for matrix bi-module problems}
\bigskip

The present subsection is devoted to introducing seven reductions
of matrix bi-module problems based on Construction 2.1.4, where the last two
do not occur in the previous papers on bocses. And finally we give a regularization
as the eighth reduction.

\medskip

{\bf Proposition 2.2.1}\, (Localization) Let $(\mf{A},\mf B)$ be a
pair with $R_X=k[x,\phi(x)^{-1}]$, and
$R'_X=k[x,\phi(x)^{-1}c(x)^{-1}]$ a finitely generated localization
of $R_X$. Define $\bar R=R$, the minimal algebra $R'=R'_X\times\prod_{Y\in\T\setminus\{X\}} R_Y$,
$L=R'$. Then $L$ is an admissible $R'$-$\bar R$-bi-module.

(i)\, There exists an induced matrix bi-module problem
$\mf{A}'=(R',\K',\M',H')$ of $\mf{A}$  and a fully faithful functor
$\vartheta: R(\mf{A}')\rightarrow R(\mf{A})$.

(ii) The induced bocs $\mf{B}'$ of $\mf B$ given by localization
\cite[4.8]{CB1} is the associated bocs of $\mf A'$.

\medskip

{\bf Proposition 2.2.2}\, (Loop mutation) Let $(\mf{A},\mf B)$
be a pair, $X\in\T_0, a_1:X\mapsto X,\dz(a_1)=0$. Define $\bar R=R[a_1]$, a minimal algebra
$R'=R'_X\times\prod_{Y\in\T\setminus\{X\}}R_Y$, with $R'_X=k[x]$, and $L=R'$. Then $L$ is an admissible $R'$-$\bar R$-bi-module.

(i)\, There exists an induced
matrix bi-module problem ${\mf A}'=(R',\K',\M',\td')$ of $\mf A$, and a
equivalent functor $\vartheta: R({\mf A}')\rightarrow
R(\mf{A})$.

(ii) The induced bocs $\mf{B}'$ of $\mf B$ given by $\theta': A'\rightarrow B'$, with
$\theta'(Y)=Y,\forall\, Y\,\in\T$, $\theta'(a_1)=x$, is the  associated bocs of $\mf A'$ by Proposition 2.1.3.

\medskip

{\bf Proposition 2.2.3} (Deletion) Let $(\mf{A},\mf B)$ be a
pair, $\T'\subset \T$. Define $\bar R=R$, $R'=\prod_{X\in \T'}R_X$, and $L=R'$.
Then $L$ is an admissible $R'$-$\bar R$-bi-module.

(i)\,  There exists an induced
matrix bi-module problem $\mf{A}'=(R',\K',\M',H')$ of
$\mf{A}$, and a fully faithful functor $\vartheta:
R(\mf{A}')\rightarrow R(\mf{A})$.

(ii) The induced bocs
$\mf{B}'$ obtained by deletion of $\T\setminus\T'$ given by \cite[4.6]{CB1}
is the associated bocs of $\mf A'$.

\medskip

{\bf Proposition 2.2.4} (Unraveling) Let $(\mf{A},\mf B)$ be a
pair with $R_X=k[x,\phi(x)^{-1}]$. Define $\bar R=R$, $R'=S\times\prod_{Z\in\T\setminus\{X\}} R_Z$,
$L=K\oplus(\oplus_{Z\in\T\setminus\{X\}} R_Z)$ with $S$ and $K$ given by Formula (2.1-3).
Then $L$ is an admissible $R'$-$\bar R$-bi-module according to Definition 2.1.1 (a1) \ding{173}.

(i)\, There exists an induced matrix bi-module
problem $\mf{A}'=(R',\K',\M',H')$ and a fully faithful
functor $\vartheta: R(\mf{A}')\rightarrow R(\mf{A})$.

(ii) The induced bocs $\mf{B}'$ given by unraveling in \cite[4.7]{CB1}
is the associated bocs of $\mf A'$.

\medskip

{\bf Proposition 2.2.5} (Edge reduction)\, Let $(\mf{A},\mf B)$
be a pair, $X,Y\in\T_0,a_1:X\mapsto Y,\dz(a_1)=0$.
Define $\bar R=R[a_1]$, $R'=S\times\prod_{Z\in\T\setminus\{X,Y\}} R_Z$,
$L=K\oplus(\oplus_{Z\in\T\setminus\{X,Y\}} R_Z)$ with $S$ and $K$ defined in Formula
(2.1-4). Then $L$ is an admissible $R'$-$\bar R$-bi-module by Definition 2.1.1 (a1) \ding{174}.

(i)\, There exists an induced matrix bi-module problem $\mf{A'}=(R',\K',\M',H')$, and an equivalence $\vartheta:
R(\mf{A}')\rightarrow R(\mf{A})$.

(ii)\, The induced bocs
$\mf{B}'$ given by edge reduction in \cite[4.9]{CB1}
is the associated bocs of $\mf A'$.

\medskip

{\bf Proposition 2.2.6}\, Let $(\mf A,\mf B)$ be a pair, $X,Y\in\T_0,a_1:X\mapsto Y,\dz(a_1)=0$.
Set $\bar R=R[a_1]$, $R'=R$, $L=K\oplus(\oplus_{U\in\T\setminus\{X,Y\}}R_U)$ with $K:R_X\stackrel{(0)}\rightarrow R_Y$.
Then $L$ is an admissible $R'$-$\bar R$-bi-module.

(i)\, There is an induced problem $\mf A'=(R',\K', \M', \td')$,
with $\K'=\K$; $\M'=\M^{(1)}$; $H'=H$. And an induced fully faithful functors $\vartheta:
R(\mf{A}')\rightarrow R(\mf{A})$.
The subcategory of $R(\mf A)$ consisting of
representations $P$ with $P(a_1)=0$ is equivalent to $R(\mf A')$.

(ii)\,  The induced bocs
$\mf{B}'$ given by the admissible functor $\theta':A'\rightarrow B'$ with
$\theta'(U)=U,\forall\, U\in\T$ and $\theta'(a_1)=0$ is the associated bocs of $\mf A'$.

\medskip

{\bf Proposition 2.2.7}\,  Let $(\mf A,\mf B)$ be a pair, $X,Y\in\T_0,a_1:X\mapsto Y,\dz(a_1)=0$.
Set $\bar R=R[a_1]$ $R'=S\times\prod_{U\in\T\setminus\{X,Y\}} R_U$, $L=K\oplus(\oplus_{U\in\T\setminus\{X,Y\}} R_U)$
with $S$ and $K$ defined in Formula (2.1-5). Then $L$ is an admissible $R'$-$\bar R$-bi-module.

(i)\, There is an induced local problem $\mf A'$, and an induced fully faithful functors $\vartheta:
R(\mf{A}')\rightarrow R(\mf{A})$.
The subcategory of $R(\mf A)$ consisting of $P$ of size vector $\m$
with $m_{_X}=m_{_Y}$, and rank$(P(a_1))$ $=m_{_X}$ for $Y\in\T_0$, or
$P(a_1)^{-1}P(x)P(a_1)=P(y)$ for $Y\in\T_1$, is equivalent to $R(\mf A')$.

(ii)\,  The induced bocs
$\mf{B}'$ given by the admissible functor $\theta':A'\rightarrow B'$ with $\theta'(X)=Z=\theta'(Y)$;
$\theta'(x)=z,\theta'(a_1)=(1)$, or in addition $\theta'(y)=z$ if $Y\in\T_1$,
is the associated bocs of $\mf A'$.

\medskip

Let $\mf{A}=(R,\K,\M,\td)$ be a matrix
bi-module problem, $\mf C=(R,\C,\N,\partial)$ be the associated
bi-co-module problem, and $\mf B$ the bocs of $\mf C$. Then
$$\begin{array}{c}\td(V_1)=A_1+ \sum_{l>1}\zeta_{1l} A_l,\,\,\td(V_j)\in
\M^{(1)}_1,j\geqslant 2\,\, \mbox{in}\,\, \mf A\\
\Longleftrightarrow  \partial (a_1)=v_1\,\, \mbox{in}\,\, \mf C
\Longleftrightarrow \delta(a_1)=v_1\,\,\mbox{in}\,\,\mf B.\end{array}\eqno{(2.2\mbox{-}1)}$$
In fact, since $\partial(a_1)=\sum_{j=1}^m\zeta_{j1}v_{j}$,
we have $\partial(a_1)=v_1$, if and only if $\zeta_{11}=1_{s(a_1)}\otimes_k1_{t(a_1)}$, and
$\zeta_{j1}=0$ for all $j\geqslant 2$, if and only if
$\td(V_1)=A_1+\sum_{l>1}\zeta_{1l} A_l$ and $\td( V_j)\in
\M^{(2)}$ for all $j\geqslant 2$, since $\td(V_{j})=\zeta_{j1}A_1+
\sum_{i>1}\zeta_{ji}A_i$. Finally, since
$\iota_1(a_1)=0, \tau_1(a_1)=0$ by triangularity of
$\mf C$, $\delta(a_1)=v_1$ in $\mf B$, if and
only if $\partial(a_1)=v_1$ in $\mf C$ by Theorem 1.4.2.

\medskip

{\bf Remark}\, Let $\mf{A}, \mf C$ be given above, with $\partial(a_1)=v_1$, then

(i)\, $\K^{(1)}=\K_0\oplus(\oplus_{j=2}^m
\bar\Delta\otimes_{_{R^{\otimes 2}}} V_j)$ is a sub-algebra of $\K$,
and $\M^{(1)}=\oplus_{i=2}^n\bar\Delta\otimes_{_{R^{\otimes 2}}}a_i$ is a
$\K^{(1)}$-$\K^{(1)}$-sub-bi-module;

(ii)\, $\C^{(1)}=\bar\Delta\otimes_{_{R^{\otimes 2}}}v_1$ is a
co-ideal of $\C$, $\C^{[1]}=\C/\C^{(1)}$ is a
quotient co-algebra, and
$\N^{(1)}=\bar\Delta\otimes_{_{R^{\otimes 2}}}a_1$ is a
$\C$-$\C$-sub-bi-co-module, thus $\N^{[1]}=\N/\N^{(1)}$ is a
$\C^{[1]}$-$\C^{[1]}$-quotient bi-co-module.

\smallskip

{\bf Proof}\, (i)\, For any $V_i,V_j\in\K_1$,
$$\begin{array}{l}\td(V_iV_j)= \td(\sum_{l=1}^m
\gamma_{ijl}\otimes_{_{R^{\otimes 2}}}V_l)
=\sum_{l=1}^m \gamma_{ijl}\otimes_{_{R^{\otimes 2}}}\td(V_l)
=\sum_{l,p} (\gamma_{ijl}\otimes_{_{R^{\otimes 2}}}\zeta_{lp})\otimes_{_{R^{\otimes 2}}}
A_p.\end{array}$$
By triangularity (2.1-1), $\td
(V_iV_j)=\td(V_i)V_j+V_i\td(V_j)\in \M^{(1)}$. Consequently, the
coefficient of $A_1$ in the formula above $\sum_{l}
(\gamma_{ijl}\otimes_{_{R^{\otimes 2}}}\zeta_{l1})=0$, where $\zeta_{11}=1$, $\zeta_{l1}=0$ for
$l>1$ by the hypothesis, so that $\gamma_{ij1}=0$ for all
$1\leqslant i,j\leqslant m$. Therefore
$V_iV_j=\sum_{l>1}\gamma_{ijl}\otimes_{_{R^{\otimes 2}}}V_l\in\K^{(1)}$
and hence $\K^{(1)}$ is a subalgebra of $\K$. Finally, $\M^{(1)}$ is a
$\K^{(1)}$-$\K^{(1)}$-bi-module still by triangularity.

(ii) Since $\mu(v_1)=\mu(\partial(a_1))
=(\partial\otimes\id)(\iota(a_1))+(\id\otimes\partial)(\tau( a_1))
\stackrel{(2.1\mbox{-}2)}=(\partial\otimes \id)(e_{s(a_1)}\otimes_{_R}a_1)
+(\id\otimes\partial)(a_1\otimes_{_R}e_{t(a_1)})=0$,
$\C^{(1)}$ is a co-ideal of $\C$, which finishes the proof.

\medskip

{\bf Proposition 2.2.8} (Regularization)\, Let
$(\mf{A},\mf B)$ be a pair with $\delta(a_1)=v_1$.

(i) There is an induced matrix bi-module problem $\mf{A}'=(R,
\K^{(1)}, \M^{(1)}, H)$ of $\mf A$, and an
equivalent functor $\vartheta: R(\mf{A}') \rightarrow
R(\mf{A})$.

(ii) The induced bocs
$\mf{B}'$ given by regularization \cite[4.2]{CB1} is the associated bocs of $\mf A'$.

\smallskip

{\bf Proof}\, (i) $\mf A'$ is a matrix bi-module problem by Formula (2.2-1) and Remark (i) above.
Note that $R'=R,\T'=\T$, for any $P\in R(\mf A)$
of size vector $\m$, let $f=\sum_{X\in\T}I_{m_X}\ast E_X+P(a_1)\ast V_1$, then
$P'=f^{-1}Pf\in R(\mf A')$. Therefore, $\vartheta$ is an equivalent functor.

(ii) $\mf C'=(R,\C^{[1]}, \M^{[1]}, \bar\partial)$ with $\bar\partial$ induced
from $\partial$ is the associated bi-co-module problem of $\mf A'$ by
Remark (ii) above. Thus the associated
bocs $\mf B'$ is given by regularization from $\mf B$.

\medskip

Let $(\mf A,\mf B)$ be a pair, with a layer
$L=(R;\omega;a_1,\cdots,a_n;v_1,\cdots,v_m)$ in $\mf B$. Suppose $a_1: X\mapsto Y$,
$\dz(a_1)=\sum_{j=1}^m f_{j}(x,y)v_{j}\ne 0$. In order to obtain
$\dz(a_1)=h(x,y)v'_1$, we make the following base change:
$$\begin{array}{c}(v_{1}',\cdots,v_{m}')=(v_{1},\cdots,v_{m})F(x,y)\end{array}\eqno{(2.2{\mbox{-}}2)}$$
with $F(x,y)\in\IM(R\otimes_kR)$ invertible. When $X\in \T_0$ or $Y\in \T_0$, $R$ is preserved;
but when $X, Y\in \T_1$, some localization $R_{X}'=R_{X}[c(x)^{-1}]$
(resp. $R_{Y}'=R_{Y}[c(y)^{-1}]$) is needed, see \cite[$\S
5$]{CB1}. Consequently, we have a base change of $\K_1$ dually given by
$$\begin{array}{c}(V'_{1},\cdots,V'_{m})=(V_{1},\cdots,V_{m})F(x,y)^{-T}.\end{array}\eqno{(2.2{\mbox{-}}3)}$$

Finally we mention a simple fact according to all the reductions defined above.
Suppose we start from a matrix bi-module problem
$\mf{A}^0=(R^0,\K^0,\M^0,H=0)$ with $\T^0$
trivial, if there is a sequence of reductions $\mf A^0,\mf A^1,\cdots,\mf A^r$
with $\mf{A}^r=(R^r,\K^r,\M^r,H^r)$, and $X\in \T^r_1$,
$H_X^r=(h_{ij}(x))$ of size $t^r$, then
$h_{ij}(x)=a_{ij}+b_{ij}x\in k[x]$ is of degree $1$.

\bigskip
\bigskip
\noindent{\bf 2.3 Canonical Forms}
\bigskip

We will give a canonical form (cf.\cite{S}) for each representation of a matrix
bi-module problem, and a sequence of reductions in the subsection.

\medskip

{\bf Convention 2.3.1}\, Suppose $\mf{A}$ is a matrix bi-module problem,
$\mf{A}'$ an induced problem and
$\vartheta: R(\mf A')\rightarrow R(\mf A)$ an induced functor.
Let $\m'$ be a size vector over $\T'$ of $\mf A'$,
define a size vector $\m=(m_1,m_2,\ldots,m_t)$ over $\T$ of $\mf A$ based on $\m'$:

(i) for regularization, loop mutation, localization, and Proposition 2.2.6, set
$\m=\m'$;

(ii) for deletion, set $m_i=m'_i$ if $i\in X, X \in \T'$, and $0$ if
$i\in X, X\in \T\setminus \T'$;

(iii) for edge reduction, set $m_i=m'_i$ if $i\in Z, Z\ne X,Y$,
$m_i=m'_{Z_1}+m'_{Z_2}$ if $i\in X$, and $m_i=m'_{Z_2}+m'_{Z_3}$ if $i\in
Y$; for proposition 2.2.7, set $m_{_X}=m_{_{Z}}'=m_{_Y}$;

(iv) for unraveling, set  $m_i=m'_i$ if $i\not\in X$, and
$m_i=\sum_{i=1}^s\sum_{j=1}^r im'_{Z_{ij}}+m'_{Z_0}$ if $i\in X$.

Then $\m$ is said to be the {\it size vector determined by $\m'$},
and is denoted by $\vartheta(\m')$.

\medskip

Let $\mf A=(R,\K,\M,H)$ be a matrix bi-module problem with trivial $\T$,
and $\m$ be a size vector. For the sake of simplicity,
write $$\begin{array}{c}H_\m(k)=\sum_{X\in \T} H_X(I_{m_{_X}}),\quad
H(k)=\sum_{X\in \T}H_X(1).\end{array}\eqno{(\ast)}$$

Note that if a size vector $\m$ over $\T$ is not sincere, let
$\T'=\{i\mid m_i\ne 0\}$, set the induced bi-module problem $\mf A'$
given by a deletion of $\T\setminus \T'$, then the size vector
$\m'=(m_i\mid m_i\ne 0)$ is sincere over $\T'$.

\medskip

{\bf Lemma 2.3.2} (cf.\cite{S})\, Let $\mf{A}=(R,\K,\M,H)$ be a matrix bi-module
problem with $\T$ being trivial. Let $P$ be a representation of $\mf A$ with a
sincere size vector $\m$, by Definition 1.3.4,
$$\begin{array}{c}P=H_{\m}(k)+\sum_{i=1}^nP(a_i)\ast A_i.\end{array}$$
Then there exists a matrix bi-module problem
$\mf{A}'=(R',\K',\M',H')$ induced by one of the following three compositions:

(i) Regularization;

(ii) Edge reduction: edge reduction + deletion;

(iii) Loop reduction: loop mutation + unraveling + deletion,

\noindent with a fully faithful functor $\vartheta:
R(\mf{A}')\rightarrow R(\mf{A})$, such that there is a representation $P'\in R(\mf A')$
having a sincere size vector $\m'$ over $\T'$ with
$P\simeq\vartheta(P')$ and $\vartheta(\m')=\m$, where
$$\begin{array}{c}P'=H_{\m}(k)+
B\ast A_1+\sum_{i=1}^{n'}P(a'_i)\ast A'_i\end{array}$$
as a $k$-matrix, $B$ is given by one of Formulae (2.3-1)-(2.3-3) below.

\smallskip

{\bf Proof}\, Let $a_1: X\mapsto Y$ in the associated bocs $\mf B$ of $\mf A$.

(i)\, If $\dz(a_1)=v_1$, denote by $\emptyset$ a distinguished zero matrix-block. Set
$$B=\emptyset_{m_{_X}\times m_{_Y}};\quad G=\emptyset_{1\times 1}.\eqno{(2.3\mbox{-}1)}$$
We proceed with a regularization for $\mf{A}$, obtain an induced problem
$\mf{A}'$, and an equivalence $\vartheta:
R(\mf{A}')\rightarrow R(\mf{A})$. Then $P\simeq\vartheta(P')$,
$P'$ is given in the proof (i) of Proposition 2.2.8.

(ii) If $\dz(a_1)=0$ and $X\ne Y$, we
proceed with an edge reduction for $\mf{A}$ and obtain an induced
problem $\mf A_1$. Let the invertible matrices
$f_X\in\IM_{m_{_X}}(k),f_Y\in \IM_{m_{_Y}}(k)$, such that
$$\begin{array}{c}
B=f_{X}^{-1}P(a_1) f_{Y}={{\,0\,\, I_r}\choose{0\,\,\
0}}_{m_{_X}\times m_{_Y}}
\quad \mbox{with}\,\, r=\rank(P(a_1));\\
G=(0),\,\,(1_{Z_2}),\,\,(0\, 1_{Z_2}),\,\,{{1_{Z_2}}\choose{0}},\,\,
{{\,0\,\, 1_{Z_2}}\choose{0\,\,\,\,\,\
0}},\end{array}\eqno{(2.3\mbox{-}2)}
$$
where the first four cases of $G$ are obtained by a deletion for $\mf A_1$:
\ding{172} $r=0$, delete $Z_2$; now suppose $r>0$, \ding{173} $m_{_X}=r=m_{_Y}$, delete
$Z_1$ and $Z_3$; \ding{174} $m_{_X}=r, m_{_Y}>r$, delete $Z_1$; \ding{175}
$m_{_X}>r, m_{_Y}=r$, delete $Z_3$. We obtain an induced problem
$\mf{A}'$ of $\mf A_1$, and a fully faithful functor
$\vartheta: R(\mf{A}')\rightarrow R(\mf{A})$.
Let $\m'=(m'_i)_{i\in T'}$ be a size vector over $\T'$, with $m'_{_Z}=m_{_Z}$ for $Z\in
\T\setminus \{X,Y\}$, $m_{_{Z_1}}'=m_{_X}-r$, $m_{_{Z_2}}'=r$ and
$m_{_{Z_3}}'=m_{_Y}-r$. Then there is some $P'=f^{-1}Pf\in R(\mf A')$ of size $\m'$
with $P\simeq\vartheta(P')$ in $R(\mf A)$, where $f=f_{X}\ast E_{X}+f_{Y}\ast E_{Y}+
\sum_{Z\in\T\setminus\{X,Y\}}I_{m_{_Z}}\ast E_Z$.

(iii)\, If $\dz(a_1)=0$ and $X=Y$, suppose $P(a_1)\simeq J=
\oplus_{i=1}^s(\oplus_{j} J_j(\lambda_i)^{e_{ij}})$,
a Jordan form over $k$ with the maximal
size $r$ of the Jordan blocks. We first proceed with a loop
mutation $a_1\mapsto (x)$, then with an unraveling for the polynomial
$g(x)=(x-\lambda_1)\cdots(x-\lambda_s)$ and the positive integer
$r$, thus obtain an induced problem $\mf A_1$ of $\mf A$. Let the invertible matrix
$f_{X}\in\IM_{m_{_X}}(k)$, such that
$$B=f_{X}^{-1}P(a_1)f_{X}=W;\quad G=\bar W,\eqno{(2.3\mbox{-}3)}$$
where $W$ is a Weyr matrix, and $G=\bar W$ being a Weyr matrix similar
to $\oplus_{e_{ij}\ne 0} J_j(\lambda_i)1_{Z_{ij}}$ over $S$.
Finally, delete a set of vertices $\{Z_0; Z_{ij}\,|\,
e_{ij}=0\}$ from $\mf A_1$. We obtain an induced problem $\mf{A}'$ of $\mf A_1$, and a fully
faithful functor $\vartheta: R(\mf{A}')\rightarrow R(\mf{A})$. Let
$\m'=(m'_i)_{i\in T'}$ be a size vector over $\T'$ with $m'_{_Z}=m_{_Z}$
for $Z\in\T\setminus\{X\}$, and $m'_{_{Z_{ij}}}=e_{ij}\ne 0$; let $f=f_{X}\ast E_{X}+
\sum_{Z\in\T\setminus\{X\}}I_{m_{_Z}}\ast E_Z$, then $P'=f^{-1}Pf\in R(\mf A')$ with
a sincere size vector $\m'$, such that $P\simeq\vartheta(P')$ in $R(\mf A)$.

In all the cases, $P'=H_{\m}(k)+B\ast A_1+\sum_{i=1}^{n'}P(a'_i)\ast A'_i$.
The proof is completed.

\medskip

Repeating the procedure of Lemma 2.3.2, we obtain the following theorem by induction.

\medskip

{\bf Theorem 2.3.3} (cf.\cite{S})\, Let $\mf{A}=(R,\K,\M,H=0)$ be a matrix bi-module
problem with $\T$ trivial. Let $\m$ be a sincere size vector over
$\T$ and $P\in R(\mf{A})$ be a representation of size $\m$.
Then there exist a unique sequence of matrix bi-module problems, and a
unique sequence of representations:
$$
\begin{array}{lllllll}\mf{A}=\mf{A}^0,& \mf{A}^1, &\ldots,
&\mf{A}^{i}, &\mf{A}^{i+1},& \ldots,
&\mf{A}^r \\G^0=0,&G^1,&\ldots,&G^{i},&G^{i+1},&\ldots,&G^r\\
B^0=0,&B^1,&\ldots,& B^{i},&B^{i+1},&\ldots,&B^r,\\
P^0=P,&P^1,&\ldots,&P^{i},&P^{i+1},&\ldots,&P^r\end{array}
$$ where $\mf A^{i+1}$ is obtained from
$\mf A^{i}$ according to one of the three reductions given in Lemma 2.3.2
for $i=0,\cdots,r-1$, such that

(i)\, $B^{i+1}$ and $G^{i+1}$ is defined by one formula of (2.3-1)-(2.3-3).

(ii)\, Let $\vartheta^{i,i+1}: R(\mf{A}^{i+1})\rightarrow R(\mf{A}^{i})$ be
the induced functor, there is a sincere size vector $\m^{i+1}$ over $\T^{i+1}$ with
$\vartheta^{i,i+1}(\m^{i+1})=\m^{i}$ and some $P^{i+1}\in R(\mf A^{i})$
of size $\m^{i+1}$ with $\vartheta^{i,i+1}(P^{i+1})\simeq P^{i}$. Denote by
$A_1^j$ the first quasi-base matrix of $\M^j_1$, then
$P^{i+1}=H_{\m^{i+1}}(k)+\sum_{j=1}^{n^{i+1}}M^{i+1}(a_j)\ast A^{i+1}_j$.

(iii)\, Denote for $i<j$ the composition of induced functors
$\vartheta^{ij}=
\vartheta^{i,i+1}\cdots\vartheta^{j-1,j}:R(\mf A^j)\rightarrow R(\mf
A^{i})$. Then $\vartheta^{0,i+1}(H_{\m^{i+1}}(k))=
\sum_{j=0}^{i} B^{j+1}\ast A^j_1\in R(\mf A)$.

\medskip

With the notation of Theorem 2.3.3, if $\mf A^r$ is minimal:
$$\begin{array}{c}\vartheta^{0r}(P^r)=\vartheta^{0r}(H_{\m^r}(k))
=\sum_{i=0}^{r-1} B^{i+1}\ast A^i_1\simeq P\in R(\mf A).
\end{array}\eqno{(2.3\mbox{-}4)}$$
The matrix $\vartheta^{0r}(P^r)$ is called the {\it canonical form} of
$P$, and denoted by $P^{\infty}$.

In the second and the third sequence of Theorem 2.3.3, if $G^{i+1}$ is obtained by an edge
or loop reduction, then $``1"$ appearing in $B^{i+1}$, which is not an
eigenvalue in the case of $B^{i+1}$ being a Weyr matrix, is called a
{\it link} of $P^\infty$. And denote by
$l(P^\infty)$ the number of links in $P^\infty$.

\medskip

{\bf Corollary 2.3.4} \cite{S,XZ}\, The canonical form of any representation $P$ over
a matrix bi-module problem $\mf{A}=(R, \K,\M,H=0)$ with $R$ trivial is uniquely determined.
Moreover,

(i)\, for any $P, Q\in R(\mf{A})$, $P\simeq Q$ if and only if $P$
and $Q$ have the same canonical form;

(ii)\, $P$ is indecomposable if and only if $l(P^\infty)=\mbox{dim}(P)-1$.

\medskip

{\bf Corollary 2.3.5}\,  Let $\mf{A}=(R,\K,\M,H=0)$ be a matrix
bi-module problem with $\T$ trivial, let $\mf{A}'=(R',
\K',\M',H')$ be an induced problem
obtained by a sequence of reductions, and $\vartheta: R(\mf{A}')\rightarrow
R(\mf{A})$ be the induced functor.  If
$\T'$ is trivial, then there is a unique reduction sequence
$\mf{A}=\mf{A}^0, \mf{A}^1, \cdots, \mf{A}^{i}, \mf{A}^{i+1}, \cdots,
\mf{A}^r=\mf A'$ in the sense of Theorem 2.3.2 performed for
$P=\vartheta(H'(k))\in R(\mf A)$ according to Theorem 2.3.3.

\medskip

Under the hypothesis of Corollary 2.3.5, let $\m^r=(1,\cdots,1)$ and $\m^{i}=\vartheta^{i,r}(\m^r)$,
we give a $R^r$-structure on $B^{i+1}$ in Theorem 2.3.3 and denoted by $G_r^{i+1}\in
\IM_{m^i_{s(a_1^i)}\times m_{t(a_1^i)}}(R^r)$: a non-zero element $g_{pq}$ of $G^{i+1}_r\ast A^r_1$ belongs to
$R^r_{X^r}$, whenever the entrance of the identity matrix $E^r\in\K^r$ at the $p$-th row,
as well as the $q$-th column by the definition of $H^r$,
is $1_{X^r}\in R^r_{X^r}$. Then
$$\begin{array}{c}H^r=\sum_{i=0}^{r-1}G^{i+1}_r\ast A^i_1.\end{array} \eqno{(2.3\mbox{-}5)}$$
It is easy to see, that $\mf A^r$ is local if and only if
$l\big(\vartheta^{0r}(H^r(k))\big)=$dim$\big(\vartheta^{0r}(H^r(k))\big)-1$.
The non-eigenvalue $1_{X^r}$ appearing in $G^i_r$ is called a link of $H^r$.

\bigskip
\bigskip
\noindent{\bf 2.4 Defining systems}
\bigskip

We introduce a concept of defining system in the subsection.

\medskip

Let $B=(b_{ij})_{t\times t}$ and $C=(c_{ij})_{t\times t}$ be two
$t\times t$ matrices over $k$. Given $1\leqslant p,q\leqslant t$, the notation
$B\equiv_{\prec(p,q)}C$ (resp. $B\equiv_{(p,q)}C$, $B\equiv_{\preccurlyeq(p,q)}C$) means that
$b_{ij}=c_{ij}$ for any $(i,j)\prec(p,q)$ (resp.
$(i, j)=(p,q),\, (i,j)\preccurlyeq(p,q)$). One can define the similar notions
for partitioned matrices.

Let $\mf A=(R,\K,\M,H=0)$ be a matrix bi-module problem with the associated bocs $\mf B$, where
$\T$ and $\V=\{V_1,\cdots,V_m\}$ are  both trivial. Suppose there exists a sequence with each reduction
being in the sense of Lemma 2.3.2:
$$
(\mf{A},\mf{B})=(\mf{A}^0,\mf{B}^0), (\mf{A}^1,\mf{B}^1) \cdots, (\mf{A}^{i},\mf{B}^i) (\mf{A}^{i+1},\mf B^{i+1}) \cdots,
(\mf{A}^r,\mf B^r),\cdots,({\mf{A}^s},\mf B^s).\eqno{(2.4\mbox{-}1)}
$$

Sometimes, it is difficult to determine the dotted arrows in the induced bocs
after some reductions. Instead, we may consider a system of equations on dotted elements
as variables (not dotted arrows), and give explicitly the linear relations on those elements
(used in section 4.5).

\medskip

{\bf Theorem 2.4.1}\, With the assumption above, for $i=0,\cdots,s$, there exists a matrix equation $\IE^i$ over
$R^i\otimes_kR^i$, such that

(i) there is a basic solution of the system $\IE^i$, which forms a $R^i\otimes_kR^i$-quasi-basis
of $\K_1^i$;

(ii) the free variables in the basic solution form
a $R^i\otimes_kR^i$-quasi-basis of $\C_1^i$.

\smallskip

{\bf Proof}\,
Let $\Phi_{\m^0}=\sum_jv_j\ast V_j$, with $\m^0=(1,\cdots,1)$,
be the formal product $\Pi$ of the pair $(\mf A,\mf B)$, and let
$\IE^0: \Phi_{\m^0} H^0\equiv_{\prec (p,q)} H^0\Phi_{\m^0}$
be a matrix equation, where the leading position of $A_1$ is $(p,q)$, $H^0=0$.
Then the basic solution of $\IE^0$ is $\{V_j\}_j$, the $R\otimes_kR$-basis of $\K_1$;
and the set of free variables is $\{v_j\}_j$, the $R\otimes_kR$-basis of $\C_1$. The assertion follows for $i=0$.

Suppose a equation system $\IE^i$ for the pair $(\mf A^i,\mf B^i)$ satisfying (i)
and (ii) has been obtained:
$$\begin{array}{c}\IE^i: \Phi_{\m^i}H^i\equiv_{\prec(p^{i},q^{i})} H^i\Phi_{\m^i},\quad
\Phi_{\m^i}=\sum_{X\in\T}\bar w_{_X}^i\ast E_X+\sum_j\bar v^i_j\ast V_j,\end{array}\eqno{(2.4\mbox{-}2)}$$
where $(p^i,q^i)$ is the leading position of $A_1^i$, $\bar v_j^i=(v^i_{jpq})$ is the split of $v_j$ of size $m^i_{s(v_j)}\times m^i_{t(v_j)}$,
$\bar w_{_X}^i=(w_{{_X}pq}^i)$ is a square matrix of the size $m^i_{_X}$, and $v^i_{jpq},w_{{_X}pq}^i$
are said to be {\it dotted elements} (not dotted arrows).
Both of them are over $R^i\otimes_k R^i$, in fact, $v^i_{jpq}$ (resp. non-zero $w^i_{_{X}pq}$) $:X^i\mapsto Y^i$,
provided the entry of the identity matrix $E^i\in\K^i$ with the same row (resp. column) index of $v^i_{jpq}$ (resp. $w^i_{_{X}pq}$)
is $1_{_{X^i}}\in R^i_{X^i}$
(resp. $1_{_{Y^i}}\in R^i_{Y^i}$). We now construct the system $\IE^{i+1}$.

(i) In the case of Regularization, we have $\m^{i+1}=\m^{i}$,
then equation $\Phi_{\m^{i}} H^{i}\equiv_{(p^{i},q^{i})}H^{i}\Phi_{\m^{i}}$,
combining with $\IE^i$ form the equation system $\IE^{i+1}:\Phi_{\m^{i+1}}H^{i+1}
\equiv_{\prec(p^{i+1},q^{i+1})} H^{i+1}\Phi_{\m^{i+1}}$.

\smallskip

(ii) In the cases of Loop or Edge reductions of Lemma 2.3.2, set the first arrow
$a^i_1:X^i\mapsto Y^i$ in $\mf B^i$, denote by $\n=\vartheta^{i,i+1}(1,1\cdots,1)$  the size
vector of $H^{i+1}(k)$ over $\T^{i}$, thus the size vector $\vartheta^{0i}(\n)=\m^{i+1}$
over $\T$ according to Convention 2.3.1.
Suppose $\bar v_j^i=(v^i_{jpq})$, define a $n_{s(v_{jpq})}\times n_{t(v_{jpq})}$-matrix block
according to Theorem 2.1.3 (ii):
$$\begin{array}{c}\bar v_j^{i+1}=\sum_{\alpha,\beta}(\f_{(s(v_{jpq}),\alpha)}
\otimes_k\textsf e_{(t(v_{jpq}),\beta)})\ast v_{jpq}.\end{array}$$
If $\bar w^i_{_X}=(w^i_{_{Xpq}})$,
let $\bar w_{_X}^{i+1}$ be given below for loop and edge reduction respectively:
$$\begin{array}{c}\bar w_{_X}^{i+1}=\sum_{\alpha,\beta}(\f_{(s(w_{_{Xpq}}),\alpha)}
\otimes_k\textsf e_{(t(w_{_{Xpq}}),\beta)})\ast w_{_{Xpq}}^i\\
+\left\{\begin{array}{l}\big(\sum_{h,j,l}(\textsf e_{_{(hj1)}}\otimes_{_R}\f_{_{(hj'l)}})
\ast F_{hjj'l}\big)\ast E^i_{X^i}\mid_{\textsf p\in X};\\
\big((\textsf e_{_{{Z_2}_{(X^i,1)}}}\otimes_{_R}\f_{_{{Z_1}_{(X^i,2)}}})\ast F_{Z_2Z_1}\big)\ast E^i_{X^i}
+\big((\textsf e_{_{{Z_3}_{(Y^i,1)}}}\otimes_{_R}\f_{_{{Z_2}_{(Y^i,2)}}})
\ast F_{Z_3Z_2}\big)\ast E^i_{Y^i})\mid_{\textsf p\in X},\end{array}\right.\end{array}$$
where $``\mid_{\textsf p\in X}"$ means that the matrix blocks are
restricted inside the $(\textsf p,\textsf p)$-th block
for any $\textsf p\in X$ partitioned by $\T$. Since the admissible bi-module depends
only on the vertex set $\T$, $\ast |_{\textsf p\in X}$ are all equal when $\textsf p$ runs over $X$.
Thus $\IE^{i+1}: \Phi_{\m^{i+1}}H^{i+1}\equiv_{\prec(p^{i+1},q^{i+1})} H^{i+1}\Phi_{\m^{i+1}}$ with
$\Phi_{\m^{i+1}}=\sum_{X\in\T}\bar w_{_X}^{i+1}\ast E_X+\sum_j\bar v^{i+1}_j\ast V_j$
in both cases, the theorem follows by induction.

\medskip

$\Phi_{\m^{i}}$ and $\IE^i$ in Formula (2.4-2) are called the {\it variable matrix} of size vector $\m^i$, and the
{\it defining system of the pair $(\mf A^i,\mf B^i)$} respectively.
It is clear that the system $\IE^i$
consisting of the equations at the
$(\textsf p_i,\textsf q_i)$-th block for some $1\leqslant i\leqslant n$, where
$(\textsf p_i,\textsf q_i)$ is the index of the leading
block determined by $A_i$ partitioned under $\T$.

Next, we give a deformed system based on the defining system.
Fix some $0<r<s$, write:
$$\begin{array}{c}H^i=H^i_1+H^i_2,\quad H_1^i=\sum_{j=1}^{r-1}G_i^{j+1}\ast A^j_1,
\quad H^i_2=\sum_{j=r}^{i-1}G_i^{j+1}\ast A^j_1.\end{array}$$
Suppose $\mf A^r=(R^r,\K^r,\M^r, H^r)$,
the size vector $\m^{ri}=(m^{ri}_1,\cdots, m^{ri}_{t^r})=\vartheta^{ri}(1,\cdots,1)$ over $\T^r$.
Take a trivially normalized quasi-basis $\{V^r_1,\ldots,V^r_{m^r}\}$ of $\K^r_1$, define a variable matrix
of size vector $\m^{ri}$ over $\T^{r}$ and a matrix
equation:
$$\begin{array}{c}\Psi_{\m^{ri}}=\sum_{X^r\in\T^r}\bar w_{_{X^r}}^{ri}\ast E_{X^r}+
\sum_j\bar v^{ri}_j\ast V^r_j,\\ \Psi_{\m^{ri}} (H^i_1+H^i_2)\equiv_{\prec(
p^i,q^i)}(H^i_1+H^i_2)\Psi_{\m^{ri}},
\end{array}$$
where $\bar v_j^{ri}=(v^{ri}_{jpq})$ is the split of $v^r_j$ of size $m_{s(v^r_j)}^{ri}\times m_{t(v^r_j)}^{ri}$,
$\bar w_{_{X^r}}^{ri}=(w_{{_{X^r}}pq}^{ri})$ is a square matrix of the size $m_{_{X^r}}^{ri}$. Furthermore, the
equation system is equivalent to:
$$\begin{array}{c}{\IF}^{ri}: \Psi_{\m^{ri}}
H^i_2\equiv_{\prec(p^i,q^i)} \Psi^0_{\m^{ri}}+H^i_2 \Psi_{\m^{ri}},\\
\Psi_{\m^{ri}}^0=H^i_1\Psi_{\m^{ri}}-\Psi_{\m^{ri}} H^i_1.\end{array} \eqno
(2.4\mbox{-}3)$$

{\bf Corollary 2.4.2}\, The matrix equations of formulae
(2.4-2) and (2.4-3) are equivalent.

\medskip

Now we give an altered Theorem
used in section 5.2-5.4, which is very easy to be proved.

\medskip

{\bf Theorem 2.4.3}\, (i) For each $0\leqslant i\leqslant s$ in the sequence (2.4-1),
there is a {\it defining system}:
$$\begin{array}{c}\mbox{\IE}^{i}: \,\,\ \Phi_{\m^{i}}
H^{i}(k)\equiv_{\prec(p^{i},q^{i})}H^{i}(k)\Phi_{\m^{i}};\quad
\Phi_{\m^i}=\sum_{X\in\T}Z_X\ast E_X +\sum_{j} Z_j\ast V_j,\end{array}\eqno
{(2.4\mbox{-}4)}$$
where $Z_X=(z_{pq}^X)_{m_{_X}\times m_{_X}}$, $\forall\, X\in \T$, and
$Z_j=(z_{rs}^j)_{m_{s(v_j)}\times m_{t(v_j)}}$ for all quasi-basis matrix $V_j$ of $\K_1$ in $\mf A$,
where $z^X_{pq},\; z^j_{pq}$ are algebraically independent variables over
$k$. $\Phi_{\m_i}$ is called a {\it variable matrix}. Then the solution space of $\IE^i$
is $\K_0^{i}\oplus \K_1^{i}$ forgotten the $(R^{i}\oplus(R^i\otimes_kR^i))$-structure.

(ii) Fix some $1<r<s$, define algebraically independent variable matrices
$Z^{ri}_{Y}$ of size $m^{ri}_{Y},\forall\,Y^r\in \T^r$, and $Z^{ri}_j$ of size $m^{ri}_{s(v^r_j)}\times
m^{ri}_{t(v^r_j)}, 1\leqslant j\leqslant m^r$. Then the equation $\IF^{ri}$ below is equivalent to $\IE^i$:
$$\begin{array}{c}{\IF}^{ri}: \Psi_{\m^{ri}}
H^i_2(k)\equiv_{\prec(p^i,q^i)} \Psi^0_{\m^{ri}}+H^i_2(k) \Psi_{\m^{ri}},\\
\Psi_{\m^{ri}}^0=H^i_1(k)\Psi_{\m^{ri}}-\Psi_{\m^{ri}} H^i_1(k),\quad
\Psi_{\m^{ri}}=\sum_{Y\in\T^r}Z^{ri}_{Y^r}\ast E_{Y^r}^r +\sum_{j=1}^{m^r}
Z^{ri}_j\ast V^r_j.\end{array} \eqno{(2.4\mbox{-}5)}$$

{\bf Corollary 2.4.4}\, $\delta(a_1^i)=0$ in ${\mf B}^i$ if and only if the equation
${\IE}^i_{(p^i,q^i)}$ is a linear combination of the equations of
${\IE}^i_{\prec(p^i,q^i)}$.

\medskip

Finally, we perform reduction procedure for the matrix bi-module
problem given in Example 1.4.5 in order to show some concrete
calculation to end the sub-section.

\medskip

{\bf Example 2.4.5}\, (i) Making an edge reduction for the first arrow $a:X\rightarrow Y$ by
$a\mapsto G^1=(1_{Z})$, we obtain an induced local bi-module problem
$\mf{A}^1$ (resp. induced bocs $\mf B^1$), with $R^1=k1_Z$;
$H^1=(1_Z)\ast A$.

(ii) Making a loop reduction for $b: Z\rightarrow Z$ by
$b\mapsto G^2=J_2(0)1_X$, we obtain
an induced local pair $(\mf{A}^2,\mf{B}^2)$ with $R^2=k1_X, H^2={{1_X\,\,\,0\,}\choose{\,\,0\,\,1_X}}\ast A
+{{\,0\,1_X}\choose{0\,\,\,0\,}}\ast B$.
Whose formal equation consists of two matrix equations:
$$\small{~_{\begin{array}{ll}  \left(\begin{array}{cc} e& v\\
0&e\end{array}\right)\left(\begin{array}{cc} c_{11}& c_{12}\\
c_{21}&c_{22}\end{array}\right)+\left(\begin{array}{cc} u^2_{11}& u^2_{12}\\
u^2_{21}&v^2_{22}\end{array}\right)\left(\begin{array}{cc} 0& 1_{X}\\
0&0\end{array}\right)\\
= \left(\begin{array}{cc} 0& 1_{X}\\
0&0\end{array}\right)\left(\begin{array}{cc} v^2_{11}& v^2_{12}\\
v^2_{21}&v^2_{22}\end{array}\right)+\left(\begin{array}{cc} c_{11}& c_{12}\\
c_{21}&c_{22}\end{array}\right)\left(\begin{array}{cc} e& v\\
0&e\end{array}\right),    \end{array}}}$$
$$\small{~_{\begin{array}{ll}  \left(\begin{array}{cc} e& v\\
0&e\end{array}\right)\left(\begin{array}{cc} d_{11}& d_{12}\\
d_{21}&d_{22}\end{array}\right)+\left(\begin{array}{cc} u^1_{11}& u^1_{12}\\
u^1_{21}&u^1_{22}\end{array}\right)\left(\begin{array}{cc} 0& 1_{X}\\
0&0\end{array}\right)+\left(\begin{array}{cc} u^2_{11}& u^2_{12}\\
u^2_{21}&u^2_{22}\end{array}\right)\left(\begin{array}{cc} 1_X& 0\\
0&1_X\end{array}\right) \\
= \left(\begin{array}{cc} 1_{X}&0 \\
0&1_X\end{array}\right)\left(\begin{array}{cc} v^2_{11}& v^2_{12}\\
v^2_{21}&v^2_{22}\end{array}\right)+\left(\begin{array}{cc} 0& 1_{X}\\
0&0\end{array}\right)\left(\begin{array}{cc} v^1_{11}& v^1_{12}\\
v^1_{21}&v^1_{22}\end{array}\right)+\left(\begin{array}{cc} d_{11}& d_{12}\\
d_{21}&d_{22}\end{array}\right)\left(\begin{array}{cc} e& v\\
0&e\end{array}\right). \end{array}}}$$
Where $(c_{pq})_{2\times 2},(d_{pq})_{2\times 2}$ are splits from $c,d$
respectively, $e=e_{_{X}}$; $v\in \C_1^2$ dual to $V={{0\,\,1_X\otimes_k1_X}\choose{0\qquad 0\quad}}
\ast E_X\in \K^2_1$, $E_X=(1_XI_{10},1_XI_{10})$.

(iii) Making a loop mutation $c_{21}\mapsto (x)$, then 3 steps of regularization, such
that $c_{22}\mapsto\emptyset, u^2_{21}= xv$;
$c_{11}\mapsto\emptyset, v^2_{21}=vx $; $c_{12}\mapsto\emptyset,
u^2_{11}=v^2_{22}$, we obtain an induced pair
$(\mf A^3,\mf B^3)$ with the differentials in $\mf B^3$:
$$
\small{\left\{\begin{array}{l} \delta(d_{21})= xv-vx\\
\delta(d_{22})=u^1_{21}+u^2_{22}-v^2_{22}-d_{21}v\\
\delta(d_{11})= u^2_{11}-v^2_{11}-v^1_{21}+vd_{21}\\
\delta(d_{12})= u^1_{11}+u^2_{12}-v^2_{12}-v^1_{22} -d_{11}v +vd_{22}.
\end{array}\right.}
$$

(iv) Note that the blocks splitting from $d_{22},d_{11},d_{12}$
will be going to $\emptyset$ by regularization for any possible reductions for $x$ and $d_{21}$.

\bigskip
\bigskip
\bigskip
\centerline{\bf  3 Classification of Minimally Wild Bocses}
\bigskip

The present section is devoted to classifying so called minimally wild
bocses, which are divided into five classes. Then we
prove the non-homogeneity for the bocses in four classes.
But the last class has been proved to be strongly homogeneous.

\bigskip
\bigskip
\noindent {\bf 3.1 Exact structures on representation categories of
bocses}
\bigskip

In this sub-section we will recall the exact structure on representation
categories of bocses. We refer to \cite{GR}  and \cite{DRSS} for the general concept of exact structure
on additive categories with Krull-Schmidt
property.

\medskip

Let $\mathfrak{B}=(\Gamma,\Omega)$ be a bocs with a layer
$L=(\Gamma';\omega;a_1,\cdots,a_n;v_1,\cdots,v_m)$. From now on
we always assume that $\mf B$ is {\it triangular on the dotted arrows},
i.e. $\delta(v_j)$ involves only $v_1, \cdots, v_{j-1}$.

The bocs $\mathfrak{B}_0=(\Gamma,\Gamma)$
is called a {\it principal bocs} of $\mf{B}$.
The representation category $R(\mf{B}_0)$ is just the
module category $\Gamma$-mod.

\medskip

{\bf Lemma 3.1.1}  \cite{O}\, Let $\mathfrak{B}=(\Gamma,\Omega)$
be a layered bocs, which is  triangular on the dotted arrows, and has a principle bocs
$\mf{B}_0$.

(i) If $\iota: M\rightarrow E$ is a morphism of $R(\mathfrak{B})$
with $\iota_0$ injective, then there exists an isomorphism $\eta$
and a commutative diagram in $R(\mathfrak{B})$, such that  the bottom row
is exact in $R(\mathfrak{B}_0)$.
Dually, if $\pi: E\rightarrow N$ is a morphism of $R(\mathfrak{B})$
with $\pi_0$ surjective, then there exists an isomorphism $\eta$ and
a commutative diagram in $R(\mathfrak{B})$, such that  the bottom row
is exact in $R(\mathfrak{B}_0)$.
$$\small{
\begin{CD}
 &  & M @>\iota>> E\\
 && @VidVV @VV\eta V\\
0@>>>M @>>\iota^{\prime}>E^{\prime}
\end{CD}}\qquad\qquad
\small{\begin{CD}
  E @>\pi>> N\\
  @V\eta VV @VVidV\\
E^{\prime}@>>\pi^{\prime}>N @>>>0
\end{CD}}
$$

(ii) If $(e):  M\stackrel{\iota}{\longrightarrow}
E\stackrel{\pi}{\longrightarrow} N$ with $\iota\pi=0$ is a pair of
composable morphisms  in $R(\mathfrak{B})$ and
$(e_0):0\longrightarrow M\stackrel{\iota_0}{\longrightarrow}
E\stackrel{\pi_0}{\longrightarrow} L\longrightarrow 0$ is exact in
the category of vector spaces, then there exists an isomorphism
$\eta$ and a commutative diagram in $R(\mathfrak{B})$:
$$\small{\begin{CD} (e) \qquad &     & M@>\iota>> E @>\pi>> N  \\
& &@VidVV @VV\eta V@VVidV \\
(e') \qquad   0 @>>> M @>>\iota^{\prime}> E^{\prime}@>>\pi^{\prime}>
N @>>>0
\end{CD}}
$$
such that $(e')$ is an exact sequence in $R(\mathfrak{B}_0)$.
Moreover, by choosing a suitable basis of $M,E',N$, we are able to obtain
$\iota'_X=(0, I)$ and $\pi'_X=(I,0)^T$ for all $X\in \T$.

\medskip

{\bf Lemma 3.1.2} Let $\mf{B}=(\Gamma,\Omega)$ be a layered
bocs, which is triangular on the dotted arrows.

(i) $\iota: M\rightarrow E$ is monic in $R(\mf B)$ if $\iota_0: M\rightarrow E$ is
injective. Dually, $\pi: E\rightarrow N$ is epic in $R(\mf B)$ if $\pi_0: E\rightarrow N$ is
surjective.

(ii)  A pair of composable morphisms $(e):
M\stackrel{\iota}{\longrightarrow} E\stackrel{\pi}{\longrightarrow}
N$ with $\iota\pi=0$ is exact in $R(\mf{B})$, if $(e_0):
0\longrightarrow M\stackrel{\iota_0}{\longrightarrow}
E\stackrel{\pi_0}{\longrightarrow} N\longrightarrow 0$ is exact as
a sequence of vector spaces.

\smallskip

{\bf Proof.} (i) If $\iota_0$ is injective, Lemma 3.1.1 (i) gives a
commutative diagram with $\iota': M\rightarrow E'$ in
$R(\mf{B}_0)$. Given any morphism $\varphi: L\rightarrow M$
with $\varphi\iota=0$, we have $\varphi\iota\eta=\varphi\iota'=0$.
Then $\varphi_0\iota'_0=0$ yields $\varphi_0=0$. And for any $v_l: X
\rightarrow Y$, $\delta(v_l)=\sum_{i,j}u_i
\otimes_{\Gamma} u_{j}$ with $u_i,u_j\in \oplus_{l'<l}\Gamma v_{l'}\Gamma$,
using induction: $0=(\varphi\iota')(v_l)=\varphi(v_l)\iota'_Y+\varphi_X\iota'(v_l)+\sum_{i,j}
\varphi(u_i)\iota'(u_j)=\varphi(v_l)\iota'_Y$, which yields $\varphi(v_l)=0$ by
the injectivity of $\iota'_Y$. Thus $\varphi=0$ and $\iota$ is monic. The second one
is proved dually.

(ii) We first prove that $\iota$ is the kernel of $\pi$.  \ding{172} If
$(e_0)$ is exact, then (i) tells that $\iota$ is monic.
\ding{173} $\iota\pi=0$. \ding{174} Lemma 3.1.1 (ii) gives a commutative diagram.
If $\varphi: L\rightarrow E $ with $\varphi\pi=0$, then
$\varphi(\eta\pi^{\prime})=0$. Let $\xi=\varphi\eta$, then
$\xi\pi^{\prime}=0$ implies that $\xi_X\pi_X^{\prime}=0$ and
$\xi(v)\pi_Y^{\prime}=0$ for any vertex $X$ and any dotted arrow $v:
X\rightarrow Y$. Let $\varphi^{\prime}:L\rightarrow M$ be given by
$\varphi^{\prime}_X\iota^{\prime}_X=\xi_X$,
$\varphi^{\prime}(v)\iota^{\prime}_Y=\xi(v)$,  then we obtain
$\varphi^{\prime}\iota^{\prime}=\xi$. Thus
$\varphi^{\prime}\iota^{\prime}\eta^{-1}=\varphi$, i.e.
$\varphi^{\prime}\iota=\varphi$. Therefore $\iota$ is a kernel of
$\pi$. It can be proved dually that $\pi$ is a cokernel of $\iota$.
The proof is finished.

\medskip

Let a layered bocs $\mathfrak{B}=(\Gamma,\Omega)$ be triangular on the dotted arrows. We define a class
${\cal E}$ of composable morphisms
in $R(\mf{B})$, such that $M\stackrel{\iota}{\longrightarrow}
E\stackrel{\pi}{\longrightarrow} L$ in ${\cal E}$, provided that
$\iota\pi=0$ and
$$0\longrightarrow M\stackrel{\iota_0}{\longrightarrow}
E\stackrel{\pi_0}{\longrightarrow} L\longrightarrow 0\eqno{(3.1\mbox{-}1)}$$ is exact as
a sequence of vector spaces. It is clear that ${\cal E}$ is closed
under isomorphisms.

\medskip

{\bf Proposition 3.1.3}\, (Theorem 4.4.1 of \cite{O} or \cite{BBP})
Suppose a layered bocs $\mf B$ is triangular on the dotted arrows,
then ${\cal E}$ defined by Formula (3.1-1) is an exact structure on
$R(\mathfrak{B})$, and $(R(\mathfrak{B}),{\cal E})$ is an exact category.

\medskip

In particular,
the bocs $\mf B$ associated to a matrix bi-module problem $\mf A=(R,\K,\M,H)$ is triangular.
In fact, the basis $\V=\{V_1,\cdots,V_m\}$ of $\K_1$ possesses a natural partial ordering:
$V_i\prec V_j$ provided their leading position
$(p_i,q_i)\prec (p_j,q_j)$. Since $\K_1\subset\mathbb N_t(R\otimes_kR)$,
$V_i\otimes_RV_j\in\sum_{V_l\prec V_i,V_j} R^{\otimes 3}\otimes_{R^{\otimes 2}}V_l$. Thus
$\dz(v_l)$ contains only $v_i\otimes_{_R}v_j$
with $v_i,v_j\prec v_l$.

\medskip

{\bf Corollary 3.1.4} (\cite{B1} and Lemma 7.1.1 of \cite{O})\,
Let $\mf{B}=(\Gamma, \Omega)$ be a layered bocs.

(i)\, For any $M\in R(\mf{B})$ with
$\dim M=\m$. If $m_X\ne 0$ for some vertex $X\in \T_1$, then $M$ is
neither projective nor injective.

(ii)\, For any positive integer $n$, there are only finitely many
iso-classes of indecomposable projectives and injectives in $R(\mf B)$
of dimension at most $n$.

\medskip

{\bf Remark 3.1.5} (\cite{BCLZ}, Def. 4.4.1 of \cite{O})\,
Let $\mf{B}=(\Gamma, \Omega)$ be a layered bocs, such that $(R(\mf{B}),
{\cal E})$ is an exact category. The almost split conflations has been defined
in a general exact category, consequently in $R(\mf B)$.

(i) An indecomposable representation $M\in R(\mf{B})$ is said to be
{\it homogeneous} if there is an almost split conflation
$M\stackrel{\iota}{\longrightarrow}E\stackrel{\pi}{\longrightarrow}M$.

(ii) The category $R(\mf{B})$ (or bocs $\mf B$) is said to be {\it homogeneous} if for
each positive integer $n$, almost all (except finitely many)
iso-classes of indecomposable representations in $R(\mf{B})$ with
size at most $n$ are homogeneous.

If $\mf B$ is of representation tame type,
then $R(\mf B)$ is homogeneous \cite{CB1}.

(iii) The category $R(\mf{B})$ (or bocs $\mf B$) is said to be {\it strongly
homogeneous} if  there exists neither projectives nor injectives,
and all indecomposable representations in $R(\mf{B})$ are
homogeneous.

If $\mf B$ is a local bocs with a layer
$(R;\omega; a; v), R=k[x,\phi(x)^{-1}]$, and the differential
$\delta(a)= xv-vx$. Then $R(\mf B)$ is strongly homogeneous and
representation wild type. In particular the induced bocs given in
Example 2.4.5 (iv) is strongly homogeneous \cite{BCLZ}.

\medskip

Note that $(R(\mathfrak{B}),{\cal E})$ may not have any almost split
conflation. For example, set quiver
$Q=\xymatrix{\cdot\ar@(ul,dl)_{a}\ar@(ur,dr)^{b}}$
and $\Gamma=kQ$. Then $R(\mf{B})$ for the principal bocs
$\mf{B}=(\Gamma,\Gamma)$ has no almost split conflations, see
\cite{V,ZL} for detail.

Recall from \cite{CB1}, let
$\mf{B}=(\Gamma,\Omega)$ be a minimal bocs. Then
for any $X\in \T_1$ with $R_X=k[x, \phi_{_X}(x)^{-1}]$, and for
any $\lambda\in k$, $\phi_{_X}(\lambda)\ne 0$, there is an almost split conflation:
$$
\begin{array}{c}
S(X,1,\lambda)\stackrel{(0\, 1)}\rightarrow S(X,2,\lambda)
\stackrel{{1}\choose{0}}\rightarrow
S(X,1,\lambda)\quad\mbox{in}\,\,\,R(\mf{B}),
\end{array}\eqno{(3.1\mbox{-}2)}
$$
where
$S(X,1,\lambda)$ (resp. $S(X,2,\lambda)$) is given by $\xymatrix{k\ar@(ur,dr)^{J_1(\lambda)}}$
(resp. $=\xymatrix{k\ar@(ur,dr)^{J_2(\lambda)}}$) at $X$, and $\{0\}$ at other vertices.

\bigskip
\bigskip
\noindent{\bf 3.2 Almost split conflations in reductions}
\bigskip

In this subsection we always assume that $\mf B=(\Gamma,\Omega)$ is a
layered bocs with triangular property on the dotted arrows, and $\mf B'=(\Gamma',\Omega')$ is
an induced bocs given by one of 8 reductions of sections 2.1-2.2.
We will study the almost split conflations during the reductions.

\medskip

{\bf Lemma 3.2.1} \cite{B1}\, Let $N'$ be an indecomposable representation in
$R(\mf B')$. If $N'$ is non-projective (resp. non-injective) in
$R(\mf B')$, then so is $\vartheta(N')$ in $R(\mf B)$.

\medskip

{\bf Lemma 3.2.2} \cite{B1}\, (i) If $\iota': M'\rightarrow E'$ is a morphism
in $R(\mf B')$ with $\vartheta(\iota'): \vartheta(M')\rightarrow
\vartheta(E')$ being a left minimal almost split inflation in
$R(\mf B)$, then so is $\iota'$.
Dually if $\pi': E'\rightarrow N'$ is a morphism in $R(\mf B')$
with $\vartheta(\pi'): \vartheta(E')\rightarrow \vartheta(N')$
being a right minimal almost split deflation in $R(\mf B)$,  then so is
$\pi'$.

(ii) If $(e'): M' \stackrel {\iota'} \rightarrow
E'\stackrel {\pi'} \rightarrow M'$ is a conflation in
$R(\mf B')$ with $\vartheta (e'): \vartheta (M')
\stackrel {\vartheta (\iota')} \rightarrow \vartheta
(E')\stackrel {\vartheta(\pi')} \rightarrow \vartheta
(M')$ being an almost split conflation in $R(\mf B)$, then
so is $(e')$.

\medskip

Let $\mf A=(R, \K, \M, H=0)$ be a matrix bi-module problem with $R$ being trivial and
$\mf B$ the corresponding bocs. Suppose that
$$(\mf A,\mf B)=(\mf A^0,\mf B^0), \cdots,(\mf A^i,\mf B^i),(\mf A^{i+1},\mf B^{i+1}),\cdots,
(\mf A^s,\mf B^s)\eqno {(3.2\mbox{-}1)}$$ is a sequence of reductions,
such that $\mf A^{i+1}$ is obtained from $\mf A^i$
in the sense of Lemma 2.3.2 for
$i=0,\cdots,s-1$; the first arrow $a^s_1\in\mf B^{s}$
is a loop at $X^s_1$ and $\delta(a_1^s)=0$.

\medskip

{\bf Theorem 3.2.3}\, With the assumption above, taken any indecomposable
$M^s\in R(\mf A^s)$ of size vector $\m^s$ with $m^s_{_{X^s_1}}=1$, $M^s(a_1^s)=(\lambda)$.
Suppose $\vartheta^{0s}(M^s)=M$ is homogeneous with an almost
split conflation $(e): M\rightarrow E\rightarrow M$ in $R(\mf A)$.
Then for $i=0,1,\cdots,s$ there exits
an almost split conflation $(e^i): M^i\rightarrow E^i\rightarrow M^i$ in $R(\mf A^i)$,
such that $\vartheta^{0i}(e^i)\simeq (e)$.

\smallskip

{\bf Proof}\, We stress that $M^s_X=k$, and according to Theorem 2.3.3 (iii), as $k$-matrices:
$$\begin{array}{c}M^s=H^s_{\m^s}(k)+\sum_{j}M^s(a^s_j)\ast A^s_j,\quad
H^s_{\m^s}(k)=\sum_{j=1}^{s-1}B^{j+1}\ast A^j_1.\end{array}$$
The theorem is obviously true for $i=0$. Suppose the assertion is valid
for some $0\leqslant i<s$, we will reach the $(i+1)$-th stage.
Set $\m^i=\vartheta^{is}(\m^s)$, Formula (3.2-2) gives
as $k$-matrices:
$$\begin{array}{c}M^i=M^s=H^i_{\m^i}(k)+B^{i+1}\ast A^i_1+
\sum_{j=2}^{n^i}M^i(a^i_j)\ast A^i_j.\end{array}\eqno {(3.2\mbox{-}2)}$$
If there exits an object
$E^i= H^i_{2\m^i}(k)+\sum_{j=1}^{m^i}E^i(a_j^i)\ast A_j^i\in R(\mf A^i)$, and
an almost split conflation $(e^i): M^i\stackrel{\iota^i}\rightarrow E^i
\stackrel{\pi^i}\rightarrow M^i$ in $R(\mf A^i)$,
such that $\vartheta^{0i}(e^i)\simeq(e)$, we will prove in (i) and (ii) below, that there exists
some isomorphism $\eta:E^i\mapsto \hat E^i$ in $R(\mf A^i)$ with
$\hat E^i(a^i_1)=B^{i+1}\oplus B^{i+1}$.

If this is the case, let $a^i_1:X\rightarrow Y$, $S_{X}$ and $S_{Y}$ are
invertible matrices determined by changing certain rows and columns of
$B^{i+1}\oplus B^{i+1}$, such that $S_X^{-1}(B^{i+1}\oplus B^{i+1})S_Y
=I_2\otimes_k B^{i+1}$, define
a matrix $S=\sum_{Z\in\T^i}P_Z\ast E_Z$ with
$S_Z=I_{m_{_Z}}$ for $Z\in\T^i\setminus \{X,Y\}$. Then as $k$-matrices:
$$\begin{array}{l}R(\mf B^i)\ni S^{-1}\hat E S=H^i_{2\m^i}(k)+
(I_2\otimes_k B^{i+1})\ast A^i_1+
\sum_{j=2}^{n^i}S_{s(a^i_j)}^{-1}\hat E(a^i_j)S_{t(a^i_j)}\ast A^i_j\\
\quad\quad\quad\quad\quad\quad\quad\quad\quad\quad\quad\stackrel{\small\mbox{define}}
=H^{i+1}_{2\m^{i+1}}(k)+\sum_{j=1}^{n^{i+1}}E^{i+1}(a^{i+1}_j)\ast A^{i+1}_j
=E^{i+1}\in R(\mf A^{i+1}).\end{array}$$
Thus we obtain a conflation $(\hat e^{i}):M^i\rightarrow
S\hat E^iS^{-1}\rightarrow M^i$ equivalent to $(e^i)$ in $R(\mf A^i)$,
and an almost split conflation $(e^{i+1})$ in $R(\mf A^{i+1})$ with
$\vartheta^{i,i+1}(e^{i+1})\simeq(\hat e^i)\simeq (e^i)$ by Lemma 3.2.2 (ii).
Consequently $\vartheta^{0,i+1}(e^{i+1})\simeq\vartheta^{0i}(e^i)\simeq(e)$.

(i)\, If $\dz(a_1^i)=v^i_1\ne 0$, after a regularization,
Set $\hat E^i=\{\hat E^i_Z\mid$ dim$(E^i_Z)=2m^i_{_Z}, Z\in\T^i\}$ be a set of vector spaces.
Let $\eta: E^i\rightarrow\hat E^i$, such that $\eta_{_Z}=I_{2m^i_{_Z}}$ for any $Z\in \T^i$,
and $\eta(v_1^i)=E^i(a_1^i), \eta(v_j)=0$ for any $j=2,\cdots,m^i$.
Let $\hat E^i=\eta^{-1}E^i\eta\in R(\mf A^i)$, then $\hat E^i(a^i_1)=E^i(a^i_1)-
\eta(v_1^i)=(0)_{2m_{_X}\times 2m_{_Y}}
=B^{i+1}\oplus B^{i+1}$ as desired.

(ii)\, If $\dz(a_1^i)=0$ in the case of edge or loop reduction, the proof is divided into three parts.

\ding{172}\, We define an object $L^s=H^s_{2\m^s}(k)+\sum_{i=1}^{n^s}L^s(a^s_i)\in R(\mf A^s)$
with $L^s(a^s_1)={{\lambda\,1}\choose{0\,\lambda}}$.
Let $\varphi^s: L^s\rightarrow M^s$ be a morphism in $R(\mf A^s)$,
such that $\varphi^s_{X^s}={I_{X^s}\choose 0},\forall\,X^s\in \T^s$,
and $\varphi^s({v^s_j})=0$
for any dotted arrow $v^s_j$ in $\mf{B}^s$. Clearly, $\varphi^s$ is
not a split epimorphism. Thus $\vartheta^{is}(\varphi^s):\vartheta^{is}(L^s)\mapsto\vartheta^{is}(M^s)=M^i$
is not a split epimorphism, since the functor $\vartheta^{is}$ is fully faithful.

\ding{173}\, Because $\vartheta^{is}(L^s)=H^i_{2\m^i}(k)+(I_2\otimes_kB^{i+1})\ast A^i_1+
\sum_{j=2}^{n^i}\vartheta^{is}(L^s)(a_i)\ast A^i_j$
by Formula (3.2-2), we are able to construct an object $L^i$ with $L^i(a^i_1)=B^{i+1}\oplus B^{i+1}$
by changing certain rows and columns in $I_2\otimes_k B^{i+1}$, then there is an isomorphism
$L^i\stackrel{\eta}\rightarrow\vartheta^{is}(L^s)$.
We obtain a lifting
$\wt \varphi^i$ of $\varphi^i=\vartheta^{is}(\varphi^s)\eta$ in $R(\mf{A}^i)$, such that
$\varphi^i=\wt \varphi^i \pi^i$, since $\pi^i: E^i\rightarrow M^i$ is
right almost split in $R(\mf{A}^i)$. The
triangle and the square below are both commutative:
$$\xymatrix {&L^i\ar[dl]_{\wt
\varphi^i}\ar[dr]^{\varphi^i}&\\E^i\ar[rr]_{\pi^i}&&
M^i}\quad\quad\quad \xymatrix {L^i_X\ar[r]^{L^i(a^i_1)}\ar[d]_{\wt
\varphi^i_X}&L^i_Y\ar[d]^{\wt \varphi^i_Y}\\
E^i_X\ar[r]_{E^i(a^i_1)}&E^i_Y}$$

\ding{174} We may assume according to Lemma 3.1.1 that in the sequence $(e^i)$,
$\iota^i_{_Z}=(0\,\, I_Z),\pi^i_{_Z}={I_Z\choose 0},\forall\,Z\in \T^i$, then
$E^i(a^i_j)={{M^i(a^i_j)\,\,\,K^i_j}\,\,\,\choose{\quad 0\quad M^i(a^i_j)}}$.
The commutative triangle
forces $\wt \varphi_Z={{I_Z\,\,C_Z}\choose {0_Z\,\, D_Z}}$ for each $Z\in \T^i$. The commutative
square yields an equality
$$\left
(\begin{array}{cc}I_X&C_X\\&D_X\end{array}\right )\left
(\begin{array}{cc}B^{i+1}&K^i_1\\& B^{i+1}\end{array}\right
)=\left (\begin{array}{cc}B^{i+1}&\\&B^{i+1}\end{array}\right
)\left (\begin{array}{cc}I_Y&C_Y\\&D_Y\end{array}\right ).$$
Let $\hat E^i=\{\hat E^i_Z\mid$ dim$(E^i_Z)=2m^i_{_Z}, Z\in\T^i\}$ be a set of vector spaces.
Define a set of isomorphisms $\eta:
E^i\rightarrow\hat E^i$, such that $\eta_{_{X}}={{I_{X}\,\,C_{X}}\choose
{\,0\,\,\,\,\,I_{X}}}$,
$\eta_{_{Y}}={{I_{Y}\,\,C_{Y}}\choose {\,0\,\,\,\,\,I_{Y}}}$,
and $\eta_{_Z}=I_{2m^i_{_Z}}$ for $Z\in \T^i\setminus\{X,Y\}$; $\eta(v_j)=0$
for any $j=1,\cdots,m^i$.
Let $\hat E^i(a^i_j)=\eta_{_{s(a^i_j)}}E^i(a^i_j)\eta_{_{t(a^i_j)}}^{-1}$ for $j=1,\cdots,n^i$, we
obtain an object $\hat E^i=\eta E\eta^{-1}\in R(\mf A^i)$ with $\hat
E^i(a^i_1)=B^{i+1}\oplus B^{i+1}$ as desired.
The proof is completed.

\medskip

Suppose $\mf B$ is a bocs, and bocs $\mf B'$ is induced from $\mf B$
by deletion of a vertex set $\T'\subset\T$.
If $M'\in R(\mf B'), M=\vartheta(M')\in R(\mf B)$ is homogeneous with
an almost split conflation $(e):M\rightarrow E\rightarrow M$, then
there exists an almost split conflation $(e')\in R(\mf B')$
with $\vartheta(e')\simeq(e)$. In fact, dim$(E_X)=2$dim$(M_X)$ for any $X\in\T$
from the definition of the exact structure (3.1-1). So that
$E_X\ne\{0\}$ if and only if $M_X\ne\{0\}$, and thus $X\notin\T'$, $E\in R(\mf B')$.

\medskip

Suppose we have the following sequence with the first part up to the $s$-th pair
is given by Formula (3.2-1):
$$(\mf A,\mf B)=(\mf A^0,\mf B^0) \cdots,
(\mf A^s,\mf B^s),(\mf A^{s+1},\mf B^{s+1})\cdots,(\mf A^i,\mf B^i),
(\mf A^{i+1},\mf B^{i+1}),\cdots,(\mf A^t,\mf B^t).\eqno {(3.2\mbox{-}3)}$$

Assume that in the sequence (3.2-3),
$\mf B^s$ is local with $\T^s=\{X\}$; $\mf B^{s+1}$ is induced from $\mf B^s$
by a loop mutation; and the reduction from $\mf B^i$ to $\mf B^{i+1}$
is given by a localization then a regularization, such that $R^{i+1}=k[x,\phi^{i+1}(x)^{-1}]$
for $s<i<t$, and $\mf B^t$ is minimal.

\medskip

{\bf Corollary 3.2.4}\, Suppose $M^t\in R(\mf B^t)$ with $M^t_{X}=k,M^t(x)=(\lambda)$ being regular,
such that $\vartheta^{0t}(M^t)=M\in R(\mf B)$ is homogeneous with an almost split conflation $(e)$. Then
there exists an almost split conflation $(e^t)$ of Formula (3.1-2) in $R(\mf B^t)$,
such that $\vartheta^{0t}(e^t_\lambda)\simeq (e)$.

\smallskip

{\bf Proof}\, Set $M^s=\vartheta^{st}(M^t)$, then $M^s(a^s_1)=(\lambda)$. Lemma 3.2.3 gives an
almost split conflation $(e^s)$ in $R(\mf B^s)$ with $\vartheta^{0s}(e^s)\simeq (e)$.
Moreover, $R(\mf A^{s+1})\simeq R(\mf A^{s})$, and for
$i\geqslant s$, $R(\mf A^{i+1})$ is equivalent to a subcategory of $R(\mf A^{i})$ consisting of the
objects $M^{i}$ with the eigenvalues of $M^i(x)$ are not the roots of $\phi^{i+1}(x)$.
The proof is finished.

\medskip

Assume that in the sequence (3.2-3), $\mf B^s$ has two vertices
$X,Y$, the first arrow $a^s_1:X\mapsto X$ with $\dz(a_1^s)=0$, $\mf B^{s+1}$ is
induced from $\mf B^s$ by a loop mutation $a^s_1\mapsto (x)$; there exists a certain index $s<l<t$,
such that $a^{s+1}_1,\cdots,a^{s+1}_{l-s}$ are either loops at $X$, or edges from
$Y$ to $X$, especially $a^{s+1}_{l-s}:Y\mapsto X$. The reduction from $\mf B^i$ to $\mf B^{i+1}$ is
given by

\ding{172} a localization then a regularization for the first loop at $X$, or a
regularization for the first edge; or
a reduction given by proposition 2.2.6 for the first edge, when $s<i<l$;

\ding{173} a reduction given by proposition 2.2.7, when $i=l$,
then $\mf B^{l+1}$ is local;

\ding{174} a localization then a regularization for the first loop, when $l<i<t$.

\medskip

{\bf Corollary 3.2.5}\,  Suppose $M^t\in R(\mf B^t)$ with $M^t_Z=k,M^t(z)=\lambda$ being regular,
such that $\vartheta^{0t}(M^t)=M\in R(\mf B)$ is homogeneous with an almost split conflation $(e)$, then
there exists an almost split conflation $(e^t_\lambda)\in R(\mf B^t)$ given by Formula (3.1-2),
such that $\vartheta^{0t}(e^t_\lambda)\simeq (e)$.

\smallskip

{\bf Proof}\, Set $M^s=\vartheta^{st}(M^t)$, then $M^s_X=k,M^s_Y=k$, and
$M^s(a^s_1)=(\lambda)$. Lemma 3.2.3 gives an
almost split conflation $(e^s):M^s\rightarrow E^s\rightarrow M^s$ in $\mf B^s$.
Since $R(\mf B^{s+1})\simeq R(\mf B^s)$, we use induction starting from $\mf B^{s+1}$.
Suppose for some $i>s$, there exists an almost split conflation
$$\begin{array}{c}(e^i):M^i\stackrel{\iota^i}\rightarrow
E^i\stackrel{\pi^i}\rightarrow M^i\in R(\mf B^i)\quad \mbox{with}\quad
M^i=\vartheta^{it}(M^t),\,\, \vartheta^{si}(e^i)\simeq(e^s), \end{array}$$
we now construct an almost split conflation $(e^{i+1})$ with $\vartheta^{i,i+1}(e^{i+1})\simeq(e^i)$.

\ding{172}, A regularization for an edge yields an equivalence; the proof of a regularization
for a loop is similar to that of Corollary 3.2.4.

Denote by $a^i_1:Y\rightarrow X$ the first edge of $\mf B^i$ with $\dz(a^i_1)=0$.
By Lemma 3.1.1 (ii), we may assume $\iota^i=(0\, 1),\pi^i={{1}\choose{0}}$.
Proposition 2.2.6 tells $M^i(a^i_1)=(0)$, therefore $E^i(a^i_1)={{0\,b}\choose{0\, 0}}$.
If $0\ne b\in k$, define $L\in R(\mf A^i)$ with $L_X=k^2=L_Y,
L(x)=J_2(\lambda), L(a^i_j)=0,\forall j$, and a morphism $g: L\rightarrow M^i$
with $g_X={{1}\choose{0}}=g_Y$, $g(v)=0$ for all dotted arrows.
Then $g$ is not a split epimorphism, there exists a lifting $\tilde g:L\rightarrow E^{i}$
with $\tilde g_X={{1\,  c}\choose{0\,d}}, \tilde g_Y={{1\,  c'}\choose{0\,d'}}$. Since $\tilde g$
is a morphism, $L^i(a^i_1)\tilde g_{X}=\tilde g_{Y}E^i(a^i_1)$, which yields
$0={{0\, b}\choose{0\,0}}$ a contradiction.
Therefore $b=0,E^i(a^i_1)={{0\,0}\choose{0\,0}}$. Set $E^{i+1}\in R(\mf B^{i+1})$ with
$E^{i+1}(x)=J_2(\lambda), E^{i+1}(a^{i+1}_{j-1})=E^i(a^i_j)$
for $j=2,\cdots,n^i$, then
$\vartheta^{i,i+1}(E^{i+1})=E^i$.

\ding{173} Proposition 2.2.7 tells $M^i(a^i_1)=(1)$, similar to proof \ding{172},
we may assume $E^i(a^i_1)={{1\, b}\choose{0\,1}}\simeq{{1\, 0}\choose{0\,1}}$.
Therefore there is a $E^{i+1}\in R(\mf B^{i+1})$ with $E^{i+1}(a^{i+1}_{j-1})=E^i(a^i_j)$
for $j=2,\cdots,n^i$, such that $\vartheta^{i,i+1}(E^{i+1})\simeq E^i$.

\ding{174}\, Similar to the proof of Corollary 3.2.4, the proof is completed.

\medskip

{\bf Lemma 3.2.6}\, (i)\, Suppose that $f(x,y)\in k[x,y]$ with $f(\lambda,
\mu)\ne 0$. Let $W_{\lambda}, W_{\mu}$ be Weyr matrices of size
$m,n$ and eigenvalues $\lambda,\mu$ respectively, and
$V=(v_{ij})_{m\times n}$ with $v_{ij}$ being $k$-linearly
independent. Let $f(W_\lambda,W_\mu)V= (u_{ij})_{m\times n}$. Then
$u_{ij}$ are also $k$-linearly independent for $1\leqslant
i\leqslant m, 1\leqslant j\leqslant n$.

(ii)\, Let $\mf B$ be a bocs with $R=R_X\times R_Y$, where
$R_X=k[x,\phi_X(x)^{-1}]$, $R_Y=k[y,\phi_Y(y)^{-1}]$,
and $a_i:X\mapsto Y$. Define {\it $\dz^0(a_i)$ to be a part
of $\delta(a_i)$} without the terms involving  $a_j, \forall
j<i$. It is possible that $X=Y$, in this case $x$ stands for the
multiplication from left and $y$ from right.
$$\left\{\begin{array}{l}\dz^0(a_1)=f_{11}(x,y)v_{1}\\
\dz^0(a_2)=f_{21}(x,y)v_1+f_{22}(x,y)v_2,\\
\quad \quad\cdots\quad \cdots\\
\dz^0(a_n)=f_{n1}(x,y)v_1+f_{n2}(x,y)v_2+\cdots
f_{nn}(x,y)v_n,\end{array}\right.$$ where $f_{ii}(x,y)\in R_X\times
R_Y$ are invertible for $i=1,2,\cdots,n$. Set $x\mapsto W_X$ of size $m$
with eigenvalues $\lambda$ satisfying $\phi_{_X}(\lambda)\ne 0$,
$y\mapsto W_Y$ of size $n$ with eigenvalues $\mu$ satisfying $\phi_{_Y}(\mu)\ne 0$.
Then the solid arrows splitting from $a_1,\cdots, a_n$ are all going
to $\emptyset$ by regularization in further reductions.

\smallskip

{\bf Proof}\, (i)\, Set $f(x,y)=\sum_{i,j\geqslant 0}\alpha_{ij}x^iy^j$,
then $f(W_\lambda,W_\mu)V=\sum_{i,j\geqslant
0}\alpha_{ij}W_\lambda^iVW_\mu^j$, we have
$u_{ij}=f(\lambda,\mu)v_{ij}+
\sum_{(i',j')\succ(i,j)}f_{i'j'}(\lambda,\mu)v_{i'j'}$ with
$f_{i'j'}(x,y)\in k[x,y]$. The conclusion follows by induction on
the ordered index set $\{(i,j)\mid 1\leqslant i\leqslant m,
1\leqslant j\leqslant n\}$.
(ii) follows by (i) inductively.

\bigskip
\bigskip
\noindent {\bf 3.3\, Minimal wild bocses}
\bigskip

In this sub-section we will define five classes of minimal wild bocses
in order to prove the main theorem.
Our classification releases on the well-known Drozd's wild
configurations, which is a refinement at some last reduction steps of those.

\medskip

{\bf Proposition 3.3.1} (\cite{CB1},Prop.3.10)\, Let
$\mathfrak{B}=(\Gamma, \Omega)$ be a bocs of representation wild type
with a layer $L=(R;\omega;a_1,\ldots,a_n;v_1,\ldots,v_m)$.
We are bound to meet one of the following configurations at
some stage of reductions:

{\bf Case 1}\, $X\in \T_1,\, Y\in \T_0$ (or dual), $\dz(a_1)=0$.

{\bf Case 2}\, $X, Y\in \T_1$ (possibly $X=Y$), $\dz(a_1)=f(x,y)v_1$ with
$f(x,y)\in k[x,y,(\phi_X(x)\phi_Y(y))^{-1}]$ non-invertible.

\medskip

We first fix some notations. Let $\mf B$ be a bocs with dotted arrows
$v_1,\cdots,v_m$. Consider a vector space $\mathcal S$ over $k(x,y;z)$,
the fractional field of the polynomial ring $k[x,y;z]$
of three indeterminates. Suppose there is a linear combination:
$$\begin{array}{c}G=f_1(x,y;z)v_1
+\cdots+f_m(x,y;z)v_m,\quad f_i(x,y;z)\in k[x,y;z].\end{array}\eqno {(3.3\mbox{-}1)}$$
Let $h(x,y;z)$ be the greatest common factor of $f_1,\cdots,f_m$, then $f_1/h,$ $\cdots,$ $f_m/h$
are co-prime. There exists some $s_i(x,y;z)
\in k[x,y;z]$ for $i=1,\cdots,m$, such that $\sum_{i=1}^ms_i(g_i/h)=c(x,y)\in k[x,y]$.
Then $G=h\sum_{i=1}^m(g_i/h)v_j$.
Since  $S=k[x,y;z,c(x,y)^{-1}]$ is a Hermite ring, there exists
some invertible $F(x,y;z)\in \IM_{m}(S)$ with the first column $(f_1/h,\cdots,f_m/h)$.
Then we make a base change of the form $(w_1,\cdots,w_m)=(v_1,\cdots,v_m)F$,
$G=h(x,y;z)w_1$.

Let $f(x,\bar x)\in k[x,\bar x]$, in the form $f(x,\bar x)v$, $x$ stands for
the left multiple and $\bar x$ for the right.

\medskip

{\bf Classification 3.3.2}\, Let a wild bocs $\mf B^0$ be
given by  Proposition 3.3.1.
Then we are bound to meet a bocs $\mf B$ with a
layer $L=(R;\omega;a_1,\ldots,a_n;v_1,\ldots,v_m)$ in one of the five classes at
some stage of reductions, which are called {\it minimally wild bocses}.

\smallskip

Suppose the bocs $\mf B$ has two vertices $\T=\{X,Y\}$, such that the
induced local bocs $\mf B_X$ is tame infinite with $R_X=k[x,\phi_{_X}(x)^{-1}]$.

\smallskip

{\bf MW1}\, $\mf B_Y$ is finite with $R_Y=k1_Y$, $\dz(a_1)=0$:
$$\xymatrix {X\ar[rr]^{a_1}\ar@(ul, dl)[]_{x}&&Y}.$$

\medskip

{\bf MW2}\, $\mf B_Y$ is tame infinite with $R_Y=k[y,\phi_{_Y}(y)^{-1}]$,
$\dz(a_1)=f(x,y)v_1$, such that $f(x,y)\in k[x,y,
\phi_{_X}(x)^{-1}\phi_{_Y}(y)^{-1}]$ is non-invertible:
$$\xymatrix {X\ar@(ul, dl)[]_{x}\ar[rr]^{a_1}&& Y\ar@(ur,
dr)[]^{y}}.$$

Suppose now we have a local bocs $\mf B$ with $R=k[x,\phi(x)^{-1}]$:
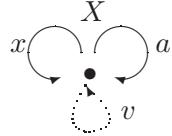
\begin{center}
\unitlength=0.5pt
\begin{picture}(60,60)
\put(-10,30){\oval(40,40)[t]} \put(-10,30){\oval(40,40)[bl]}
\put(-10,10){\vector(3,1){5}} \put(40,30){\oval(40,40)[t]}
\put(40,30){\oval(40,40)[br]} \put(40,10){\vector(-3,1){5}}

\put(11,9){$\bullet$} \qbezier[10](13,5)(-8,-28)(16,-30)
\qbezier[10](13,5)(50,-28)(13,-30) \put(17,-1){\vector(-1,2){3}}

\put(9,54){\makebox{$X$}}\put(40,-20){$v$}
\put(-44,30){\makebox{$x$}} \put(66,30){\makebox{$a$}}
\end{picture}
\end{center}\vskip 3mm

{\bf MW3}\, The differential $\dz^0$ of the solid arrows of $\mf B$ are:
$$\begin{array}{c}
\left\{
\begin{array}{ccl}
\dz^0(a_{1}) &=& f_{11}(x,\bar x)w_{1},\\
\cdots   &  & \qquad\cdots \\
\dz^0(a_{n}) &=& f_{n1}(x,\bar x)w_{1}
+\cdots+f_{nn}(x,\bar x)w_{n},
\end{array} \right.
\end{array}\eqno {(3.3\mbox{-}2)}$$
where $w_1,\cdots,w_n$ are given by base changes,
$f_{ii}(x,x)\in k[x,\phi(x)^{-1}]$ are all invertible;
$f_{11}(x,\bar x)\in k[x,\bar x,\phi(x)^{-1}\phi(\bar x)^{-1}]$ is non-invertible.

\smallskip

Now suppose there exists some $1\leqslant n_1\leqslant n$, such that:
\begin{equation*}
\left\{
\begin{array}{cllll}
\dz^0(a_{1}) &=& f_{11}(x,\bar x)w_{1},\\
\cdots   &  & \qquad\cdots &\\
\dz^0(a_{n_1-1}) &=& f_{n_1-1,1}(x,\bar x)w_{1}
&+\cdots&+f_{n_1-1,1}(x,\bar x)w_{n_1-1},\\
\dz^0(a_{n_1}) &=& f_{n_1,1}(x,\bar x)w_{1}
&+\cdots& +f_{n_1,n_1-1}(x,\bar x)w_{n_1-1}+f_{n_1,n_1}(x,\bar x)\bar w,
\end{array} \right. \eqno {(3.3\mbox{-}3)}
\end{equation*}
where $f_{ii}(x,x)\in k[x,\phi(x)^{-1}],1\leqslant i<n_1,$ are invertible;
$\bar w=0$, or $\bar w\ne 0$ but $f_{n_1,n_1}(x,x)=0$.
Denote by $x_1$ the solid arrow $a_{n_1}$, there exist a polynomial $\psi(x,x_1)$
being divided by $\phi(x)$. Write $\delta^1$ the part of differential $\delta$
by deleting all the terms involving some solid arrow except $x,x_1$, and the further
unraveling on $x$ is restricted to $x\mapsto (\lambda)$ for $\psi(\lambda,x_1)\ne 0$,
such that
$$\left\{
\begin{array}{llll}
\dz^1(a_{n_1+1}) &=K_{n_1+1}&+f_{n_1+1,n_1+1}(x,x_1,\bar x_1)w_{n_1+1},&\\
&\cdots \cdots&\cdots \\
\dz^1(a_n) &=K_n&+f_{n,n_1+1}(x,x_1,\bar x_1)w_{n_1+1}&+ \cdots+
f_{nn}(x,x_1,\bar x_1)w_{n},
\end{array} \right.\eqno {(3.3\mbox{-}4)}$$
with $K_{i}=\sum_{j=1}^{n_1-1}f_{ij}(x,x_1,\bar x_1)w_j$,
where $w_{n_1+1},\cdots,w_n$ are given by the base changes (3.3-1) step
by step, $f_{ii}(x,x_1,\bar x_1)\in k[x,x_1,\bar x_1,
\phi(x,x_1)^{-1}\phi(x,\bar x_1)^{-1}]$
are all invertible for $n_1<i\leqslant n$.

\smallskip

{\bf MW4}\, $\bar w=0$, or $\bar w\ne 0$ and $(x-\bar x)^2\mid f_{n_1n_1}(x,\bar x)$ in Formula (3.3-2).

\smallskip

{\bf MW5}\, $\bar w\ne 0$ and $(x-\bar x)^2\nmid f_{n_1n_1}(x,\bar x)$ in Formula (3.3-2).

\medskip

The proof of the classification depends on the classification of local bocses, whose proof is
based on Formulae (3.3-2)-(3.3-9) and two Lemmas below.

Let $\mf B$ be a local bocs having a layer
$L=(R;\omega; a_1, \cdots, a_n; v_1,\cdots,v_m)$. If $R=k1_X$ is trivial,
then the differentials of the solid arrows have two possibilities. First:
\begin{equation*}
\left\{
\begin{array}{ccl}
\dz^0(a_1) &=& f_{11}w_1,\\
\cdots   &  & \qquad\cdots \\
\dz^0(a_{n}) &=& f_{n1}w_1 +\cdots+f_{nn}w_{n},
\end{array} \right. \eqno {(3.3\mbox{-}5)}
\end{equation*}
with $f_{ij}\in k,h_{ii}\ne 0$ for $1\leqslant i\leqslant n$,
and $w_1,\cdots,w_n$ given by base changes. Second, there exists
some $1\leqslant n_0\leqslant n$ such that:
\begin{equation*}
\left\{
\begin{array}{ccl}
\dz^0(a_1) &=& f_{11}w_1,\\
\cdots   &  & \qquad\cdots \\
\dz^0(a_{n_0-1}) &=& f_{n_0-1,1}w_1 +\cdots +f_{n_0-1,n_0-1}w_{n_0-1},\\
\dz^0(a_{n_0}) &=& f_{n_01}w_1 \quad +\cdots +f_{n_0,n_0-1}w_{n_0-1},
\end{array} \right. \eqno {(3.3\mbox{-}6)}
\end{equation*}
with $f_{ij}\in k,f_{ii}\ne 0$ for $1\leqslant i<n_0$.
Set $a_i\mapsto \emptyset,i=1,\cdots, n_0-1$
by a series of regularization, then $a_{n_0}\mapsto (x)$ by a loop mutation,
we obtain an induced local bocs $\mf B'$.

\medskip

Without loss of generality, we may still denote $\mf B'$ by $\mf B$ with the layer $L$, but $R=k[x]$.
The differentials $\dz^0$
have again two possibilities. First one is given by Formula (3.3-2),
such that $f_{ii}(x,x)\ne 0$ for $i=1,\cdots,n$. Define a polynomial:
$$\begin{array}{c}\phi(x)=\prod_{i=1}^{n}c_i(x)f_{ii}(x,x)\end{array}\eqno {(3.3\mbox{-}7)}$$
with $c_i(x)$ appearing at the localization in order to do a base
change before the $i$-th step of regularization, thus $f_{ii}(x,x)$
are invertible in $k[x,\phi(x)^{-1}]$.

\medskip

{\bf Lemma 3.3.3}\, Let $\mf B$ be a bocs given by Formula (3.3-2) with a polynomial $\phi(x)$
(3.3-7). There exist two cases:

(i)\, $f_{ii}(x,\bar x)\in k[x,\bar x,\phi(x)^{-1}\phi(\bar x)^{-1}]$
are all invertible for $1\leqslant i\leqslant n$;

(ii)\, There exists some minimal $1\leqslant s\leqslant n$, such that
$f_{ss}(x,\bar x)\in k[x,\bar x,\phi(x)^{-1}\phi(\bar x)^{-1}]$
is non-invertible.

\medskip

The second possibility of the differential $\dz^0$ in the case $R=k[x$] is given by Formula (3.3-3)
for some fixed $1\leqslant n_1\leqslant n$,
where $f_{ii}(x,x)\ne 0$ for $1\leqslant i<n_1$; $\bar w=0$, or
$\bar w\ne 0$ but $f_{n_1,n_1}(x,x)=0$. Define
$$\phi(x)=\left\{\begin{array}{ll}\prod_{i=1}^{n_1-1}c_i(x)f_{ii}(x,x),&\bar w=0;\\
c_{n_1}(x)\prod_{i=1}^{n_1-1}c_i(x)f_{ii}(x,x),&\bar w\ne 0,\end{array}\right.\eqno {(3.3\mbox{-}8)}$$
then $f_{ii}(x,x),1\leqslant i<n_1,$ are invertible in $k[x,\phi(x)^{-1}]$.

There are two possibilities in the further reductions for the third time. First
possibility is given by Formula (3.3-4),
such that $f_{ii}(x,x_1,x_1)\ne 0$ for $n_1<i\leqslant n$. There is a sequence
of localizations given by polynomials $c_{i}(x,x_1)$
in order to do base changes before each regularizations for $i=n_1+1,\cdots,n$.
Define a polynomial
$$\begin{array}{c}\psi(x,x_1)=\phi(x)
\prod_{i=n_1+1}^n c_i(x,x_1)f_{ii}(x,x_1,x_1).\end{array}\eqno {(3.3\mbox{-}9)}$$

\medskip

{\bf Lemma 3.3.4}\, Let the bocs $\mf B$ be given by Formulae (3.3-3)-(3.3-4)
with polynomials $\phi(x)$ in (3.3-8), and $\psi(x,x_1)$ in (3.3-9). We obtain two cases.

(i)\, There exists some $\lambda\in k$ with $\psi(\lambda,x_1)\ne0$, and
a minimal $n_1+1\leqslant s\leqslant n$, such that
$f_{ss}(\lambda,x_1,\bar x_1)\in k[x_1,\bar x_1,
\psi(\lambda,x_1)^{-1}\psi(\lambda,\bar x_1)^{-1}]$
being non-invertible, i.e.,
after making a unraveling $x\mapsto (\lambda)$, and then a series of
regularization $a_i\mapsto\emptyset,w_i=0$
for $i=1,\cdots,n_1-1$, the induced local bocs $\mf B_{(\lambda)}$
with $R_{(\lambda)}=k[x_1,\psi(\lambda,x_1)^{-1}]$ satisfies Lemma 3.3.3 (ii).

(ii)\, For any $\lambda\in k$ with $\psi(\lambda,x_1)\ne0$,
$f_{ii}(\lambda,x_1,\bar x_1)\in k[x_1,\bar x_1,
\psi(\lambda,x_1)^{-1}\psi(\lambda,\bar x_1)^{-1}]$
are invertible for $n_1<i\leqslant n$, i.e.,
the induced bocs $\mf B_{(\lambda)}$
with $R_{(\lambda)}=k[x_1,\psi(\lambda,x_1)^{-1}]$ has 3.3.3 (i).

(iii)\, The case (ii) is equivalent to $f_{ii}(x,x_1,\bar x_1)\in k[x,x_1,\bar x_1,
\phi(x,x_1)^{-1}\phi(x,\bar x_1)^{-1}]$
being invertible for all $n_1<i\leqslant n$.

\smallskip

{\bf Proof}\, (ii)$\Longrightarrow$(iii)\, If there exists some $n_1< s\leqslant n$ with
$f_{ss}(x,x_1,\bar x_1)$ non-invertible, then $f_{ss}$ contains a non-trivial factor
$g(x,x_1,\bar x_1)$ co-prime to $\psi(x,x_1)\psi(x,\bar x_1)$. Consider the
variety $V=\{(\alpha,\beta,\gamma)\in k^3\mid g(\alpha,\beta,\gamma)=0,
\psi(\alpha,\beta)\phi(\alpha,\gamma)=0\}$.
Since dim$(V)\leqslant 1$, there exists a co-finite subset $\mathscr L\subset k$, such that
$\forall\,\lambda\in\mathscr L$, the plane $x=\lambda$
of $k^3$ intersects $V$ at only finitely many points. Thus $g(\lambda,x,\bar x_1)$
and $\psi(\lambda,x_1)\psi(\lambda,\bar x_1)$ are co-primes, consequently $g(\lambda,x,\bar x_1)$,
thus $f_{ss}(\lambda,x_1,\bar x_1)\in k[x_1,\bar x_1,
\psi(\lambda,x_1)^{-1}\psi(\lambda,\bar x_1)^{-1}]$
is not invertible.

(iii)$\Longrightarrow$(ii)\, If $f_{ii}(x,x_1,\bar x_1)\in k[x,x_1,\bar x_1,
\psi(x,x_1)^{-1}\psi(x,\bar x_1)^{-1}]$ is invertible, then for any
$\lambda\in k,\psi(\lambda,x_1)\ne 0$, $f_{ii}(\lambda,x_1,\bar x_1)\in k[x_1,\bar x_1,
\psi(\lambda,x_1)^{-1}\psi(\lambda,\bar x_1)^{-1}]$ is
invertible. The proof is finished.

\medskip

The second possibility of $\dz^1$ is that there is some $n_2$ with $n_1<n_2\leqslant n$,
such that
\begin{equation*}
\small\left\{
\begin{array}{l}
\dz^1(a_{n_1+1})= K_{n_1+1}+f_{n_1+1,n_1+1}(x,x_1,\bar x_1)w_{n_1+1},\\
\cdots \cdots \\
\dz^1(a_{n_2-1}) =K_{n_2-1}+f_{n_2-1,n_1+1}(x,x_1,\bar x_1)w_{n_1+1}
+\cdots+f_{n_2-1,n_2-1}(x,x_1,\bar x_1)w_{n_2-1},\\
\dz^1(a_{n_2})\quad =K_{n_2}\quad+f_{n_2,n_1+1}(x,x_1,\bar x_1)w_{n_1+1}
+\cdots +f_{n_2,n_2-1}(x,x_1,\bar x_1)w_{n_2-1}+f_{n_2,n_2}(x,x_1,\bar x_1)\bar w'
\end{array} \right. \eqno {(3.3\mbox{-}10)}
\end{equation*}
with $K_{i}=\sum_{j=1}^{n_1-1}f_{ij}(x,x_1,\bar x_1)w_j$ for $i=n_1+1,\cdots,n_2$,
where $f_{ii}(x,x_1,x_1)\ne 0$, for $n_1< i<n_2$; $\bar w'=0$, or
$\bar w'\ne 0$ but $f_{n_2,n_2}(x,x_1,x_1)=0$. Define a polynomial
$$\psi_1(x,x_1)=\left\{\begin{array}{ll}\phi(x)
\prod_{i=n_1+1}^{n_2-1}c_i(x,x_1)f_{ii}(x,x_1,x_1), &\mbox{if}\,\,\bar w'=0;\\
c_{n_2}(x,x_1)\phi(x)
\prod_{i=n_1+1}^{n_2-1}c_i(x,x_1)f_{ii}(x,x_1,x_1), &\mbox{if}\,\,
\bar w'\ne 0.\end{array}\right.\eqno {(3.3\mbox{-}11)}$$
$f_{ii}(x,x_1,x_1)$ are invertible in $k[x,x_1,\psi_1(x,x_1)^{-1}]$.

\medskip

Suppose we meet a bocs $\mf B$, whose differential is given by Formula (3.3-3) and (3.3-10),
with a polynomial $\psi_1(x,x_1)$ below (3.3-11).
Fix any $\lambda_0\in k$
with $\psi_1(\lambda_0,x_1)\ne 0$, there is an induced
bocs $\mf B_{(\lambda_0)}$ with $R_{(\lambda_0)}=k[x_1,\psi_1(\lambda_0,x_1)^{-1}]$
given by an unraveling $x\mapsto(\lambda_0)$, and then
a series of regularization $a_i\mapsto\emptyset,w_i=0$ for $i=1,\cdots,n_1-1$.
We obtain three cases:

\ding{172} $\mf B_{(\lambda_0)}$ is in the case of
Lemma 3.3.4 (i), then there exists some $\lambda_1$,
after sending $x_1\mapsto (\lambda_1)$ and a series of regularization, the induced bocs
$\mf B_{(\lambda_0,\lambda_1)}$ satisfies Lemma 3.3.3 (ii);

\ding{173} $\mf B_{(\lambda_0)}$ is in the case of Lemma 3.3.4 (iii);

\ding{174} $\mf B_{(\lambda_0)}$ is in the case
of Formulae (3.3-3) and (3.3-10).

\noindent In the case \ding{174}, we repeat the above procedure once again
for $\mf B_{(\lambda_0)}$. By induction on the indices of the finitely many solid arrows,
we finally reach case \ding{172} or \ding{173}.

\medskip

{\bf Classification 3.3.5}\, Let $\mf B$ be a local bocs with $R$ trivial, there exists four cases:

(i)  $\mf B$ satisfies formula (3.3-5).

(ii)  $\mf B$ has an induced bocs $\mf B'$ satisfying Lemma 3.3.3 (i).

(iii) $\mf B$ has an induced local bocs
$\mf B_{(\lambda_0,\lambda_1,\cdots,\lambda_{l})}$ for some $l<n$ satisfying Lemma 3.3.3 (ii).

(iv) $\mf B$ has an induced local bocs
$\mf B_{(\lambda_0,\lambda_1,\cdots,\lambda_{l-1})}$ for some $l<n$, satisfying Lemma 3.3.4 (ii).

\medskip

{\bf Proof of 3.3.2}\, (i)\, Suppose we meet a two-point wild bocs,
if $\mf B_X$ or $\mf B_Y$ is in the case of Classification 3.3.5
(iii) or (iv), we may consider the wild induced local bocs given by
deletion. Therefore we assume that one of them satisfies Formula (3.3-5) and
another satisfies Lemma 3.3.1 (i), or both are in the case of Lemma 3.3.3 (i).
MW1 or MW2 follows by Proposition 3.3.1.

(ii)\, If we meet a local wild bocs in the case of Classification 3.3.5 (iii),
then there is an induced bocs satisfying Lemma 3.3.3 (ii), we reach MW3.

(iii)\, If we meet a local wild bocs in the case of 3.3.5 (iv),
then there is an induced bocs satisfying Lemma 3.3.4 (ii), we reach MW4 or MW5.
The classification is completed.

\bigskip
\bigskip
\noindent{\bf 3.4 Non-homogeneity of MW1-4}
\bigskip

Through out the sub-section, let $\mf A^0=(R^0,\K^0,\M^0,H^0=0)$ be
a matrix bi-module problem with trivial $R^0$, and associated bocs $\mf B^0$.

\medskip

{\bf Proposition 3.4.1} \cite{B1}\, If $\mf B^0$ has an induced
bocs $\mf{B}$ in the case of MW1,
then $\mf B^0$ is non-homogeneous.

\smallskip

{\bf Proof}\,  (i)\, Let $\mf{B}_X$ be the induced local bocs,
$\vartheta_1: R(\mf{B}_X)\rightarrow R(\mf{B})$,
$\vartheta_2:R(\mf{B})\rightarrow R(\mf{B}^0)$ be two induced functors,
$\vartheta=\vartheta_2\vartheta_1$.
Set $\mathscr L_X=k\setminus\{$ the roots of $\phi_{_X}(x)\}$, for
any $\lambda\in \mathscr L_X$ define a
representation $S'_{\lambda}\in R(\mf{B}_X)$ given by
$(S'_{\lambda})_X=k$, $S'_{\lambda}(x)=(\lambda)$.
If $\mf B^0$ is homogeneous, then there is a co-finite subset
$\mathscr L\subseteq \mathscr L_X$ such that $\{\vartheta(S'_{\lambda})\in
R(\mf{B}^0)\mid \lambda\in\mathscr L\}$ is a family of pairwise
non-isomorphic homogeneous objects of
$R(\mf{B})$. By Corollary 3.2.4, there is an almost split conflation
$(e'_{\lambda}):S'_{\lambda}\stackrel{\iota'}{\longrightarrow}
E'_{\lambda}\stackrel{\pi'}{\longrightarrow}
S'_{\lambda}$ in $R(\mf B_X')$ with $E'(x)=J_2(\lambda)$ and $\vartheta(e'_\lambda)$
is an almost split conflation in $R(\mf B^0)$.
Fix any $\lambda\in \mathscr L^0$, and the conflation
$(e_{\lambda})=\vartheta_1(e'_\lambda)$ in $R(\mf{B})$:
$$
S_{\lambda}\stackrel{\iota}{\longrightarrow}
E_{\lambda}\stackrel{\pi}{\longrightarrow}
S_{\lambda},\,\,\mbox{with}\,\, \left\{\begin{array}{llll}
(S_\lambda)_X=k, &(S_\lambda)_Y=0,&S_\lambda(x)=(\lambda),&\mbox{others zero};\\
(E_\lambda)_X=k^2,&(E_\lambda)_Y=0,&E_\lambda(x)=J_2(\lambda),
&\mbox{others zero}.\end{array}\right.
$$
Since $\vartheta_2(e_\lambda)=\vartheta(e_\lambda')$, $(e_\lambda)$
is almost split by Lemma 3.2.2 (ii).

(ii)\, Define a representation $L\in R(\mf{B})$ given by $L_X=L_Y=k$,
$L(x)=(\lambda)$ and $L(a_1)=(1)$. Let $g: L\rightarrow S_{\lz}$ be a morphism
with $g_{_X}=(1)$, $g_{_Y}=(0)$ and
$g(v)=0$ for all dotted arrow $v$'s. We assert that $g$ is not a retraction.
Otherwise, if there is a morphism $h:S_{\lz}\rightarrow L$ such
that $hg=id_{S_{\lz}}$, then $h_X=(1)$ and $h_Y=(0)$.
But $h$ is a morphism implies that
$(1)(1)=h_XL(a)=S_{\lz}(a)h_Y =(0)(0)$, a contradiction.

(iii)\, There exists a lifting $\tilde{g}: L\rightarrow E_{\lambda}$ with
$\tilde{g}\pi=g$. If $\tilde{g}_{_{X}}=(a,b)$, then
$\tilde{g}_{_{X}}\pi_{_{X}}=g_{_X}$ yields $(a,b){1\choose 0}=(1),a=1$.
$\tilde{g}$ being a morphism from $L
$ to $E_{\lz}$ implies that $\tilde{g}_{_X} E_{\lz}(x)=L
(x)\tilde{g}_{_X}$, i.e., $(1,b){{\lz\,\, 1}\choose {0\,\,
\lz}}=(\lz)(1,b)$, $(\lambda,1+b\lambda)=(\lambda,\lambda b)$, a
contradiction. Thus $\mf{B}^0$ is not homogeneous, the proof is finished .

\medskip

{\bf Proposition 3.4.2} \cite{B1}\, If $\mf B^0$ has an induced
bocs $\mf{B}$ in the case of MW2, then $\mf B^0$ is
non-homogeneous.

\smallskip

{\bf Proof}\, Since
$f(x,y)\in k[x,y,\phi_{_X}(x)^{-1}\phi_{_Y}(y)^{-1}]$
is non-invertible, after dividing the dotted arrows $v_j$ by some powers
of $\phi_{_X}(x)$ and $\phi_{_Y}(y)$, we may
assume that $f(x,y)=l(x,y)\alpha(x)\beta(y)$, where $l(x,y)\in
k[x,y]$ and $\alpha(x)$ (resp. $\beta(y)$) is a product of some
factors of $\phi_{_X}(x)$ (resp. $\phi_{_Y}(y)$),
such that $(l(x,y), \phi_{_X}(x)\phi_{_Y}(y))=1$. By
Bezout's theorem there is an infinite set
$$\mathscr L'=\{(\lambda,\mu)\in k\times k\mid l(\lz,\mu)=0,
\phi_{_X}(\lz)\phi_{_Y}(\mu)\ne 0\}.$$ Without loss of generality we may
assume that $\mathscr L_X=\{\lambda\in k\mid
(\lambda,\mu)\in\mathscr L'\}$ is an infinite set.

(i)\, as the same as the proof (i) of Theorem 3.4.1.

(ii)\, Define a representation $L\in R(\mf{B})$ given by $L_X=k=L_Y$,
$L(x)=(\lambda),\lambda\in\mathscr L_X; L(y)=(\mu),(\lambda,\mu)\in \mathscr L'$, and $L(a_1)=(1)$.
Let $g: L\rightarrow S_{\lz}$ be a morphism in
$R(\mf{B})$ with $g_X=(1)$, $g_{_Y}=(0)$ and
$g(v)=0$ for all dotted arrow $v$'s, then $g$ is not a retraction.
Otherwise, if there is a morphism $h:S_{\lz}\rightarrow L$ such
that $hg=id_{S_{\lz}}$, then $h_{_X}=(1)$ and $h_{_Y}=(0)$.
But $h$ is a morphism implies that
$0-1=S_{\lz}(a)h_{_Y}-h_{_X}L (a)=h(\dz(a))
=f(\lambda,\mu)h(v)=0$, a contradiction.

(iii)\, There exists a lifting $\tilde{g}: L\rightarrow E_{\lambda}$
with $\tilde{g}\pi=g$. A contradiction appears as the same as
3.4.1 (iii), which finishes the proof.

\medskip

{\bf Proposition 3.4.3} \cite{B1}\, If $\mf B^0$ has an induced bocs $\mf{B}$
in the case of MW3, then $\mf B^0$ is
non-homogeneous.

\smallskip

{\bf Proof}\, Let $f_{11}(x,\bar x)=l(x,\bar x)\alpha(x)\beta(\bar x)$, where
$l(x,\bar x)$, $\alpha(x)$, $\beta(\bar x)$ and $\mathscr L'$ are given in
the beginning of the proof of Proposition 3.4.2. Suppose $\mathscr L'=\{\lambda\in k
\mid (\lambda,\mu)\in\mathscr L'\subseteq k\}$ is an infinite
set, and $\vartheta:R(\mf B)\rightarrow R(\mf B^0)$
the induced functor.

Taken $S_{\lz}\in R(\mf{B})$ given by $(S_{\lz})_X=k$,
$S_{\lz}(x)=(\lambda),\forall\,\lz\in\mathscr L'$, then
$S_{\lz}(a_i)=(0)$ for $1\leqslant i\leqslant n$, since
$f_{ii}(\lambda,\lambda)\ne 0$.
If $\mf B^0$ is homogeneous, then there is a co-finite subset
$\mathscr L\subseteq \mathscr L'$ such that $\{\vartheta(S_{\lambda})\in
R(\mf{B}^0)\mid \lambda\in\mathscr L\}$ is a family of pairwise
non-isomorphic homogeneous objects of
$\mf{B}^0$. By Corollary 3.2.4, there is an almost split conflation
$(e_{\lambda}):S_{\lambda}\stackrel{\iota}{\longrightarrow}
E_{\lambda}\stackrel{\pi}{\longrightarrow}
S_{\lambda}$ in $R(\mf B)$ with $\vartheta(e_\lambda)$
is an almost split conflation in $R(\mf B^0)$, where $(E_{\lambda})_X=k^2$,
$E_{\lambda}(x)=J_2(\lambda)$,
$E_{\lz}(a_i)=0$ for $1\leqslant i\leqslant n$ by Lemma 3.2.6 (ii).

Fix any $\lambda\in \mathscr L$
with $(\lambda,\mu)\in \mathscr L'$, then $f_{11}(\lambda,\mu)=0,\phi(\lambda)\phi(\mu)\ne 0$.
Define $L\in R(\mf{B})$ with $L(X)=k^2$; $L(x)={{\lambda\, 0}\choose{0\,\mu}}$,
$L(a_1)=J_2(0)$; and $L(a_i)=0$ for
$2\leqslant i\leqslant n$. Let $g: L \rightarrow S_{\lz}$ be a
morphism with $g_{_X}={1\choose 0}$ and $g(v_j)={0\choose 0}$ for all $j$, then
$g$ is not a retraction. Otherwise, there is a morphism $h:S_{\lz}\rightarrow L $ with
$hg=id_{S_{\lz}}$, thus $h_X=(1,b)$. Set $h(v_1)=(c,d)$,
$$\begin{array}{c}S_{\lz}(a_1)h_X-h_XL(a_1)=f_{11}
(S_{\lz}(x), L(x))h(v_1)\Longrightarrow\\
0(1,b)-(1,b){{0\, 1}\choose{0\, 0}}=f_{11}\big(\lambda,
{{\lambda\,0}\choose{0\,\mu}}\big)(c,d)
=(f_{11}(\lambda,\lambda)c, f_{11}(\lambda,\mu)d),\end{array}$$
which leads $-(0,1)=(\ast,0)$, a contradiction.
$$\xymatrix{L\ar@(ul,ur)^{{\lambda\,0}\choose{0\,\mu}}
\ar@(ul,dl)_{{0\,1}\choose{0\,0}}\ar@/^/[rr]^{g}&&
{S_\lambda}\ar@(ul,ur)^{(\lambda)}
\ar@(ur,dr)^{(0)}\ar@/^/[ll]^{h}}$$

Therefore there exists a lifting $\tilde{g}: L \rightarrow E_{\lambda}$ such that
$\tilde{g}\pi=g$. Set $\tilde{g}_X={{a\,\, b}\choose {c\,\, d}}$,
$\tilde{g}_{_X}\pi_{_X}=g_{_X}$, yields
$\tilde{g}_{_X}={{1\,\, b}\choose {0\,\, d}}$. On the other hand,
$\tilde{g}:L\mapsto E_{\lz}$ is a morphism,
$\tilde{g}_{_X} E_{\lz}(x)=L (x)\tilde{g}_{_X}$, i.e.,
${{1\,\,b}\choose {0\,\,d}}{{\lz\,\,1}\choose {0\,\,\lambda}}
={{\lz\,\, 0}\choose {0\,\, \mu}}{{1\,\, b}\choose{0\,\, d}}$, which
leads ${{\lz\,\, 1+\lz b}\choose {0\,\,\,\,\,\,\lz d}}={{\lz\,\,
\lz b}\choose {0\,\, \mu d}}$, a contradiction. Therefore $\mf{B}^0$
is not homogeneous. The proof is finished.

\medskip

{\bf Proposition 3.4.4}\, If $\mf B^0$ has an induced bocs $\mf{B}$
in the case of MW4, then $\mf B^0$ is
non-homogeneous.

\smallskip

{\bf Proof}\, Fix some $\lambda\in k$
with $\psi(\lambda, x_1)\ne 0$, then $\mathscr L'=\{\mu\mid\psi(\lambda,\mu)\ne 0\}\subseteq k$
is a co-finite subset.
Let $x\mapsto (\lambda)$, since $f_{ii}(\lambda,\lambda)\mid
\phi(\lambda)\mid\psi(\lambda,x_1),f_{ii}(\lambda,\lambda)\ne 0$
for $1\leqslant i\leqslant n_1-1$.
After a series of regularizations $a_i\mapsto\emptyset,w_i=0$,
we obtain an induced bocs of $\mf B$.
Furthermore, since $f_{ii}(\lambda, x_1,\bar x_1)\in k[x_1,\psi(\lambda,x_1)^{-1}
\psi(\lambda,\bar x_1)^{-1}]$ are invertible for $n_1+1\leqslant i\leqslant n$,
after regularizations
$a_i\mapsto\emptyset, w_i=0$ for $n_1<i\leqslant n$,
we obtain an induced minimal local bocs $\mf B_\lambda$
and an induced functor $\vartheta_1$:
$$\mf B_\lambda:\xymatrix {X\ar@(ul,dl)[]_{x_1}},\quad
R_\lambda=k[x_1,\phi(\lambda,x_1)^{-1}],\quad
\vartheta_1: R(\mf{B}_\lambda)\rightarrow R(\mf{B}).
$$
Set $\vartheta_2:R(\mf B)\rightarrow R(\mf{B}^0)$
and $\vartheta=\vartheta_2\vartheta_1$.
Let $S'_\mu\in R(\mf{B}_\lambda)$ be given by $(S'_\mu)_X=k$ and $S'_\mu(x_1)=(\mu)$
for any $\mu\in\mathscr L'$.

If $\mf B^0$ is homogeneous, there exists a co-finite subset
$\mathscr L\subseteq\mathscr L'$, such that $\{\vartheta(S'_{\mu})\in
R(\mf{B}^0)\mid \mu\in\mathscr L\}$ is a family of pairwise
non-isomorphic homogeneous objects of
$\mf{B}^0$. By Corollary 3.2.4, there is an almost split conflation
$(e'_{\mu})$ in $R(\mf B_\lambda)$, such that $\vartheta(e'_\mu)$
is so in $R(\mf B^0)$. Fix any $\mu\in\mathscr L$, consider the
conflation $(e_\mu)=\vartheta_1(e_\mu')$ in $R(\mf B)$, which
is almost split by Lemma 3.2.2 (ii), denote $n_1$ by $s$ for simple:
$$
(e_\mu):S_\mu\stackrel{\iota}{\longrightarrow}
E_\mu\stackrel{\pi}{\longrightarrow} S_\mu, \qquad
(E_\mu)_X=k^2, E_\mu(x)=\lambda I_2, E_\mu(a_s)=J_2(\mu).
$$

Define a representation $L$ of $\mf{B}$ given by $L_X=k^2$,
$L(x)=J_2(\lambda)$, $L(a_s)=\mu I_2$ and others zero. $L$ is well defined,
in fact if we make an unraveling
for $x\mapsto J_2(\lambda)$, then by Lemma 3.2.6 (ii), after a sequence of
regularizations ${{a_{i11}\,a_{i12}}\choose{a_{i21}\,a_{i22}}}$ (splitting from $a_i$)
$\mapsto\emptyset$; and ${{w_{i11}\,w_{i12}}\choose{w_{i21}\,w_{i22}}}$ (splitting from $w_i$)
$=0$ for  $i=1,\cdots,s-1$, we obtain an induced bocs with
$\dz^0\big({{a_{s11}\,a_{s12}}\choose{a_{s21}\,a_{s22}}}\big)=0$.

Let $g: L \rightarrow S_{\mu}$ be a morphism with $g_{_X}={1\choose 0}$ and
$g(v_j)={0\choose 0}$ for all possible $j$.
It is obvious that $g$ is not a retraction.
$$\xymatrix{L\ar@(ul,ur)^{J_2(\lambda)}
\ar@(ul,dl)_{\mu I_2}\ar[rr]^{\tilde g}&&
{E_\lambda}\ar@(ul,ur)^{\lambda I_2}
\ar@(ur,dr)^{J_2(\mu)}}$$
Thus there exists a lifting $\tilde{g}: L \rightarrow E_{\lambda}$ with
$\tilde{g}\pi=g$. $\tilde{g}_{_X}\pi_{_X}=g_{_X}$
leads to $\tilde{g}_{_X}={{1\,b}\choose {0\,d}}$. On the other
hand, since $(x-y)^2\mid f_{11}(x,y)$ if $\bar w\ne 0$,
$(J_2(\lz)-\lz I_2)^2=0$ and hence $f_{11}(J_2(\lambda), \lz
I_2)=0$. $\tilde{g}:L\rightarrow E_{\lambda}$ being a morphism implies that
$$\begin{array}{c}L(a_1)\tilde{g}_{_X}-\tilde{g}_{_X} E_{\lambda}(a_1)=\left\{\begin{array}{cc}
f_{11}(L(x),E_{\lambda}(x))\tilde{g}(\bar w)=0,&\bar w\ne0;\\
0,&\bar w=0,\end{array}\right.\\
\Longrightarrow \mu{{1\,\, b}\choose{0\,\,d}}-{{1\,\, b}\choose{0\,\,d}}J_2(\mu)=0,
\quad \mbox{i.e.} -{{0\, 1}\choose{0\,0}}=0.\end{array}$$
The contradiction tells that $\mf B^0$ is non-homogeneous. The proof is completed.

\medskip

{\bf Proposition 3.4.5}\, Let $\mf B^0$ be a bocs and $\mf{B}$
be an induced bocs with $\T=\{X,Y\};R_X=k[x,\phi(x)^{-1}], R_{Y}=k1_Y$; and a layer
$L=(R;\omega;a_1,\cdots,a_n;
v_1,\cdots,v_{m}),\dz(a_1)=0$:
$$\begin{array}{c}
\mf B:\quad \xymatrix {X\ar[rr]^{a_1}\ar@(ul, dl)[]_{x}
&& Y}\end{array}
$$

Making a reduction given by proposition 2.2.7,
we obtain an induced local bocs. Suppose all the loops $a_2,\cdots,a_n$
in the induced bocs are going to $\emptyset$ by a sequence of regularizations, and some
localizations before each base changes, the
induced bocs $\mf B'$ is minimal with $R'=k[x,\phi'(x)^{-1}]$. Then $\mf B^0$ is non-homogeneous.

\smallskip

{\bf Proof}\, (i) The infinite subset $\mathscr L'=\{\lambda\mid\phi'(\lambda)\ne 0\}\subseteq k$
gives a set of pairwise non-isomorphic objects $S'_\lambda\in R(\mf B')$
with $S'_\lambda(X)=k,S'_\lambda(x)=(\lambda), \forall\,\lambda\in\mathscr L'$.
Write the induced functor $\vartheta_1:R(\mf B')\rightarrow R(\mf B)$, then
$$S_\lambda=\vartheta_1(S'_\lambda):\,\,\,
\xymatrix {k\ar[rrr]^{S_\lambda(a_1)=(1)}\ar@(ul,dl)[]_{S_\lambda(x)=(\lambda)} &&& k}\,\,\,\in R(\mf B),$$
Let $\vartheta_2:R(\mf B)\rightarrow R(\mf B^0)$ be the induced functor,
$\vartheta=\vartheta_2\vartheta_1$.
If $\mf B^0$ is homogeneous, then there is a co-finite subset
$\mathscr L\subseteq \mathscr L'$ such that $\{\vartheta(S'_{\lambda})\in
R(\mf{B}^0)\mid \lambda\in\mathscr L\}$ is a family of homogeneous objects.
Using Corollary 3.2.5,
there is an almost split conflation
$(e_{\lambda}'):S'_{\lambda}\stackrel{\iota}{\longrightarrow}
E'_{\lambda}\stackrel{\pi}{\longrightarrow}
S'_{\lambda}$ in $R(\mf B')$ with $E'_\lambda(x)=J_2(\lambda)$,
such that $\vartheta(e_\lambda')$ is an almost split conflation in $R(\mf B^0)$.

(ii)\, Fix any $\lambda$, there is an almost split conflation:
$(e_\lambda)=\vartheta_1(e_\lambda'): S_\lambda\rightarrow
E_\lambda\rightarrow S_\lambda$ by lemma 3.2.2 (ii). Define an
object $L'\in R(\mf B_X)$ with $L'_X=k^2$,
$L'(x)=\lambda I_2$, $L'(a_i)=(0)_{2\times 2}$
for all $a_i:X\mapsto X$. Set the induced functor
$\vartheta_3:R(\mf B_X)\rightarrow R(\mf B)$
and $L=\vartheta_3(L')$, then $L(a_i)=0,1\leqslant i\leqslant n$. Define a morphism
$g: S_\lambda\rightarrow L$ with $g_{_X}=(0\,1)$, $g_{_Y}=0$, and
$g(v)=0$ for any dotted arrows. We claim that $g$ is not a
retraction. Otherwise, if there is a morphism
$h:L\rightarrow S_{\lz}$ with $gh=id_{S_{\lz}}$,
then $h_{_X}={{1}\choose{c}}$. Since $h_X S(a_1)=L(a_1)h_Y$, i.e. ${{1}\choose{c}}(1)=0$, a
contradiction
$$
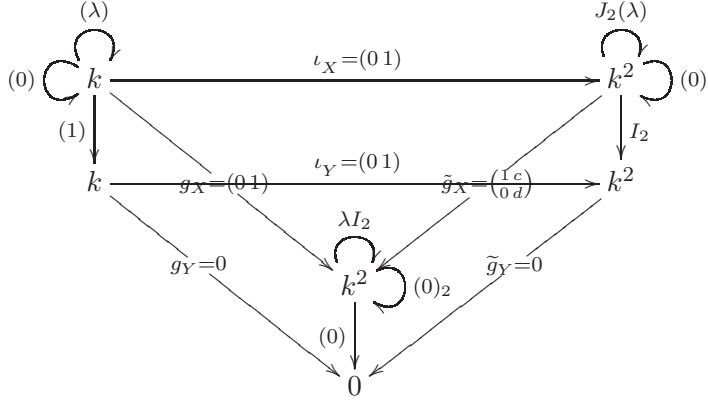
\xymatrix {k\ar@(ul,dl)[]_{(0)}\ar@(ul,ur)[]^{(\lambda)}
\ar[rrrrrr]^{\iota_{_X}={(0\,1)}}\ar[rrrdd]|(.5){g_X=
{(0\,1)}}\ar[d]_{(1)}
&&&&&&k^2\ar@(ul,ur)[]^{J_2(\lambda)}\ar@(ur,dr)[]^{(0)}
\ar[llldd]|(.5){\tilde{g}_X={{1\, c}\choose{0\,d}}}\ar[d]^{I_2}\\
k\ar[rrrrrr]^{\iota_{_Y}={(0\,1)}}\ar[rrrdd]|(.4){g_Y=0}
&&&&&& k^2\ar[llldd]|(.4){\wt{g}_Y=0}
\\
&&& k^2\ar@(ul,ur)[]^{\lambda I_2}\ar@(ur,dr)[]^{(0)_2}\ar[d]_{(0)}&&&
\\&&&0&&&}$$
Therefore, there exists a lifting $\tilde
g: L\rightarrow E_\lambda$ with $\iota\tilde g=g$,
$\tilde g_X={{a\,b}\choose{0\,1}}$. Since $\tilde g$ is a
morphism, we have $E(x)\tilde g_{_X}=\tilde g_{_X}
L(x)$, i.e. ${{0\, 1}\choose{0\, 0}}={{0\, 1}\choose{0\, 0}}{{a\, b}\choose {0\, 1}}
={{a\, b}\choose {0\, 1}}{{0\, 0}\choose{0\,
0}}$, a contradiction.
Therefore $\mf{B}^0$ is not homogeneous. The proof is completed.

\bigskip
\bigskip
\bigskip
\centerline {\bf 4. One-sided pairs}
\bigskip

Through out the present section, we always assume that $\mf A=(R,\K,\M,H=0)$ is a bipartite
matrix bi-module problem with RDCC condition given by Remark 1.4.4.
We will define some special quotient problem of some induced bi-module
problem of $\mf A$.

\bigskip
\bigskip
\noindent{\bf 4.1 Definition of one-sided pairs}
\bigskip

Let $\mf A^0=(R^0,\K^0,\M^0,H^0=0)$ be a bipartite matrix bi-module problem
with RDCC condition and the corresponding bi-co-module problem $\mf C^0$ and bocs $\mf B^0$. Let a sequence of pairs
$$(\mf A^0,\mf B^0),(\mf A^1,\mf B^1),\cdots,(\mf A^r,\mf B^r)$$
be given by reductions in the sense of Lemma 2.3.2
at each step.
Suppose that the leading position of
the first base matrix $A^r_1$ of $\M^r_1$ is $(p^r,q^r+1)$ over $\T^r$, which
is sitting at the leading block $(\textsf{p},\textsf{q})$
of a certain base matrix of $\M^0_1$ partitioned by $\T^0$.
We further assume that $d_1, \cdots, d_m$ are the first $m$ solid arrows of $\mf B^r$,
which locate at the $p$-th row in the formal product $\Theta^r$, such that $d_m$ is
sitting at the last column of the $(\textsf{p},\textsf{q})$-block, see the picture below:

\begin{center}
\setlength{\unitlength}{1mm}
\begin{picture}(40,25)
\put(0,0){\line(1,0){40}}  \put(0,0){\line(0,1){25}}
\put(0,25){\line(1,0){40}} \put(40,0){\line(0,1){25}}
\put(15,10){\line(1,0){25}} \put(15,10){\line(0,1){5}}
\put(15,15){\line(1,0){25}}
\put(20,10){\line(0,1){5}} \put(35,10){\line(0,1){5}}
\put(25,11){$\cdots$}\put(16, 11){$d_1$}
\put(35.5, 11){$d_m$}\put(86,12){$(4.1\mbox{-}1)$}
\end{picture}
\end{center}

Recall the notation below Formulae (2.1-2): we have the quotient problem of the matrix bi-module problem
${(\mf A^r)}^{[m]}=(R^r, \K^r,(\M^r)^{[m]},H^r)$ of $\mf A^r$ ; the sub-co-bi-module problem
$(\mf C^r)^{(m)}=(R^r, \C^r,(\N^r)^{(m)},\partial\mid_{(\N^r)^{(m)}})$ of $\mf C^r$
with the quasi-basis $d_1,\cdots,d_m$ of $(\N^r)^{(m)}$; and the sub-bocs
$(\mf{B}^r)^{(m)}$ of $\mf B^r$. We obtain a {\it quotient-sub-pair $((\mf
A^r)^{[m]}, (\mf B^r)^{(m)})$}.
Denote by $\bar A^r_1, \ldots,\bar A_m^r$
the quasi-basis of $(\M^r)^{[m]}$ over $R^r\otimes_k R^r$,
we have the formal equation of the pair:
$$\begin{array}{c}(\sum_{Y\in \T^r}
e_{_{Y}}\ast E_{Y}+\sum_{j=1}^{t^r}v^r_j\ast V^r_j)
(\sum_{i=1}^m d_i\ast\bar A_i^r+H^r)\\
=(\sum_{i=1}^m d_i\ast\bar A^r_i+H^r)(\sum_{Y\in
\T^r}e_{_{Y}}\ast E_{Y}+\sum_{j=1}^{t^r}v^r_j\ast V^r_j).\end{array}$$
Let $d_j: X\rightarrow Y_j$ (possibly $Y_j=X$), which can be
rewritten as a {\it reduced formal equation}:
$$e_{_{X}}(d_1,d_2,\cdots,d_m)= (w_1,w_2,\cdots,w_m)+(d_1,d_2,\cdots,d_m)
\left(\begin{array}{cccc}e_{_{Y_1}}&w_{12}&\cdots & w_{1m}\\
&e_{_{Y_2}}&\cdots& w_{2m}\\&&\ddots &\vdots\\&&&e_{_{Y_m}}
\end{array}\right) \eqno (4.1\mbox{-}2)$$

{\bf Remark 4.1.1}\, (i)\, In the formula (4.1-2), $e_{_{X}}$ is the $(p, p)$-th entry of the formal
product $\sum_{Y\in {\T^r}}e_{_{Y}}\ast E_{Y}$; and $e_{_{Y_\xi}}$ the $(q+\xi,q+\xi)$-th entry
of that for $\xi=1,\cdots,m$.

(ii)\, $(w_1,\cdots,w_m)$ is the $(p,q+1)$-th, $\cdots$, $(p,q+m)$-th entries of
$(\sum_{j=1}^{t^r}v^r_j\ast V^r_j)H^r-
H^r(\sum_{j=1}^{t^r}v^r_j\ast V^r_j)$.
$w_\xi=\sum_j\alpha_{\xi}^jv^r_j$, where $j$ runs over
$s(v^r_j)\ni p,e(v^r_j)\ni q+\xi,\alpha_{\xi}^j\in k,1\leqslant\xi\leqslant m$.

(iii)\, For $1\leqslant\eta<\xi\leqslant m$, $w_{\eta\xi}$ is the
$(q+\eta,q+\xi)$-th entry of $\sum_{j=1}^{t^r}v^r_j\ast V^r_j$. And
$w_{\eta\xi}=\sum_j\beta_{\eta\xi}^jv^r_j$ where $j$ runs over
$s(v^r_j)\ni q+\eta,e(v^r_j)\ni q+\xi,\beta_{\eta\xi}^j\in k$.

(iv)\, The differential of any solid arrows can be read off from the reduced formal
product, we notice that the solid arrows appear in each monomial only once from the left
to a dotted arrow:
$$\begin{array}{c}
-\delta(d_i)=w_{i} +\sum_{j<i} d_jw_{ij},\quad 1\leqslant i\leqslant
m.\end{array}\eqno{(4.1\mbox{-}3)}
$$

{\bf Definition 4.1.2}\, A bocs $\mf{B}$  with a layer $L=(R;\omega;
a_1,\ldots,a_m,b_1,\cdots,b_n;\underline u_j,\underline v_j,\bar u_j,\bar v_j)$, see the picture below,
is called {\it one-sided}, provided
$R$ is trivial; $\dz(a_i)=\sum_j\alpha_{ij}\underline v_j
+\sum_{b_{i'}\prec a_i,j}\beta_{ii'j}b_{i'}\underline v_j
+\sum_{i'<i,j}\gamma_{ii'}a_{i'}\underline u_j$,
$\dz(b_i)=\sum_j\lambda_{ij}\bar u_j+\sum_{i'<i,j}\mu_{ii'}b_{i'}\bar u_j
+\sum_{a_{i'}\prec b_i,j}\nu_{ii'j}a_{i'}\bar v_j$ with constant coefficients.
{\unitlength=1mm
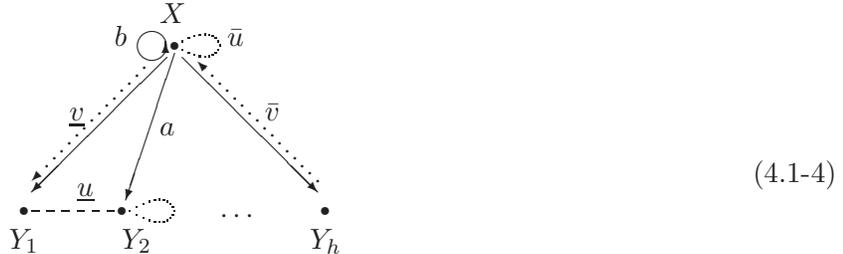
$$\begin{array}{c}
\begin{picture}(80,45)
\put(29,39){\circle{4}} \put(30.8,39){\vector(0,1){0.5}}
\put(32,39){\circle*{1}} \qbezier[10](32,39)(37,42)(38,39)
\qbezier[10](32,39)(37,36)(38,39)

\put(32,38){\vector(-1,-3){6.5}} \put(31,37.5){\vector(-1,-1){18}}
\put(33,37.5){\vector(1,-1){18}} \put(25,17){\circle*{1}}
\put(12,17){\circle*{1}}  \put(52,17){\circle*{1}}
\put(10,12){$Y_1$}  \put(50,12){$Y_h$} \put(25,12){$Y_2$}
\put(30,42){$X$}

\qbezier[10](26,17)(31,20)(32,17) \qbezier[10](26,17)(31,14)(32,17)

\mput(13,17)(2,0){6}{\line(1,0){1}}
\mput(28,36)(-1,-1){14}{\circle*{0.5}}
\put(15,23){\vector(-1,-1){2}}

\mput(36,36)(1,-1){16}{\circle*{0.5}} \put(37,35){\vector(-1,1){2}}

\put(39,39){\makebox{$\bar u$}} \put(18,29){\makebox{$\underline v$}}
\put(44,29){\makebox{$\bar v$}} \put(19,19){\makebox{$\underline u$}}
\put(24,39){\makebox{$b$}} \put(30,27){\makebox{$a$}}
\put(38,15.5){\makebox{$\cdots$}}
\end{picture}
\end{array}\eqno{(4.1\mbox{-}4)}
$$}\vskip -15mm

\noindent The associated bocs $(\mf B^r)^{(m)}$ of
$(\mf A^r)^{[m]}$ is one sided by Formula (4.1-3).
We call this quotient-sub-pair a {\it one-sided pair},
and denoted by $(\bar{\mf A},\bar{\mf B})$ for simple.

\medskip

Write $\bar{\mf A}=(R, \K,\bar\M,F)$, where $R=k1_X\times
\prod_{j=1}^hk1_{Y_j}$ is a trivial sub-algebra of $R^r$;
$\{X;Y_1,\cdots,Y_h\}\subseteq \T^r$; $\bar\M_1$ has an $R$-$R$-quasi-basis
$(1_{_{XY_{t(d_1)}}},0,\cdots,0)$, $\cdots$, $(0,0,\cdots,1_{_{XY_{t(d_m)}}})$;
$F=(0)$; let $E_X=(1_{_X})\in\IM_{1\times 1}(R)$, and $E_{Y_j}=\sum_{p\in
Y_j}1_{_{Y_j}}E_{pp}$ with $E_{pq}$ the matrix units of
$\IM_{m\times m}(R)$, then $\{E_X;E_{Y_j}\}_j$ forms a part of the
$R$-quasi-basis of $\K_0$. However, there exist some linear relations between
$w_i$ and $w_{ij}$ in the formula (4.1-2).

\medskip

Now we start the reduction procedure in the sense of
Lemma 2.3.2. Let $(\mf A,\mf B)$ be any pair, $(\mf A^{[p]},\mf B^{(p)})$ be a quotient-sub-pair.
Since the reduction for any pair is made with respect to
an admissible $\bar R$-$R'$-bi-module $L$ by Proposition 2.2.1-2.2.7,
or a regularization in Proposition 2.2.8,
which is completely as the same as to make a
reduction for the quotient-sub-pair.
Therefore there are two sequences of reductions:
$$\begin{array}{ccccccc}&(\bar{\mf A},\bar{\mf B}),&
(\bar{\mf A}^1,\bar{\mf B}^1),&\cdots,&(\bar{\mf A}^i,\bar{\mf B}^i)
&\cdots,&(\bar{\mf A}^s,\bar{\mf B}^s);\\
&(\mf A^r,{\mf B}^r), &({\mf A}^{r+1},{\mf B}^{r+1}),&\cdots,&({\mf A}^{r+i},
{\mf B}^{r+i}),&\cdots,&({\mf A}^{r+s},
{\mf B}^{r+s}).\end{array} \eqno{(4.1\mbox{-}5)}$$
From now on, we perform reductions according to Formula (2.4-3). More precisely, the system $\IF^{ri}$
in (2.4-3) for the pair
$(\mf A^{r+i},\mf B^{r+i})$ can be written as
a reduced form $\bar\IF^i$ for the pair $(\bar{\mf A}^i,\bar{\mf B}^i)$:
$$\begin{array}{c}{\bar\IF}^i: \bar\Psi_{\n^i}^{l}F^i
\equiv_{\prec(\bar p^i,\bar q^i)} \bar\Psi_{\n_i}^{m}+F^i\bar\Psi_{\n^i}^{r},
\end{array} \eqno{(4.1\mbox{-}6)}$$
where the upper indices $l,m,r$ on $\bar\Psi$ indicate left, middle and right parts
of the matrices of dotted elements.
Since it is difficult to determine the dotted arrows after a reduction,
instead we describe the linear relation of the dotted elements in $\bar\Psi$
appearing at the reduction. In order to do so, we define a
{\it pseudo reduced formal equation}
of the pair $(\bar{\mf A}^i,\bar{\mf B}^i)$:
$${e}^i_{_X}(F^i+\bar\Theta^i)=(W_1,\cdots, W_m)+(F^i+\bar\Theta^i)
\left( \begin{array}{cccc}e_{_{Y_1}}^i&W_{12}&\cdots & W_{1m}\\
& e_{_{Y_2}}^i&\cdots&W_{2m}\\&&\ddots &\vdots\\&&& e_{_{Y_m}}^i
\end{array}\right),\eqno{(4.1\mbox{-}7)}$$
where $W_h$ and $W_{hl}$ split from
$w_h,w_{hl}$; $e^i_{_Z}$ is given by the proof (ii) of Theorem 2.4.1.
$\bar\Theta^i$ is the reduced formal product of $(\bar\M^i_1,\bar\N^i_1)$.
$F^i+\bar\Theta^i$ is an $(1\times\m)$-partitioned
matrix under $\bar\T$ with size
vector $\n^i=(n^i_0;n^i_1,\cdots,n^i_m)$ as shown below, where $F^i$ is sitting in the blank part:
\unitlength=1mm
$$F^i+\bar\Theta^i=\begin{picture}(85,10)
\put(0,-10){\line(1,0){20}} \put(0,-10){\line(0,1){20}}
\put(0,10){\line(1,0){20}} \put(20,-10){\line(1,0){65}}
\put(15,-10){\line(0,1){20}} \put(20,10){\line(1,0){65}}
\put(85,-10){\line(0,1){20}} \put(24,-10){\line(0,1){20}}
\put(25,7){$d_q^i\;\cdots$} \put(24,6){\line(1,0){61}}
\put(39,-10){\line(0,1){20}} \put(40,3){$d_p^i\; d_{p+1}^i \;
\cdots\;  $} \put(42,7){$\cdots\quad \cdots$} \put(50,-1){$d_1^i\;
d_2^i\cdots$} \put(39,2){\line(1,0){46}} \put(49,-2){\line(0,1){4}}
\put(49,-2){\line(1,0){36}} \put(63,-10){\line(0,1){20}}
\put(65,7){$\cdots$} \put(65,-1){$\cdots$} \put(65,3){$\cdots $}
\put(71,-10){\line(0,1){20}}\put(72,7){$\cdots\; \cdots$}
\put(72,3){$\cdots\; \cdots$} \put(72,-1){$\cdots\;d_{p-1}^i$}
\put(77,-7){$ $} \put(50,-7){$ $} 
\end{picture}
$$
\vskip 10mm

{\bf Remark 4.1.3}\, Let $(\bar{\mf A}',\bar{\mf B}')$ be any induced pair of $(\bar{\mf A},\bar{\mf B})$
after several reductions in the sense of Lemma 2.3.2. And $(\bar{\mf A}'',\bar{\mf B}'')$
is an induced pair of $(\bar{\mf A}',\bar{\mf B}')$ by one reduction of 2.3.2.

(i) If we have a linear relation of dotted elements $\sum_ju_j=0$ in $\bar\IF'$,
then $\sum_j\bar u_j=0$ in ${\bar\IF}^{''}$ with $\bar u_j$ being the split of $u_j$.

(ii) Suppose the first arrow $\dz(a_1')=v+\sum_j\alpha_ju_j$ in $\mf B'$, where
$v, u_j$ are dotted elements of ${\bar\Psi}_{\n'}$. If $v$ is a dotted arrow, and $v\notin\{u_j\}$, then $\dz(a_1')\ne 0$.

(iii) Set $a_1'\mapsto\emptyset$ given above, we sometimes say that $v$ is {\it replaced by $-\sum_ju_j$
in $\bar\Psi_{\n''}$}.

(iv) If we are able to determine, that a dotted element $v$ of $\bar\Psi_{\n''}$ is
linearly independent of all the others, then $v$ is said to be a {\it dotted arrow preserved in
$\bar{\mf B''}$}.

We are able to read off the differential $\dz$ of the solid arrows
from the pseudo reduced formal equation by Theorem 1.6.4.
Denote the $(p,q)$-th entry of $F^i$ splitting from $d_l:X\mapsto Y$ by $f^i_{lpq}$:
$$\begin{array}{c}\dz(d^i_{l\bar p\bar q})=-w^i_{lpq}-
\sum_{j>l,q<\bar q}f^i_{l\bar pq}w^i_{lj,q\bar q}-
\sum_{q<\bar q}f^i_{l\bar pq}w^i_{Yq\bar q}+\sum_{q,j<l}w^i_{X\bar pq}f^i_{lq\bar q}.
\end{array}\eqno{(4.1\mbox{-}8)}$$

\bigskip
\bigskip
\noindent{\bf 4.2\, The differentials in one sided pairs}
\bigskip

Let $(\bar{\mf A}, \bar{\mf B})$ be a one-sided pair given in formulae (4.1-2) and
(4.1-4). This subsection is devoted to calculating the differentials of
the solid arrows of $\bar{\mf B}$.

\medskip

Suppose first that $\bar{\mf B}$ is a one sided local bocs having a layer
$L=(R;\omega;b_1,\dots,b_n;v_1,\cdots,v_m)$ with $R=k1_X$. Recall Classification 3.3.2,
which is simpler in one sided case.

\medskip

{\bf Classification 4.2.1}\, Let $\bar{\mf B}$ be a one sided local bocs
given above.

(i)\, $\bar{\mf B}$ satisfies the triangular formula (3.3-5) (letter $a$ is changed to letter $b$).
After making an easy base change, (3.3-5) can be written as:
$\delta^0(b_1)= \bar u_1,
\cdots,\delta^0(b_{n})= \bar u_{n}$.

Formulae (3.3-6) can be written as following given by a base change,
$$\left\{\begin{array}{l}
\delta^0(b_1)=\bar u_1, \cdots\,\, \cdots
\delta^0(b_{n_0-1})= \bar u_{n_0-1},\\
\delta^0(b_{n_0})= \alpha_{n_0,1}\bar u_1+\cdots+ \alpha_{n_0,n_0-1}\bar u_{n_0-1},
\end{array}\right.\eqno{(4.2\mbox{-}1)}$$
with $\alpha_{n_0,j}\in k$. After a series of regularization, Formula (3.3-2) with $x=b_{n_0}$
can be written as:
$$\left\{\begin{array}{ll}\dz^0(b_{n_0+1})&=h_{n_0+1,1}(x)\bar u_{n_0+1},\\
\dz^0(b_{n_0+2})&=h_{n_0+2,1}(x)\bar u_{n_0+1}+h_{n_0+2,2}(x)\bar u_{n_0+2},\\
&\cdots\quad \cdots\\ \dz^0(b_n)&=h_{n,n_0+1}(x)\bar u_{n_0+1}
+h_{n,n_0+2}(x)\bar u_{n_0+2}+\cdots+h_{n,n}(x)\bar u_n,
\end{array}\right.\eqno{(4.2\mbox{-}2)}$$
with the polynomial $\phi=1$, and
$h_{n_0+i,j}(x)=\alpha_{n_0+i,j}^0+\alpha_{n_0+i,j}^1x$ of degree $1$.

(ii)\, $\bar{\mf B}$ satisfies Lemma 3.3.3 (i), if $h_{ii}(x)\in k\setminus\{0\}$ for $n_0<i\leqslant n$.

(iii)\, $\bar{\mf B}$ has an induced bocs $\bar{\mf B}_{(\lambda_0,\cdots,\lambda_{l})}$
satisfies Lemma 3.3.3 (ii), if $h_{ss}(x)\in k[x]\setminus k$ non-invertible for some $n_0<s\leqslant n$,
which is in the case of MW3.

(iv)\, $\bar{\mf B}$ has an induced bocs $\bar{\mf B}_{(\lambda_0,\cdots,\lambda_{l-1})}$
for some $l< n$, which satisfies Lemma 3.3.4 (ii) with a triangular formula similar to
Formula (4.2-2) for $b_1,\cdots,b_{n_1}$ with a polynomial $\phi(x)=\prod_{i=1}^{n_1-1}h_{ii}(x)\ne 0$,
and $h_{n_1,n_1}(x)=0$, i.e. $\bar w=0$.
Moreover the differential $\dz^1$ with respect to
$x,x_1=b_{n_1}$, and a polynomial $\psi(x)=\phi(x)\prod_{i=n_1+1}^{n}c_{i}(x)\ne 0$
given by Formula (3.3-4):
$$\left\{\begin{array}{lll}\dz^1(b_{n_1+1})&=K_{n_1+1}+h_{n_1+1,n_1+1}(x)u_{n_1+1},&\\
&\cdots\\ \dz^1(b_n)&=K_{n}\quad\,\,\,+h_{n,n_1+1}(x,x_1)u_{n_1+1}&
+\cdots+h_{nn}(x)u_n,\end{array}\right.\eqno{(4.2\mbox{-}3)}$$
where $h_{ii}(x)\mid\psi(x)$ for $i=n_1+1,\cdots,n$.
Then $\bar{\mf B}_{(\lambda_0,\cdots,\lambda_{l-1})}$ is in the case of MW4 with
$\bar w=0$. In particular the case of MW5 can not occur.

\medskip

Now we concern the general one sided pairs $(\bar{\mf A},\bar{\mf B})$. Suppose $\bar{\mf B}_X$,
the induced local bocs of $\bar{\mf B}$, is in the case of Classification 4.2.1 (iii) or (iv),
then $\bar{\mf B}_X$ is wild and non-homogeneous, so is $\bar{\mf B}$.
Since $\bar{\mf B}_X$ satisfies 4.2.1 (i) is relatively simple, we now
concentrate on the case of 4.2.1 (ii). Denote the solid edges of $\bar{\mf B}$
before $b_{n_0}$ by $a_1,\cdots a_h$, the
differential $\dz^0$ on $a$'s can be written as:
$$\begin{array}{c}
\delta^0(a_1)=\underline v_1,
\,\,\cdots\cdots,\,\,
\delta^0(a_h)=\underline v_{h};
\end{array}\eqno(4.2\mbox{-}4)$$
or $\delta^0(a_1)=\underline v_1,
\cdots,\delta^0(a_{h_1-1})=\underline v_{h_1-1},
\delta^0(a_{h_1})=\sum_{j=1}^{h_1-1}\beta_{h_1j}\underline v_j$;
then $\delta^0(a_{h_1+1})=\underline v_{h_1+1},\cdots$,
$\delta^0(a_{h_2-1})=\underline v_{h_2-1}$,
$\delta^0(a_{h_2})=\sum_{j\ne n_1,j=1}^{h_2-1}\beta_{h_1j}\underline v_j$;
continuously we obtain say $s$ formulae with $\beta\in k$ and $s$ edges
$a_{h_1}, a_{h_2}, \cdots, a_{h_s}$.
Set $\Lambda=\{1,\cdots,h\}\setminus\{h_1,\cdots,h_s\}$, then
$$\left\{\begin{array}{lll}\dz^0(a_i)&=\quad\quad\underline v_i,&i\in\Lambda;\\
\dz^0(a_{h_l})&=\sum_{j\in\Lambda,j<h_l}\beta_{h_l,j}\underline v_j,&l=1,\cdots,s.
\end{array}\right. \eqno(4.2\mbox{-}5)$$

\medskip

{\bf Convention 4.2.2}\, Suppose $\bar{\mf B}$ is
a one-sided bocs with $\bar{\mf B}_X$ satisfying 4.2.1 (ii). \ding{172}
All the loops $b_1,\cdots,b_n$ at $X$ are called {\it $b$-class arrows},
where the loop $\bar b=b_{n_0}$
is said to be {\it effective} or $\bar b$-class; the others
are {\it non-effective}.  \ding{173}
The edges $a_1,\cdots,a_h$ before $\bar b$ are called {\it $a$-class arrows},
where $\{\bar a_i=a_{h_i}\mid 1\leqslant i\leqslant s\}$
are said to be {\it effective} or $\bar a$-class;
the others are {\it non-effective}.
\ding{174} Let $c_1,c_2,\ldots,c_t$ be the solid edges after $\bar b$, which are called
$c$-class arrows and effective.

An solid arrow splitting from one of the classes $a,\bar a;c;b,\bar b$, and a dotted element from $\underline u,
\underline v,\bar u,\bar v$ is said to be the same class.

\medskip

Next we give a special case of the differentials $\dz^1$ with respect to $\bar b$ on $c$-class arrows,
where $\{\underline v_j\}_{j\in \Lambda}\cup\{\underline v_{h+j}\}_{1\leqslant j\leqslant t}$
are dotted arrows; $\gamma_{ij}(\bar b)=\gamma^0_{ij}+\gamma^1_{ij}\bar b,\gamma^0_{ij},\gamma^1_{ij}\in k$;
$\gamma_{i,h+i}(\bar b)\ne 0,1\leqslant i\leqslant t$:
$$\left\{\begin{array}{l}\dz^1(c_1)=\sum_{j\in\Lambda}\gamma_{1j}(\bar b)
\underline v_j+\gamma_{1,h+1}(\bar b)\underline v_{h+1},\\
\quad\cdots\quad\cdots\quad\cdots\\
\dz^1(c_t)=\sum_{j\in\Lambda}\gamma_{tj}(\bar b)
\underline v_j+\gamma_{t,h+1}(\bar b)\underline v_{h+1}+\cdots
+\gamma_{t,h+t}(\bar b)\underline v_{h+t},\\
\end{array}\right.\eqno{(4.2\mbox{-}6)}$$

{\bf Lemma 4.2.3}\, Let $\bar{\mf B}$ be
a one-sided bocs with $\bar{\mf B}_X$
satisfying 4.2.1 (ii). If Formula (4.2-6) fails,
i.e. there exists some $1\leqslant t_1\leqslant t$
with $\gamma_{t_1,h+t_1}=0$, then $\bar{\mf B}$ is not homogeneous.

\smallskip

{\bf Proof}\, We make a sequence of regularization
$a_j\mapsto\emptyset,\underline v_j=0$ for
$j\in\Lambda$, $b_j\mapsto\emptyset,\bar u_j=0$ for $1\leqslant j<n_0$,
and edge reduction $\bar a_i\mapsto (0)$ for $i=1,\cdots, s$, then obtain
an induced pair $(\mf A',\mf B')$. Without loss of generality,
suppose $\T'=\{X,Y\}$. Make a loop mutation $\bar b\mapsto (x)$, a
localization $\phi(x)=\prod_{i=1}^{t_1-1}g_{i,h+i}(x)$,
and regularizations $c_i\mapsto\emptyset,\underline v'_{h+i}=0$ for $1\leqslant i<t_1$,
we obtain an induced bocs $\xymatrix{\cdot\ar@(ul,dl)[]|{x}\ar[r]^{c_{t_1}}&\cdot}$ of
two vertices with $\dz(c_{t_1})=0$, which is in the case of MW1.
Thus $\bar{\mf B}'$, consequently $\bar{\mf B}$, is wild and non-homogeneous, this finishes the proof.

\medskip

Let $\bar\dz$ be obtained from the differential $\dz$ by removing all
the monomial involving some non-effective $a,b$-class solid arrows.
Now we write $\bar\dz$ acting on all $a,b,c$-class arrows in the following
three formulae (where the third one is given under the assumption that the differential
$\dz^1$ of $c$-class arrows satisfies Formula (4.2-6)):
$$\begin{array}{c}\bar\dz(a_i)=\underline v_i
+\sum_{h_l<i}\bar a_{l}(\sum_{j}\epsilon_{ilj}\underline u_{j}),
\quad\epsilon_{ilj}\in k, i\in\Lambda;\\
\bar\dz(\bar a_\tau)=\sum_{j<h_\tau}\bar\beta_{\tau j}\underline v_j
+\sum_{l<\tau}\bar a_{l}(\sum_{j}\bar\epsilon_{\tau lj}\underline u_{j}),
\quad\bar\beta_{\tau j},\bar\epsilon_{\tau lj}\in k,
1\leqslant \tau\leqslant s.\end{array}\eqno{(4.2\mbox{-}7)}$$
$$\begin{array}{c}\bar\dz(b_i)=\bar u_i
+\sum_{\bar a_{l}\prec b_i}\bar a_{l}(\sum_{j}\varepsilon_{ilj}\bar v_{j}),\qquad
\varepsilon_{ilj}\in k,\quad i<n_0;\\
\bar\dz(\bar b)=\sum_{j=1}^{n_0-1}\bar\alpha_{j} \bar u_{j}
+\sum_{l=1}^s\bar a_{l}(\sum_{j}\bar\varepsilon_{lj}\bar v_{j}),\quad
\bar\alpha_{j}, \bar\varepsilon_{lj}\in k,\quad i=n_0;\\
\bar\dz(b_i)=\bar u_i+\sum_{j=1,j\ne n_0}^{i-1}\alpha_{ij}(\bar b)\bar
u_{j}+\sum_{l=1}^s\bar a_l(\sum_{j}\varepsilon_{ilj}\bar v_j)
+\sum_{c_l\prec b_i}c_l(\sum_j\varepsilon'_{ilj}\bar
v_{j}),\\\alpha_{ij}(\bar b)\in k[\bar b],\quad \varepsilon_{ilj},\varepsilon'_{ilj}\in k,
\quad i>n_0.\end{array}\eqno{(4.2\mbox{-}8)}$$
$$\begin{array}{c}\bar\dz(c_\tau)=\sum_{l=1}^s\bar a_{l}
(\sum_{j}\zeta_{\tau lj}\underline u_{j})+
\sum_{j\leqslant h+\tau}\gamma_{\tau j}(\bar b)\underline v_j+
\sum_{l=1}^{i-1}c_l(\sum_j\xi_{\tau lj}
\underline u_{j}),\,\,\zeta_{\tau lj},\xi_{\tau lj}\in k.
\end{array}\eqno{(4.2\mbox{-}9)}$$

\medskip
\bigskip
\noindent {\bf 4.3. Reduction sequences of one-sided pairs}
\bigskip

In the present subsection, we will construct a reduction sequence starting
from a one-sided pair $(\bar{\mf A},\bar{\mf B})$ with at least two
vertices, where $\bar{\mf B}_X$ satisfies Classification 4.2.1 (ii).

\medskip

{\bf Condition 4.3.1 (BRC)}\, Let $({\mf A},{\mf B})$ with trivial $\T$ be any pair of a matrix bi-module
problem and the associated bocs.

(i)\, Suppose the solid arrows $\mathcal D=\{d_1,\cdots,d_{q}\}$ and
$\mathcal E=\{e_1,\cdots,e_{p}\}$
locate at the bottom row of $\Theta$ form the first
$p+q$ arrows of $\mf B$ (not necessarily full of the row), such that $e_1,\cdots, e_{p-1}$ are edges starting from $X$,
$e=e_{p}$ is a loop at $X$, and $d_i\prec e_p,1\leqslant i\leqslant q$.
There exists a set of dotted arrows $\mathcal U=\{u_1,\cdots,u_q\}$,
write $\mathcal W=\V^\ast\setminus\U=\{w_1,\cdots w_t\}$.

(ii)\, Denote by $\bar\dz$ the part of the differential of any solid arrow by
removing all the monomials containing any solid arrow besides of
$\{e_1,\cdots,e_p\}$, such that
$$\begin{array}{c}\bar\dz(d_i)=u_i+
\sum_{j=1}^t(\sum_{e_l\prec d_i}\lambda_{ijl}e_l)w_j,\quad
\lambda_{ijl}\in k,i=1,\cdots,q;\\
\bar\dz(e_i)=\sum_{d_j\prec e_i}\alpha_{ij}
u_j+\sum_{j=1}^{t}(\sum_{l=1}^{i-1}\mu_{ijl}
e_l)w_j,\quad  \alpha_{ij},\mu_{ijl}\in k,i=1,\cdots,p.\end{array}$$
Then $({\mf A},{\mf B})$ is said to have the
{\it bottom row condition} (BRC) with respect to $(\mathcal D,\mathcal U)$
and $(\mathcal E,\mathcal W)$.

\medskip

Suppose the pair $({\mf A},{\mf B})$ satisfies (BRC) with $p>1$ or $q>1$
and the first arrow $a_1:X\mapsto Y$, where $X\ne Y$ if $a_1=e_1$. Now we
discuss the condition (BRC) on the induced pair $({\mf A}',{\mf B}')$.

(i)\, If $a_1=d_1$, then $d_1\mapsto\emptyset$, let $\mathcal D'=\{d_2,\cdots,d_q\},
\mathcal U'=\{u_2,\cdots,u_q\}$ and $\mathcal E'=\mathcal E,\mathcal W'=\mathcal W$.

(ii)\, If $a_1=e_1$, and $e_1\mapsto(0)$, let $\mathcal D'=\mathcal D,\mathcal U'=\mathcal U$;
$\mathcal E'=\{e_2,\cdots,e_p\},\mathcal W'=\mathcal W$.

(iii)\, If $a_1=e_1$, $e_1\mapsto{1\choose 0}$, let
$\mathcal D'=\{d_{i21},d_{i22} \mid\, t(d_i)=X\}\cup\{d_{j2} \mid\,
t(d_j)\ne X\}$, correspondingly $\mathcal U'=\{u_{i21},u_{i22},u_{j2}\}$;
$\mathcal E'=\{e_{22},\cdots,e_{p-1,2},e_{p21},e_{p22}\}$, where $\bullet_{i2}$, resp.
$\bullet_{i2l}$, stands for the arrow of $\mf B'$ locating at the second row of the
$2\times 1$, resp. $2\times 2$, matrix splitting from $\bullet_i$ of $\mf B$.

(iv)\, If $a_1=e_1,e_1\mapsto {{0\, 1}\choose {0\, 0}}$, set
$\mathcal D'=\{d_{i21},d_{i22} \mid\, t(d_i)=X\,$ or $Y\}\cup\{d_{j2} \mid\,
t(d_j)\ne X,Y\}$, $\mathcal U=\{u_{i21},u_{i22},u_{j2}\}$;
$\mathcal E'=\{e_{i21},e_{i22} \mid\, t(e_i)=Y\}\cup\{e_{j2} \mid\,
t(e_j)\ne Y\}\cup\{ e_{p21}, e_{p22}\}$.

\medskip

{\bf Lemma 4.3.2}\, Suppose the pair $({\mf A},{\mf B})$ satisfies (BRC) with $p>1$ or $q>1$
and the first arrow $a_1:X\mapsto Y$. Then after making a reduction $a_1\mapsto G$ as above (i)-(iv),
the induced pair $({\mf A}',{\mf B}')$ satisfies (BRC) with respect to
$(\mathcal D',\mathcal U')$ and $(\mathcal E',\mathcal W')$ defined as above (i)-(iv).

\smallskip

{\bf Proof}\, The case (i) and (ii) are trivial. Suppose $G={{1}\choose {0}}$ in (iii) resp.
${0\,1}\choose {0\, 0}$ in (iv), then $X,Y$ split into two vertices $X', Y'$ resp. three vertices $X', Y', Y''$:
$$e_{_X}\mapsto e'_{_X}=\left (\begin{array}{cc}e_{_{Y'}}&w\\
&e_{_{X'}}\end{array}\right ),\quad e_{_Y}\mapsto e_{_Y}'=e{_{_{Y'}}}\,\,\, \mbox{resp.}\,\,\,
e_{_Y}\mapsto e_{_Y}'=\left (\begin{array}{cc}e_{_{Y''}}&w'\\
&e_{_{Y'}}\end{array}\right ).$$
In the two cases, $e_{i2}$ or $e_{i21},e_{i22}$ for $1<i<p$
starting at $X'$ do not end at $X'$, since $t(e_i)\ne X$ by (RBC) (i) on $({\mf A},{\mf B})$; the edge
$e_{p21}: X'\mapsto Y'$, and the loop $e_{p22}:X'\mapsto X'$.
By (BRC)(ii),
$$\begin{array}{c}\bar \dz(D_i)=U_i+\sum_{j=1}^t\lambda_{ij1}GW_j
+\sum_{j=1}^t(\sum_{e_l\prec d_i}\lambda_{ijl}E_l)w_j,\\
\bar\dz(E_i)=\sum_{d_j\prec e_i}\alpha_{ij}U_j+\sum_{j=1}^t\mu_{ij1}GW_j
+\sum_{j=1}^t(\sum_{l=2}^{i-1}\mu_{ijl}E_l)W_j+e'_{_X}E_i-E_ie'_{_Y},\end{array}$$
where $\bar\dz(M)=(\bar\dz(a_{ij}))$ for $M=(a_{ij})$. Since the bottom row of $G$ is $(0)$
or $(0\,0)$ and $e_{_X}',e'_{_{Y}}$ are upper triangular, (BRC) on
$({\mf A}',{\mf B}')$ follows, the proof is finished.

\medskip

{\bf Lemma 4.3.3}\, Let $(\bar{\mf A},\bar{\mf B})$ be a one-sided pair with $\bar\T$ trivial,
$|\bar\T|\geqslant 2$, and $\bar{\mf B}_X$ having 4.2.1 (ii).
Then the pair satisfies (BRC) with respect to the sets $(\mathcal D,\mathcal U)$ and $(\mathcal E,\mathcal W)$,
where
$$\begin{array}{c}\mathcal D=\{a_i,i\in\Lambda\}\cup\{b_j,j<n_0\},
\mathcal U=\{\underline v_i,i\in\Lambda\}\cup\{\bar u_j,j<n_0\};
\mathcal E=\{\bar a_\tau,1\leqslant\tau\leqslant s\}\cup\{\bar b\}.\end{array}$$


{\bf Theorem 4.3.4}\, Let $(\bar{\mf A},\bar{\mf B})$ be a one-sided pair with $\bar\T$ trivial having at
least two vertices, such that $\bar{\mf B}_X$ satisfies Classification 4.2.1 (ii).
Then there exists a reduction sequence:
$$(\bar{\mf A},\bar{\mf B})=(\bar{\mf A^0},\bar{\mf B}^0),(\bar{\mf A}^1,\bar{\mf B}^1),\cdots,
(\bar{\mf A}^\gamma,\bar{\mf B}^\gamma),(\bar{\mf A}^{\gamma+1},\bar{\mf B}^{\gamma+1}),
\cdots,(\bar{\mf A}^{\kappa-1},\bar{\mf B}^{\kappa-1}),
(\bar{\mf A}^\kappa,\bar{\mf B}^\kappa)
\eqno{(4.3\mbox{-}1)}$$ in the sense of Lemma 2.3.2. Where the $\kappa$-th pair possesses the minimal property that:
if a row of $\Theta^\kappa$ contains some $\bar b$-class arrows of $\bar {\mf B}^\kappa$,
then the same row of $F^\kappa$ contains one and only one nonzero entry which is a link in some $G^i_\kappa$
given by an edge reduction.

\begin{itemize}
\item[(i)] For $i=0, 1,\cdots, (\kappa-2)$,  the reduction from $\bar{\mf A}^{i}$ to
$\bar{\mf A}^{i+1}$  is a composition of a series of reductions $\bar{\mf
A}^{i}=\bar{\mf A}^{i,0}, \bar{\mf A}^{i,1}, \cdots,\bar{\mf A}^{i,r_{i}}$,
$\bar{\mf A}^{i,r_{i}+1}=\bar{\mf A}^{i+1}$, where.
\begin{itemize}
\item[\ding{172}] The reduction from $\bar{\mf A}^{i,j}$ to $\bar{\mf
A}^{i,j+1}$ is a sequence of regularization for non-effective
$a,b$-class arrows and finally an edge reduction of the form $(0)$
for an effective $a$ or $b$-class arrow, $0\leqslant j< r_{i}-1$. The reduction form
$\bar{\mf A}^{i,r_{i}-1}$ to $\bar{\mf A}^{i,r_i}$ is a sequence
of regularization for non-effective $a,b$-class arrows.
\item[\ding{173}] The first arrow $a_1^i:X^i\mapsto Y^i$ of $\bar{\mf B}^{i,r_i}$
is an effective $a$ or $b$-class edge with $\dz(a_1^{i})=0$. Making an edge reduction $a_1^{i}\mapsto
{{1}\choose {0}}$ or ${0\, 1}\choose {0\,0}$,
we obtain $\bar{\mf A}^{i,r_{i}+1}=\bar{\mf A}^{i+1}$.
\end{itemize}
\item[(ii)] It is possible that there exist a minimal integer $\gamma$, and an
index $1\leqslant j\leqslant r_{\gamma}+1$,
such that the first arrow of $\bar{\mf B}^{\gamma,j}$ locates outside the matrix block coming from $\bar b$,
but the first arrow of $\bar{\mf B}^{\gamma,j+1}$ locates at the first column of
the block.
\item[(iii)] The reduction from $\bar{\mf A}^{\kappa-1,0}$ to
$\bar{\mf A}^{\kappa-1,r_{\kappa-1}}$  is a composition of a series of reductions
given by (i) \ding{172}. $\bar{\mf B}^{\kappa-1,r_{\kappa-1}}$ is non-local. There are two possibilities:
\begin{itemize}
\item[\ding{172}] The first arrow
$a^{\kappa-1}_1$ of $\bar{\mf B}^{\kappa-1,r_{\kappa-1}}$ is an effective $a$ or $b$-class solid
edge with $\dz(a^{\kappa-1}_1)=0$. Making an edge reduction
$a_1^{\kappa-1}\mapsto (1)$ or $(0\; 1)$, we obtain $\bar{\mf A}^{\kappa-1,r_{\kappa-1}+1}=\bar{\mf A}^\kappa$;
\item[\ding{173}] The first arrow
$a^{\kappa-1}_1$ is an effective $b$-class loop at the
down-right corner of the matrix block
splitting from $\bar b$ with $\dz(a^{\kappa-1}_1)=0$. Making a loop reduction
$a_1^{\kappa-1}\mapsto W$, a trivial Weyr matrix over $R^{\kappa}$, we obtain
$\bar{\mf A}^{\kappa-1,r_{\kappa-1}+1}=\bar{\mf A}^\kappa$.
\end{itemize}
\end{itemize}

{\bf Proof}\, 
If the number of $\bar a$-class edges $s=0$, after a series of regularization,
we reach the unique effective loop with $\dz(\bar b)=0$. Let $\bar b\mapsto W$, the final pair
$(\bar{\mf A}^1,\bar{\mf B}^1)$ satisfies (iii) \ding{173} with $\kappa=1$.
Suppose $s>0$, we make regularization for $a_i,b_j$ before $\bar a_1$,
the corresponding $\underline v_i,\bar u_j=0$. Thus $\dz(\bar a_1)=0$ by Formula
(4.2-7), if $\bar a_1\mapsto (0)$, $\bar{\mf A}^{0,1}$ of (i) \ding{172} follows.
If $r_0>1$, repeating
the procedure in (i) \ding{172}, we finally reach $\bar{\mf A}^{0,r_0}$
with the first arrow $a^{0}_1$ and $\dz(a^{0}_1)=0$.
If $a^{0}_1$ is $\bar a$-class and $a^{0}_1\mapsto(1),(0\, 1)$, we obtain (iii) \ding{172};
If $a^{0}_1$ is $\bar b$-class and $a^{0}_1\mapsto W$, we obtain (iii) \ding{173} with $\kappa=1$.
Otherwise, if
$a^{0}_1\mapsto G^1, G^1={{1}\choose {0}}$ or ${0\, 1}\choose {0\,0}$ in the case of (i) \ding{173},
we obtain the induced pair $(\bar{\mf A}^1,\bar{\mf B}^1)$.

Suppose we have reached $(\bar{\mf A}^{i},\bar{\mf B}^{i})$ for some $i<\kappa-1$ given in (i), now
continue the reductions up to the induced pair $(\bar{\mf A}^{i+1},\bar{\mf B}^{i+1})$.
$(\bar{\mf A}^{i},\bar{\mf B}^{i})$ satisfies (BRC) by Lemma 4.3.2-4.3.3.
Suppose the first arrow of $\bar{\mf B}^{i,0}$, $a_1^{i,0}=a_{\tau n^il}$ or $b_{\tau' n^il'}$, splits from a non-effective
$a_\tau$ or $b_{\tau'}$ with $n^i$ the lowest row index, and $l$ or $l'$ the column index of the splitting block,
then $\dz(a_1^{i,0})=\underline v_{\tau n^il}$,
or $\bar u_{\tau' n^il'}$ by Formulae (4.2-7)-(4.2-8) and (4.1-8).
Thus $a_1^{i,0}\mapsto\emptyset,\underline v_{\tau n^il}=0$ or
$\bar u_{\tau' n^il'}=0$ by Remark 4.1.3 (ii). We continue regularization
for the non-effective arrows inductively, and finally
send an effective one to $(0)$ by Lemme 4.1.3 (i), then
obtain $\bar{\mf A}^{i,1}$ at (i) \ding{172}. With the similar argument as above
we reach $\bar{\mf A}^{i,r_i}$ with the first arrow
$\dz(a^{i}_1)=0$ still by 4.1.3 (i). Let $a^i_1\mapsto{{1}\choose{0}}$ or
${0\,1}\choose{0\,0}$, we obtain the $(i+1)$-the pair.

$\bar{\mf A}^{\kappa-1,r_{\kappa-1}}$ is not local, since $\bar{\mf A}^{\kappa-1,0}$
is not by Condition (BRC). If $a^{\kappa-1}_1$ is an edge, then $a^{\kappa-1}_1\mapsto(1)$ or $(0\, 1)$
gives the case (iii) \ding{172};
If $a^{\kappa-1}_1$ is a loop, then $a^{\kappa-1}_1\mapsto W$ gives the case (iii) \ding{173}.
In both cases, the induced pair $(\bar{\mf A}^\kappa,\bar{\mf B^\kappa})$
possesses the minimal property, the proof is finished.

\medskip

Suppose $s(a_1^{i-1})=X^{i-1}$ in the case of Theorem 4.3.4 (i) \ding{173},
the reduction on $a^{i-1}_1$ gives $e_{_{X^{i-1}}}\mapsto{{e_{_{Y^{i}}}\,\,
\bar w^{i}\,\,}\choose{\,\,0\,\,\,\, e_{_{X^{i}}}}}$ for $1\leqslant i<\kappa$.
Denote by $\bar W^i_\kappa$ the split of $\bar w^i$ in $e_{_{X}}^\kappa$
for $1\leqslant i<\kappa$, which can be divided into $(\kappa-i)$ blocks: denote by $n^\kappa_i$
the size of $e_{_{Y^i}}^\kappa$, and by $n^\kappa_\kappa$ that
of $e^\kappa_{_{X^\kappa}}$, which is $1$ in the case of Theorem 4.3.4 (iii) \ding{172};
or as the same as that of $W$ in the case of (iii) \ding{173}.
Thus $\bar W^i_{\kappa j}$ has the size $n^\kappa_i\times n^\kappa_{j+1}$.
Write $n^\kappa=\sum_{i=1}^\kappa n^\kappa_i$, the number of rows of $e_{_X}^\kappa$.
$$ e_{_X}^\kappa=\left (\begin{array}{ccccc} e^{\kappa}_{_{Y^1}}&\bar W^1_{\kappa1}
&\cdots &\cdots &\bar W^1_{\kappa,\kappa-1}\\
& e^{\kappa}_{_{Y^2}}&\cdots &\cdots &\bar W^2_{\kappa,\kappa-1}\\
&&\cdots&\cdots &\\
&&& e^{\kappa}_{_{Y^{\kappa-1}}}&\bar W^{\kappa-1}_{\kappa,\kappa-1}\\
&&&&e^\kappa_{_{X^\kappa}}\end{array}\right
)\eqno{(4.3\mbox{-}2)}$$
When we make a reduction, the dotted arrows appearing in $J'$,
see Lemma 2.1.2, are said to be {\it $w$-class}.
Where the dotted elements in $\bar W^i_{\kappa}$ of Formula (4.3-2) for $1\leqslant i<\kappa$,
and their splits are said to be {\it $\bar w$-class};
those in $e_{_{Y^i}}^\kappa, 1\leqslant i<\kappa$,
and $e_{_{X^\kappa}}^\kappa$ are still said to be {\it $w$-class}.
Moreover the dotted arrows first appearing
at the $\varsigma$-th step of reduction for $\varsigma>\kappa$ and
their splits are also called $w$-class.

The following facts are already implied in the proof of Theorem 4.3.4.

\medskip

{\bf Corollary 4.3.5}\, The elements in $\bar W^i_\kappa$ of (4.3-2) for $1\leqslant i<\kappa$
are dotted arrows of $\mf B^\kappa$.

\bigskip
\bigskip
\noindent {\bf 4.4\, Major pairs}
\bigskip

We will show in this sub-section, that under
some further assumption the one sided pairs given in Theorem
4.3.4 are not homogeneous.

\medskip

Let $(\bar{\mf A},\bar{\mf B})$ be a one-sided pair with
$\bar{\mf B}_X$ satisfying Classification 4.2.1 (ii), whose
number of the $\bar a$-class arrows $s\geqslant 1$.
According to the coefficients of the first two Formulae of (4.2-8),
we define $s$ linear combinations of the $\bar v$-class arrows
in $\bar{\mf B}$:
$$\begin{array}{c}\hat{v}_\tau=\sum_{j}(\bar\varepsilon_{\tau j}
-\sum_{\bar a_\tau\prec b_i\prec\bar b}\bar\alpha_i\varepsilon_{i\tau j})
\bar v_{j},\quad \tau=1,\cdots,s.\end{array}\eqno{(4.4\mbox{-}1)}$$
Fix any $1\leqslant\tau\leqslant s$, making reductions according to Theorem 4.3.4
(iii) \ding{172} for $\kappa=1$, such that $a^0_1=\bar a_\tau\mapsto(1)$, we
reach the induced pair $(\bar{\mf A}^1,\bar{\mf B}^1)$. Then we continue to do further reductions
based on formulae (4.2-7)-(4.2-8) inductively
for $\bar a_{\eta}\prec a_i,b_i\prec\bar a_{\eta+1}, \eta=\tau,\cdots,s$, and
$\bar a_s\prec a_i,b_i\prec\bar b$ by Remark 4.1.3 (ii):
$$\begin{array}{c} a_i'\mapsto\emptyset,\, \underline v_i'                                                                                                                                       +(1)\sum_j\epsilon_{i\tau j}\underline u_j'=0,\quad
b_i'\mapsto\emptyset,\, \bar u_i'+(1)\sum_j\varepsilon_{i\tau j}\bar v_j'=0;
\end{array}\eqno{(4.4\mbox{-}2)}$$
on the other hand, $\bar a_\eta\mapsto\emptyset$ or $(0)$ for $\tau<\eta\leqslant s$ according to
$\dz(\bar a_{\tau+\eta}')\ne 0$ or $=0$. Replacing $\bar u_i$ by $\bar v$-class arrows
inductively by Remark 4.1.3 (iii), the second formula of (4.2-8) gives at the induced pair $(\mf{\bar A}',\mf{\bar B}')$:
$$\begin{array}{l}\dz(\bar b')=\sum_{\bar a_\tau\prec b_i\prec\bar b}\bar\alpha_{i} \bar u_{i}'
+(1)(\sum_{j}\bar\varepsilon_{\tau j}\bar v_{j}')
=\sum_{j}(\bar \varepsilon_{\tau j}-\sum_{\bar a_\tau\prec b_i\prec\bar b}
\bar\alpha_{i}\varepsilon_{i\tau j})\bar v_{j}'=\hat v_\tau'.
\end{array}\eqno{(4.4\mbox{-}3)}$$

{\bf Lemma 4.4.1}\, Let $(\bar{\mf A},\bar{\mf B})$ be a one-sided pair with $\bar\T$ trivial and
$s\geqslant 1$, such that $\bar{\mf B}_X$ satisfies
Classification 4.2.1 (ii) and the $c$-class arrows have Formula (4.2-6).
If there exists some $1\leqslant\tau\leqslant s$, with $\hat v_\tau=0$ in
Formula (4.4-1), then $\bar{\mf B}$ is wild and non-homogeneous.

\smallskip

{\bf Proof}\, (i)\, Since $\bar{\mf B}_X$ is minimal local with $R_X=k[x,\phi(x)^{-1}]$,
set $\mathscr L'=k\setminus\{$roots of $\phi(x)\}$, there is an almost split
conflation $(e'_\lambda)$
for any $\lambda\in\mathscr L'$ in $R(\mf B_X)$.
Let $\vartheta:R(\mf B_X)\rightarrow R(\mf{\bar B})$ be the induced functor, if
$\mf{\bar B}$ is homogeneous, then there is a co-finite subset $\mathscr L\subseteq\mathscr L'$,
and a set of almost split conflations $\{(e_{\lambda})=\vartheta(e_\lambda'):
S'_\lambda\rightarrow E'_\lambda\rightarrow S'_\lambda\mid\lambda\in\mathscr L\}$
with $S_\lambda(\bar b)=(\lambda),E_\lambda(\bar b)=J_2(\lambda)$.

(ii)\, According to Formula (4.4-1)-(4.4-3), set $\bar a_\tau\mapsto (1)$, we have
$\dz(\bar b')=0$ at the induced pair $(\bar{\mf A}',\bar{\mf B}')$. Thus we are
able to construct an object $L\in R(\bar{\mf B})$ with $L_X=k,L_Y=k,L(\bar a_\tau)=(1),L(\bar b)=(\lambda)$
and others zero. Similar to the proof of Lemma 3.4.1, we obtain a contradiction, which shows that $\bar{\mf B}$
is not homogeneous, the proof is finished.

\medskip

{\bf Theorem 4.4.2}\, Let $(\bar{\mf A},\bar{\mf B})$ be a one-sided pair
with $R$ trivial and $s>1$, such that
$\bar{\mf B}_X$ satisfies Lemma 4.2.1 (ii) and the $c$-class arrows have Formula (4.2-6).
If the elements $\{\hat v_1,\hat v_2,\cdots,\hat v_s\}$ defined
in Formula (4.4-1) are linearly dependent, then $\bar{\mf B}$ is wild and non-homogeneous.

\smallskip

{\bf Proof}\, Without loos of generality, we may assume $\bar\T=\{X,Y\}$.
Suppose there is a minimally linearly dependent subset
of $\{\hat v_1,\hat{v}_2,\cdots,\hat{v}_s\}$ with $l$ vectors.
Since the case of $l=1$ has been proved by Lemma 4.4.1, we assume here that $l>1$:
$$\begin{array}{c}\{\hat v_{\tau_1}, \hat v_{\tau_2}, \cdots, \hat v_{\tau_l}\},\quad
\tau_1<\tau_2<\cdots<\tau_l,\\
\hat v_{\tau_1}=\beta_2\hat v_{\tau_2}+\cdots+\beta_r\hat v_{\tau_l}, \quad \beta_2,
\cdots,\beta_l\in k\setminus \{0\}.
\end{array}\eqno{(4.4\mbox{-}5)}$$

(i) Making reductions according to Theorem 4.3.4 (i) and (iii) \ding{172} for $\kappa=l$,
such that $a^{p}_1=\bar a_{\tau_{p},p}\mapsto{{1}\choose{0}}$
for $p=1,\cdots,l-1$, $a^l_1=\bar a_{\tau_l,l}\mapsto (1)$, we obtain an
induced pair $(\bar{\mf A}^l,\bar{\mf B}^l)$. The reduced formal product $F^l+\bar\Theta^l$
looks like (with only $\bar a,\bar b, c$-class arrows):
$$
\begin{array}{|c|c|c|c|c|c|c|c|c|c|c|} \hline
& 1& \cdots& \bar{a}_{\tau_2,1}&\cdots& \bar{a}_{\tau_3,1}&\cdots &
\bar{a}_{\tau_l,1}&\bar{a}_{\tau_l+1,1}\cdots \bar{a}_{s1}&
\bar b_{11}\ \bar b_{12}\ \cdots \ \bar b_{1l} & c_{11}\ \cdots\ c_{t1}\\
\cline{3-11}
&&&1&\cdots& \bar{a}_{\tau_3,2}&\cdots &
\bar{a}_{\tau_l,2}&\bar{a}_{\tau_l+1,2}\cdots \bar{a}_{s2}&\bar b_{21}\
\bar b_{22}\ \cdots \ \bar b_{2l} & c_{12}\ \cdots\ c_{t2}\\
\cline{5-11}
0&0&&0&&1&\cdots & \bar{a}_{\tau_l,3}&\bar{a}_{\tau_l+1,3}\cdots
\bar{a}_{s3}&\bar b_{31}\
\bar b_{32}\ \cdots \ \bar b_{3l} & c_{13}\ \cdots\ c_{t3}\\
\cline{7-11}
&&&&&&\cdots && &\cdots&\\ \cline{8-11}
 &&&&&0&&1&\bar{a}_{\tau_l+1,l}\cdots \bar{a}_{sl}&\bar b_{l1}\
\bar b_{l2}\ \cdots \ \bar b_{ll} & c_{1l}\ \cdots\ c_{tl}\\ \hline
\end{array}\eqno{(4.4\mbox{-}6)}$$
We claim that the pair $(\bar{\mf A}^l,\bar{\mf B}^l)$ is local:
$\bar{\mf A}$ has two vertices, the dimension of
$\vartheta^{0l}(F^l)$ in $R(\bar{\mf A})$ equals $l+1$, and the number of links
in $F^l$ equals $l$, the assertion follows Corollary 2.3.4 (ii).

(ii) We make further reductions from $\bar{\mf B}^l$ inductively
for the $\bar p$-th row ordered by $\bar p=l,\cdots,2$ in the reduced formal product $\Theta^l$.
For $\bar p=l$, similar to Formulae (4.4-2)-(4.4-3):
$a_{il}\mapsto\emptyset,i\in\Lambda,
b_{ilq}\mapsto\emptyset,\bar u_{ilq}+(1)\sum_j
\varepsilon_{i\tau_l j}\bar v_{jlq}=0,i<n_0;
\bar a_{\eta l}\mapsto(0)$ or $\emptyset,
\tau_l<\eta\leqslant s,
\bar b_{lq}\mapsto\emptyset, \hat v_{lq}=0,
q=1,\cdots,l$.
Next, $b_{ilq}\mapsto\emptyset$
for $i>n_0, 1\leqslant q\leqslant l$ by Remark 4.1.3 (ii); and
$c_{il}\mapsto (0)$ or $\emptyset$ inductively.
The dotted arrows $\underline v_{ip}$ and $\bar u_{ipq}$
for all $i$ and $p<l$ are preserved by 4.1.3 (iv); note that $\bar v_i: Y\mapsto X$,
the size of the split of $\bar v_i$ in $\bar{\Psi}^l$ is $l\times 1$,
the $\bar v_{ipq}$ for all $i,q$ and $p<l$ are also preserved.
The induced bocs $\bar{\mf B}^{l+1}$ follows.

(iii)\, Suppose we have reached an induced bocs $\bar{\mf B}^{2l-\bar p+1}$ for some
$\bar p\leqslant l$, which satisfies the following two conditions.
\ding{172} Denote the entrances of $F^{2l-\bar p}$ which do not belong to $\cup_{j=1}^lG_{2l-\bar p}^j$
by $\bullet^0$ coming from $\bullet$, one of the $a,b,c$-class solid arrows.
Then $a_{ip}^0=\emptyset,i\in\Lambda,\bar a_{ip}^0=\emptyset$ or $(0)$, $b_{ip}^0=\emptyset,i\ne n_0,
\bar b_{pq}=\emptyset,c_{ip}=\emptyset$ or $(0)$ for any $p>\bar p$;
\ding{173} The dotted arrows $\underline v_{ip}$,
$\underline u_{ipq}$; $\bar v_{ipq}$ are preserved for all $i,q$ and $p\leqslant\bar p$.
Now we continue to make reductions for the solid arrows at the $\bar p$-th row of
$\Theta^{2l-\bar p+1}$. We have
$a_{i\bar p}^0=\emptyset,i\in\Lambda$ and
$b_{i\bar pq}=\emptyset,\bar u_{ipq}+\bar a_{\tau_p p}^0\sum_j
\varepsilon_{i\tau_p j}\bar v_{jpq}=0,i<n_0$ inductively
according to the assumptions \ding{172}-\ding{173} and Remark 4.1.3 (ii);
while $\bar a_{\eta,\bar p}\mapsto(0)$ or $\emptyset$ for $\tau_{\bar p}<\eta\leqslant s$;
$\dz(\bar b_{\bar pq})=\hat v_{\bar pq}$,
$\bar b_{\bar pq}\mapsto\emptyset,\hat v_{\bar pq}=0,q=1,\cdots,l$. Moreover $b_{i\bar pq}\mapsto\emptyset$,
$c_{i\bar p}\mapsto(0)$ or $\emptyset$ similarly to the discussion in (ii). We finally reach the $(2l-\bar p)$-th
pair with the assumptions \ding{172}-\ding{173}.
By induction, we obtain a pair $(\bar{\mf A}^{2l-1},\bar{\mf B}^{2l-1})$.

(iv) The non-effective $a_{i1},b_{j1q}\mapsto\emptyset$ for $i\in\Lambda,j<n_0,q=1,\cdots,l$, and
$\bar a_{i1}\mapsto\emptyset$ or $(0)$ for $\tau_1<i\leqslant s$
at the first row of (4.4-6), we obtain an induced bocs $\bar{\mf B}^{2l}$.
Formula (4.4-5) gives
$\hat v_{\tau_1 1q}=\beta_2\hat v_{\tau_2 2q}+\cdots+\beta_r\hat v_{\tau_l lq}=0$
thus $\dz^0(\bar b_{1q})=0$ for $q=1,\cdots,l$. Since $l\geqslant 2$, the bocs
$\bar{\mf B}^{2l}$ is wild and non-homogeneous
by Classification 4.2.1 (iii) or (iv). And hence so is $\bar{\mf B}$.

\medskip

{\bf Definition 4.4.3}\, A one-sided pair $(\bar{\mf A},\bar{\mf B})$ with $\bar{\mf
B}_X$ satisfying Classification 4.2.1 (ii) is said to be a {\it
major pair}, provided $\{\hat v_1,\hat v_2,\cdots,\hat v_s\}$
in Formula (4.4-1) are linearly independent.

\bigskip
\bigskip
\noindent {\bf 4.5 Further reductions}
\bigskip

Throughout the subsection let $(\bar{\mf A},\bar{\mf B})$ be a one-sided major pair
having at least two vertices, such that $\bar{\mf B}_X$ satisfies Classification 4.2.1 (ii)
and the $c$-class arrows satisfy Formula (4.2-6). Suppose
$(\bar{\mf A}^\kappa,\bar{\mf B}^\kappa)$ is an induced pair
given by Theorem 4.3.4 (iii). Let $(\bar{\mf A}^\varsigma,\bar{\mf B}^\varsigma)$
be an induced pair for some
$\varsigma\geqslant\kappa$ by a series of reductions in the sense of Lemma 2.3.2.
This subsection is devoted to discussing the
reduction for $\bar{\mf B}^\varsigma$ via calculating the differential of the first arrow.

\medskip

Let $(\bar{\mf A}^\varsigma,\bar{\mf B}^\varsigma)$ be given above.
Put a solid or dotted arrow in a square box; and
a matrix block in a rectangle with four boundaries. The block $G^j_\varsigma$ for
$1\leqslant j\leqslant\kappa$ is defined in Formula (2.3-5), whose
upper boundary is that of $I^j_\varsigma$ and denoted by $m_\varsigma^{j-1}$,
the lower one is that of $F^\varsigma$, the left and right boundaries are
given by the dotted lines $l^j_\varsigma$ and $r^j_\varsigma$.
Denote by $A^j_\varsigma,B^j_\varsigma,C^j_\varsigma
\subset F^\varsigma+\Theta^\varsigma$, the sets of $a,b,c$-class entrances and solid arrows
in the $j$-th block row, on the right hand side of $I_\varsigma^j$, with
the upper (resp. lower) boundary $m_\varsigma^{j-1}$ (resp. $m_\varsigma^{j}$).
\vspace{-5mm}
\begin{center}
\begin{equation*}
\begin{array}{c} {\unitlength=0.75mm
\begin{picture}(180,75)
\put(10,5){\framebox(160,70)}
\put(20,65){\framebox(10,10)}\multiput(20,5)(0,3){20}{\line(0,1){2}}
\multiput(30,5)(0,3){20}{\line(0,1){2}}\put(20,65){\line(1,0){150}}
\put(23,67){$I^1_\varsigma$}\put(70,67){$A^1_\varsigma\cup B^1_\varsigma$}
\put(145,67){$B^1_\varsigma\cup C^1_\varsigma$}\put(22,0){$G^1_\varsigma$}
\put(40,55){\framebox(10,10)}\multiput(35,5)(0,3){20}{\line(0,1){2}}
\multiput(50,5)(0,3){17}{\line(0,1){2}}\put(40,55){\line(1,0){130}}
\put(41,57){$I^2_\varsigma$}\put(85,57){$A^2_\varsigma\cup B^2_\varsigma$}
\put(145,57){$B^2_\varsigma\cup C^2_\varsigma$}\put(38,0){$G^2_\varsigma$}
\put(80,35){\framebox(10,10)}\multiput(75,5)(0,3){13}{\line(0,1){2}}
\multiput(90,5)(0,3){11}{\line(0,1){2}}\put(90,35){\line(1,0){80}}
\put(80,37){$I^{\gamma+1}_\varsigma$}\put(110,37){$B^{\gamma+1}_\varsigma$}
\put(141,37){$B^{\gamma+1}_\varsigma\cup C^{\gamma+1}_\varsigma$}\put(78,0){$G^{\gamma+1}_\varsigma$}
\put(100,20){\framebox(10,8)}\multiput(100,5)(0,3){5}{\line(0,1){2}}
\multiput(110,5)(0,3){5}{\line(0,1){2}}\put(110,20){\line(1,0){60}}\put(110,28){\line (1,0){60}}
\put(100,23){$I^{\kappa-1}_\varsigma$}\put(120,22){$B^{\kappa-1}_\varsigma$}
\put(141,22){$B^{\kappa-1}_\varsigma\cup C^{\kappa-1}_\varsigma$}\put(100,0){$G^{\kappa-1}_\varsigma$}
{\thicklines\put(10,45){\line(1,0){160}}}{\thicklines\put(10,20){\line(1,0){160}}}
\put(125,5){\line(0,1){15}}\put(140,5){\line(0,1){70}}
\put(129, 8){$W^{\kappa}_\varsigma$}\put(145,8){$B^{\kappa}_\varsigma\cup C^{\kappa}_\varsigma$}
\put(129,0){$G^{\kappa}_\varsigma$}
\put(65,48){$\cdots\cdots\qquad\qquad\quad\cdots\cdots$}\put(150,48){$\cdots$}
\put(115,30){$\cdots$}\put(150,30){$\cdots$}
\end{picture}}\end{array}\eqno{(4.5\mbox{-}1)}\end{equation*}\end{center}

In $A^j_{\varsigma},B^j_{\varsigma},C^j_{\varsigma}$, a solid arrow is denoted by
$a_{\varsigma i, pq}^j,b_{\varsigma i,pq}^j,c_{\varsigma i,pq}^j$ splits from $a_i,b_i,c_i$
respectively, an entry of $F^\varsigma$ by $\bullet_{\varsigma i,pq}^{j,0}$,
where $(p,q)$ is the index in the
$n^\varsigma\times n^\varsigma_{t(a_i)},n^\varsigma\times n^\varsigma,
n^\varsigma\times n^\varsigma_{t(c_i)}$-block matrix respectively.
If there is no confusion, we write $\bullet^j_{ipq}$ and
$\bullet_{ipq}^{j,0}$ for simple.
For the sake of convenience, $\Phi_{\m^\varsigma}^{m}$ is partitioned also
by the lines $m,l,r$ as the same as in $F^\varsigma+\Theta^\varsigma$.

\medskip

{\bf Remark}\, From now on, we consider the pseudo reduced formal equation (4.1-7)
at the $\varsigma$-th step, in order to determine the linear relation on the dotted elements.
Keep Remark 4.1.3 in mind. Since loop or edge reduction may add some $w$-class dotted arrows,
but does not cause any new linear relations among the splits of dotter elements,
we will describe the relationship of $\bar u,\bar v, \underline u,\underline v,\bar w,w$-class
elements during the regularization from $\bar{\mf B}^\varsigma$ to $\bar{\mf B}^{\varsigma+1}$.
In the following three Lemmas, suppose the first arrow
$a^\varsigma_1=\bullet_{\tau\bar p\bar q}^\iota$ of $\bar{\mf B}^\varsigma$ belongs to $A^\iota_\varsigma
\cup B^\iota_\varsigma\cup C^\iota_\varsigma$ in Picture (4.5-1),
Write an element $\bullet^{j,0}_{ipq}$ of $F^\varsigma$
with $j\geqslant\iota$, $p>\bar p$, or $p=\bar p$ but $q<\bar q$; while write the solid arrow
$\bullet^{j'}_{i'p'q'}$ with $j'\leqslant\iota$ and $p'<\bar p$,
or $p'=\bar p$ but $q'\geqslant\bar q$.

\medskip

{\bf Lemma 4.5.1}\, Let $(\mf A^\varsigma,\mf B^\varsigma)$ be an induced pair of
$(\mf A^{\kappa},\mf B^{\kappa})$, the latter is given by
Theorem 4.3.4 (iii) \ding{173}, and $\iota=\kappa$ with
the first arrow $a^\varsigma_1=\bullet_{\tau\bar p\bar q}^{\kappa}\in
B^{\kappa}_\varsigma\cup C^{\kappa}_\varsigma$, (see below the second thick line of Picture
(4.5-1) for example). Assume that

(i)\, Any $b^{\kappa,0}_{ipq}=\emptyset,i>n_0$, the corr. element
$\bar u_{ipq}^{\kappa}$ is replaced by a combination of some $\bar v$-class arrows
in $\bar{\mf B}^\varsigma$. While the dotted arrows $\bar u_{i'p'q'}^{j'}$
are preserved.

(ii)\, If $c^{\kappa,0}_{ipq}=\emptyset$, there is a linear relation among some elements
$\underline v_{i_1p_1q_1}^\kappa,h<i_1\leqslant h+i, p_1\geqslant p$ or $p_1=p,q_1<q$ and some $\underline u,w$-class.
While all the dotted arrows $\underline v_{i'p'q'}^{j'}$ are preserved.

Then after a regularization, the induced pair $(\mf A^{\varsigma+1},\mf B^{\varsigma+1})$
still satisfies (i)-(ii). In particular all the dotted arrows $\underline v_{i'p'q'}^{j'},
\bar u^{j'}_{i'p'q'}$ except $\underline v_{\tau\bar p\bar q}^{\kappa},\bar u^{\kappa}_{\tau\bar p\bar q}$;
and all the $\bar w,\bar v$-class arrows are preserved.

\smallskip

{\bf Proof}\, The assumption is valid for $\varsigma=\kappa$ by Theorem 4.3.4 and Corollary 4.3.5.

(i)\, If $a^\varsigma_1=b^\kappa_{\tau\bar  p\bar q},\tau>n_0$,
then according to Formula (4.2-8), (4.1-8) and Remark 4.1.3 (i):
$$\begin{array}{l}
\dz(b_{\tau\bar p\bar q}^\kappa)=\bar u_{\tau\bar p\bar q}^\kappa
+\sum_{n_0<i<\tau}(\alpha^0_{\tau i}\bar u_{i\bar p\bar q}^\kappa
+\sum_{q}\alpha^1_{\tau i}\bar b^{\kappa,0}_{\bar pq}\bar u_{iq\bar q}^\kappa)
+\sum_{c_i\prec b_\tau,q}c^{\kappa,0}_{i\bar pq}(\sum_{l}\varepsilon'_{\tau il}
\bar v_{lq\bar q}).\end{array}$$
Since $W$ is upper triangular, $\bar p\leqslant q$ in $\bar b^{\kappa,0}_{\bar p q}$. By the assumption (i),
$\bar u_{\tau\bar p\bar q}^\kappa$ is a dotted arrow,
thus $b_{\tau\bar p\bar q}^\kappa\mapsto\emptyset$,
$\bar u_{\tau\bar p\bar q}^\kappa$ is replaced by a combination of some $\bar v$-class arrows,
since $\bar u_{i\bar p\bar q}^\kappa,
\bar u_{iq\bar q}^\kappa$ are already so by the assumption (i).

(ii)\, If $a^\varsigma_1=c_{\tau\bar p\bar q}^\kappa$, then according to Formula (4.2-9), (4.1-8) and Remark 4.1.3 (i):
$$\small\begin{array}{l}\dz(c_{\tau\bar p\bar q}^\kappa)=
\sum_{h<i\leqslant h+\tau}(\gamma^0_{\tau i}\underline v_{i\bar p\bar q}^\kappa
+\sum_{q}\gamma_{\tau i}^1\bar b^{\kappa,0}_{\bar pq}\underline v_{iq\bar q}^\kappa)
+\sum_{i,q}
c^{\kappa,0}_{i\bar pq}(\sum_{l}\xi_{\tau il}\underline u_{lq\bar q})
-\sum_{q<\bar q}c^{\kappa,0}_{\tau\bar pq} w_{q\bar q}+\sum_{p>\bar p} w_{\bar pp}c^0_{\tau\bar pq}.
\end{array}$$
In the case of $\dz(c_{\tau\bar p\bar q})\ne 0$,
$c_{\tau\bar p\bar q}\mapsto\emptyset$, which yields an additional linear relation
among elements $\underline v_{i\bar p\bar q},\underline v_{iq\bar q},h<i\leqslant h+\tau, q\geqslant\bar p$,
and some $\underline u,w$-class elements.

The required $\underline v,\bar u$-class and all the
$\bar v,\bar w$-class dotted arrows are preserved,
the pair $(\mf A^{\varsigma+1},\mf B^{\varsigma+1})$ still satisfies assumption
(i)-(ii), which finishes the proof.

\medskip

Suppose in Formula (4.5-1), $I^j_\varsigma,j>\gamma,$
intersects the $p$-th row of $F^\varsigma$ at the $q_p^j$-th column with $\bar b^{j,0}_{pq_p^j}=(1)$.
Denote by $\bar w_{q_{p}^{j} q}^j$ for any possible $q$ the dotted element with row index $q_p^j$
at $e_{_X}^\varsigma$. If $j>\iota$, or $j=\iota$ but $p>\bar p$, then $q_{p}^{j}>q_{\bar p}^{\iota}$,
$\bar w_{q_{p}^{j} q}^j$ is sitting below $\bar w_{q_{\bar p}^{\iota}\bar q}^\iota$,
we refer to Picture (4.5-2).

\medskip

{\bf Lemma 4.5.2}\, Let $(\bar{\mf A}^\varsigma,\bar{\mf B}^\varsigma)$ be an induced pair of
$(\bar{\mf A}^\kappa,\bar{\mf B}^\kappa)$ with $\gamma$ existing in Theorem 4.3.4 (ii).
Suppose the first arrow $a^\varsigma_1=\bullet_{\tau\bar p\bar q}^\iota\in
B^\iota_\varsigma\cup C^\iota_\varsigma$ with
$\gamma<\iota\leqslant\kappa$ in 4.3.4 (iii) \ding{172}, or $\gamma<\iota<\kappa$, in \ding{173},
(see between the two thick lines of Picture (4.5-1) for example). Assume that

(i)\, Any $\bar b^{j,0}_{pq}=\emptyset$, the
corr. $\bar w_{q_{p}^{j} q}^j$ is replaced by a combination of some
$\bar w_{lq},l>q_p^j$ and $w$-class elements at $\bar{\mf B}^\varsigma$.
While the dotted arrows $\bar w_{p'q'}$ with
$p'\leqslant q_{\bar p}^{\iota}$ are preserved.

(ii)\, Any $b^{j,0}_{ipq}=\emptyset,i>n_0$, the corr.
$\bar u_{ipq}^{j}$ is replaced by a combination of
some $\bar u_{i_1p_1q_1}^{j_1}\succ\bar u_{ipq}^{j},n_0<i_1<i$, and $\bar v$-class elements
in $\bar{\mf B}^\varsigma$. While the dotted arrows $\bar u_{i'p'q'}^{j'}$
are preserved.

(iii)\, If $c^{j,0}_{ipq}=\emptyset$, there is a relation among elements
$\underline v_{i_1p_1q_1}^{j_1},h<i_1\leqslant h+i$, and $\underline u,w$-class.
While the dotted arrow $\underline v_{i'p'q'}^{j'}$, $i'\in\Lambda$,
are preserved.

Then after a regularization, the induced pair $(\mf A^{\varsigma+1},\mf B^{\varsigma+1})$
still satisfies (i)-(iii). In particular all the dotted arrows $\underline v_{i'p'q'}^{j'},i'\in\Lambda$;
$\bar u^{j'}_{i'p'q'}$ except $\bar u^{j}_{\tau\bar p\bar q}$;
$\bar w^{j'}_{p'q'},p'<q_{\bar p}^j$; and $\bar v$-class are preserved.

\smallskip

{\bf Proof}\, The assumption (i)-(iii) are valid, if the block of $a_1^\varsigma$ has
the bottom and right boundaries $(m_\varsigma^\kappa,r_{\varsigma}^\kappa)$
in the case of Theorem 4.3.4 (iii) \ding{172}; or $(m_\varsigma^{\kappa-1},r_{\varsigma}^{\kappa-1})$
in (iii) \ding{173} by Lemma 4.5.1.

(i)\, If $a^\varsigma_1=\bar b_{\bar p\bar q}^\iota$, by the notation
$\bar b^{\iota,0}_{\bar pq_{\bar p}^\iota}=(1)$, and by the assumption (i)
$\bar b^{\iota,0}_{\bar pq}=\emptyset,q<\bar q$,
$$\begin{array}{ll}\dz(\bar b^\iota_{\bar p\bar q})&=
(1)\bar w^\iota_{q^{\iota}_{\bar p} \bar q}
+\sum_{j\geqslant\iota,p>\bar p}w_{\bar pp}\bar b^{j,0}_{p\bar q},
\quad w_{\bar pp}\,\,\mbox{belongs to}\,\,\bar w\,\,\mbox{or}
\,\,w\,\,\mbox{class}.\end{array}$$
Since $\bar w_{q^{\iota}_{\bar p} \bar q}^\iota$ is a dotted arrow,
$\bar b^\iota_{\bar p\bar q}\mapsto\emptyset$ at the induced
pair $(\mf A^{\varsigma+1},\mf B^{\varsigma+1})$. As an example, we refer to
$\bar W^3_\varsigma$ and $\bar W^4_\varsigma$ in the picture (4.5-2) below.

(ii)\, If $a^\varsigma_1=b^\iota_{\tau\bar p\bar q},\tau>n_1$, then
$\bar b^{\iota,0}_{\bar pq_{\bar p}^\iota}=(1)$,
$\bar b^{\iota,0}_{\bar pq}=\emptyset,\forall\,q>q_{\bar p}^\iota$ by (i) above;
and for some $\iota_1\leqslant\iota$:
$$\begin{array}{l}
\dz(b_{\tau\bar p\bar q}^\iota)=\bar u_{\tau\bar p\bar q}^\iota
+\sum_{n_0<i<\tau}(\alpha^0_{\tau i}\bar u_{i\bar p\bar q}^\iota
+\alpha^1_{\tau i}\bar b^{\iota,0}_{\bar pq^\iota_{\bar p}}\bar u_{iq^\iota_{\bar p}\bar q}^{\iota_1})
+\sum_{c_i\prec b_\tau,q}c^0_{i\bar pq}(\sum_{l}\varepsilon'_{\tau il}
\bar v_{lq\bar q}),\,\,b_{\tau\bar p\bar q}^\iota\mapsto\emptyset,\end{array}$$
since $\bar u_{\tau\bar p\bar q}^\iota$ is a dotted arrow
by the assumption (ii), which is replaced by some $\bar u,\bar v$-class elements.

(iii) If $a^\varsigma_1=c^\iota_{\tau\bar p\bar q}$, in case of $c^\iota_{\tau\bar p\bar q}\mapsto\emptyset$,
a relation among $\underline u,\underline v,i>h,w$-class elements is added:
$$\small\begin{array}{l}\dz(c_{\tau\bar p\bar q}^\iota)=
\sum_{h<i\leqslant h+\tau}(\gamma^0_{\tau i}\underline v_{i\bar p\bar q}^\iota
+\gamma_{\tau i}^1\bar b^{\iota,0}_{\bar pq^\iota_{\bar p}}\underline v_{iq^\iota_{\bar p}\bar q}^{\iota_1})
+\sum_{i,q}c^{\iota,0}_{i\bar pq}(\sum_{l}\xi_{\tau il}\underline u_{lq\bar q})
-\sum_{q<\bar q}c^{\iota,0}_{\tau\bar pq} w_{q\bar q}+\sum_{p>\bar p} w_{\bar pp}c^{j,0}_{\tau\bar pq}.
\end{array}$$

The required $\underline v,\bar u,\bar w$-class and all the
$\bar v$-class dotted arrows are preserved,
the pair $(\mf A^{\varsigma+1},\mf B^{\varsigma+1})$ still satisfies assumption
(i)-(iii), which finishes the proof.

\medskip

Suppose in the picture (4.5-1), $I^j_\varsigma$ comes from an effective $\bar a_{i^j}$
for $j\leqslant\gamma$ if $\gamma$ exists, otherwise $j\leqslant\kappa$
in Theorem 4.3.4 (iii) \ding{172}, or $j<\kappa$ in (iii) \ding{173}, which intersects the $p$-th row at the
$q_p^j$-th column in the $n^\varsigma_{_X}\times n^\varsigma_{t(a_{i^j})}$-block
coming from $\bar a_{i^j}$ partitioned by $\bar\T$ with $\bar a^{j,0}_{i^j,pq_p^j}=(1)$.
Denote by $\hat v_{i^j,q_{p}^{j} q}^j$ splitting from $\hat v_{i^j}$ for any possible $q$, whose row index is $q_p^j$
in the $n^\varsigma_{t(a_{i^j})}\times n^\varsigma_{_X}$-block.
We write $\hat V^j_\varsigma$ at the $(n_0,h_i)$-th block
partitioned by $\bar\T$ in $\bar\Psi^{r}_{\n^\varsigma}$,
and refer to Picture (4.5-2).

\medskip

{\bf Lemma 4.5.3}\, Let $(\bar{\mf A}^\varsigma,\bar{\mf B}^\varsigma)$ be an induced bocs of
$(\bar{\mf A}^\kappa,\bar{\mf B}^\kappa)$. Suppose
the first arrow $a^\varsigma_1=\bullet_{\tau\bar p\bar q}^\iota\in
A^\iota_\varsigma\cup B^\iota_\varsigma\cup C^\iota_\varsigma$ with
$\iota\leqslant\gamma$ if $\gamma$ exists; otherwise $j\leqslant\kappa$
in Theorem 4.3.4 (iii) \ding{172}, or $j<\kappa$ in (iii) \ding{173},
(see above the first thick line of Picture (4.5-1) for example). Assume that

(i)\, Any $a^{j,0}_{ipq}=\emptyset,i\in\Lambda$,
the corr. $\underline v_{ipq}^j$ is replaced by a
combination of some $\underline u$-class elements at
$\bar{\mf B}^\varsigma$. While all the dotted arrows
$\underline v_{i'p'q'}^{j'}$ are preserved.

(ii)\, If an effective $\bar a^{j,0}_{ipq}=\emptyset$, there is a linear relation among
$\underline u,\bar w,w$-class elements at $\bar{\mf B}^{\varsigma}$. While all the
dotted arrows $\underline v_{i'p'q'}^{j'},i'\in\Lambda$, are preserved.

(iii)\, Any $b^{j,0}_{ipq}=\emptyset,i<n_0$, the corresponding $\bar u_{ipq}^j$ is replaced by a
combination of some $\bar v$-class elements at $\bar{\mf B}^\varsigma$. While all the dotted
arrows $\bar u_{i'p'q'}^{j'}$ are preserved.

(iv)\, Any $\bar b^{j,0}_{pq}=\emptyset$, $\bar v_{i^jq_{p}^{j} q}$ corr. to
$\bar a^{j,0}_{i^jpq_p^j}$ is replaced by a combination of some $\bar v$-class elements
below and some $\bar w, w$-class at $\bar{\mf B}^\varsigma$.
While the dotted arrows $\bar v_{i'p'q'}^{j'}$ with $p'\leqslant q_{\bar p}^{\iota}$ are preserved.

(v)\, Any $b^{j,0}_{ipq}=\emptyset,i>n_0$, the corr. element $\bar u_{ipq}^j$
is replaced by a combination of some $\bar u_{i_1pq}^{j}$, $n_0<i_1<i$, and $\bar v$-elements at
$\bar{\mf B}^\varsigma$. While all the dotted arrows $\bar u_{i'p'q'}^{j'}$ are preserved.

(vi)\, If $c^{j,0}_{ipq}=\emptyset$, there is
a linear relation among some elements
$\underline v_{i_1pq}^{j},h<i_1<h+\tau$,
and $\underline u,w,\bar w$-class at $\bar{\mf B}^\varsigma$. While all the dotted arrows
$\underline v_{i'p'q'}^{j'},i'\in\Lambda$, are preserved.

Then after a regularization, the induced pair $(\mf A^{\varsigma+1},\mf B^{\varsigma+1})$
still satisfies (i)-(vi). In particular, all the dotted arrows $\underline v_{i'p'q'}^{j'},i'\in\Lambda$,
except $\underline v_{\tau\bar p\bar q}^{j}$; $\bar u^{j'}_{i'p'q'}$ except $\bar u^{j}_{\tau\bar p\bar q}$;
and $\hat v^{j'}_{p'q'},p'<q_{\bar p}^\iota$ are preserved.

\smallskip

{\bf Proof}\, The assumption (i)-(vi) are valid, if $a_1^\varsigma$ has
the bottom and right boundaries $(m_\varsigma^\gamma,r_{\varsigma}^\gamma)$
when $\gamma$ exists by Lemma 4.5.2, otherwise $(m_\varsigma^\kappa,r_{\varsigma}^\kappa)$ in
Theorem 4.3.4 (iii) \ding{172}, or $(m_\varsigma^{\kappa-1},r_{\varsigma}^{\kappa-1})$
in (iii) \ding{173} by Lemma 4.5.1.

(i)\, If $a^\varsigma_1=a^\iota_{\tau\bar p\bar q}$ with $\tau\in\Lambda$, since
$\underline v_{\tau\bar p\bar q}$ is a dotted arrow by the assumption (i):
$$\begin{array}{c}\dz(a^\iota_{\tau\bar p\bar q})=\underline v^\iota_{\tau\bar p\bar q}+
\sum_{i,q}\bar a^{\iota,0}_{i\bar pq}(\sum_{l}\epsilon_{\tau il}\underline u_{lq\bar q})\,\,\Longrightarrow
\,\,a_{\tau\bar p\bar q}^\iota\mapsto\emptyset,\,\, v_{\tau\bar p\bar q}^\iota=
-\sum_{i,q}\bar a^{\iota,0}_{i\bar pq}
(\sum_{l}\epsilon_{\tau il}\underline u_{lq\bar q}),\end{array}$$

(ii)\, If $a^\varsigma_1=\bar a^\iota_{\tau pq}$ effective, by
substituting $\underline v^\iota_{i'\bar p\bar q}$ given in (i):
$$\begin{array}{c}\dz(\bar a_{\tau\bar p\bar q}^\iota)=\sum_{\bar a_{\tau^\iota}\prec a_{i'}
\prec\bar a_{\tau}}\bar\beta_{\tau i'}\underline v^\iota_{i'\bar p\bar q}
+\sum_{i,q}\bar a^{\iota,0}_{i\bar pq}(\sum_{l}\bar\epsilon_{\tau il}\underline u_{lq\bar q})
-\sum_{q<\bar q}\bar a^{\iota,0}_{\tau\bar pq} w_{q\bar q}
+\sum_{j\geqslant\iota,p>\bar p} w_{\bar pp}\bar a^{j,0}_{\tau\bar pq}\\
=\sum_i\bar a^{\iota,0}_{i\bar pq}\big(\sum_{l,q}(\bar\epsilon_{\tau il}-
\sum_{\bar a_{\tau^\iota}\prec a_{i'}\prec\bar a_i\preccurlyeq\bar a_{\tau}}
\bar\beta_{\tau i'}\epsilon_{i'il})\underline u_{lq\bar q}\big)
-\sum_{q<\bar q}\bar a^{\iota,0}_{\tau\bar pq} w_{q\bar q}+\sum_{j\geqslant\iota,p>\bar p}
w_{\bar pp}\bar a^{j,0}_{\tau\bar pq}.\end{array}$$
If $\dz(\bar a^\iota_{\tau\bar p\bar q})\ne 0$,
$\bar a^\iota_{\tau\bar p\bar q}\mapsto\emptyset$
yields a linear relation among some $\underline u$, $w$ and $\bar w$-class elements.

(iii)\, If $a^\varsigma_1=b^\iota_{\tau\bar p\bar q},\tau<n_0$,
since $\bar u_{\tau\bar p\bar q}^\iota$ is a dotted arrow by the assumption (iii):
$$\begin{array}{c}\dz(b^\iota_{\tau\bar p\bar q})=\bar u^\iota_{\tau\bar p\bar q}+
\sum_{i,q} \bar a^{\iota,0}_{i\bar pq}(\sum_{l}\varepsilon_{\tau il}\bar v_{lq\bar q})\,\,\Longrightarrow\,\,
b^\iota_{\tau\bar p\bar q}\mapsto\emptyset,\,\, \bar u^\iota_{\tau\bar p\bar q}=-
\sum_{i,q}\bar a^{\iota,0}_{i\bar pq}(\sum_{l}\varepsilon_{\tau il}\bar v_{lq\bar q}).\end{array}$$

(iv)\, If $a^\varsigma_1=\bar b^\iota_{\bar p\bar q}$ effective,
by substituting $\bar u^\iota_{i'\bar p\bar q}$ given in (iii), and Formula (4.4-1):
$$\begin{array}{ll}\bar\dz(\bar b^\iota_{\bar p\bar q})&=\sum_{i'<n_0}
\bar\alpha_{i'} \bar u_{i'\bar p\bar q}^\iota+\sum_{i,q}\bar a^{\iota,0}_{i\bar pq}
(\sum_{l}\bar\varepsilon_{il}\bar v_{lq\bar q})-\sum_{q<\bar q}\bar b^{\iota,0}_{\bar pq}
w_{q\bar q}+\sum_{j\geqslant\iota,p>\bar p} w_{\bar pp}\bar b^{j,0}_{p\bar q}\\
&=\sum_{i,q}\bar a_{i\bar pq}^{\iota,0}\big(\sum_{l}
(\bar\varepsilon_{il}-\sum_{\bar a_{\tau^\iota}
\preccurlyeq\bar a_{i'}\prec b_i\prec\bar b}\bar\alpha_{i'}\varepsilon_{i'il})\bar v_{lq\bar q}\big)
-\sum_{q<\bar q}\bar b^{\iota,0}_{\bar pq}
w_{q\bar q}+\sum_{j\geqslant\iota,p>\bar p} w_{\bar pp}\bar b^{j,0}_{p\bar q}\\
&=\hat v_{\tau^\iota, q^\iota_{\bar p}\bar q}+\sum_{\bar a_{i\bar pq}^{\iota,0}
\succ\bar a^{\iota,0}_{\tau^\iota,\bar p q^\iota_{\bar p}};q}\bar a_{i\bar pq}^{\iota,0}
\hat v_{iq\bar q}-\sum_{q<\bar q}\bar b^{\iota,0}_{\bar pq}
w_{q\bar q}+\sum_{j\geqslant\iota,p>\bar p} w_{\bar pp}\bar b^{j,0}_{p\bar q},\end{array}$$
since $\bar a^{\iota,0}_{\tau^\iota,\bar p q^\iota_{\bar p}}=1$,
$\hat v_{\tau^\iota q^\iota_{\bar p}\bar q}$ is a dotted arrow by the assumption (iv),
$\bar b^\iota_{\bar p\bar q}\mapsto\emptyset$,
$\hat v_{\tau^\iota q^\iota_{\bar p}\bar q}$ is replaced by some
$\bar w,w$-class; and some $\hat v$-class elements
below $q_{\bar p}^\iota$.
See $\hat V^1_\varsigma,\hat V^2_\varsigma$ in Picture (4.5-2).

(v)\, If $a^\varsigma_1=b^\iota_{\tau\bar p\bar q},\tau>n_0$,
since $\bar b_{\bar pq}^0=\emptyset$ for all possible $q$ by (iv) above:
$$\begin{array}{l}\dz(b^\iota_{\tau\bar p\bar q})=\bar u_{\tau\bar p\bar q}^\iota
+\sum_{i\ne n_0,i<\tau}\alpha^0_{\tau i}\bar u^\iota_{i\bar p\bar q}
+\sum_{i,q}\bar a^{\iota,0}_{i\bar pq}
(\sum_{l}\varepsilon_{\tau il}\bar v_{lq\bar q})+\sum_{c_i\prec b_\tau;q}\bar c^{\iota,0}_{i\bar pq}
(\sum_{l}\varepsilon'_{\tau il}\bar v_{lq\bar q}).\end{array}$$
Since $\bar u^\iota_{\tau\bar p\bar q}$ is a dotted arrow by the assumption (v),
$b_{\tau\bar p\bar q}\mapsto\emptyset$, $\bar u^\iota_{\tau\bar p\bar q}$ is replaced
by some $\bar u_{i\bar p\bar q},n_0<i<\tau$, and $\bar v$-elements by a replacement in (iii).

(vi)\, If $a_1^\varsigma=c_{\tau\bar p\bar q}^\iota$, since $\bar b_{\bar pq}^0=\emptyset$
for all possible $q$, by (iv):
$$\begin{array}{ll}\dz(c^\iota_{\tau\bar p\bar q})&=\sum_{i\leqslant h+\tau}\gamma^0_{\tau i}
\underline v_{i\bar p\bar q}^\iota+\sum_{i,q}\bar a^{\iota,0}_{i\bar pq}
(\sum_{l}\zeta_{\tau il}\underline u_{lq\bar q})\\
&+\sum_{i<\tau;q}\bar c^{\iota,0}_{i\bar pq}(\sum_{l}\xi_{\tau il}\underline u_{lq\bar q})
-\sum_{q<\bar q}c^{\iota,0}_{\tau\bar pq} w_{q\bar q}
+\sum_{j\geqslant\iota,p>\bar p} w_{\bar pp}c^{j,0}_{\tau p\bar q}.\end{array}$$
Suppose $\dz(c^\iota_{\tau\bar p\bar q})\ne 0$,
$c^\iota_{\tau\bar p\bar q}\mapsto\emptyset$ causes a linear relation among
$\underline v^\iota_{i\bar p\bar q},h<i\leqslant h+\tau$, and $\underline u,w,\bar w$-class elements
by a replacement in (i).

The required $\underline v,\bar u,\bar w,\bar v$-class dotted arrows are preserved,
the pair $(\bar{\mf A}^{\varsigma+1},\bar{\mf B}^{\varsigma+1})$ still satisfies assumption
(i)-(vi), which finishes the proof.

\medskip

The following picture shows a pseudo formal equation $\bar\Theta^\varsigma$ of $(\bar{\mf
A}^\varsigma,\bar{\mf B}^\varsigma)$ for $\kappa=5, \gamma=2$ in the case of
Theorem 4.3.4 (iii) \ding{173}  with only effective arrows.
From this we are able to see the correspondence
of $(\bar B^1_\varsigma,\hat V^1_\varsigma), (\bar B^2_\varsigma,\hat V^2_\varsigma)$ and
$(\bar B^3_\varsigma,\bar W^3_\varsigma),(\bar B^4_\varsigma,\bar W^4_\varsigma)$:

\medskip

\vspace{7mm}
\begin{equation*}\begin{array}{c} {\unitlength=0.8mm
\begin{picture}(185,60)\thicklines
\put(0,16){\framebox(30,30){}} \put(0,40){\line(1,0){30}}
\put(6,34){\line(1,0){24}} \put(12,28){\line(1,0){18}}
\put(18,22){\line(1,0){12}}

\put(6,34){\line(0,1){12}} \put(12,28){\line(0,1){18}}
\put(18,22){\line(0,1){24}} \put(24,16){\line(0,1){30}}

\put(0,42){\makebox{\small$e^\varsigma_{_{Y^1}}$}} \put(6,36){\makebox{\small$e^\varsigma_{_{Y^2}}$}}
\put(12,30){\makebox{\small$e^\varsigma_{_{Y^3}}$}} \put(18,24){\makebox{\small$e^\varsigma_{_{Y^4}}$}}
\put(24,18){\makebox{\small$e^\varsigma_{_{X^5}}$}} \put(22,32){\makebox{$\bar W$}}

\thicklines \put(40,16){\framebox(70,30){}}
\put(52,40){\line(0,1){6}} \put(58,34){\line(0,1){6}}
\put(70,28){\line(0,1){6}} \put(82,22){\line(0,1){6}}
\put(94,16){\line(0,1){6}} \put(52,40){\line(1,0){6}}
\put(58,34){\line(1,0){12}}\put(70,28){\line(1,0){12}}
\put(82,22){\line(1,0){12}}

\thinlines \put(52,16){\line(0,1){30}} \put(64,16){\line(0,1){30}}
\put(94,16){\line(0,1){30}} \put(58,16){\line(0,1){30}}
\put(46,16){\line(0,1){30}}

\put(40,40){\line(1,0){70}} \put(52,34){\line(1,0){58}}
\put(64,28){\line(1,0){46}} \put(76,22){\line(1,0){34}}

\put(70,16){\line(0,1){18}} \put(76,16){\line(0,1){12}}
\put(82,16){\line(0,1){12}} \put(88,16){\line(0,1){6}}

\put(47,42){\makebox{$I^1$}}\put(42,42){\makebox{$0$}}
\put(53,36){\makebox{$I^2$}} \put(65,30){\makebox{$I^3$}}
\put(77,24){\makebox{$I^4$}}
\put(88,17.5){\makebox{$\scriptstyle{W}$}}
\put(56,41){\makebox{\small$\bar A^1_2$}} \put(58,35){\makebox{\small$\bar A_{2}^2$}}
\put(78,42){\makebox{\small$\bar B^1$}}
\put(78,36){\makebox{\small$\bar B^2$}}
\put(78,30){\makebox{\small$\bar B^3$}}
\put(87,24){\makebox{\small$\bar B^4$}} \put(95,41){\makebox{\small$C_{1}^1$}}
\put(95,35){\makebox{\small$C_{1}^2$}} \put(95,29){\makebox{\small$C_{1}^3$}}
\put(95,23){\makebox{\small$C_{1}^4$}} \put(95,17){\makebox{\small$C_{1}^5$}}
\put(103,41){\makebox{\small$C_{2}^1$}} \put(103,35){\makebox{\small$C_{2}^2$}}
\put(103,29){\makebox{\small$C_{2}^3$}} \put(103,23){\makebox{\small$C_{2}^4$}}
\put(102,17){\makebox{\small$C_{2}^5$}}

\mput(64,41.5)(0,1.5){4}{\line(1,0){30}}
\mput(65.5,40)(1.5,0){20}{\line(0,-1){6}}
\mput(70,29.5)(0,1.5){3}{\line(1,0){24}}
\mput(83.5,28)(1.5,0){8}{\line(0,-1){6}}
\put(102,16){\line(0,1){30}}

\thicklines \put(120,0){\framebox(70,70){}}
\mput(120,58)(0,6){2}{\line(1,0){70}} \put(132,52){\line(1,0){58}}
\put(132,46){\line(1,0){58}} \put(144,16){\line(1,0){46}}
\put(174,8){\line(1,0){16}} \put(132,70){\line(0,-1){24}}
\put(144,70){\line(0,-1){54}} \put(174,70){\line(0,-1){62}}
\put(182,70){\line(0,-1){70}}

\put(144,40){\line(1,0){30}} \put(150,34){\line(1,0){24}}
\put(156,28){\line(1,0){18}} \put(162,22){\line(1,0){12}}

\put(150,34){\line(0,1){12}} \put(156,28){\line(0,1){18}}
\put(162,22){\line(0,1){24}} \put(168,16){\line(0,1){30}}

\thinlines\mput(144,59.5)(0,1.5){4}{\line(1,0){30}}
\mput(145.5,52)(1.5,0){20}{\line(0,1){6}}
\mput(150,41.5)(0,1.5){4}{\line(1,0){24}}
\mput(163.5,28)(1.5,0){7}{\line(0,1){6}}

\put(126,70){\line(0,-1){12}} \put(138,58){\line(0,-1){12}}

\put(144,41){\makebox{\small$e^\varsigma_{_{Y^1}}$}} \put(150,35){\makebox{\small$e^\varsigma_{_{Y^2}}$}}
\put(156,29){\makebox{\small$e^\varsigma_{_{Y^3}}$}} \put(162,23){\makebox{\small$e^\varsigma_{_{Y^4}}$}}
\put(168,17){\makebox{\small$e^\varsigma_{_{X^5}}$}} \put(159,40){\makebox{\small$\bar W^3$}}
\put(165,28){\makebox{\small$\bar W^4$}}

\put(123,62){\makebox{$T^\varsigma_1$}} \put(135,50){\makebox{$T^\varsigma_2$}}
\put(176,10){\makebox{$T^\varsigma_3$}} \put(183,1){\makebox{$T^\varsigma_4$}}
\put(155.5,59){\makebox{\small$\hat V^1$}} \put(155.5,47.5){\makebox{$\,$}}
\put(156,53){\makebox{\small$\hat V^2$}} \put(134.5,62){\makebox{$\underline U$}}
\put(179.5,53){\makebox{$\underline U$}} \put(179.5,27){\makebox{$\underline V$}}
\put(185,10){\makebox{\small$\underline U$}}

\put(63,14){\makebox{$\underbrace{\hskip 25mm}$}}
\put(77,7){\makebox{$\,$}} \put(42,9){\makebox{$\,$}}
\put(55,9){\makebox{$\,$}} \put(95,9){\makebox{$\,$}}
\put(105,9){\makebox{$\,$}} \put(13,9){\makebox{$\,$}}
\put(153,-5){\makebox{$\,$}}
\put(-0.5,46.5){\makebox{$\overbrace{\hskip 24mm}$}}
\put(52,46.5){\makebox{$\overbrace{\hskip 9mm}$}}
\put(64,46.5){\makebox{$\overbrace{\hskip 24mm}$}}
\put(94.5,46.5){\makebox{$\overbrace{\hskip 12.5mm}$}}
\put(119.5,70.5){\makebox{$\overbrace{\hskip 56mm}$}}

\put(13,50){\makebox{$e_{_X}^\varsigma$}} \put(56,50){\makebox{$A^\varsigma$}}
\put(77,51){\makebox{$B^\varsigma$}} \put(101,50){\makebox{$C^\varsigma$}}
\put(153,74){\makebox{$T^\varsigma$}}

\put(77,-4){\makebox{$\mbox{Picture}\,(4.5\mbox{-}2)$}}
\end{picture}}
\end{array}
\end{equation*}

\medskip
\bigskip
\bigskip
\noindent {\bf 4.6 Regularization for non-effective $a$ class and all the $b$ class arrows}
\bigskip

Let a one-sided pair $(\bar{\mf A},\bar{\mf B})$ and an induced pair
$(\bar{\mf A}^\varsigma,\bar{\mf B}^\varsigma)$ be given in the beginning of section 4.5.
It is clear by Lemma 4.5.1-4.5.3, that all the $b$-class and
non-effective $a$-class solid arrows go to $\emptyset$ by regularization.
Now assume that we meet a local pair $(\bar{\mf A}',\bar{\mf B}')$ induced after the $\varsigma$-th
step still in the sense of Lemma 2.3.2,
which has an induced pair in Classification 3.3.5 (iv):
$$\begin{array}{c}(\bar{\mf A}'_{(\lambda_0,\lambda_1,\cdots,\lambda_{l-1})},
\bar{\mf B}'_{(\lambda_0,\lambda_1,\cdots,\lambda_{l-1})})\end{array}\eqno{(4.6\mbox{-}1)}$$
Denote by $(\bar{\mf A}^{s},\bar{\mf B}^{s})$ an induced pair of the pair (4.6-1) after
a series of regularization with the first loop
$a^s_1$ and $\dz(a^s_1)=0$. Make a loop mutation,
we obtain an induced pair $(\bar{\mf A}^{s+1},\bar{\mf B}^{s+1})$ with $R^{s+1}=k[x]$.
We claim in particular that $W$ appearing in Theorem 4.3.4 (iii) \ding{173} must be trivial,
since $\bar{\mf B}^{\kappa-1,r_{\kappa-1}}$ is not local, but $x$
appears only in a local bocs in the case of MW5.

Suppose $\mf B^{s+1}$ satisfies Formulae (3.3-3)-(3.3-4), an induced pair $(\bar{\mf A}^t,\bar{\mf A}^t)$
is in the case of MW5, then the non-effective $a$-class
and all the $b$-class solid arrows are regularized during the reductions. In particular,
the parameter $x$ and the first arrow $a^t_1$ of $\bar{\mf B}^t$ belongs to $\bar a$ or $c$-class.
In fact, the discussion of 4.5.1-4.5.3 is still valid if we describe the linear relationship
over the fractional field $k(x)$ of the polynomial ring $k[x]$, or over the field $k(x,x_1)$
of two indeterminants instead of the base field $k$.

\medskip

{\bf Theorem 4.6.1}\, Let $(\bar{\mf A},\bar{\mf B})$ be a one sided pair
having at least two vertices, such that the induced local bocs $\bar{\mf B}_X$ satisfies
Classification 4.3.1 (ii), where the pair is major and the $c$-class arrows satisfy Formula (4.2-6).
If $(\bar{\mf A},\bar{\mf B})$ has an induced pair $(\bar{\mf A}^t,\bar{\mf B}^t)$ in the case of MW5,
then the parameter $x$ and the first arrow $a^t_1$ must split from some $\bar a$ or $c$-class arrows.

\medskip

Finally, let $(\bar{\mf A},\bar{\mf B})$ be a one-sided pair
having at least two vertices, such that $\bar{\mf B}_X$ satisfies Classification 4.2.1 (i).
Then $\bar{\mf B}$ has only $a,b$-class solid arrows, where $b_1,\cdots,b_n$
are all non-effective; and $a_i,i\in\Lambda$ satisfying the first formula of
(4.2-5) are non-effective, while $\bar a_i=a_{h_i}, i=1,\cdots,s$, satisfying the second
one are effective. If there is an induced pair $(\bar{\mf A}^t,\bar{\mf B}^t)$
in the case of MW5, we have the following observation.

(i)\, Let $(\bar{\mf A}',\bar{\mf B}')$ be an induced pair of $(\bar{\mf A},\bar{\mf B})$,
with $\bar\T'$ being trivial and the reduced formal product
$\bar\Theta'$ given by Formula (4.1-6). Similarly to condition 4.3.1, let
$\mathcal D=\{d_1,\cdots,d_r\}$ be a set of solid arrows, $\mathcal E=\{e_1,\cdots,e_s\}$
that of edges, such that $\mathcal D\cup\mathcal E$ form the lowest
row of $\bar\Theta'$. Let $\mathcal U=\{u_1,\cdots,u_r\}$ be a set of dotted arrows,
with $\mathcal W=\V'\setminus\U$, such that
$\bar\dz(d_i)$ and $\bar\dz(e_i)$ satisfy the formulae in 4.3.1 (ii),
we obtain (BRC)$'$. Then after a reduction given by Lemma 4.3.2, the induced pair still
satisfies (BRC)$'$; and the original pair $(\bar{\mf A},\bar{\mf B})$ satisfies (BRC)$'$
similar to Lemma 4.3.3, but the proofs are much easier than those.

(ii)\, For constructing a reduction sequence of $(\bar{\mf A},\bar{\mf B})$ up to
$(\bar{\mf A}^\kappa,\bar{\mf B}^\kappa)$,
we need only the part (i) and (iii) \ding{172} of Theorem 4.3.4.

(iii)\, For the further reductions, we need only Theorem 4.5.3 (i)-(iii),
then reach an induced pair $(\bar{\mf A}^s,\bar{\mf B}^s)$, where all the $b$-class,
non-effective $a$-class arrows are regularized step by step.

\medskip

{\bf Corollary 4.6.2}\, Let $(\bar{\mf A},\bar{\mf B})$ be a one sided pair
having at least two vertices, such that the induced local bocs $\bar{\mf B}_X$ satisfies
Classification 4.3.1 (i). If $(\bar{\mf A},\bar{\mf B})$ has an induced pair $(\bar{\mf A}^t,\bar{\mf B}^t)$ satisfying MW5,
then the parameter $x$ and the first arrow $a^t_1$ must split from some $\bar a$-class arrows.

\bigskip
\bigskip
\bigskip

\centerline {\bf 5. Non-homogeneity of wild bipartite problems}

\bigskip

This section is devoted to proving the main Theorem 0.3. The way to do
so is based on proving the non-homogeneous property for a wild
bipartite matrix bi-module problem having BDCC condition in the case of MW5.

\bigskip
\bigskip
\noindent{\bf 5.1  An inspiring example}
\bigskip

Let $\mf A=(R,\K,\M,H=0)$ with $R$ trivial be a bipartite matrix
bi-module problem having BDCC condition of representation wild type. Suppose
$\mf A'=(R',\K',\M',$ $H')$ is an induced bi-module problem satisfying MW5.
We classify the position of the first arrow $a'_1$ in
the formal product $\Theta'$ of $(\mf A',\mf B')$:
$$\begin{array}{c}H'+\Theta'=H'+ \sum_{i=1}^{n'}
{a'}_i\ast{A'}_i.\end{array}\eqno{(5.1\mbox{-}1)}$$
Denote by $(p,q')$ the leading position of $A'_1$ over $\T'$, which
locates in the $(\textsf p,\textsf q)$-th leading block of some base
matrix of $\M_1$ partitioned by $\T$, with $\textsf q=\textsf q_Z$ a main block
column over $Z$.

\medskip

{\bf Classification 5.1.1}\, Let the pairs $(\mf A,\mf B),(\mf A',\mf B')$
be given above with $\mf B'$ satisfying MW5. Then there are two possible position
relations between $p$ and the row indices of the links of ${H'}$ in Formula (5.1-1):

{\bf case (I)} $p <$ the row indices of all the links in the
$(\textsf p,\textsf q)$-block of ${H'}$;

\smallskip

{\bf case (II)} $p \geqslant$ some row index of at least one link in the
$(\textsf p,\textsf q)$-block of ${H'}$.

\smallskip

It is clear that there is no link above the $(\textsf p,\textsf
q)$-th block, since $\mf A'$ is already local.

\medskip

{\bf Lemma 5.1.2} Let $p_x$ be the row index of $x$ in $H'$, then
$p_x>p$ in Classification 5.1.1.

\smallskip

{\bf Proof.}  Since $x$ appears before the first arrow $a_1'$ of $\mf B'$,
$p_x\geqslant p$ by the ordering of the reductions.
If $p_x=p$, then the parameter $x$ locates at the left side of $a_1'$
in $H'+\Theta'$, $\delta(a_1')$ contains only the terms of the form
$\alpha xv,\alpha\in k$, which contradicts
to the assumption that $\mf{B}'$ is in the case of MW5. Thus $p_x>p$, the proof is finished.

\medskip

{\bf Example 5.1.3}\, Let $(\mf A,\mf B)$ be a pair of matrix
bi-module problem associated to the algebra defined in Example 1.4.5.
There is a reduction sequence $\mf A=\mf A^0,\mf A^1,\mf A^2,\mf A^3$ given in
Examples 2.4.5, such that the corresponding bocs $\mf B^3$ of $\mf A^3$ is strongly
homogeneous in the case of MW5 described in 3.1.5 (iii).
In order to prove that $(\mf A,\mf B)$ is not homogeneous, we must find
another way different from the proof of MW1-MW4. More
precisely, we will reconstruct a new reduction sequence based on the matrix
$\tilde M$ over $k[x]$ with the size vector $\tilde{\m}=(2,2,2,2,2,3,3,3,3,3)$:
$$\tilde M
=\small\left(
\begin{array}{c|c|c|c|c}
\quad0\quad& \quad0\quad& \begin{array}{ccc}0& 1 &0\\
0&0&1\end{array} &
\begin{array}{ccc} 0&0& 1 \\ 0&0&0\end{array}&
0\\ \hline
& 0&\begin{array}{ccc}  0&0& 1 \\ 0&0&0\end{array}& 0&
\begin{array}{ccc} 0& 0&0 \\ 0&x&0\end{array}\\ \hline
&& 0&0& \begin{array}{ccc}0& 0&1 \\ 0&0&0\end{array}\\ \hline
&&&0&\begin{array}{ccc} 0&1 &0\\ 0&0&1\end{array}\\ \hline
&&&\begin{array}{ccc}&&\\&&\end{array}&0
\end{array}\right).
$$
Corresponding to steps (i)-(iii) of 2.4.5, there is a reduction sequence
$\mf{A}=\tilde{\mf{A}}^0,\tilde{\mf{A}}^1,\tilde{\mf{A}}^2,\tilde{\mf{A}}^{3}$,
where the reductions from $\tilde{\mf A}^0$ to $\tilde{\mf A}^1$ is given by $a\mapsto (0\,1)$
in the sense of Lemma 2.3.2. Thus $b$ splits into $b_1, b_2$ in $\tilde{\mf B}^1$, and set
$b_1\mapsto(0), b_2\mapsto {{0\,1}\choose{0\,0}}$ from $\tilde{\mf{A}}^1$ to $\tilde{\mf{A}}^2$.
$\tilde{\mf A}^3$ is obtained from  $\tilde{\mf A}^2$
by an edge reduction going to $(0)$, a loop mutation, then a series of regularization:
$$\begin{array}{c}
\tilde{H}^3=\left(\begin{array}{ccc}0&1_{X}& 0\\
0&0&1_X\end{array}\right)\ast A +\left(\begin{array}{ccc}0&0& 1_X\\
0&0&0\end{array}\right)\ast B+\left(\begin{array}{ccc}\emptyset&\emptyset& \emptyset\\
0&x1_X&\emptyset\end{array}\right)\ast C.\end{array}$$
The $(1,5)$-th block partitioned under $\T$ in the formal equation of $(\tilde{\mf
A}^3,\tilde{\mf B}^3)$ is of the form:
$$\small{~_{\left(\begin{array}{cc} e& v\\
0&e\end{array}\right)\left(\begin{array}{ccc} d_{10}&d_{11}& d_{12}\\
d_{20}&d_{21}&d_{22}\end{array}\right)+\left(\begin{array}{cc} u^1_{11}& u^1_{12}\\
u^1_{21}&u^1_{22}\end{array}\right)\left(\begin{array}{ccc}0& 0& 1\\
0&0&0\end{array}\right)+\left(\begin{array}{cc} u^2_{11}& u^2_{12}\\
xv&u^2_{22}\end{array}\right)\left(\begin{array}{ccc}0& 1& 0\\
0&0&1\end{array}\right)=}}$$
$$\small{~_{\left(\begin{array}{ccc}0& 1&0 \\
0&0&1\end{array}\right)\left(\begin{array}{ccc}
v^{_2}_{_{00}}& v^{_2}_{_{01}}& v^{_2}_{_{02}}\\
v^{_2}_{_{10}}&v^{_2}_{_{11}}&v^{_2}_{_{12}}\\ 0&vx&v^{_2}_{_{11}}\end{array}\right)+
\left(\begin{array}{ccc}0& 0& 1\\
0&0&0\end{array}\right)\left(\begin{array}{ccc}
v^{_1}_{_{00}}& v^{_1}_{_{01}}& v^{_1}_{_{02}}\\
v^{_1}_{_{10}}& v^{_1}_{_{11}}& v^{_1}_{_{12}}\\
v^{_1}_{_{20}}&v^{_1}_{_{21}}&v^{_1}_{_{22}}\end{array}\right)
+\left(\begin{array}{ccc}d_{_{10}}& d_{_{11}}& d_{_{12}}\\
d_{_{20}}&d_{_{21}}&d_{_{22}}\end{array}\right)\left(\begin{array}{ccc}s_{00}&s_{01}&s_{02} \\0&e& v\\
0&0&e\end{array}\right)} }$$ with $e=e_{_{X}}, s_{00}=e_{_{Y}}$. Then the
differentials of the solid arrows of $\tilde{\mf B}^3$ can be read off:
$$d_{20}: X\mapsto Y,\,\,\delta(d_{20})=0;\quad
d_{21}:X\mapsto X,\,\,\dz(d_{21})=xv-vx-d_{20}s_{01}.$$
And $\bar\dz(d_{22})=u_{21}^1-v^2_{11}-d_{20}s_{02}-d_{21}v,
\bar\dz(d_{10})=-v_{10}^2-v_{20}^1+vd_{20}, \bar\dz(d_{11})=u_{11}^2-v_{11}^2-v_{21}^1
-d_{10}s_{01}, \bar\dz(d_{12})=u_{11}^1+u_{11}^2-v_{12}^1-v_{22}^1$, where $\bar\dz$ is respect to $d_{20},d_{21}$.
It is clear that the bocs $\wt{\mf B}^3$ is in the case of Proposition 3.4.5, since
for $d_{20}\mapsto(1)$, the solid loops
$d_{21},d_{22},d_{10},d_{11},d_{12}$ will be regularized, because
$s_{01},u_{21}^1,v_{10}^2,u_{11}^2,u_{11}^1$ are pairwise
different dotted arrows.
Therefore $(\mf A,\mf B)$ is not homogeneous.

\medskip

Motivated by Example 5.1.3, we consider the general cases.
Since the example satisfies Case (I) of Classification 5.1.1,
we start from Case (I) in the subsection 5.1-5.3.

Let $\mf{A}=(R,\K,\M,H=0)$ be a bipartite matrix bi-module problem
having RDCC condition. Let $\mf{A}'$ be an induced matrix
bi-module problem with $R'$ trivial. Let
$\vartheta: R(\mf A')\rightarrow R(\mf A)$ be the induced
functor, $ M= \vartheta^{0i} (H'(k))=\sum_j{M}_j\ast A_j\in R(\mf A)$.
Suppose that the size vector of $M$ is $\l\times \n$ over
$\T$. Let $\textsf q=\textsf q_Z\in T_2$ for some $Z\in \T$.
Define a size vector $\l\times\tilde{\n}$ over $\T$, and construct a representation
$\tilde{M}=\sum_j
\tilde{M}_j\ast A_j\in R(\mf{A})$ with $0$ a zero column as follows:
$$\tilde{n}_j=\left\{\begin{array}{ll}n_j, &\mbox{if}\,
j\not\in Z,\\n_j+1,& \mbox{if}\, j\in Z.\end{array}\right.\quad
\tilde{M}_j=
\left\{\begin{array}{ll}M_j,& \mbox{if}\, A_j1_Z=0,\\(0\,M_j),&
\mbox{if} A_j1_Z=A_j.\end{array}\right.\eqno(5.1\mbox{-}2)$$

Denote by $(p,q+1)$ the leading
position of $M$, and suppose that $q+1$ is the index of the
first column of the $\textsf q$-th block-column, denote by
$\tilde q$ the index of the added column in the
$\textsf{q}_Z$-th block-column of $\tilde M$. Consider the defining
system $\IE$ of $\K'_0\oplus\K'_1$ given by Formula (2.4-4), and a matrix
equation $\tilde{\IE}$:
$${\IE}:\,\,
\Phi^1_{\l}M\equiv_{_{\prec(p, q+1)}}
M\Phi^2_{\n^{i}},\quad
\tilde{\IE}:\,\, \Phi^1_{\l}\tilde M\equiv_{_{\prec({p},
\tilde{q})}}\tilde{M}\Phi^2_{\wt{\n}}.
\eqno{(5.1\mbox{-}3)}$$
Using the right hand side of ${\IE}$ and
$\tilde{\IE}$, one defines two new matrix equations respectively:
$$ {\IE}_\tau:\,\,0\equiv_{_{\prec(p, q+1)}}
M\Phi^2_{\n},\quad \tilde{\IE}_\tau:\,\,
0\equiv_{_{\prec({p},
\tilde{q})}}\tilde{M}\Phi^2_{\wt{\n}}.\eqno{(5.1\mbox{-}4)}$$
Since $\mf A$ has BDCC condition, the main block-column
in $\Phi^2_{\n^{i}}$ determined by $Z\in \T$ can be written as $\Phi^2_{\n^i,Z}=
(\Phi_{1}^2,\cdots,\Phi_{n^i}^2)^T$, such that either $\Phi_{l}^2=0$ or
$\Phi^2_l=(z^l_{pq})\ne 0$, where $z^l_{pq}$ are algebraically
independent variables over $k$, Denote the pairwise different
non-zero blocks by $\Phi^2_{l_1},\cdots,\Phi^2_{l_u}$.

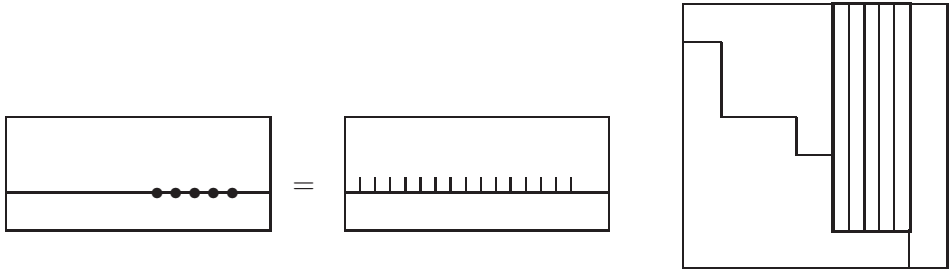
\begin{center}
\setlength{\unitlength}{1mm}
\begin{picture}(145,38)\thinlines{\linethickness{0.1mm}}
\put(30,15){\circle*{1.5}}\put(32.5,15){\circle*{1.5}}
\put(35,15){\circle*{1.5}}
\put(37.5,15){\circle*{1.5}}\put(40,15){\circle*{1.5}}
\put(10,10){\framebox(35,15)}\put(10,15){\line(1,0){35}}
\put(48, 15){$=$}

\put(55,10){\framebox(35,15)}\put(55,15){\line(1,0){35}}
\mput(55,15)(2,0){16}{\line(0,1){2}}

\put(100,5){\framebox(35,35)}
\put(100,35){\line(1,0){5}}
\put(105,25){\line(1,0){10}}\put(105,25){\line(0,1){10}}
\put(115,20){\line(1,0){5}}\put(115,20){\line(0,1){5}}
\mput(122,10)(2,0){4}{\line(0,1){30}}
\put(130,5){\line(0,1){5}}

\thicklines{\linethickness{1mm}}
\put(120,10){\framebox(10,30)}
\end{picture}
\end{center}
Taken any integer $p'\geqslant p$ and $1\leqslant j\leqslant
n_{_Z}$, the $(p',q+j)$-th entry of the right side of $\IE_\tau$ is:
$$\begin{array}{c}\sum_{\Phi_l\ne
0}\sum_{q'}\alpha^l_{p'q'}z_{q',q+j}^l.
\end{array}\eqno{(5.1\mbox{-}5)}$$ It is easy to see that
$z_{q',q+j_1}^l$ and $z_{q',q+j_2}^l$ have the same coefficient
$\alpha_{p',q'}^l$, the $(p',q')$-th entry of $H(k)$, in Formula (5.1-5).
In the picture above, $n_{_Z}=5, p'=p$, the five equations locating
at five circles have the same coefficients of each variables.

\medskip

{\bf Remark 5.1.4}\, (i)\, The $(p,q+j_1)$-the equation is a linear combination of
the previous equations in $\IE_\tau$, if and only if so is the
$(p,q+j_2)$-th equation by Formula (5.1-5). Similarly, we have
the same result in $\tilde\IE_\tau$.

(ii)\, The equations in the system $\tilde{\IE}$ (resp. $\tilde{\IE}_{\tau}$) and
those in ${\IE}$ (resp. ${\IE}_{\tau}$)  are the same everywhere except at the
$\tilde q$-column, since the entries of $\tilde M$ in the $\tilde q'$-th, the first column
of the $\textsf q'$-th block column with
$\textsf q'\sim\textsf q$ under the partition $\T$, are all zero.

(iii)\, If there exists some $1\leqslant j\leqslant n_{_Z}$, such
that the $(p,q+j)$-th equation is a linear
combination of the previous equations in $\IE_\tau$, then so is
the $(p,\tilde q)$-th equation in $\tilde\IE$, by (i)-(ii).

(iv)\, If the $(p,q+j)$-th equation is a linear combination of the previous
equations in $\IE$, then so is it in $\IE_\tau$, since the variables in $\Phi_{\l}$ and $\Phi_{\n}$
are algebraically independent.

\bigskip
\bigskip
\noindent{\bf 5.2 Bordered trivial matrices in the bipartite case}
\bigskip

This sub-section is devoted to constructing a reduction sequence
based on a given sequence and a bordered matrix, which generalizes Example 5.1.3.

\medskip

Let $\mf{A}=(R,\K,\M,H=0)$ be a bipartite matrix bi-module problem
having RDCC condition. Let
$\mf{A}^s=({R}^s,\K^s,\M^s,H^s)$ be an induced problem with $R^s$
trivial and local. Then Corollary 2.3.5 gives a unique reduction
sequence with each reduction being in the sense of Lemma 2.3.2:
$$
\mf{A}=\mf A^0, \mf{A}^1, \cdots,  \mf{A}^i,\mf{A}^{i+1}, \cdots,
\mf{A}^s, \eqno {({\ast})}
$$ Let $\vartheta^{0s}: R({\mf A}^s)\rightarrow
R(\mf A)$ be the induced functor, $M=\vartheta^{0s} (H^s(k))=\sum_{j} {M}_j\ast A_i\in
R(\mf A)$ with the size vector $\l\times \n$, and $M_j=G^j_s(k)$ given by Formula (2.3-5).
Define a size vector $\l\times\tilde{\n}$ and a representation $\tilde{M}=\sum_{j}
\tilde{M}_j\ast A_j\in R(\mf{A})$ given by Formula (5.1-2).

\medskip

{\bf Theorem 5.2.1}\, There exists a unique reduction sequence
based on the sequence $(\ast)$:
$$\begin{array}{c}\mf{A}=\tilde{\mf{A}}^0,
\tilde{\mf{A}}^1, \cdots, \tilde{\mf{A}}^{i},
\tilde{\mf{A}}^{i+1}, \cdots, \tilde{\mf{A}}^s\end{array} \eqno
{(\tilde\ast)}
$$
where $\tilde{\mf{A}}^i=(\tilde{R}^i, \tilde{\K}^i,
\tilde{\M}^i,\tilde{H}^i)$, the reduction from $\tilde{\mf{A}}^{i}$
to $\tilde{\mf{A}}^{i+1}$ is a reduction or a composition of two
reductions in the sense of Lemma 2.3.2. Moreover, $\tilde\T^s$ has two
vertices, and $$\tilde{\vartheta}^{0s}(\tilde{H}^s(k))=\tilde{M}.$$

{\bf Proof}\,  We may assume that $\l\times \n$ is
sincere over $\T$. Otherwise, after a suitable deletion, we are able to
obtain an induced problem $\mf{A}'$, which is still
bipartite having BDCC condition,
such that $M$ is a representation of $\mf{A}'$ of size vector
$\l'\times \n'$, which is sincere over $\T'$.

We will construct a sequence $(\tilde\ast)$ inductively. The original
term in the sequence is $\tilde{\mf{A}}^0=\mf{A}^0$. Suppose that
we have constructed a sequence $\tilde{\mf{A}}^0,
\tilde{\mf{A}}^1,\cdots, \tilde{\mf{A}}^{i}$ for some $0\leqslant i<s$, and $\tilde
\vartheta^{0i}: R(\tilde{\mf{A}}^{i})\mapsto R(\tilde{\mf A}^0)$ is the induced functor, such
that there exists a representation
$$\begin{array}{c}\tilde M^{i}=\tilde H^{i}_{\tilde
\m^{i}}(k)+\sum_{j=1}^{n^i} \tilde M^{i}_j\ast \tilde A^{i}_j\in R(\tilde{\mf{A}}^{i}),\quad
\mbox{with}\quad\vartheta^{0i}(\tilde{M}^{i})\simeq \tilde M\in
R(\tilde{\mf A}^0).\end{array}\eqno{(5.2\mbox{-}1)}$$

Write $M^i=\vartheta^{is}({H^s}(k))=H^{i}_{\m^{i}}(k)+\sum_{j=1}^{n^{i}}M_j^{i}\ast
A_j^{i}\in R(\mf{A}^{i})$, where $M_1^{i}=G^{i+1}_s(k)$
by Corollary 2.4.5, denoted $G^{i+1}_s(k)$ by $B$ for simplicity. Denote the
first column of the $\textsf q$-th block-column in the formal
product $\Theta^i$, which $B$ belongs to, by $\beta$. Now we are constructing $\tilde{\mf
A}^{i+1}$.

\smallskip

{\bf Case 1}\, $\tilde{\T}^{i}=\T^{i}$ and $\tilde{\mf A}^i=\mf A^i$.

\smallskip

{\bf 1.1} $B\cap \beta$ is empty. Then
$\tilde{G}^{i+1}=G^{i+1}$, $\tilde{H}^{i+1}=H^{i+1}$ and
$\tilde{\mf{A}}^{i+1}=\tilde{\mf A}^{i+1}$.

\smallskip

Before giving the following cases, we claim that if $B\cap \beta$ is non-empty, $B$ thus
$G^{i+1}$ can not be Weyr matrices. Otherwise, the first arrow $a_1^{i}$ of $\mf
B^{i}$ will be a loop. Since $\tilde{\T}^{i}=\T^{i}$, $\tilde{a}_1^{i}$
will also be a loop and hence the numbers of rows and columns of
$\tilde{B}$ are the same. When the matrix $B$ is enlarged by one
column, then also enlarged by one row, a contradiction to the construction
of $\tilde M$. So $B$ is either a regularization
block $\emptyset$ or an edge reduction block ${{0\,\, I_r}\choose
{0\,\,\,0}}$.

\smallskip

{\bf 1.2} $B\cap \beta$ is non-empty, and $B=\emptyset$ is a
zero block. Then $\tilde{B}=(\emptyset \, B)$ with $\emptyset$ a zero
column, $\tilde{H}^{i+1}=H^{i+1}$ and $\tilde{\mf{A}}^{i+1}=\mf{A}^{i+1}$
by a regularization.

\smallskip

{\bf 1.3} $B\cap \beta$ is non-empty, $B={{0\,\, I_r}\choose
{0\,\,\,0}}$ and $r<$ the number of columns of $B$. Then
$\tilde{B}=(0\, B)$ with $0$ a zero column, $\tilde{H}^{i+1}=H^{i+1}$ and
$\tilde{\mf{A}}^{i+1}=\mf{A}^{i+1}$ by an edge reduction.

\smallskip

{\bf 1.4} $B\cap \beta$ is non-empty, $B={{I_r}\choose 0}$. Then
$\tilde{B}=(0\, B)$ with $0$ a
zero column. Recall the formula (2.3-5) and Corollary 2.3.5, we define
$$
\tilde{G}^{i+1}=\left\{\begin{array}{ll}  \left(\begin{array}{cc} 0&
1_{Z_2}\end{array}\right),&
\mbox{if } G^{i+1}=(1_{Z_2});\\
&\\
\left(\begin{array}{cc} 0& 1_{Z_2}\\
0&0\end{array}\right),& \mbox{if } G^{i+1}=
\left(\begin{array}{c}  1_{Z_2}\\
0\end{array}\right),\end{array}\right. \quad \tilde{H}^{i+1}=\sum_{X\in
\T^{i}}I_{X}\ast H^{i}_{X}+\tilde{G}^{i+1}\ast A_1^{i}.
$$
Then $\tilde{\mf{A}}^{i+1}$ is induced from $\tilde{\mf A}^i$ by an edge
reduction in the sense of Lemma 2.3.2.

We stress, that after the edge reduction 1.4, $\tilde{\T}^{i+1}=\T^{i+1}\cup
\{Y\}$ with $Y$ an equivalent class consisting of the indices of the added columns
in the formal product $\tilde\Theta^{i+1}$ of the pair $(\tilde{\mf A}^{i+1},\tilde{\mf B}^{i+1})$,
and $(\tilde{\mf{A}}^{i+1},\tilde{\mf B}^{i+1})\ne (\mf A^{i+1},\mf B^{i+1})$ from this stage.

\smallskip

We show the above $B$ as a small block in the corr. leading block partitioned by $\T$:

\bcen \unitlength 1mm
\begin{picture}(110,20)

\put(0,0){\framebox(30,20)[l]{\dashbox{1}(3,20)}}
\put(10,10){\framebox(10,8)}

\put(40,0){\framebox(30,20)[l]{\dashbox{1}(3,20)}}
\put(40,10){\framebox(10,8)}

\put(80,0){\framebox(30,20)} \put(83, 0){\line(0,1){20}}
\put(80,10){\framebox(10,8)}


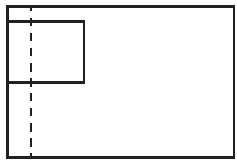
\put(9, -5){Case 1.1} \put(40, -5){Cases 1.2 and 1.3} \put(87,
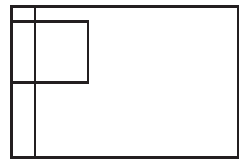
-5){Case 1.4}
\end{picture}

\end{center}
\vskip 3mm

\smallskip

{\bf Case 2.} $\tilde{\T}^{i}=\T^{i}\cup \{Y\}$.

\smallskip

{\bf 2.1} $B\cap \beta$ is empty. Then $\tilde{B}=B$,
$\tilde{G}^{i+1}=G^{i+1}$, and $\tilde{H}^{i+1}=\sum_{\tilde{X}\in
\tilde\T^{i}}\tilde{I}_{\tilde{X}}\ast\tilde H^{i}_{\tilde
{X}}+\tilde{G}^{i+1}\ast\tilde{A}_1$.

\smallskip

If $B\cap \beta$ is non-empty. Denote by $\tilde a^{i}_0$ and $\tilde
a^{i}_1$ the first and the second solid arrows of $\tilde {\mf B}^{i}$, which
locate at $(p^{i}, \tilde q_0^{i})$ and $(p^{i},\tilde q^{i}_0+1)$
in the formal product $\tilde\Theta^i$ respectively.

\smallskip

{\bf 2.2} $B\cap \beta$ is non-empty, and there exists some
$1\leqslant j\leqslant n^i_Z$, such that the
$(p^i,q^i+j)$-th equation is a linear combination of previous equations in
$\IE_\tau^i$. Then $\delta(\tilde{a}_0^{i})=0$ by Remark 5.1.4 (iii) and Corollary 2.4.4.
We make two reductions: the first is an
edge reduction given by $\tilde{a}_0^{i}\mapsto (0)$; the second
one for $\tilde{a}_1^i$ is as the same as that for $a_1^{i}$ by 5.1.4 (ii),
then we obtain an induced
problem $\tilde{\mf{A}}^{i+1}$ with $\tilde{B}=(0 \, B)$, where
$0$ is a zero column.

\smallskip

{\bf 2.3} $B\cap \beta$ is non-empty, and for all $1\leqslant
j\leqslant m^i_Z$, the $(p^i,q^i+j)$-th equation
is not a linear combination of previous equations in
$\IE^i_\tau$. Thus $\dz(\tilde a^i_0)\ne 0$ still by 5.1.4 (iii) and 2.4.4.
And $(p^i,q^i+j)$-th equation neither is in $\IE^i$ by 5.1.4 (iv),
$\dz(a_1^i)\ne 0$ again by 2.4.4.
Then we make two regularization
$\tilde{a}_0^{i}\mapsto \emptyset,\; \tilde{a}_1^{i}\mapsto \emptyset$,
and $\tilde B=(\emptyset\, B)$ with $\emptyset$ being a zero
column.

\smallskip

In the cases 2.2-2.3, we define $\tilde{G}^{i+1,0}=(0)$ or $\emptyset$,
$\tilde{G}^{i+1,1}=G^{i+1}$, and set
$\tilde{H}^{i+1}=\sum_{\tilde X\in\tilde{\T}^{i}}
\tilde{I}_{\tilde X}\ast \tilde {H}^{i}_{{\tilde
X}}+\tilde{G}^{i+1,0}\ast\tilde{A}_0^{i}
+\tilde{G}^{i+1,1}\ast\tilde{A}_1^{i}.$

By summary up all the cases, we obtain an induced pair
$(\tilde{\mf{A}}^{i+1},\tilde{\mf{B}}^{i+1})$ and a representation
$\tilde M^{i+1}$ satisfying Formula (5.2-1).
The theorem follows by induction.

\medskip

{\bf Corollary 5.2.2} With the notations as in Theorem 5.2.1. The
main diagonal block $\tilde{e}^i_Z,Z\in\T,$ of $\tilde{\K}^i_0\oplus\tilde{\K}^i_1$
has the form with $m=n^i_Z$:
$$
\left(\begin{array}{ccccc} s_{00}& s_{01}& s_{02}&\cdots& s_{0m}\\
&s_{11} & s_{12}& \cdots & s_{1m}\\
&&s_{22}& \cdots &s_{2m}\\
&&&\ddots &\vdots\\
&&&&s_{mm}
\end{array}\right).
$$
Where $s_{01}, s_{02},\ldots,s_{0m}$ are dotted arrows of $\tilde{\mf B}^i$.

\smallskip

{\bf Proof}\, By the construction of $\tilde{H}^i$, the added
``$0$-column'' can be only $0$ or $\emptyset$. Therefore except
$s_{00}$, the elements at the $0$-th row: $s_{01}, s_{02},
\ldots,s_{0m}$ do not appear in the defining equation system of
$\tilde{\mf{A}}^i$, and thus they are free. The proof is finished.

\bigskip
\bigskip
\noindent{\bf 5.3  Bordered matrices with a parameter $x$ in the bipartite case}
\bigskip

This sub-section is devoted to proving that the bipartite pair
with an induced minimal wild pair of MW5 and satisfying
Classification 5.1.1 (I) is not homogeneous.

Suppose we have a reduction sequence:
$$\mf{A}=\mf A^0, \mf{A}^1,\cdots,\mf{A}^s, \mf{A}^{s+1},  \cdots,\mf{A}^\varepsilon,
\cdots,\mf{A}^\epsilon,\cdots,\mf A^t=\mf A' \eqno{(\ast')}$$
where the reduction from $\mf A^{i}$ to $\mf A^{i+1}$ is given by
Lemma 2.3.2 for $1\leqslant i=1<s$, and $\mf A^{s}$ is local; the
reduction from $\mf A^{s}$ to $\mf A^{s+1}$ is a loop mutation, then
we obtain a parameter $x$ at the $(p_x, q_x)$-position and $R^{s+1}=k[x]$; the reduction
from $\mf A^{i}$ to $\mf A^{i+1}$ is a regularization for
$s<i<t$. The pair $(\mf A^t,\mf B^t)$ is
in the case of MW5 and Classification 5.1.1 (I).

Note that the set $T^i$ of integers and its partition $\T^i$ of $\mf A^i$
are all the same for $i=s,\cdots,t$, we may write $\n^i$ uniformly by $\n$;
$t^i$ by $\bar t$. Suppose the first arrow $a^t_1$ of $\mf B^t$ locates at $(p,q')$-position
in the formal product $\Theta^t$ with $q'=q+j$ for
some $1\leqslant j\leqslant n^t_Z$; the first arrow
$a^\epsilon_1$ of ${\mf B}^\epsilon$ locates at $(p,q+1)$ in $\Theta^\epsilon$,
where $q+1$ is the index of the
first column in the $\textsf q$-th block-column.
The picture below shows the position of the first solid arrows in
the formal products $\Theta^i$ of $({\mf A}^i,{\mf B}^i)$ for $i=s,
\epsilon, t$ (when we ignore the added $\tilde q_0$-th column):

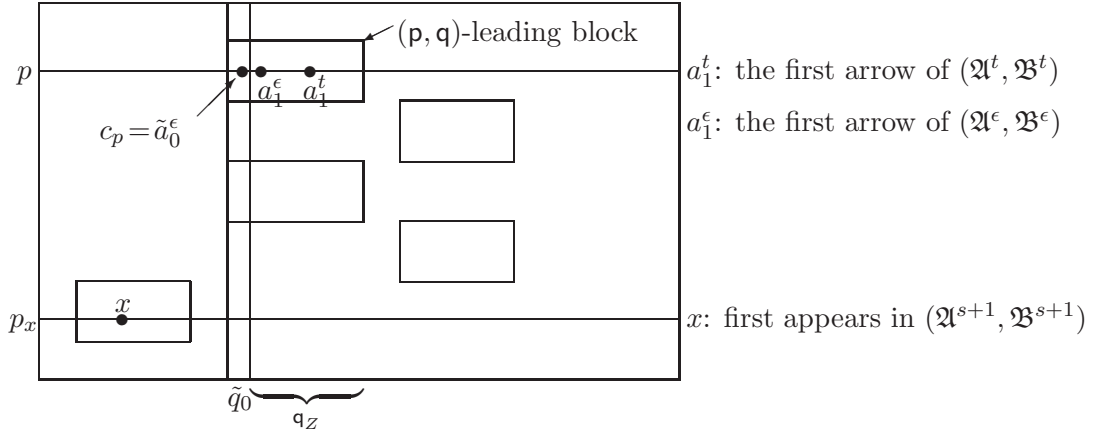
\begin{center}
\unitlength=1mm
\begin{picture}(120,50)
\put(0,0){\framebox(85,50){}} \put(5,5){\framebox(15,8){}}
\put(25,21){\framebox(18,8){}} \put(25,37){\framebox(18,8){}}
\put(48,13){\framebox(15,8){}} \put(48,29){\framebox(15,8){}}
\put(25,0){\line(0,1){50}} \put(28,0){\line(0,1){50}}
\put(0,8){\line(1,0){85}} \put(0,41){\line(1,0){85}}
\put(-4,7){$p_x$} \put(-3,40){$p$} \put(10,7){$\bullet$}
\put(26,40){$\bullet$}\put(35,40){$\bullet$}
\put(10,9){$x$}\put(28.5,40){$\bullet$}
\put(8,32){$c_{p}\!=\!\tilde{a}_0^\epsilon$}
\put(35,37.3){$a^t_1$}\put(29,37.3){${a}_1^\epsilon$}
\put(86,7){$x$: first appears in $(\mf{A}^{s+1}, \mf{B}^{s+1})$}
\put(86,40){$a^t_1$: the first arrow of $(\mf{A}^t,
\mf{B}^t)$} \put(86,33){${a}_1^\epsilon$: the first arrow of
$(\mf{A}^\epsilon, \mf{B}^\epsilon)$} \put(21,35){\vector(1,1){5}}
\put(25,-3.5){$\tilde{q}_0$} \put(47,47){\vector(-2,-1){4}}
\put(47.5,45){$(\textsf p,\textsf q)$-leading block}
\put(28,-1){\makebox{$\underbrace{\hskip 15mm}_{\textsf q_Z}$}}
\end{picture}\vskip 8mm
\end{center}
For any $s<i\leqslant t$, assume that $R^i=k[x,\phi^i(x)^{-1}]$
in ${\mf A}^i$, $H^i\in\IM_{t^i}(k[x,\phi^i(x)^{-1}])$.
Then $H^i$ are all in the same
form but with different $\phi^i(x)$. Since $k(x)$ is a $k[x,\phi^i(x)^{-1}]$-bi-module,
$$H^i\otimes_{R^i}1_{k(x)}\in \IM_{t^i}(k[x,\phi^i(x)^{-1}])
\otimes_{R^i}k(x)\simeq\IM_{t^i}(k(x)).$$
Remark 5.1.4 (i) is still valid when we consider the equations over $k(x)$
instead of $k$ in $\IE_\tau$ and $\tilde{\IE}_\tau$. Define a variable matrix
$\Phi_{\tilde\n}^{(\bullet)}$ of size $\tilde\n$ with the $\tilde q$-th column
as the same as that of $\Phi_{\tilde\n}^2$ and others zero.
In the systems below, see Formula (5.1-4), we only need to consider the second one:
$$\begin{array}{c}
\tilde{\IE}_\tau: 0\equiv
\tilde M\Phi_{\tilde{\n}}^2; \quad \tilde{\IE}_{\tau}^{(\bullet)}: 0
\equiv \tilde M\Phi_{\tilde\n}^{(\bullet)},\end{array} \eqno{(5.3\mbox{-}1)}$$

Suppose $x$ locates at the $p_x$-th row of the $(\textsf p_x,\textsf q_{Z'})$-th
main block partitioned under $\T$, see the picture above for the
case of $Z'\ne Z$; and the example 5.1.3 for $Z'=Z$.
Thus the equation system $\tilde\IE_\tau^{(>p_x)}$ consisting of the
equations below the $p_x$-th row is over the base field $k$.
Denote by $\K_\tau^{(> p_x)}\subset$ IM$_{\bar t\times\bar t}(k)$ the solution space
of $\tilde\IE^{(>p_x)}_\tau$. Since $\tilde\K_1^s\oplus\tilde\K_1^s$ is local and upper triangular,
one base matrix $E$ of $\K_\tau^{(> p_x-1)}$ with entries all zero, except $1_Y$ at the
position $(\tilde q,\tilde q)$, $\{E\}$ forms the basis of $\K_{\tau,0}^{(> p_x-1)}$; and other
base matrix $V_j,j\geqslant 1$, with non-zero entries above the $\tilde q$-row at the $\tilde q$-th column,
form a basis of $\K_{\tau,1}^{(>p_x-1)}$.

Suppose we have solved the
equations $\tilde\IE^{(>h)}_\tau$ for some $p< h\leqslant p_x$, and obtained
$R^{(>h)}_\tau=k[x,\prod_{\eta=p_x}^{h-1}d^{\eta}(x)^{-1}]\times k1_Y$;
and a quasi-free module
$\K_{\tau,1}^{(>h)}$ with a quasi-basis $\{U_1,\cdots,U_\kappa\}$ over
$R^{(>h)}_\tau\otimes_k R^{(>h)}_\tau$. Consider the formal product and the equation system,
see Theorem 2.4.1:
$$\begin{array}{c}\Pi^{(>h)}_\tau=\sum_{\zeta=1}^{\kappa}
u_\zeta\ast U_\zeta,\quad\quad
\tilde\IE^{(>h)}_\tau:\,0\equiv H(k(x))\Pi^{(> h)}_\tau.\end{array}\eqno{(5.3\mbox{-}2)}$$
Where the $h$-th equation of $\tilde\IE^{(>h)}_\tau$ is $\sum_{\zeta=1}^{\kappa} f_\zeta(x)u_\zeta,
f_\zeta(x)\in\,R^{(>h)}_\tau$. There are two possibilities:

(i)\, $f_\zeta(x)=0$ for $\zeta=1,\cdots,\kappa$,
then $\K^{(> h-1)}_{\tau,1}=\K^{(>h)}_{\tau,1}$,
set $d^{h}(x)=1$ and the quasi-basis of $\K^{(>h)}_{\tau,1}$
are preserved in $\K^{(> h-1)}_{\tau,1}$.

(ii)\, There exists some $f_{\zeta}(x)\ne 0$, without loss of generality
we may assume that $f_\kappa(x)\ne 0$. Choose a new basis of Hom$_{k(x)}(\K^{(>h)}_{\tau,1}\otimes_{R^{(>h)}}k(x),k(x))$
at the first line of Formula (5.3-3) below, we have the corresponding base change
of $\K^{(>h)}_{\tau,1}$ over $R^{(>h)}_\tau$ at the second line:
$$\left\{\begin{array}{l}u_\zeta'=u_\zeta,\\
U'_\zeta=U_\zeta-f_\zeta(x)/f_\kappa(x)U_\zeta;
\end{array}\right.\,\,\mbox{for}\,\,1\leqslant\zeta<\kappa;\quad\left\{\begin{array}{l}
u'_\kappa=\sum_{\zeta=1}^\kappa f_\zeta(x)u_\zeta,\\U'_\kappa=1/f_\kappa(x)U_\kappa.\end{array}\right.
\eqno{(5.3\mbox{-}3)}$$
Where $u_\kappa=0$ is the solution of the $(h,\tilde q)$-th equation in the system (5.3-2),
thus $\K^{(>h-1)}_{\tau,1}$ possesses the quasi-free-basis
$\{U_\zeta\mid \zeta=1,\cdots,\kappa-1\}$. Let $d^{h}(x)\in k[x]$ be the numerator
of $f_{\kappa}(x)$, and $R^{(>h-1)}_\tau=k[x,\prod_{\eta=1}^hd^\eta(x)]\times k1_Y$.

By induction, we finally reach the equation system $\tilde\IE^{(> p-1)}$
with the solution space $\K_{\tau,1}^{(>p-1)}$ and a polynomial $d^{p}(x)$.
Consider the pair $(\mf A^i,\mf B^i)$ for some $s+1\leqslant i\leqslant \epsilon$,
suppose $R^i=k[x,\phi^i(x)^{-1}]$, and the row index of the first arrow is
$p^i, p_x\leqslant p^i\leqslant p$, in the formal product $\Theta^i$. Define
$$\begin{array}{c}\tilde{\phi}^{i}(x)=\phi^i(x)\,\Pi_{\eta=p_x}^{p^i} d^\eta(x)\in k[x],\quad
\mbox{in particular}\quad \tilde{\phi}^{t}(x)=\phi^t(x)\,\Pi_{\eta=p_x}^{p} d^\eta(x).
\end{array}\eqno{(5.3\mbox{-}4)}$$
Denote $M^i=H^i(k[x,\tilde{\phi}^i(x)^{-1}])=\sum_j{M}_j^i\ast A_j$
having the size vector $\l\times \n$ over $\T$. Then
we are able to construct a size vector $\l\times\tilde{\n}$ over $\T$;
and a matrix $\tilde M^{i}=\sum_{j}\tilde{M}_j^i\ast A_j$ by Formula (5.1-2).
Write the matrix equations for $i\geqslant s$:
$$\begin{array}{ll}\IE^i:
\Phi_{\l} M^i\equiv_{\prec (p^i,q^i)}
M^i\Phi_{{\n}},&\tilde\IE^i:
\Phi_{\l}\tilde M^i\equiv_{\prec (p^i,q^i)}
\tilde M^i\Phi_{{\tilde\n}};\\ \IE^i_\tau:0\equiv_{\prec(p^i,q^i)}
M^i\Phi_{{\n}},&\tilde\IE^i_\tau:0\equiv_{\prec(p^i,q^i)}
\tilde M^i\Phi_{{\tilde\n}}.\end{array}\eqno{(5.3\mbox{-}5)}$$

{\bf Remark 5.3.1}\, Remark 5.1.4 (i)-(iv) are still valid if we consider the matrices
$M^i$ and $\tilde M^i$ over $k[x,\tilde\phi^{i}(x)^{-1}]$ instead of over $k$.

\medskip

{\bf Theorem 5.3.2}\, There exists a unique reduction sequence
based on the sequence $(\ast')$:
$$\begin{array}{c}\mf{A}=\tilde{\mf{A}}^0,
\tilde{\mf{A}}^1, \cdots, \tilde{\mf{A}}^{s},
\tilde{\mf{A}}^{s+1},\cdots, \tilde{\mf{A}}^\varepsilon,\cdots,
\tilde{\mf{A}}^\epsilon,\cdots,
\tilde{\mf A}^t=\tilde{\mf A}' \end{array}\eqno{(\tilde\ast')}
$$
where the first part of the sequence till to $\tilde{\mf A}^s$ is given
by Theorem 5.2.1; the reduction from $\tilde{\mf A}^s$ to $\tilde{\mf
A}^{s+1}$ is given by a loop mutation $a^{s+1}_1\mapsto(x)$, or
first an edge reduction $(0)$, then a loop mutation $(x)$; the
reduction from $\tilde{\mf A}^i$ to $\tilde{\mf A}^{i+1}$ for $s+1\leqslant
i<t$ is given by a regularization, or two regularization,
or an edge reduction $(0)$ then a regularization.

\smallskip

{\bf Proof}\, We make reduction from $\tilde{\mf{A}}^{s}$ to
$\tilde{\mf{A}}^{s+1}$ by a loop mutation when $x$ does not locate at
the $\textsf q_Z$-the block column of the formal product $\Theta^{s+1}$, or by an edge reduction $(0)$ then a loop
mutation when $x$ locates at the $\textsf q_{_Z}$-column by Corollary 2.4.4.

Now suppose we have an induced bi-module problem $\tilde{\mf A}^{i}$ for some $i>s$.
If the first arrow $a^{i}_1$ of ${\mf A}^{i}$ does not locate
at the $(q+1)$-th column of $\Theta^{i}$, make a
regularization $\tilde a^{i}_1\mapsto \emptyset$.
Otherwise, there are two possibilities:
\ding{172} there exists some $1\leqslant j\leqslant n_Z$, the $(p^i,q+j)$-th equation
is a linear combination of the previous equations in $\IE^i_\tau$,
then $\dz(\tilde a_0^i)=0$ by Remark 5.3.1 and Corollary 2.4.4, set $\tilde a^i_0\mapsto(0), \tilde
a^i_1\mapsto \emptyset$; \ding{173} otherwise
$\dz(\tilde a_0^i)\ne0$, set $\tilde a^i_0\mapsto\emptyset$ and $\tilde a^i_1\mapsto
\emptyset$. The sequence $(\tilde\ast')$ is completed by induction as desired.

\medskip

{\bf Corollary 5.3.3}\, $\mf B^t$ at $(\ast')$ satisfying MW5
implies $\dz(\tilde a^\epsilon_{0})=0$ in $\tilde{\mf B}^\epsilon$ at $(\tilde
\ast')$.

\smallskip

{\bf Proof}\, The first arrow $a_1^t$ of $\mf B^t$ locates at
the $(p,q+j)$-th position, therefore $\dz^0(a^\epsilon_{j})=0$ in $\Theta^\epsilon$.
By Remark 5.3.1, the $(p,l)$-th equation is a linear combination of previous equations
in $\tilde\IE_\tau^\epsilon$ for $0\leqslant l\leqslant n_{_Z}$.
Thus $\dz(\tilde a_0^\epsilon)=0$ in
$\tilde{\mf B}^\epsilon$ by Corollary 2.4.4, the proof is finished.

\medskip

{\bf Proposition 5.3.4}\, Let $\mf{A}=(R,\K,\M,H=0)$ with $R$ trivial be a bipartite
matrix bi-module problem having RDCC condition. If there exists an induced
pair $(\mf A',\mf B')$ of $(\mf A,\mf B)$ in the case of MW5, and
$H'+\Theta'$ satisfies Classification 5.1.1 (I), then $\mf{A}$ is not
homogeneous.

\smallskip

{\bf Proof}\, Suppose we have a sequence $(\ast')$ with $\mf B'=\mf
B^t$, then there is a sequence
$(\tilde\ast')$ based on $(\ast')$ by Theorem 5.3.2. Corollary 5.3.3 tells that the
first arrow $\tilde{a}_0^\epsilon$ of $\tilde{\mf B}^\epsilon$ is an edge
with $\dz(\tilde a_0^\epsilon)=0$, and hence we may set
$\tilde{a}_0^\epsilon\mapsto (1)$ according to Proposition 2.2.7.
Denote by $(\hat{\mf A},\hat{\mf B})$ the induced pair, which is
obviously local. Thus we are able to use the triangular formulae
given in the section 3.3, and obtain an induced pair in one
case of Classification 3.3.3 (ii)-(iv).

{\bf Case 1.} If we meet 3.3.3 (ii), then $\tilde{\mf
B}^\epsilon$ is in the case of Proposition 3.4.5. We are done.

{\bf Case 2.} If we meet 3.3.3 (iii), then there exists an
induced bocs of $\hat{\mf B}$ satisfying MW3,
we are done by Proposition 3.4.3.

{\bf Case 3.}  If we meet 3.3.3 (iv), then there exists
an induced bocs satisfying MW4, we are done by Proposition 3.4.4.

{\bf Case 4.}  If we meet of 3.3.3 (iv), and any induced
minimal wild bocs satisfying MW5, then we choose one of them, say
$(\hat{\mf A}^1, \hat{\mf B}^1)$, where the first arrow $\hat a_1$ in the sense of MW5
locates at the $p^1$-th row in the formal product $\hat\Theta^1$. We claim that
$p^1<p$. In fact, the arrows
$\tilde{a}_j$ for $j=1,\cdots,n_{_Z}$ at the $p$-th row in $\hat\Theta$ of the pair
$(\hat{\mf A},\hat{\mf B})$ has
the differentials $\delta^0 (\tilde a_j)=s_{0j}+\cdots$, and hence will
be regularized according to Corollary 5.2.2.

Repeating the above process for $(\hat{\mf A}^1, \hat{\mf B}^1)$, if
we meet one of the cases 1-3, the procedure stops. Otherwise if we
meet the case 4 repeatedly, there is a sequence of local pairs
and a decreasing sequence of the row indices:
$$\begin{array}{ccccc}(\hat{\mf A}, \hat{\mf B}),&
(\hat{\mf A}^1, \hat{\mf B}^1),& (\hat{\mf A}^2, \hat{\mf B}^2),&
\cdots, &(\hat{\mf A}^\beta, \hat{\mf B}^\beta),\\
\qquad p\quad>&\qquad p^1\quad>&\qquad p^2\quad>&\cdots>&p^\beta.\end{array}$$
Since the number of the rows of $\hat{H}^i$ for $i=1,\cdots,v$ is
fixed, the procedure must stop at some stage $\beta$, such that one of the
Cases 1-3 appears. The proposition is proved by induction.

\bigskip
\bigskip
\noindent{\bf 5.4  Bordered matrices in one-sided case}
\bigskip

This subsection is devoted to constructing a reduction sequence
starting from a one sided pair based on some bordered matrices.

\medskip

Let $(\mf{A},\mf B)$ be a bipartite pair having RDCC condition, let
$(\mf{A}^r,\mf B^r)$ be an induced matrix
bi-module problem with $R^r$ trivial given by
Formula (4.1-1). Which gives a quotient-sub pair
$((\mf A^r)^{[m]},(\mf B^r)^{(m)})$ denoted by $(\bar{\mf A},\bar{\mf B})$.
$\bar{\mf B}$ has a layer $L=(R; \omega; d_1,\cdots, d_m;
v_1,\cdots,v_t)$ by Definition 4.1.2.
Denote by $\bar T_R=\{0\}$ and $\bar T_C=\{1,,2\cdots,m\}$, the row and column
indices of $(d_1,d_2,\cdots,d_m)$ in the reduced formal product
$\bar\Theta$ in Formula (4.1-2),
then the vertex set $\bar T=\bar T_R\times \bar T_C$.

Let $(\bar{\mf A}',\bar{\mf B}')$ be an induced pair of $(\bar{\mf A},\bar{\mf B})$,
with the induced functor $\bar\vartheta:R(\bar {\mf A}')\rightarrow R(\bar{\mf A})$.
Write $\bar M=\bar\vartheta(F(k))=\sum_{j=1}^m\bar M_j\ast E_j\in R(\bar{\mf A})$ with the size vector
$\n=(n_0;n_1,\cdots,n_m)$ over $\bar\T$.

\medskip

{\bf Remark 5.4.1}\, Based on Formula (2.4-5) and using the reduced form
parallel to Formula (4.1-7), we will construct the defining system $\bar \IF$ of the quotient-sub-pair
$(\bar {\mf A}',\bar {\mf B}')$.

(i)\, According to Remark 4.1.1 (i), denote by $\bar Z_0$ the
$(p^r,p^r)$-the square block of $\Psi_{\m^{r\prime}}$ of size $n_0$ with $p^r\in X^r$; by $\bar
Z_{\xi\xi}$ the $(q^r+\xi,q^r+\xi)$-th square block of size $n_\xi$
with $q^r+\xi\in Y^r_i$, then
$$\begin{array}{c}\bar Z_0=Z_{X^r}=(z^{X^r}_{pq})_{n_0\times n_0},
\quad \bar Z_{\xi\xi}=Z_{Y^r_i}=(z^{Y^r_i}_{pq})_{n_\xi\times n_\xi}.\end{array}
$$

(ii)\, Denote by $\bar Z_\xi$ the $(p^r,q^r+\xi)$-block of $\Psi_{\m^{r\prime}}$ of size
$n_0\times n_\xi$, and by $\bar Z_{\eta\xi}$ the $(q^r+\eta,q^r+\xi)$-block of size
$n_\eta\times n_\xi,\eta<\xi$. Write the variable matrix
$Z_j=(z^{j}_{pq})$ for $j=1,\cdots, t^r$, then
$$\begin{array}{lll} \bar Z_\xi=\sum_{j}\alpha_{\xi}^{j} Z_{j},&\mbox{where $j$
runs over}\,\, s(V^r_j)\ni p^r,\,\,e(V^r_{j})\ni q^r+\xi,& \alpha_\xi^j\in k;\\
\bar Z_{\eta\xi}=\sum_j\beta_{\eta\xi}^j Z_j,&\mbox{where $j$ runs over}\,\,
s(V^r_j)\ni q^r+\eta,\,\, e(V^r_j)\ni q^r+\xi, &\beta_{\eta\xi}^j\in k.
\end{array}$$
Fix an integer
$l\in \{1, \cdots,m\}$ with $l\in Y\ne X$ in Definition 4.1.2,
thus $d_l:X\rightarrow Y$ is a solid edge.
Let $\rho\in Y$, such that $(p_\rho,q_\rho+1)$ is
the leading position of $\bar H'(k)$ in $\bar{\mf A}'$, and suppose that
$(q_\rho+1)$ is the index of the
first column of the $\rho$-th block-column over $\bar\T$. Write
$$\begin{array}{c}{\IF}:\,\,(\bar Z_0\ast E_0)\bar M
\equiv_{\prec (p_\rho^i,q_\rho^i+1)}\sum_{\xi=1}^m\bar Z_\xi\ast E_{\xi}+
\bar M(\sum_{1\leqslant\eta\leqslant\xi\leqslant m}\bar Z_{\eta\xi}\ast
E_{\eta\xi}).\end{array}\eqno{(5.4\mbox{-}1)}$$
Similar as in Equation (5.1-3), we have the right hand side of $\bar\IF$:
$$\begin{array}{c}{\IF}_\tau:\,\,0
\equiv_{\prec (p_\rho,q_\rho+1)}\sum_{\xi=1}^m\bar Z_\xi\ast E_{\xi}+\bar M
(\sum_{1\leqslant \eta\leqslant\xi\leqslant m}\bar Z_{\eta\xi}\ast
E_{\eta\xi}).\end{array}\eqno{(5.4\mbox{-}2)}$$

Define a size vector $\tilde {\underline n}=(\tilde n_0;\tilde n_1, \cdots,\tilde n_m)$
over $\bar\T$ as follows: $\tilde n_\xi=n_\xi$ if $\xi\notin Y$; $\tilde
n_\xi=n_\xi+1$ if $\xi\in Y$. Construct a representation based on
$\bar M$:
$$\tilde {M}=\sum_{j=1}^m \tilde{M}_j\ast E_j\in R(\bar{\mf A}),\quad
\tilde{M}_j=\left\{\begin{array}{ll}\bar{M}_j, & \mbox{if}\,\, E_j1_{Y}=0;\\
(0\,\bar{M}_j),&\mbox{if}\,\,
E_j1_{Y}=E_j,\end{array}\right.\eqno{(5.4\mbox{-}3)}$$ with
$0$ a column vector. Write $\tilde Z_0, \tilde Z_{\xi\xi}$ the variable matrices of size
$\tilde n_0\times \tilde n_0, \tilde n_\xi\times \tilde n_\xi$;
and $\tilde Z_\xi=\sum_j\alpha_{\xi}^j \tilde Z_j$ of size $\tilde n_0\times\tilde n_\xi$,
$\tilde Z_{\eta\xi}=\sum_{\eta\xi}^j\beta_{\eta\xi}^j \tilde Z_j$ of size $\tilde n_\eta\times\tilde n_\xi$
according to Remark 5.4.1 respectively. Thus we
obtain the following matrix equation with $\tilde q_{\rho}$ being the index of the
first column of the $\rho$-th block-column of $\tilde M$:
$$\begin{array}{ccc}
{\tilde{\IF}}:&(\tilde Z_0\ast E_0)\tilde{M}
&\equiv_{\prec (p_\rho,\tilde q_{\rho})}\sum_{\xi=1}^m\tilde Z_{\xi}\ast
E_{\xi}+\tilde{M} (\sum_{1\leqslant \eta\leqslant\xi\leqslant
m}\tilde Z_{\eta\xi}\ast E_{\eta\xi}),\\
{\tilde{\IF}}_\tau:&0&\equiv_{\prec (p_\rho,\tilde
q_{\rho})}\sum_{\xi=1}^m\tilde Z_{\xi}\ast E_{\xi}+\tilde{M}
(\sum_{1\leqslant \eta\leqslant\xi\leqslant m}
\tilde Z_{\eta\xi}\ast E_{\eta\xi}).\end{array}\eqno{(5.4\mbox{-}4)}$$

\medskip

Taken any integer $p'\geqslant p_\rho$ and $1\leqslant
h\leqslant n_{\rho}$, the $(p',q_\rho+h)$-th entry of
$\IF_\tau$ equals
$$\begin{array}{c} \sum_{p'}\gamma_{p,p'}z^{Y^r}_{p',q_\rho^i+h} +\sum_{p'}\nu_{p,p'}z^{j}_{p',
q_\rho^i+h},\quad \gamma_{p,p'},\nu_{p,p'}\in k.\end{array}\eqno{(5.4\mbox{-}5)}$$
It is clear that $z^{Y}_{p',q_\rho+h_1}$ and $z^{Y}_{p',q_\rho+h_2}$, $z^{j}_{p',q_\rho+h_1}$
and $z^{j}_{p',q_\rho+h_2}$ in Formula (5.4-5) have the same coefficients
respectively for all $1\leqslant h_1, h_2\leqslant n_\rho, h_1\ne
h_2$. The assertion is also valid for $\tilde\IF_\tau$.

The picture below shows four equations (abridged by four circles)
of $\IF_\tau$ in Formula (5.4-2). There are three solid
edges ending at $Y$, where $|Y|=3$, $n_\rho=4$, the equations
locate at the $(p_2,q_2+h)$-th positions for $h=1,2,3,4$ have
the same coefficients as shown in Formula (5.4-5).

\medskip

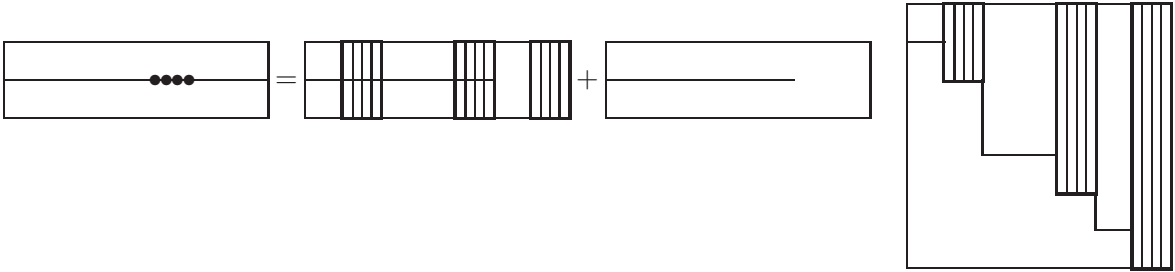
\begin{center}
\setlength{\unitlength}{1mm}
\begin{picture}(170,33)
\thinlines{\linethickness{0.1mm}}
\put(5,25){\framebox(35,10)}\put(5,30){\line(1,0){35}}
\put(25,30){\circle*{1.5}}\put(26.5,30){\circle*{1.5}}\put(28,30)
{\circle*{1.5}}\put(29.5,30){\circle*{1.5}}
\put(41,29){$=$}
\put(81,29){$+$}
\put(45,25){\framebox(35,10)}
\mput(51.25,25)(1.25,0){3}{\line(0,1){10}}\mput(66.25,25)(1.25,0){3}{\line(0,1){10}}
\mput(76.25,25)(1.25,0){3}{\line(0,1){10}}\put(45,30){\line(1,0){25}}

\put(85,25){\framebox(35,10)}\put(85,30){\line(1,0){25}}

\put(125,5){\framebox(35,35)}
\mput(131.25,30)(1.25,0){3}{\line(0,1){10}}\mput(146.25,15)(1.25,0){3}{\line(0,1){25}}
\mput(156.25,5)(1.25,0){3}{\line(0,1){35}}
\put(125,35){\line(1,0){5}}
\put(135,20){\line(1,0){10}}\put(135,20){\line(0,1){10}}
\put(150,10){\line(1,0){5}}\put(150,10){\line(0,1){5}}

\thicklines{\linethickness{1mm}}
\put(50,25){\framebox(5,10)}\put(65,25){\framebox(5,10)}\put(75,25){\framebox(5,10)}
\put(130,30){\framebox(5,10)}\put(145,15){\framebox(5,25)}\put(155,5){\framebox(5,35)}
\end{picture}
\end{center}

{\bf Remark 5.4.2}\,  Parallel to Remark 5.1.4, we have the
following facts.

(i)\, For any $1\leqslant h_1,h_2\leqslant n_\rho$, the $(p_\rho,q_\rho+h_1)$-th
equation is a linear combination of
the previous equations in $\IF_\tau$, if and only if so is the
$(p_\rho,q_\rho+h_2)$-th equation by Formula (5.4-5). Similarly, we have
the same result in $\tilde{\IF}_\tau$.

(ii)\, The equations in the system $\tilde{\IF}$ (resp. $\tilde{\IF}_{\tau}$) and
those in ${\IF}$ (resp. ${\IF}_{\tau}$)  are the same everywhere except at the $\tilde
q_{\rho}$-column, since the entries of $\tilde{M}$
at the $\tilde q_{\rho}$-th column is zero for $\rho=1,\cdots,|Y|$ respectively.

(iii)\, If there exists some $1\leqslant h\leqslant n_\rho$, such
that the $(p_\rho,q_\rho+h)$-th equation is a linear
combination of the previous equations in $\IF_\tau$, then so is
the $(p_\rho,\tilde q_{\rho})$-th equation in $\tilde\IF$, by (ii) and (i) above.

(iv)\, If the $(p_\rho,q_\rho+h)$-th equation is a linear combination of the previous
equations in $\IF$, then so is in $\IF_\tau$, since the variables in
$\{\bar Z_0,\bar Z_{\rho\rho}, Z_j\}_j$ are algebraically independent.

\medskip

Suppose
there is a reduction sequence given by Formula (4.1-5) with $\bar{\mf A}^0=\bar{\mf A}$
and each reduction being in the sense of Lemma 2.3.2:
$$\begin{array}{llllllll}\mf A, \mf A^1,\cdots,\mf A^{r-1},
&\mf A^r,&\mf A^{r+1},&\cdots,&\mf A^{r+i},&\mf A^{r+i+1},&\cdots,&
\mf{A}^{r+s};\\
&\bar{\mf{A}}^0,&\bar{\mf{A}}^1,&\cdots,&\bar{\mf{A}}^{i},&\bar{\mf{A}}^{i+1},&\cdots,&
\bar{\mf{A}}^s,\end{array}\eqno {(\bar \ast)}
$$

{\bf Proposition 5.4.3}\, Parallel to Theorem 5.2.1, there exists a unique reduction sequence based
on the sequence $(\bar \ast)$:
$$\begin{array}{llllllll}\tilde{\mf A}, \tilde{\mf A}^1,\cdots,\tilde{\mf A}^{r-1},
&\tilde{\mf A}^r,&\tilde{\mf A}^{r+1},&\cdots,&\tilde{\mf A}^{r+i},&\tilde{\mf
A}^{r+i+1},&\cdots,&\tilde{\mf A}^{r+s};\\
&\tilde{\bar{\mf{A}}},&\tilde{\bar{\mf{A}}}^1,&\cdots,&\tilde{\bar{\mf{A}}}^{i},
&\tilde{\bar{\mf{A}}}^{i+1},&\cdots,&
\tilde{\bar{\mf{A}}}^s.\end{array}\eqno {(\tilde{\bar\ast})}
$$

(i)\, $\tilde{\mf{A}}^i=\mf{A}^i$ for $i=0,1,\cdots,r$.

(ii)\, The reduction from $\tilde{\bar{\mf{A}}}^{i}$ to
$\tilde{\bar{\mf{A}}}^{i+1}$ is a reduction or a composition of two
reductions in the sense of Lemma 2.3.2 for $i=0,\cdots,s-1$, such
that $\tilde{\bar\T}^{s}$ has two vertices, and
$\tilde{\vartheta}^{0s}(\tilde{F}^s)=\tilde{M}^s$.

(iii)\, The reduction from $\tilde{\mf{A}}^{r+i}$ to $\tilde{\mf{A}}^{r+i+1}$ is as the same as
that from $\tilde{\bar{\mf{A}}}^{i}$ to $\tilde{\bar{\mf{A}}}^{i+1}$.

(iv)\, The diagonal block $\tilde e_{X}$ in $\tilde{\bar\K}^{s}_0$ of
$\tilde{\bar {\mf A}}^s$ over $\bar\T$ is of the form given in Corollary 5.2.2.

\smallskip

{\bf Proof}\, (i) is clear. The proof of (ii) is parallel to that of
Theorem 5.2.1, the only difference
is that for each $\rho\in Y$, we need to add a column into the $\rho$-th
block column for each $\rho=1,\cdots, |\bar Y|$ step by step.
(iii)\, follows from Formula (4.1-5). The proof of
(iv)\, is parallel to Corollary 5.2.2. The proof of the theorem is completed.

\medskip

Parallel to $(\ast')$ of the subsection 5.3, suppose we have the following sequences:
$$\begin{array}{llllllllll}\mf A,\mf A^1,\cdots,\mf A^{r-1},
&\mf{A}^{r}, &\mf{A}^{r+1},&\cdots,&{\mf{A}}^{r+s},
&{\mf{A}}^{r+s+1},  &\cdots, &{\mf{A}}^{r+\epsilon},&\cdots,&{\mf
A}^{r+t}; \\
&\bar{\mf{A}},&\bar{\mf{A}}^1,&\cdots,&\bar{\mf{A}}^s,
&\bar{\mf{A}}^{s+1},
&\cdots,&\bar{\mf{A}}^\epsilon,&\cdots,&\bar{\mf A}^t,\end{array} \eqno{(\bar\ast')}
$$
the reduction from $\bar{\mf A}$ (resp. ${\mf A}^{r}$)
to $\bar{\mf A}^{s}$ (resp. ${\mf A}^{r+s}$) is given by $(\bar\ast)$; from $\bar{\mf A}^{s}$ (resp. $\mf A^{r+s}$)
to $\bar{\mf A}^{s+1}$ (resp. $\mf A^{r+s+1}$) is a loop mutation and
we obtain a parameter $x$; the reduction from $\bar{\mf A}^{i}$ (resp. ${\mf A}^{r+i}$)
to $\bar{\mf A}^{i+1}$ (resp. ${\mf A}^{r+i+1}$) is a
regularization for $i=s+1,\cdots,t-1$.  The pair $(\mf A^{r+t},\mf B^{r+t})$ is
minimally wild in the case of MW5 and Classification 5.1.1 (II).

\medskip

{\bf Remark 5.4.4}\, (i)\, If the first arrow $a^t_1$ of $\bar{\mf B}^t$ splits from $d_l$ of
the one sided sub-bocs $\bar{\mf B}$,
then $d_l:X\mapsto Y$ is an edge by Theorem 4.6.1 and Corollary 4.6.2. Consequently
we are able to add some columns according to Theorem 5.4.3.

(ii)\, We will discuss how to determine
$\mf A^r$ in $(\bar\ast')$ in the next subsection.

(iii)\, Suppose $a_1^t$ locates at the $(p,q')$-th position
in the reduced formal product $\bar\Theta^t$ with $q'=q+j$ for some $1\leqslant j\leqslant n_l$,
and the first arrow $a^\epsilon_1$ of $\bar{\mf B}^\epsilon$ locates at the $(p,q+1)$-th
position in $\bar\Theta^\epsilon$.
Let integers $\rho\in Y$, parallel to Formula (5.3-1) we write an equation
system $\tilde{\IF}^{(\bullet,\rho)}_\tau$. Thus we obtain
a polynomial $d^{\eta\rho}(x)$ inductively for $\rho=1,\cdots,|Y|,\eta=p_x,\cdots,p$,
according to the discussion given
in the first part of the subsection 5.3. If $\bar R^t=k[x,\phi^t(x)^{-1}]$, define
$$\begin{array}{c}\tilde\phi^t(x)=\phi^t(x)\prod_{\eta=p_x}^p\prod_{\rho=1}^{|Y|}d^{\eta\rho}(x),\end{array}$$

\medskip

{\bf Proposition 5.4.5}\, Parallel to Theorem 5.3.4, there exists a unique reduction sequence based
on the sequence $(\bar \ast')$ satisfying the proposition (i)-(iv) below:
$$\begin{array}{llllllllll}\tilde{\mf A},\tilde{\mf A}^1,\cdots,\tilde{\mf A}^{r-1}
&\tilde{\mf{A}}^{r}, &\tilde{\mf A}^{r+1},&\cdots,&\tilde{\mf{A}}^{r+s},
&\tilde{\mf{A}}^{r+s+1},  &\cdots,
&\tilde{\mf{A}}^{r+\epsilon},&\cdots,&\tilde{\mf
A}^{r+t}; \\
&\tilde{\bar{\mf{A}}},&\tilde{\bar{\mf{A}}}^1,&\cdots,&\tilde{\bar{\mf{A}}}^s,
&\tilde{\bar{\mf{A}}}^{s+1},
&\cdots,&\tilde{\bar{\mf{A}}}^\epsilon,&\cdots,&\tilde{\bar{\mf
A}}^t.\end{array} \eqno{(\tilde{\bar\ast}')}
$$

(i) The first parts of the two sequences up to $r+s$ and $s$ respectively
are given by $(\tilde{\bar\ast})$.

(ii)\, $\tilde{\bar{\mf A}}^{s+1}$ is induced from
$\tilde{\bar{\mf A}}^s$ by a loop mutation $a^{s+1}_1\mapsto(x)$,
or an edge reduction $(0)$, then a loop mutation $(x)$; the reduction from
$\tilde{\bar{\mf A}}^{s+i}$ to $\tilde{\bar{\mf A}}^{s+i+1}$ is given by a
regularization, or a composition of two regularizations, or an edge
reduction $(0)$ then a regularization for $i=1,\cdots,t-1$.

(iii) The reduction from $\tilde{{\mf A}}^{r+s+i}$ to $\tilde{{\mf
A}}^{r+s+i+1}$ is as the same as that from $\tilde{\bar{\mf A}}^{s+i}$ to
$\tilde{\bar{\mf A}}^{s+i+1}$ for $i=1,\cdots,(t-s-1)$.

(iv) Denote by $a_1^t\in \bar{\mf B^t}$ the first solid arrow,
if $\dz(a_1^t)=vx-xv$, then the first solid edge $\tilde{a}^\epsilon_0$
of $\tilde{\bar{\mf B}}^\epsilon$ with the differential $\dz(\tilde{a}^\epsilon_{0})=0$.

\smallskip

{\bf Proof}\, (i) is obvious. The proof of (ii) is parallel to that of Theorem 5.3.2. (iii)
follows from Formula (4.1-5). (iv) is parallel to Corollary 5.3.3. The proof is finished.

\bigskip
\bigskip
\noindent{\bf 5.5 Non-homogeneity in the case of MW5 and classification (II)}
\bigskip

Suppose a bipartite pair $(\mf A,\mf B)$ has an induced pair $(\mf A',\mf B')$ in the case
of MW5 and Classification 5.1.1 (II). This subsection is devoted
to determining the one sided quotient-sub pair according to the position of the first arrow $a_1'$ in
the formal product $H'+\Theta'$, and then proving that $(\mf A,\mf B)$ is not homogeneous.

\medskip

Let $\mf A$ be a bipartite matrix bi-module problem having RDCC condition.
The sequence
$$
(\mf{A},\mf B),(\mf{A}^1,\mf B^1),\cdots, (\mf{A}^{\varsigma},\mf{A}^{\varsigma}),
(\mf{A}^{\varsigma+1},\mf{B}^{\varsigma+1}),\cdots,
(\mf{A}^{\tau},\mf B^\tau)=(\mf A',\mf B')\eqno{(5.5\mbox{-}1)}
$$
satisfies the following condition:
$R^i$ is trivial for $i\leqslant\varsigma$, the reduction from $\mf{A}^{i}$ to $\mf{A}^{i+1}$ is in the
sense of Lemma 2.3.2 for $i<\varsigma$; $\mf{A}^{\varsigma}$ is
local with $\dz(a^\varsigma_1)=0$, after a loop mutation, we have
$R^{\varsigma+1}=k[x]$ in $\mf B^{\varsigma+1}$;
finally we make regularization for $i>\varsigma$, and $\mf B^\tau=\mf B'$ is in the case of MW5 and
Classification 5.1.1 (II). Suppose the first arrow  $a_1^\tau$  of
$\mf B^\tau$ locates at the $(p^\tau,q^\tau)$-th position of the
$(\textsf p,\textsf q)$-th block in the formal product $\Theta^\tau$.
Since the size of $H^\varsigma$ coincides with that of $H^\tau$, and since we make regularization from
$\mf B^{\varsigma+1}$ to $\mf B^\tau$, according to Formula (2.3-5):
$$\begin{array}{c}H^\tau(k[x,\phi^\tau(x)^{-1}])=\sum_{i=1}^\varsigma G_\tau^i\ast A^{i-1}_1
+(x)\ast A^{\varsigma}_1.
\end{array}\eqno{(5.5\mbox{-}2)}$$
From now on, we call $G^i_\tau$, the leading block of $G^i_\tau\ast A_1^{i-1}$,
a {\it $G$-type matrix} of $H^\tau(k[x,\phi^\tau(x)^{-1}])$, and
sometimes do not distinguish $G^i$ over $R^i$ or $G^i(k)$ over $k$.

Let $i<\varsigma$, $\vartheta^{i\tau}: R(\mf A^\tau)\rightarrow R(\mf A^i)$
be the induced functor, and $\n^{i\tau}=\vartheta^{i\tau}(1,1,\cdots,1)$.
There is a simple fact, that any row (column)
index $\rho$ of $H^i+\Theta^i$ in the pair $(\mf A^i,\mf B^i)$ determines a row (column) index
$n^{i\tau}_1+\cdots+n^{i\tau}_\rho$ of $H^\tau+\Theta^\tau$ in the pair $(\mf A^\tau,\mf B^\tau)$.
Consequently, if the upper (resp. lower, left or right) boundaries of two $G$-type matrices
$G^i(1),G^i(2)$ of $H^i(k)$ are colinear, then
the same boundaries of two splitting blocks $G_\tau^i(1),G^i_\tau(2)$
of $H^\tau(k[x,\phi(x)^{-1}])$ are still colinear.

\medskip

{\bf Remark 5.5.1}\, Consider the $G$-type matrices inside the
$(\textsf p,\textsf q)$-th block. If $G^i_\tau$ and $G^{i+1}_\tau$
are both in the $(\textsf p,\textsf q)$-th block, the relative position of their
upper boundaries has three possibilities
according to Formulae (2.3-1)-(2.3-3).

(i) The upper boundaries of $G^i_\tau$ and $G^{i+1}_\tau$ are co-linear,
if and only if the reduction from $\mf A^{i-1}$ to $\mf A^{i}$ is
given by one of the following: $G^{i}=(0), (1),(0\,1)$;
$G^i=(\lambda)$; $G^{i}=\emptyset$; and the right boundary of $G^i_\tau$ is not that of
the $(\textsf p,\textsf q)$-th block. In this case their lower boundaries are also colinear.

(ii)\, The upper boundary of $G^{i+1}_\tau$ is strictly lower than that of $G^i_\tau$,
if and only if $G^{i}={1\choose 0}$ or ${{0\,1}\choose{0\,0}}$; or $G^{i}=W$
of size being strictly bigger than $1$, and the right boundary of $G^i_\tau$ is not that of
the $(\textsf p,\textsf q)$-th block. In this case,
the lower boundaries of $G^i_\tau$ and $G^{i+1}_\tau$ are also colinear.

(iii)\, The lower boundary of $G^{i+1}_\tau$ is the upper boundary of $G^i_\tau$,
if and only if the right boundary of $G^i_\tau$ coincides with that
of the $(\textsf p,\textsf q)$-th block.

\medskip

Collect all the $G$-type matrices of $H^\tau$ inside the $(\textsf p,\textsf q)$-th
block, such that their upper boundaries are above or at that of $a^\tau_1$:
$$G_\tau^{q_1},G_\tau^{q_2},\ldots, G_\tau^{q_u},\quad \mbox {with}\quad
q_1<q_2<\cdots<q_u<\varsigma.\eqno{(5.5\mbox{-}3)}$$
The $G$-type matrices $G_\tau^{q_i} (1\leqslant i\leqslant u)$ in (5.5-3) are
grouped into $h$ groups according to their upper boundaries
are colinear or not, and denoted by $\rho_j$ the common upper boundary of the $j$-th group
for $j=1,\cdots,h$, where $\rho_{j+1}$ is
strictly lower than $\rho_{j}$:
$$\{G_\tau^{q_{11}},\dots,G_\tau^{q_{1,u_1}}\},
\cdots\cdots, \{G_\tau^{q_{h,1}},\cdots, G_\tau^{q_{h,u_{h}}}\}, \quad u_1+\cdots+u_{h}=u.
\eqno{(5.5\mbox{-}4)}$$
The matrices $G_\tau^{q_{j,l}}$ and $G_\tau^{q_{j,l+1}}$ in
the $j$-th group have two possibilities:
\ding{172} If $G_\tau^{q_{j,l}}$ is in the case of Remark 5.5.1 (i), then $G_\tau^{q_{j,l+1}}$ comes from
the next reduction with $q_{j,l+1}=q_{j,l}+1$. \ding{173} If $G_\tau^{q_{j,l}}$ is in the case of
Remark 5.5.1 (ii), then $G_\tau^{q_{j,l+1}}$ follows by
a sequence of reductions with the upper boundary of the $G$-type matrices lower
than that of the $p^\tau$-row, and including at least one reduction in the case of Remark 5.5.1 (iii). At last
the sequence reaches $G_\tau^{q_{j,l+1}}$ with the upper boundary
$\rho_j$ as a neighbor of $G_\tau^{q_{j,l}}$, thus $q_{j,l+1}>q_{j,l}+1$.

\medskip

{\bf Lemma 5.5.2}\, $G_\tau^{q_{j,u_j}}$ must be in the case of Remark 5.5.1 (ii) for $j=1,\cdots,h$.

\smallskip

{\bf Proof}\, If $G_\tau^{q_{j,u_j}}$ is in the case of 5.5.1 (iii), then $\rho_j$
is lower than $\rho_{j+1}$, a contradiction to the grouping of Formula (5.5-4);
and $a_1^\tau$ is sitting upper $\rho_{h}$ for $j=h$, a contradiction to the choice
of the sequence (5.5-3).

Suppose $G_\tau^{q_{j,u_j}}$ satisfies 5.5.1 (i).  Then for $j<h$,
the upper boundaries of $G_\tau^{q_{j,u_j}}$ and $G_\tau^{q_{j,u_j}+1}$ coincide,
a contradiction to the grouping in (5.5-4). For $j=h$, $G^{q_{h,u_h}}_\tau=0, I,(0\,I)$, or
$\lambda I$, or $\emptyset$ with hight $d\geqslant 1$.
Suppose the next reduction gives $G^{q_{h,u_h}+1}_{\tau}$ and denoted by $G'_\tau$ for simplicity.
If $G'_\tau$ satisfies 5.5.1 (i) and (ii), then $G'_{\tau}$ and $G^{q_{h,u_h}}_\tau$
have the same upper boundary, a contradiction to the grouping of (5.5-4); if $G'_\tau$ satisfies 5.5.1 (iii),
then $a_1^\tau$ locates above $\rho_h$, a contradiction to the choice of (5.5-3).
Therefore there is no any further reduction in the sense of Lemma 2.3.2.
If the hight $d>1$, $\mf B^\tau$ is not local, so $d=1$. Since $a^\tau_1$
locates between the lower and the upper boundary of $G^{q_{h,u_h}}_\tau$,
which forces $G^{q_{h,u_h}}_\tau$ sitting
at the $p^\tau$-th row. But the parameter $x$ appears after
$G_\tau^{q_{h,u_h}}$ and before $a_1^\tau$, thus
locates at the $p^\tau$-th row, a contradiction to Lemma 5.1.2.
Therefore $G^{q_{j,u_j}}_\tau$ is in the case of 5.5.1 (ii), the proof is completed.

\medskip

{\bf Definition 5.5.3}\, We define $h$ rectangles in $\Theta^\tau$: for $j<h$,
the $j$-th rectangle has the upper boundary $\rho_j$, lower boundary
$\rho_{j+1}$, and the left boundary is the right boundary of $G^{q_{j,u_j}}$,
the right boundary is that of the $(\textsf p,\textsf q)$-th block.
While the $h$-th rectangle has the upper boundary $\rho_h$, lower boundary
is that of $G^{q_{h,u_h}}$. The rectangle with the upper boundary $\rho_j$ is
said to be the $j$-th {\it lader}, there are altogether $h$ laders.

\medskip

The picture below shows an example, in which $h=3$, the three groups given by sawtooth patterns with some dots,
and the last $G$-type matrix in each group is given by a rectangle without dots.
The upper boundaries of the three laders are shown by dotted lines.

\medskip

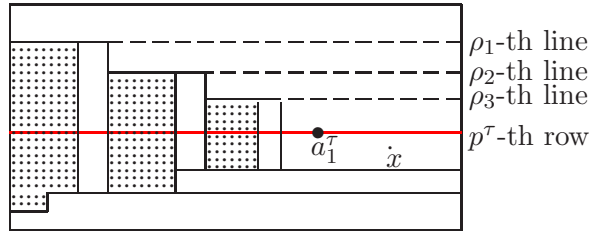
\begin{center}
\setlength{\unitlength}{1mm}
\begin{picture}(80,30)
\put(0,0){\line(1,0){60}}  \put(0,0){\line(0,1){30}}
\put(0,30){\line(1,0){60}} \put(60,0){\line(0,1){30}}

\put(0,2.5){\line(1,0){5}}\put(5,2.5){\line(0,1){2.5}}
\put(5,5){\line(1,0){55}} {\color{red}\put(0,13){\line(1,0){60}}}
\put(9,5){\line(0,1){20}} \put(13,5){\line(0,1){20}}
\put(0,25){\line(1,0){13}}
\multiput(13,25)(3,0){16}{\line(1,0){2}}
\put(22,5){\line(0,1){16}}
\put(26,8){\line(0,1){13}} \put(33,8){\line(0,1){9}}
\put(36,8){\line(0,1){9}} \put(13,21){\line(1,0){13}}
\multiput(25,21)(3,0){12}{\line(1,0){2}}
\put(26,17.5){\line(1,0){10}}\multiput(37,17.5)(3,0){8}{\line(1,0){2}}

\put(22, 8){\line(1,0){38}} \put(40,12){$\bullet$}\put(40, 10){$a^\tau_1$}
\put(50, 10){$\cdot$}\put(50,8.5){$x$}

\put(61, 11.3){$p^\tau$-th row}
\put(61, 16.8){$\rho_3$-th line}
\put(61, 20){$\rho_2$-th line}
\put(61, 24){$\rho_1$-th line}

\multiput(0,-2.5)(1,0){5}{\begin{picture}(0,0)\multiput(0,5)(0,1){2}{$\cdot$}\end{picture}}
\multiput(0,0)(1,0){9}{\begin{picture}(0,0)\multiput(0,5)(0,1){19}{$\cdot$}\end{picture}}
\multiput(13,0)(1,0){9}{\begin{picture}(0,0)\multiput(0,5)(0,1){15}{$\cdot$}\end{picture}}
\multiput(26,2.5)(1,0){7}{\begin{picture}(0,0)\multiput(0,5)(0,1){9}{$\cdot$}\end{picture}}
\end{picture}
\end{center}

\medskip

{\bf Lemma 5.5.4}\, Let $r=q_{h,u_h}-1$ in the formula (5.5-1).
We define a one sided quotient-sub pair $(\bar{\mf A},\bar{\mf B})=((\mf A^r)^{[m]},(\mf B^r)^{(m)})$
of the pair $(\mf A^r,\mf B^r)$ consisting of the solid arrows
$d_1,\cdots,d_m$ sitting at the $p^r$-row with the right boundary of $d_m$ is that
of the $(\textsf p,\textsf q)$-th block as shown in Formula (4.1-1). Then

(i)\, $m>1$;

(ii)\, all the $G$-type matrices in $H^\tau$
coming from $d_2,\cdots,d_m$ locate below the $p^\tau$-th row;

(iii)\, $a_1^\tau$ is split from $d_l$ with $l>1$. If $(\bar{\mf A},\bar{\mf B})$ satisfies
Theorem 4.6.1 or Corollary 4.6.2, then $d_l$ is a solid edge.

(iv)\, $\varsigma=r+s$ and $\tau=r+t$ in the formula (5.5-1), which coincides with
sequence $(\bar\ast')$ given above Remark 5.4.4.

\smallskip

{\bf Proof}\, (i)\, follows from Lemma 5.5.2.
(ii)\, comes from the choice of the sequence (5.5-3).
(iii) and (iv) are obvious. The proof is finished.

\medskip

{\bf Lemma 5.5.5}\, (i) Let $(\mf A_{X}^r,\mf B^r_{X})$ be the induced local pair at $X$
of $(\mf A^r,\mf B^r)$ defined above, denote by $h_{_{X}}$ the number
of the inheriting ladders in $H^r_{X}+\Theta^r_{X}$, then $h_{_{X}}\leqslant h$.

(ii)\, $H^{r}+\Theta^{r}$ of the pair $(\mf A^{r},\mf B^{r})$
has also $h$ ladders in the $(\textsf p,\textsf q)$-th block.
The boundaries of the $j$-th ladder of $H^\tau+\Theta^\tau$ comes from
that of the $j$-th ladder of $H^{r}+\Theta^{r}$ for $j=1,\cdots,h$,
according to the simple fact above Remark 5.5.1.

(iii)\, Return to the formulae $(\bar\ast')$ and $(\tilde{\bar\ast}')$ of Proposition 5.4.5,
then $\tilde H^{r+t}+\tilde\Theta^{r+t}$ has also $h$ ladders.
Furthermore, the number of rows in the $h$-th (resp. in the $j$-th, for $j=1,\cdots,h-1$)
ladder in $\tilde H^{r+t}+\tilde\Theta^{r+t}$ is as the same as (resp. as the same as or more than) that in
$H^{r+t}+\Theta^{r+t}$.

\medskip

{\bf  Proposition 5.5.6}\, Let $\mf{A}=(R,\K,\M,H=0)$ be a bipartite
matrix bi-module problem having RDCC condition. If there exists an induced pair $(\mf
A',\mf B')$ of $(\mf A,\mf B)$, which satisfies MW5 and
$H'+\Theta'$ of $(\mf A',\mf B')$ is in the case of
Classification 5.1.1 (II), then $\mf{A}$ is not homogeneous.

\smallskip

{\bf Proof}\, Suppose the induced pair $(\mf A',\mf B')$ is the last term
$(\mf A^{r+t},\mf B^{r+t})$ of the sequence $(\bar\ast')$ above Remark 5.4.4.
We assume in addition that the number of the laders in $H^{r+t}+\Theta^{r+t}$
is minimal with the property of MW5 and Classification 5.1.1 (II).

\smallskip

(I)\, Let $X$ be given by Definition 4.1.2. If $(\mf A^r_{X},\mf B^r_{X})$
is wild, denote by $((\mf A^r_{X})',(\mf B^r_{X})')$ the induced minimally wild
local pair obtained by using the triangular Formulae of subsection 3.3
with the parameter $x'$ and the first arrow $a_1'$.

(I-1)\, If $(\mf B^r_{X})'$ is of MW3, MW4, or MW5 with $H_{X}'+\Theta_{X}'$ being
in the case of Classification 5.1.1 (I), then it is not homogeneous by Proposition
3.5.3-3.5.5, we are done.

(I-2)\, If $(\mf B^r_{X})'$ is in the case of MW5 and Classification 5.1.1 (II),
then, the number of the inheriting ladders $h_{_X}$ in $H^r_{X}+\Theta^r_{X}$ with
$h_{_X}\leqslant h$ by Lemma 5.5.5 (i). Suppose $a_1'$ locates at the $h'$-lader. If $h'=h_{_X}=h$,
since this ladder contains only one row by Lemma 5.5.4, $x'$ must locate at the same row,
a contradiction to Lemma 5.1.2,; if
$h'<h_X=h$, or $h'=h_X<h$, then it contradicts to the minimality assumption on the number of ladders.

\smallskip

(II)\, Suppose $(\mf A^r_{X},\mf B^r_{X})$ is tame infinite, the quotient-sub-pair
$(\bar {\mf A}_{X},\bar{\mf B}_{X})$ is in the case of Classification 4.2.1 (ii).

(II-1) If the one sided pair $(\bar{\mf A},\bar{\mf B})$ satisfies the hypothesis of Lemma 4.2.3
or 4.4.1, since the unique effective loop $\bar b$ of $\bar{\mf B}_{X}$ is as the same as that of $\mf B^r_{X}$,
$\mf B^r$ is not homogeneous.

(II-2)\, If $(\bar{\mf A},\bar{\mf B})$ satisfies Theorem 4.4.2, then we use
triangular formulae of the subsection 3.3 for the local wild pair $(\mf A^{r+2l},\mf B^{r+2l})$.
If we reach the cases of MW3, MW4, or MW5 and Classification 5.1.1 (I), we are done. If we meet
again MW5 and Classification 5.1.1 (II), the first arrow must be outside of the $h$-th ladder,
a contradiction to the minimality assumption.

\smallskip

(III)\, Now we consider the following two cases: (i) $\mf B^r_{X}$ is tame infinite,
$\bar{\mf B}_{X}$ satisfies Classification 4.2.1 (ii),
the pair $(\bar {\mf A},\bar{\mf B})$ is major and satisfies Formula (4.2-6);
(ii) $\mf B^r_{X}$ is tame infinite or finite, and $\bar{\mf B}_{X}$ is finite.
Then in both cases $d_l$ of $\bar{\mf B}$, from which $a_1^t$ split, is a solid edge by Lemma 5.5.4 (iii).
Consequently, Formula $(\tilde{\bar\ast})$ of Proposition 5.4.3
can be used with respect to $d_l$. Furthermore by Formula $(\tilde{\bar\ast}')$ of Proposition 5.4.5,
$\dz(\tilde{a}^\epsilon_0)=0$ in $\tilde{\bar{\mf A}}^{\epsilon}$ by 5.4.5 (iv).
Set the edge $\tilde{a}^\epsilon_0\mapsto (1)$ then all the other arrows split from
$d_l$ at the same row are mapped to $\emptyset$ by Proposition 5.4.3 (iv).
The induced pair is obviously local of tame infinite or wild type, then we are able to
use the triangular formulae once again, and obtain an induced pair $(\hat {\mf{A}}^1,\hat {\mf{B}}^1)$
in the cases (ii)-(iv) of Classification 3.3.3.

\smallskip

(III-1)\, If the induced local pair $(\hat {\mf{A}}^1,\hat {\mf{B}}^1)$ is tame
infinite, then Proposition 3.5.5 ensures that the two-point pair
$(\tilde{\mf A}^{r+\epsilon},\tilde{\mf B}^{r+\epsilon})$ is wild and non-homogeneous,
we are done.

(III-2)\, If $(\hat {\mf{A}}^{1},\hat {\mf{B}}^{1})$ is in the case of MW3,
or MW4, or MW5 and Classification 5.1.1 (I), then it is non-homogeneous, we are done.

(III-3)\, If $(\hat {\mf{A}}^{1},\hat {\mf{B}}^{1})$ is in the case of MW5
and classification 5.1.1 (II), and suppose in addition, whose first arrow locates at the
$h_1$-th ladder with $h_1<h$, which contradicts to the
minimality number assumption of the ladders.

(III-4)\, If $(\hat {\mf{A}}^{1},\hat {\mf{B}}^{1})$ is in the case of MW5
and classification 5.1.1 (II), and suppose in addition, whose first arrow locates still at the
$h$-th ladder. We need to do induction on some pairs of integers.
Denote by $\sigma$ the number of the rows in the $h$-th ladder of $H^{r+t}+\Theta^{r+t}$,
which is a constant after making some bordered matrices by Lemma 5.5.5 (iii);
and $m$, the number of the solid arrows in the pair $(\bar{\mf A},\bar{\mf B})$,
is also a constant.
Define a finite set with $\sigma m$ pairs:
$$\begin{array}{c}\mathcal S=\{(\varrho, \zeta)\,|\,1\leqslant \varrho\leqslant \sigma,
\zeta=1,\cdots, m\},\\
\mbox{ordered by}\quad(\varrho_1, \zeta_1)\prec (\varrho_2, \zeta_2)
\Longleftrightarrow \varrho_1>\varrho_2, \quad \mbox{or}\quad
\varrho_1=\varrho_2,\,\,\zeta_1<\zeta_2.\end{array}$$
Denote the induced minimally wild local pair $(\mf A^{r+t},\mf B^{r+t})$ in $(\bar\ast')$ by
$(\hat {\mf{A}}^{0},\hat {\mf{B}}^{0})$ for unifying the notations.
Let $(\varrho^0,\zeta^0)\in \mathcal S$, such that $\varrho^0=p$ is
the row-index of the first arrow $\hat a_1^0=a_1'$
in the $h$-th lader by Remark 5.4.4 (iii); $\zeta^0=l$,
since $a_1'$ splits from the edge $d_l$ by Lemma 5.5.4 (iii). Similarly denote by
$(\varrho^1,\zeta^1)\in \mathcal S$ determined by the first arrow
$\hat{a}^1_1\in \hat{\mf{B}}^1$. Proposition 5.4.3 (iv) ensures
$(\varrho^0, \zeta^0)\prec (\varrho^1, \zeta^1)$.

Now we start the procedure (III) from the pair $(\hat{\mf A}^{1},\hat{\mf B}^{1})$
instead of $(\hat{\mf A}^{0},\hat{\mf B}^{0})$.
If (III-4) appears repeatedly, then after finitely
many steps, we reach an induced pair of (III-1)-(III-3)
by induction on $\mathcal S$. The proof is completed.

\bigskip
\bigskip

{\bf 5.6 The proof of the main Theorem}

\bigskip

We are ready to prove the main theorem 3.

\medskip

{\bf  Theorem 5.6.1} Let $\mf{A}=(R,\K,\M,H=0)$ be a bipartite
matrix bi-module problem having RDCC condition. If $\mf A$ is of wild type, then
$R(\mf{A})$ is not homogeneous.

\smallskip

{\bf Proof}\, For any wild bocs, there exists an
induced $\mf{B}'$ satisfying one of MW1-MW5 according to
Classification 3.3.2. Then Proposition 3.4.1-3.4.4 proved that $\mf B'$
in the case of MW1-MW4 is not homogeneous. When $\mf B$ is the associated
bocs of a bipartite matrix bi-module problem having RDCC condition, Proposition 5.3.4 and 5.5.6
proved that the induced pair $(\mf A',\mf B')$ in the case of MW5 is not
homogeneous. Therefore $(\mf A,\mf B)$ is not homogeneous, the proof is
completed.

\medskip

{\bf Proof of Main Theorem 3}\,  Let $\Lambda$ be a finite-dimensional basic
algebra over an algebraically closed field $k$. If $\Lambda$ is of
wild representation type, then mod-$\Lambda$ is not homogeneous.

\smallskip

In fact, let $\mf{A}$ be the matrix bi-module problem
associated to $\Lambda$. Then $\mf{A}$ is bipartite, having RDCC condition by Remark
1.4.4, and is representation wild type. So the pair $(\mf{A},\mf{B})$ is not
homogeneous by Theorem 5.6.1. Note that there is an almost
one-to-one correspondence between almost split sequences in
$\mod\Lambda$ and almost split conflations in $R(\mf{B})$, see
\cite{B2} and \cite{ZZ}, therefore $\mod \Lambda$ is not
homogeneous. The proof is finished.

\bigskip
\bigskip

{\bf Acknowledgement}

\bigskip

We would like to express our sincere thanks to S.Liu for his
proposal of $\Delta$-algebra, to Y.Han for his
suggestion of the concept on co-bi-module problem. Y.B.Zhang
is indebted to R.Bautista for giving the open problem in 1991; and
to W.W.Crawley-Boevey, D.Simson, C.M.Ringel for some discussions
during the long procedure of solving the problem. In particular,
she is grateful to Y.A.Drozd for his invitation to visit
the Institute of Mathematics, National Academy of Sciences of
Ukraine, and helpful discussions.

\end{document}